\pdfoutput=1

\documentclass[oneside,reqno,11pt]{amsart}

\usepackage{preamble}

\myexternaldocument[I:]{corank1_ASW_I}
\myexternaldocument[II:]{corank1_ASW_II}
\myexternaldocument[III:]{corank1_ASW_III}

\title{Co-rank \texorpdfstring{$1$}{1} Arithmetic Siegel--Weil IV: \\ Analytic local-to-global}
\author{Ryan C. Chen}
\date{March 18, 2026}
\address{Department of Mathematics,
	Fine Hall - Washington Road.
        Princeton, NJ 08544, USA}
\email{rcchen.math@gmail.com}

\begin{document}
    
    \begin{abstract}
This is the fourth in a sequence of four papers, where we prove the arithmetic Siegel--Weil formula in co-rank $1$ for Kudla--Rapoport special cycles on exotic smooth integral models of unitary Shimura varieties of arbitrarily large even arithmetic dimension. 
Our arithmetic Siegel--Weil formula implies that degrees of Kudla--Rapoport arithmetic special $1$-cycles are encoded in near-central first derivatives of unitary Eisenstein series Fourier coefficients.

In this paper, we pin down precise normalizations for some $U(m,m)$ Siegel Eisenstein series, give local Siegel--Weil special value formulas with explicit constants, and record a geometric Siegel--Weil result for degrees of complex $0$-cycles. Using this, we complete the proof of our arithmetic Siegel--Weil results by patching together the local main theorems from our companion papers.
\end{abstract}

    \maketitle
    
    \tableofcontents
    
    \clearpage


        \section{Introduction}
        \label{sec:part_IV:intro}
            This paper is a continuation of our companion papers \cite{corank1_ASW_I.pdf,corank1_ASW_II.pdf,corank1_ASW_III.pdf}. In this paper, we complete the proof of our global arithmetic Siegel--Weil results using our local main theorems from \cite{corank1_ASW_I.pdf,corank1_ASW_II.pdf} and the geometric local-to-global reduction procedure from \cite{corank1_ASW_III.pdf}.

For the reader's convenience, we recall the statements of our global arithmetic Siegel--Weil results (\cref{ssec:part_IV:intro:arith_Siegel--Weil});
these were also sketched in \crefext{I:sec:intro}.  

In the present paper, the main new ingredients are (1) precise normalizations for Eisenstein series and their relation with singular Fourier coefficients of co-rank $1$, (2) classical local Siegel--Weil formulas with precise constants, and (3) a geometric Siegel--Weil formula for complex $0$-cycles, which will be treated in \cref{sec:geometric_Siegel-Weil} (along with an observation about complex volumes of unitary Shimura varieties, which may be of independent interest). 

For introductory purposes, \cref{ssec:part_IV:intro:classical_and_geometric_Siegel-Weil} contains some background on classical and geometric Siegel--Weil. This is included for comparison with arithmetic Siegel--Weil, and also helps fix some notation. The material in \cref{ssec:part_IV:intro:classical_and_geometric_Siegel-Weil} is mostly expository, but some of our formulations may be new, particularly in our normalizations for Eisenstein series. The same normalization choices play an amplified role in our main arithmetic Siegel--Weil results. We also mention some results (comparison of complex volume and degrees of complex $0$-cycles to Eisenstein series) which seem to be new or at least not explicit in the literature; see discussion following \cref{equation:intro:arith_Siegel-Weil:geometric_Siegel_mass,equation:intro:arith_Siegel_Weil:geometric_Siegel-Weil}.

In \cref{ssec:part_IV:intro:outline}, we outline the structure of this paper and its relation with our companion papers \cite{corank1_ASW_I.pdf,corank1_ASW_II.pdf,corank1_ASW_III.pdf}.

            \subsection{Eisenstein series}
            \label{ssec:part_IV:intro:Eisenstein}
                In our work, we focus on the unitary/Hermitian case.
Consider an imaginary quadratic field $F / \Q$ with ring of integers $\mc{O}_F$ and odd discriminant $\Delta$. Given $m \in \Z_{\geq 0}$ and an even integer $n \in \Z$, we consider the (normalized) \emph{Siegel Eisenstein series}
    \begin{align}\label{equation:part_IV:intro:classical_and_geometric_Siegel-Weil:Eisenstein}
    E^*(z,s)^{\circ}_n & \coloneqq \Lambda_m(s)^{\circ}_n \sum_{\begin{psmallmatrix} a & b \\ c & d \end{psmallmatrix} \in P_1(\Z) \backslash SU(m,m)(\Z)} \frac{\det(y)^{s - s_0}}{\det(c z + d)^n |\det (c z + d)|^{2(s - s_0)}}
    \end{align}
for the group
    \begin{equation}\label{equation:part_IV:intro:classical_and_geometric_Siegel-Weil:results:U(m,m)}
    U(m,m) \coloneqq \left \{ h \in \Res_{\mc{O}_F / \Z} \GL_{2m} : {}^t \overline{h} \begin{pmatrix} 0 & 1_m \\ -1_m & 0 \end{pmatrix} h = \begin{pmatrix} 0 & 1_m \\ -1_m & 0 \end{pmatrix} \right \}
    \end{equation}
where $\Lambda_m(s)^{\circ}_n$ is the normalizing factor
    \begin{align}
    \Lambda_m(s)^{\circ}_n &\coloneqq \frac{(2 \pi)^{m(m-1)/2}}{(-2 \pi i)^{nm}} \pi^{m(- s + s_0)} |\Delta|^{m(m-1)/4 + \lfloor m / 2 \rfloor (s + s_0)} 
    \\
    & \mathrel{\phantom{\coloneqq}} \cdot \left ( \prod_{j = 0}^{m - 1} \Gamma(s - s_0 + n - j) \cdot L(2s + m - j, \eta^{j + n})  \right ). \notag
    \end{align}
In \eqref{equation:part_IV:intro:classical_and_geometric_Siegel-Weil:results:U(m,m)}, the notation $1_m$ stands for the $m \times m$ identity matrix, we wrote $SU(m,m) \subseteq U(m,m)$ for the determinant $1$ subgroup, and we set $P_1 \coloneqq P \cap SU(m,m)$ for the Siegel parabolic $P \subseteq U(m,m)$ (consisting of $m \times m$ block upper triangular matrices). The variable $s \in \C$ is a complex parameter, we set $s_0 = (n - m)/2$, and the element $z = x + i y$ lies in Hermitian upper half-space (i.e. $x \in \mrm{Herm}_m(\R)$ and $y \in \mrm{Herm}_m(\R)_{>0}$; the latter means that $y$ is positive definite).\footnote{Here, the notation $\mrm{Herm}_m$ denotes a scheme over $\Spec \Z$, e.g. $\mrm{Herm}_m(\R)$ denotes $m \times m$ complex Hermitian matrices, and $\mrm{Herm}_m(\Q)$ denotes $m \times m$ Hermitian matrices with entries in $F$.} The symbol $\eta$ denotes the quadratic character associated to $F / \Q$ (via class field theory).
The sum in \eqref{equation:part_IV:intro:classical_and_geometric_Siegel-Weil:Eisenstein} is convergent for $\mrm{Re}(s) > m / 2$, and admits meromorphic continuation to all $s \in \C$. When $m = 1$, the expression in \eqref{equation:part_IV:intro:classical_and_geometric_Siegel-Weil:Eisenstein} is a classical Eisenstein series on the usual upper half-plane. 

The normalized Eisenstein series has a symmetric functional equation 
    \begin{equation}
    E^*(z,s)^{\circ}_n = (-1)^{m (m - 1)(n - m - 1)/2} E^*(z,-s)^{\circ}_n
    \end{equation}
as in \crefext{ssec:Eisenstein:global_normalized_Fourier:global_normalization}. Our definition of the normalizing factor $\Lambda_m(s)^{\circ}_n$ is motivated by symmetry of global and local functional equations, along with certain local special value formulas; see \crefext{sec:setup,sec:Eisenstein:Weil_rep,sec:local_Whittaker,sec:Eisenstein:local_functional_equations,sec:Eisenstein:global_normalized_Fourier} for further discussion. The function $\Lambda_m(s)^{\circ}_n$ should be closely related with the $L$-function of an Artin--Tate motive attached to the group $U(m,m)$, in the sense of Gross \cite{Gross97} (see \cite[{Remark 1.1.1}]{BH21}).

Given $T \in \mrm{Herm}_m(\Q)$, the Eisenstein series $E^*(z,s)^{\circ}_n$ has \emph{$T$-th Fourier coefficient}
    \begin{equation}\label{equation:part_IV:intro:classical_and_geometric_Siegel-Weil:Fourier_coeff}
    E^*_T(y,s)^{\circ}_n \coloneqq 2^{m(m-1)/2} |\Delta|^{-m(m-1)/4} \int_{\mrm{Herm}_m(\Z) \backslash \mrm{Herm}_m(\R)} E^*(z, s)^{\circ}_n e^{-2 \pi i \mrm{tr}(Tz)} ~dx
    \end{equation}
for $z = x + i y$ in Hermitian upper half-space, where this integral is taken with respect to the Euclidean measure\footnote{The factor $2^{m(m-1)/2} |\Delta|^{-m(m-1)/4}$ disappears in the (usual) equivalent ad\`elic formulation, upon taking a certain self-dual Haar measure. The ad\`elic formulations of \eqref{equation:part_IV:intro:classical_and_geometric_Siegel-Weil:Eisenstein} and \eqref{equation:part_IV:intro:classical_and_geometric_Siegel-Weil:Fourier_coeff} are used in \crefext{sec:setup}.} 
on $\mrm{Herm}_m(\R)$. The integral is convergent for $\mrm{Re}(s) > m / 2$, and admits meromorphic continuation to all $s \in \C$. When $\det T \neq 0$, there is a factorization into normalized \emph{local Whittaker functions}
    \begin{equation}
    E^*_T(y,s)^{\circ}_n = W^*_{T,\infty}(y,s)^{\circ}_n \prod_p W^*_{T,p}(s)^{\circ}_n
    \end{equation}
over all places (\cref{part:part_IV:Eisenstein}). The preceding setup is as in \crefext{I:ssec:intro:Eisenstein}, and is taken from loc. cit. essentially verbatim.

For example, if $n = 2$ and $m = 1$, we have
    \begin{align} \label{equation:intro:arith_Siegel-Weil:weight_2_Eisenstein}
    & W^*_{T,p}(s)^{\circ}_2 = p^{v_p(T)(s+1/2)} \s_{-2s}(p^{v_p(T)}) \quad \quad \s_{s}(r) \coloneqq \sum_{d \mid r} d^s
    \\
    & W^*_{T,\infty}(y,s)^{\circ}_2 = \Gamma(s - 1/2)^{-1} |4 \pi T y|^{s - 1/2} \int_a^{\infty} e^{- 4 \pi |T| y u} (u \pm 1)^{s + 1/2} u^{s - 3/2} ~ du \notag
    \end{align}
for any nonzero $T \in \Z$, where
the lower limit of integration is $a = 0$ and the sign $\pm$ is $+$ (resp. $a = 1$ and the sign $\pm$ is $-$) if $T > 0$ (resp. if $T < 0$).

            \subsection{Classical and geometric Siegel--Weil}
            \label{ssec:part_IV:intro:classical_and_geometric_Siegel-Weil}
                Let $V$ be an $n$-dimensional $F$-vector space, equipped with a non-degenerate Hermitian pairing $(-,-)$. Set $G = U(V)$ and assume $n > 0$. Fix a full-rank $\mc{O}_F$-lattice $L \subseteq V$. For simplicity, we assume in the introduction that $L$ is self-dual.\footnote{We always mean self-dual for the bilinear \emph{trace pairing} $\mrm{tr}_{F / \Q}(v,w)$ unless otherwise specified; see conventions in \crefext{I:ssec:Hermitian_conventions:lattices}.} Write $K_{L,f} \subseteq G(\A_f)$ for the stabilizer of $L \otimes_{\Z} \hat{\Z}$, where $\A_f$ denotes the finite ad\`ele ring of $\Q$.

First consider the case where $V$ is positive definite. Since we assumed $L$ is self-dual, this forces $n \equiv 0 \pmod{4}$ (by the global product formula for local invariants of Hermitian spaces).
Given any positive definite Hermitian $\mc{O}_F$-lattice $\mc{L}$, we set 
    \begin{equation}
    \mc{Z}_{T,\mc{L}} \coloneqq \{\underline{x} \in \mc{L}^m : (\underline{x}, \underline{x}) = T \},
    \end{equation}
where $(\underline{x}, \underline{x})$ denotes the Gram matrix\footnote{If $\underline{x}$ is the $m$-tuple $[x_1, \ldots, x_m]$, the notation $(\underline{x},\underline{x})$ will mean the matrix with $i,j$-th entry $(x_i, x_j)$. We often write e.g. $[x_1, \ldots, x_m]$ instead of $(x_1, \ldots, x_m)$ for tuples, to avoid confusion with Hermitian pairings $(-,-)$.} of $\underline{x}$. When $m \leq n$, we have
    \begin{align}
    2 \Lambda_n(0)^{\circ}_n & = \sum_{\substack{\text{$\mc{O}_F$-lattices } \mc{L} \\ \text{self-dual, rank $n$,} \\ \text{positive definite}}} \frac{1}{|\Aut(\mc{L})|} = \# [G(\Q) \backslash ( G(\A_f) / K_{L,f}) ] \label{equation:intro:arith_Siegel-Weil:Siegel_mass}
    \\
    \frac{2 \Lambda_n(0)^{\circ}_n}{\kappa \Lambda_m(s_0)^{\circ}_n} E^*_T(y,s_0)^{\circ}_n & = \sum_{\substack{\text{$\mc{O}_F$-lattices } \mc{L} \\ \text{self-dual, rank $n$,} \\ \text{positive definite}}} \frac{|\mc{Z}_{T,\mc{L}}|}{|\Aut(\mc{L})|} 
    \quad \quad \text{for any $T \in \mrm{Herm}_m(\Q)$,} \label{equation:intro:arith_Siegel-Weil:classical_Siegel-Weil}
    \end{align}
where $\kappa = 2$ (resp. $\kappa = 1$) if $m = n$ (resp. if $m < n$). The sums run over isomorphism classes of positive definite rank $n$ self-dual Hermitian $\mc{O}_F$-lattices, the notation $\Aut(\mc{L})$ means the (unitary) automorphism group of $\mc{L}$. The symbols $\# [-]$ and $|-|$ mean groupoid and set cardinality, respectively. That is, we have
    \begin{equation}\label{equation:intro:arith_Siegel-Weil:theta_series}
     \frac{2 \Lambda_n(0)^{\circ}_n}{\kappa \Lambda_m(s_0)^{\circ}_n} E^*(z,s_0)^{\circ}_n = \sum_{\substack{\text{$\mc{O}_F$-lattices } \mc{L} \\ \text{self-dual, rank $n$,} \\ \text{positive definite}}} \frac{1}{|\Aut(\mc{L})|} \Theta_{\mc{L}}(z)
    \end{equation}
which re-expresses the Eisenstein series at $s = s_0$ as a weighted sum of theta series for the lattices $\mc{L}$.

\Cref{equation:intro:arith_Siegel-Weil:Siegel_mass,equation:intro:arith_Siegel-Weil:classical_Siegel-Weil} are special cases of (unitary analogues of) the classical \emph{Siegel mass formula} and \emph{Siegel--Weil formula} respectively. For \cref{equation:intro:arith_Siegel-Weil:Siegel_mass}, see \crefext{proposition:geometric_Siegel-Weil:complex_volumes:positive_definite}. \Cref{equation:intro:arith_Siegel-Weil:classical_Siegel-Weil} follows from \cite[{Proposition 6.2}]{Ichino04}, \cite[{Theorem 1.1}]{Ichino07}, and \cite[{Theorem 2.2}]{Yamana11} (in combination with \cref{equation:intro:arith_Siegel-Weil:Siegel_mass}).

Next, consider the case where $V$ has arbitrary signature $(n - r, r)$. Since $L$ was assumed self-dual, this forces $n \equiv 2 r \pmod{4}$,
i.e. we must have $\varepsilon(V_{\R}) = (-1)^{n(n-1)/2 + r} = + 1$ in the notation of \crefext{I:ssec:Hermitian_conventions:lattices}.
We will mostly be interested in the cases $r = 1$ and $r = 0$.
There is an associated Hermitian symmetric domain $\mc{D}$ which parameterizes maximal negative definite subspaces of the complex Hermitian space $V_{\R}$. For sufficiently small open compact $K_f \subseteq K_{L,f}$ (so that we have manifolds instead of orbifolds, for simplicity), there is an associated complex Shimura variety
    \begin{equation}
    \mrm{Sh}_{K_{f},\C}(G) = [G(\Q) \backslash (\mc{D} \times G(\A_f) / K_{f})]
    \end{equation}
of dimension $(n - r) r$ (analytification suppressed from notation). In the signature $(n, 0)$ and $(n - 1, 1)$ cases respectively, we have ``geometric Siegel mass formulas''
    \begin{equation}\label{equation:intro:arith_Siegel-Weil:geometric_Siegel_mass}
    2 \Lambda_n(0)^{\circ}_n = \frac{\mrm{vol}(\mrm{Sh}_{K_{f},\C}(G))}{[K_{L,f} : K_f]} \quad \quad - 2 \Lambda_n(0)^{\circ}_n = \frac{\mrm{vol}(\mrm{Sh}_{K_{f},\C}(G))}{[K_{L,f} : K_f]}
    \end{equation}
where $\mrm{vol}(\mrm{Sh}_{K_{L,f},\C}(G))$ is the volume with respect to the Chern form of a certain dual tautological bundle. The case of signature $(n - 1, 1)$ may be extracted from \cite[{Theorem A}]{BH21}, see \crefext{proposition:geometric_Siegel-Weil:complex_volumes:n-1_1}.\footnote{In \crefext{proposition:geometric_Siegel-Weil:complex_volumes:n-1_1}, note that we took $\mrm{vol}(-)$ with respect to the tautological bundle. Here we are taking volume with respect to the dual bundle, which produces the minus sign in \cref{equation:intro:arith_Siegel-Weil:geometric_Siegel_mass}.} The case of signature $(n, 0)$ is an equivalent reformulation of the classical Siegel mass formula \cref{equation:intro:arith_Siegel-Weil:Siegel_mass}: if we allow the (stacky) level $K_f = K_{L,f}$, then there is a canonical equivalence of groupoids
    \begin{equation}
    \mrm{Sh}_{K_L,\C}(G) \cong \left \{ \begin{array}{l} \text{Hermitian $\mc{O}_F$-lattices $\mc{L}$ which are} \\ \text{self-dual and signature $(n,0)$} \end{array} \right \}
    \end{equation}
in that case.

In \emph{geometric Siegel--Weil formulas}, the sets $\mc{Z}_{T,\mc{L}}$ (from classical Siegel--Weil) are replaced by \emph{special cycles} $\mc{Z}_{T,\C}$ over the Shimura variety, and the theta series $\Theta_{\mc{L}}(z)$ become \emph{generating series of special cycles}. One can define $\mc{Z}_{T,\C}$ by the complex uniformization
    \begin{equation}\label{equation:intro:arith_Siegel_Weil:complex_cycle_uniformization}
    \mc{Z}_{T,\C} \coloneqq \Biggl [ G(\Q) \backslash \Biggl ( \coprod_{\substack{\underline{x} \in V^m \\ (\underline{x}, \underline{x}) = T}} \mc{D}(\underline{x}_{\infty}) \times \mc{D}(\underline{x}_f) \Biggr ) \Biggr ]
    \end{equation}
where $\mc{D}(\underline{x}_{\infty}) \subseteq \mc{D}$ is the closed complex submanifold consisting of those complex lines perpendicular to all elements of the $m$-tuple $\underline{x}$, and
    \begin{equation}
    \mc{D}(\underline{x}_f) \coloneqq \{ g \in G(\A_f) / K_f : g^{-1} x_i \in L \otimes_{\Z} \hat{\Z} \text{ for all $x_i \in \underline{x}_f$} \} \subseteq G(\A_f)/K_f,
    \end{equation}
with $\underline{x}_{\infty}$ and $\underline{x}_f$ denoting the images of $\underline{x}$ in $V(\R)^m$ and $V(\A_f)^m$, respectively.
The definition in \cref{equation:intro:arith_Siegel_Weil:complex_cycle_uniformization} is (a reformulation of) a definition due to Kudla \cite{Kudla04} (there for $\mrm{GSpin}$), with unitary analogue as in \cite[{\S 3}]{Liu11}.  We call $\mc{D}(\underline{x}_{\infty})$ an Archimedean \emph{local special cycle} and $\mc{D}(\underline{x}_{f})$ an ``away-from-$\infty$'' local special cycle. There is a natural map $\mc{Z}_{T,\C} \ra \mrm{Sh}_{K_f,\C}$, which is a disjoint union of closed immersions of complex manifolds after possibly shrinking $K_f$.

A \emph{geometric Siegel--Weil formula} for signature $(n - 1, 1)$ is an identity of the shape
    \begin{equation}\label{equation:intro:arith_Siegel_Weil:geometric_Siegel-Weil}
    -\frac{2 \Lambda_n(0)^{\circ}_n}{\kappa \Lambda_m(s_0)^{\circ}_n} E^*_T(y,s_0)^{\circ}_n = \frac{\mrm{vol}(\mc{Z}_{T,\C})}{[K_{L,f} : K_f]}
    \end{equation}
for $m \leq n - 1$ (so $\kappa = 1$). In the case of signature $(n,0)$, the expression in \cref{equation:intro:arith_Siegel_Weil:geometric_Siegel-Weil} (without the minus sign on the left) is an equivalent reformulation of the classical Siegel--Weil formula \cref{equation:intro:arith_Siegel-Weil:classical_Siegel-Weil}: if we allow the (stacky) level $K_f = K_{L,f}$, there is a canonical equivalence of groupoids
    \begin{equation}
    \mc{Z}_{T,\C} \cong \left \{ \begin{array}{l} \text{pairs $(\mc{L}, \underline{x})$, where $\mc{L}$ is a self-dual Hermitian $\mc{O}_F$-lattice} \\ \text{of signature $(n,0)$ and $\underline{x} \in \mc{L}^m$ is an $m$-tuple with $(\underline{x},\underline{x}) = T$} \end{array} \right \}.
    \end{equation}
Our presentation of the geometric Siegel--Weil formula in \cref{equation:intro:arith_Siegel_Weil:geometric_Siegel-Weil} may be nonstandard. Its appearance is intended to highlight the similarity with our formulation of arithmetic Siegel--Weil in \cref{equation:part_IV:intro:arith_Siegel_Weil:arithmetic_Siegel-Weil} below.

Strictly speaking, geometric Siegel--Weil formulas in the literature typically restrict to $V$ satisfying Weil's convergence condition (meaning $V$ anisotropic or $m < n - 1$ in the signature $(n - 1, 1)$ Hermitian setup), see remarks following \cite[{Theorem 4.1}]{Kudla04} and \cite[{Theorem 3.6.1}]{Li22IHES}. It is also typical to phrase geometric Siegel--Weil formulas in terms of ``coherent'' Eisenstein series, while our $E^*(z,s)^{\circ}_n$ is described in terms of an incoherent ad\`elic Hermitian space (positive definite at $\infty$ and self-dual at all finite places), see \cref{part:part_IV:Eisenstein}. Outside those cases available in the literature, geometric Siegel--Weil formulas may need additional care. For example, when $m = 1$ and $n = 2$ and $T = 0$ (which is essentially about ``complex volume of modular curve''), the formula in \cref{equation:intro:arith_Siegel_Weil:geometric_Siegel-Weil} is only valid up to a non-holomorphic correction term $\frac{2 h_F}{w_F} \cdot \frac{1}{8 \pi y}$ on the left, where $h_F$ (resp. $w_F$) is the class number of (resp. number of roots of unity in) $\mc{O}_F$. In this case, the right-hand side is $\frac{2 h_F}{w_F} \cdot \frac{\zeta(-1)}{2} = \frac{-h_F}{12 w_F}$.
Such non-holomorphic correction terms do not appear for the compact Shimura curves considered in \cite[{(1.0.14)}]{KRY06}. 

We will need the following geometric Siegel--Weil result which does not seem to be covered by the literature discussed in the previous paragraph: we prove \cref{equation:intro:arith_Siegel_Weil:geometric_Siegel-Weil} when $T$ is nonsingular of rank $m = n - 1$, see \cref{proposition:geometric_Siegel-Weil:degrees} (also complex uniformization from \crefext{III:ssec:Arch_uniformization:uniformization}, as well as \cref{equation:arithmetic_Siegel-Weil:main_results:height_constant}); in that case, $\mc{Z}_{T,\C}$ is $0$-dimensional. For example, when $n = 2$ and $\mc{O}_F^{\times} = \{ \pm 1 \}$, the special cycle $\mc{Z}_{T,\C}$ can be described in terms of Hecke translates of CM elliptic curves (\cref{ssec:arithmetic_Siegel-Weil:Serre_tensor}), and \cref{equation:intro:arith_Siegel_Weil:geometric_Siegel-Weil} is then the (well-known) statement that the $T$-th Hecke correspondence (over the modular curve) has bidegree
    \begin{equation}
    - \frac{1}{h_F} \frac{2 \Lambda_2(0)^{\circ}_2}{\Lambda_1(1/2)^{\circ}_2} E^*_T(y,1/2)^{\circ}_2 = \s_1(T)
    \end{equation}
for $T \in \Z_{>0}$.
The extra factor of $h_F$ accounts for multiple connected components in the Shimura variety, see \crefext{ssec:arithmetic_Siegel-Weil:Serre_tensor}.

We remark that our proof of \cref{equation:intro:arith_Siegel_Weil:geometric_Siegel-Weil} (for $T$ nonsingular of rank $m = n - 1$) is inspired by \cite[{Remark 4.6.2}]{LZ22unitary}, and may be carried out using either complex or non-Archimedean (Rapoport--Zink) uniformization. We need that case of \cref{equation:intro:arith_Siegel_Weil:geometric_Siegel-Weil} as an ingredient for our main arithmetic Siegel--Weil results.

            \subsection{Arithmetic Siegel--Weil}
            \label{ssec:part_IV:intro:arith_Siegel--Weil}
                Since the work of Kudla--Rapoport \cite{KR14} (also Rapoport--Smithling--Zhang \cite{RSZ21}), it has been customary to define special cycles $\mc{Z}(T) \ra \mc{M}$ over (stacky) integral models $\mc{M} \ra \Spec \mc{O}_F$ for Shimura varieties associated to $G' \coloneqq \Res_{F / \Q} \G_m \times G$. After adding enough level $K'_f \subseteq G'(\A_f)$ to $\mc{M}_{\C}$, we have a finite covering map $\mc{M}_{K'_f,\C} \ra \mrm{Sh}_{K_f, \C}(G)$. In this paper, we mainly take $\mc{M} \ra \Spec \mc{O}_F$ to be the ``exotic smooth'' Rapoport--Smithling--Zhang (RSZ) integral model of odd relative dimension $n - 1$ \cite[{\S 6}]{RSZ21},
for $G = U(V)$ with $V$ a Hermitian space of signature $(n - 1, 1)$
which contains an everywhere self-dual lattice of full rank.
These conditions imply $n \equiv 2 \pmod{4}$,
but see also remarks after Theorem \ref{theorem:part_IV:intro:results:main}
for the case of $n \equiv 0 \pmod{4}$, and Remark \ref{remark:arithmetic_Siegel-Weil:main_results:other_levels}
for variants involving Hermitian spaces of possibly odd dimension.
When $n = 2$, the stack $\mc{M}$ is essentially a disjoint union of (stacky) modular curves (\cref{ssec:arithmetic_Siegel-Weil:Serre_tensor}, also \crefext{I:example:integral_models:Serre_tensor_global}).

The stack $\mc{M}$ admits a moduli description: it parameterizes tuples $(A_0, \iota_0, \lambda_0, A, \iota, \lambda)$ where $A_0$ and $A$ are abelian schemes (dimensions $1$ and $n$ respectively) with $\mc{O}_F$-actions $\iota_0$ and $\iota$, and with compatible quasi-polarizations $\lambda_0$ and $\lambda$. The datum $(A_0, \iota_0, \lambda_0, A, \iota, \lambda)$ satisfies a few additional conditions, which we suppress in the introduction (see \crefext{I:ssec:part_I:arith_intersections:integral_models} and \crefext{III:ssec:ab_var:integral_models}). We are able to prove versions of our main global results for more general $\mc{M}$ (including odd arithmetic dimension $n$) at the price of discarding finitely many primes (particularly ramified primes for odd $n$); see \cref{remark:arithmetic_Siegel-Weil:main_results:other_levels}.

The moduli stack $\mc{M}$ carries a natural family of Hermitian $\mc{O}_F$-lattices 
    \begin{equation}
    \mc{L}at \ra \mc{M} \quad \quad \mc{L}at \coloneqq \underline{\Hom}_{\mc{O}_F}(A_0, A).
    \end{equation}
Given any $T \in \mrm{Herm}_m(\Q)$, the associated \emph{Kudla--Rapoport} special cycle $\mc{Z}(T) \ra \mc{M}$ is defined as the substack
    \begin{equation}\label{equation:part_IV:intro:arith_Siegel-Weil:KR_cycle}
    \mc{Z}(T) \coloneqq \{\underline{x} \in \mc{L}at^m : (\underline{x}, \underline{x}) = T\} \subseteq \mc{L}at^m
    \end{equation}
consisting of $m$-tuples with Gram matrix $T$. More precisely, see \crefext{III:ssec:ab_var:integral_models}. This is in close analogy with classical Siegel--Weil: there we considered $\mc{O}_F$-lattices varying in a given \emph{ad\`elic isomorphism class},\footnote{
    Non-degenerate Hermitian lattices $L$ and $L'$ lie in the same \emph{ad\`elic isomorphism class}
    (classical terminology: \emph{genus})
    if $L \otimes_{\Z} \Z_p$ and $L' \otimes_{\Z} \Z_p$ are isomorphic Hermitian lattices
    for all primes $p$, and if $L$ and $L'$ moreover have the same signature.
    A genus may contain multiple isomorphism classes.
    The self-dual Hermitian lattices of a given signature form a single genus,
    so this recovers our previous setup with isomorphism classes of positive definite rank $n$ self-dual Hermitian lattices.
} and here we are considering $\mc{O}_F$-lattices varying over the moduli stack $\mc{M}$. In the complex fiber, the special cycles $\mc{Z}(T)_{\C}$ recover the special cycles $\mc{Z}_{T,\C}$ appearing in \cref{equation:intro:arith_Siegel_Weil:complex_cycle_uniformization}, up to $\mc{M}_{K'_f,\C}$ being a finite cover of $\mrm{Sh}_{K_f,\C}(G)$ (for suitable $K'_f$); see \crefext{III:ssec:Arch_uniformization:uniformization}. The morphism $\mc{Z}(T) \ra \Spec \mc{O}_F$ is smooth of relative dimension $n - 1 - \rank(T)$ in the generic fiber over $\Spec F$. If $T$ is not positive semi-definite, then $\mc{Z}(T)$ is empty.

An \emph{arithmetic Siegel--Weil formula} is an (in general conjectural) identity roughly of the shape
    \begin{equation}\label{equation:part_IV:intro:arith_Siegel_Weil:arithmetic_Siegel-Weil}
    \frac{h_F}{w_F} \frac{d}{ds} \bigg|_{s = s_0} \frac{2 \Lambda_n(s-s_0)^{\circ}_n}{\kappa \Lambda_m(s)^{\circ}_n} E^*_T(y,s)^{\circ}_n \overset{?}{=} \widehat{\mrm{vol}}_{\widehat{\mc{E}}^{\vee}}([\widehat{\mc{Z}}(T)])
    \end{equation}
with $\kappa$ as in \cref{ssec:part_IV:intro:classical_and_geometric_Siegel-Weil}. The right-hand side of \cref{equation:part_IV:intro:arith_Siegel_Weil:arithmetic_Siegel-Weil} denotes an \emph{arithmetic volume}, which is a real number ``defined'' by an arithmetic intersection product
    \begin{equation}
    \widehat{\mrm{vol}}_{\widehat{\mc{E}}^{\vee}}([\widehat{\mc{Z}}(T)]) ``\coloneqq" \widehat{\deg}([\widehat{\mc{Z}}(T)] \cdot \widehat{c}_1(\widehat{\mc{E}}^{\vee})^{n - m})
    \end{equation}
in an arithmetic Chow ring $\arithCh^*(\mc{M})_{\Q}$ (roughly in the sense of Gillet--Soul\'e \cite{GS87})
for a certain dual metrized tautological bundle $\widehat{\mc{E}}^{\vee}$ on $\mc{M}$ (the bundle $\widehat{\mc{E}}^{\vee}$ is discussed in 
\crefext{I:ssec:part_I:arith_intersections:tautological_bundle,I:ssec:part_I:arith_intersections:metrized_taut_bundle}). 
The notation $[\widehat{\mc{Z}}(T)]$ indicates a class in $\arithCh^m(\mc{M})_{\Q}$, which is expected to involve $\mc{Z}(T)$ and some additional Archimedean data (e.g. from a Green current on the complex Shimura variety), as appearing in arithmetic intersection theory.
As mentioned in \crefext{I:ssec:intro:arith_Siegel-Weil}, a precise formulation of arithmetic Siegel--Weil has not been proposed in full generality (on both the analytic and geometric sides). Our normalization on the analytic side is already nonstandard; we would like to highlight the similarity with normalizations for geometric and classical Siegel--Weil, as presented in \cref{equation:intro:arith_Siegel-Weil:classical_Siegel-Weil,equation:intro:arith_Siegel_Weil:geometric_Siegel-Weil}. The normalization is in general delicate for arithmetic Siegel--Weil formulas, where both the derivative and special value at $s = s_0$ may have meaning (we do observe this for our main theorem, see \crefext{I:remark:intro:results:main:value_and_derivative_both} and discussion below).

In \crefext{III:sec:arith_cycle_classes} and \crefext{I:ssec:part_I:arith_intersections:vertical_classes}, we proposed a new candidate definition of arithmetic cycle classes
    \begin{equation}
    [\widehat{\mc{Z}}(T)] \coloneqq [\widehat{\mc{Z}}(T)_{\ms{H}}] + \sum_{p \text{ prime}} [{}^{\mathbb{L}}\mc{Z}(T)_{\ms{V},p}] \in \arithCh^m(\mc{M})_{\Q}
    \end{equation}
associated to arbitrary (possibly singular) $T$, where $\arithCh^m(\mc{M})_{\Q}$ is an arithmetic Chow group associated to $\mc{M}$. 
Here, $[\widehat{\mc{Z}}(T)_{\ms{H}}]$ should describe ``horizontal'' contributions and ${}^{\mathbb{L}}\mc{Z}(T)_{\ms{V},p}$ should describe ``vertical'' contributions.

Due to non-properness of $\mc{M} \ra \Spec \mc{O}_F$ in general, one should likely modify $[\widehat{\mc{Z}}(T)]$ on a suitable compactification of $\mc{M}$. 
If $\mc{Z}(T) \ra \Spec \mc{O}_F$ is proper, however, we consider certain ``arithmetic degrees without boundary contributions'' (a real number)
    \begin{align}\label{equation:part_IV:intro:results:if_proper}
        \widehat{\deg}([\widehat{\mc{Z}}(T)] \cdot \widehat{c}_1(\widehat{\mc{E}}^{\vee})^{n-m}) & \coloneqq \left ( \int_{\mc{M}_{\C}} g_{T,y} \wedge c_1(\widehat{\mc{E}}^{\vee}_{\C})^{n-m} \right ) 
        \\
        & \hphantom{\coloneqq} + \widehat{\deg}((\widehat{\mc{E}}^{\vee})^{n - \rank(T)}|_{\mc{Z}(T)_{\ms{H}}}) \notag
        \\
        & \hphantom{\coloneqq} + \sum_{p \text{ prime}} \deg_{\F_p}({}^{\mathbb{L}}\mc{Z}(T)_{\ms{V},p} \cdot (\mc{E}^{\vee})^{n-m}) \log p \notag
    \end{align}
conditional on convergence of the integral
involving a certain current $g_{T,y}$
associated to $\mathcal{Z}(T)$, satisfying a modified Green current equation
\crefext{III:equation:part_III:global_theorem_statement:arithmetic_Siegel-Weil:modified_current_equation}.
Convergence of the integral in the settings of our arithmetic Siegel--Weil results
was explained after \crefext{III:equation:Arch_uniformization:Archimedean}.
Here we set $\mc{M}_{\C} \coloneqq \mc{M} \times_{\Spec \mc{O}_F} \Spec \C$ for any choice of embedding $F \ra \C$. The middle term is mixed characteristic in nature: for $\rank T = n - 1$, it is (essentially) a weighted sum of Faltings heights of abelian varieties (\cref{remark:arithmetic_Siegel-Weil:main_results:Faltings_height}).
For proper $\mc{Z}(T) \ra \mc{O}_F$,
the quantity in \eqref{equation:part_IV:intro:results:if_proper} should coincide with the arithmetic degree (without boundary contributions) of a version of $[\widehat{\mc{Z}}(T)]$ on any reasonable compactification of $\mc{M}$.

The following is the main global theorem for our four-paper sequence.

\begin{theoremLetter}[Co-rank $1$ arithmetic Siegel--Weil]\label{theorem:part_IV:intro:results:main}
Assume the prime $2$ splits in $\mc{O}_F$.
    \begin{enumerate}[(1)]
        \item For any $T \in \mrm{Herm}_n(\Q)$ with $\rank(T) = n - 1$ and any $y \in \mrm{Herm}_n(\R)_{>0}$, we have
            \begin{equation}\label{equation:part_IV:intro:results:singular}
            \frac{h_F}{w_F} \frac{d}{d s} \bigg |_{s = 0} E^*_{T}(y,s)^{\circ}_n = \widehat{\deg}([\widehat{\mc{Z}}(T)]).
            \end{equation}
        \item For any $T^{\flat} \in \mrm{Herm}_{n - 1}(\Q)$ with $\det T^{\flat} \neq 0$ and any $y^{\flat} \in \mrm{Herm}_{n-1}(\R)_{>0}$, we have
            \begin{equation}\label{equation:part_IV:intro:results:nonsingular}
            2 \frac{h_F}{w_F} \frac{d}{d s} \bigg|_{s = 0} \left ( \frac{\Lambda_n(s)_n^{\circ}}{\Lambda_{n - 1}(s + 1/2)^{\circ}_n} E^*_{T^{\flat}}(y^{\flat}, s + 1/2)^{\circ}_n \right ) = \widehat{\deg}([\widehat{\mc{Z}}(T^{\flat})] \cdot \widehat{c}_1(\widehat{\mc{E}}^{\vee})).
            \end{equation}
    \end{enumerate}
\end{theoremLetter}

This appears below as Theorem \ref{theorem:arithmetic_Siegel-Weil:main_results:main}. Note that part (1) concerns the central derivative of a $U(n,n)$ Eisenstein series, while part (2) concerns a non-central derivative of a $U(n - 1, n - 1)$ Eisenstein series. For $n \equiv 0 \pmod {4}$, Theorem \ref{theorem:part_IV:intro:results:main}(1) also holds in the sense that there is no self-dual $\mc{O}_F$-lattice of signature $(n - 1, 1)$ and the left-hand side is $0$ (Remark \ref{remark:arithmetic_Siegel-Weil:main_results:0_mod_4}).

\begin{remark}
Our result in \cref{theorem:part_IV:intro:results:main}(2) combined with our geometric Siegel--Weil result (\cref{equation:intro:arith_Siegel_Weil:geometric_Siegel-Weil} and surrounding discussion, also \cref{proposition:geometric_Siegel-Weil:degrees}) shows that, for the relevant (normalized) $U(n - 1, n - 1)$ Eisenstein series, both the derivative and special value at $s = 1/2$ \emph{simultaneously} have arithmetic-geometric meaning (also mentioned in \crefext{I:remark:intro:results:main:value_and_derivative_both}). This phenomenon amplifies the sensitivity of \cref{theorem:part_IV:intro:results:main}(2) to the choice of normalization for the Eisenstein series, and is one reason for the importance of \cref{part:part_IV:Eisenstein} in this paper (where we discuss our normalizations for Eisenstein series and local Whittaker functions).
\end{remark}

We highlight the simplicity of the analytic side in Theorem \ref{theorem:part_IV:intro:results:main}(1). It is expected that arithmetic Siegel--Weil for integral models with bad reduction should be corrected on the analytic side, e.g. by special values of other Eisenstein series. See for example \cite{HSY22a,HLSY22} for bad reduction in the nonsingular case $\det T \neq 0$ for the central derivative at $s = 0$ (i.e. $T$ is $n \times n$), or \cite{KRY06} for quaternionic Shimura curves. 
We do not know whether the analytic formulation \cite{HSY22a,HLSY22} is expected to hold for singular $T$. 

We argue that arithmetic Siegel--Weil formulas should be simplest to formulate on integral models with everywhere good reduction, as in our case. We thus propose a precise formulation of the analytic side of the central derivative arithmetic Siegel--Weil formula in our setup.

\begin{question*}[Arithmetic Siegel--Weil]
Let $T \in \mrm{Herm}_n(\Q)$ be arbitrary. For a suitable current $g_{T,y}$, a suitable compactification of $\mc{M}$, and a possibly modified class $[\widehat{\mc{Z}}(T)]$ on the compactification, do we have
    \begin{equation}\label{equation:part_IV:intro:results:conjecture}
    \frac{h_F}{w_F} \frac{d}{d s} \bigg|_{s = 0} E^*_T(y,s)^{\circ}_n \overset{?}{=} \widehat{\deg}([\widehat{\mc{Z}}(T)]).
    \end{equation}
\end{question*}

Our theorem verifies this proposed arithmetic Siegel--Weil formula for all singular $T \in \mrm{Herm}_n(\Q)$ of rank $n - 1$, in the sense of ``arithmetic degrees without boundary contributions''. The formula also holds (in the same sense) for all nonsingular $T \in \mrm{Herm}_n(\Q)$. This latter case (``central derivative nonsingular arithmetic Siegel--Weil'') is possibly considered known to experts up to a volume constant by collecting the local theorems in \cite{Liu11,LZ22unitary,LL22II}. This particular global statement does not appear in the literature, though other variants are available (e.g. for unramified CM fields $F / F_0$ with $F_0 \neq \Q$ \cite{LZ22unitary} or on integral models with bad reduction and correction terms by special values of other Eisenstein series \cite{HLSY22}). In our setup, we will compute the volume constant and explain how to extract the $\det T \neq 0$ case of \eqref{equation:part_IV:intro:results:conjecture} from the literature (Remark \ref{remark:arithmetic_Siegel-Weil:main_results:nonsingular}).

For a discussion on difficulties in approaching
arithmetic Siegel--Weil for higher-dimensional cycles
(arbitrary $T$)
from our point of view,
we refer the reader to \crefext{I:remark:intro:strategy_overview:higher_corank_difficulties}
or the end of \crefext{III:ssec:part_III:intro:local-to-global}.

\begin{remark}
Part (2) of Theorem \ref{theorem:part_IV:intro:results:main} is the special case of part (1) when $T = \mrm{diag}(0,T^{\flat})$ and $y = \mrm{diag}(1, y^{\flat})$. The geometric sides agree essentially by definition \eqref{equation:part_IV:intro:results:if_proper}. On the analytic side, the relation is provided by the formula
    \begin{align}\label{equation:part_IV:intro:arith_Siegel--Weil:Eisenstein_corank1_Fourier_unfold}
    E_T^*(y, s)^{\circ}_n & = \frac{\Lambda_n(s)^{\circ}_n}{\Lambda_{n - 1}(s + 1/2)^{\circ}_n} E^*_{T^{\flat}}(y^{\flat}, s + 1/2)^{\circ}_n -\frac{\Lambda_n(- s)^{\circ}_n}{\Lambda_{n - 1}(- s + 1/2)^{\circ}_n} E^*_{T^{\flat}}(y^{\flat}, s - 1/2)^{\circ}_n
    \end{align}
from Corollary \ref{corollary:Eisenstein:singular_Fourier:corank_1}, along with the functional equation $E^*_{T^{\flat}}(y^{\flat},s)^{\circ}_n = E^*_{T^{\flat}}(y^{\flat},-s)^{\circ}_n$. The general case of Theorem \ref{theorem:part_IV:intro:results:main} is proved in a similar way as the special case $T = \mrm{diag}(0,T^{\flat})$, with an additional ``local diagonalizability argument'' (Proof of Theorem \ref{theorem:arithmetic_Siegel-Weil:main_results:main}) where the identity is proved modulo $\sum_{\ell \neq p} \Q \cdot \log \ell$ for any given $p$ (varying $p$ removes the ambiguity).

Unlike $E^*_{T^{\flat}}(y,s)^{\circ}_n$, the Fourier coefficient $E^*_T(y,s)^{\circ}_n$ does not admit an obvious Euler product decomposition into local Whittaker functions as $T$ is singular. Since our proof of Theorem \ref{theorem:part_IV:intro:results:main} is local (via local Whittaker functions and local special cycles), the decomposition \cref{equation:part_IV:intro:arith_Siegel--Weil:Eisenstein_corank1_Fourier_unfold} is crucial for our method (to use local Whittaker functions to describe $E^*_T(y,s)^{\circ}_n$). This decomposition is sensitive to the normalization used to define $E^*(z,s)^{\circ}_n$, and is another reason for the importance of \cref{part:part_IV:Eisenstein} in this paper, where we treat normalized Eisenstein series and local Whittaker functions.
\end{remark}

\begin{remark}
An expected application of arithmetic Siegel--Weil formulas is in the theory of \emph{arithmetic theta lifting}. One expects to form automorphic \emph{arithmetic theta series} as generating series
    \begin{equation}\label{equation:intro:arith_Siegel-Weil:arith_theta}
    \widehat{\Theta} = \sum_T [\widehat{\mc{Z}}(T)] q^T
    \end{equation}
with ``Fourier coefficients'' $[\widehat{\mc{Z}}(T)]$ valued in the arithmetic Chow group $\arithCh^m(\mc{M})_{\Q}$. These should be analogous to (weighted averages of) classical theta series, as in the classical Siegel--Weil formula. In analogy with classical theta lifting, one expects to use $\widehat{\Theta}$ as an integral kernel to lift $U(m,m)$ automorphic forms to elements of $\arithCh^m(\mc{M})_{\Q}$. In analogy with the classical Rallis inner product formula, one expects to use the doubling method and arithmetic Siegel--Weil formulas to relate the derivative of an $L$-function with the arithmetic inner product of this arithmetic theta lift \cite[{Part III}]{Kudla04}. We refer to \cite{KRY06,BHKRY20II,LL22I,LL22II} for some cases where versions of this have been realized, with applications to Beilinson--Bloch. For modularity results on generating series of arithmetic divisors, see \cite{KRY06,BBK07,BHKRY20,Qiu22}.
\end{remark}

It is also possible to formulate and prove Theorem \ref{theorem:part_IV:intro:results:main} in terms of Faltings heights (i.e. replacing the middle term in \eqref{equation:part_IV:intro:results:if_proper} with the degree of the metrized Hodge bundle). The formulation in Theorem \ref{theorem:part_IV:intro:results:main} seems more natural to us, but the version with Faltings heights is in Remark \ref{remark:arithmetic_Siegel-Weil:main_results:Faltings_height}.

Since our proof of Theorem \ref{theorem:part_IV:intro:results:main} will be local in nature, we also have a version for more general moduli stacks $\mc{M}$ (including odd arithmetic dimension $n$) at the price of discarding finitely many primes (particularly the ramified ones). This is explained in \cref{remark:arithmetic_Siegel-Weil:main_results:other_levels}.

The simplest case of Theorem \ref{theorem:part_IV:intro:results:main} is the case $n = 2$. When $\mc{O}_F^{\times} = \{ \pm 1 \}$, the Serre tensor construction gives an open and closed embedding $\ms{M}_0 \times_{\Spec \mc{O}_F} \ms{M}_{\text{ell}} \ra \mc{M}$, where $\ms{M}_0$ is the moduli stack of elliptic curves with signature $(1,0)$ action by $\mc{O}_F$ and $\ms{M}_{\text{ell}}$ is the moduli stack of all elliptic curves, base-changed to $\mc{O}_F$ (Section \ref{ssec:arithmetic_Siegel-Weil:Serre_tensor}). In this case, the special cycle $\mc{Z}(j) \ra \mc{M}$ for $j \in \Z_{>0}$ pulls back to the $j$-th Hecke correspondence. Then the proof of Theorem \ref{theorem:arithmetic_Siegel-Weil:main_results:main}(2) gives the following corollary (appearing below as Corollary \ref{corollary:arithmetic_Siegel-Weil:Serre_tensor}). One might think of this corollary as reformulating a result of Nakkajima--Taguchi \cite{NT91} (they compute Faltings heights of elliptic curves with CM by possibly non-maximal orders) by averaging over Hecke translates and expressing the result in terms of Eisenstein series Fourier coefficients. 

\begin{corollary}\label{corollary:part_IV:intro:results:n=2}
Assume $2$ is split in $\mc{O}_F$. Fix any elliptic curve $E_0$ over $\C$ with $\mc{O}_F$-action. For any integer $j > 0$, we have
    \begin{equation}
    \sum_{E \xra{w} E_0} (h_{\mrm{Fal}}(E) - h_{\mrm{Fal}}(E_0)) = \frac{1}{2} \frac{d}{ds} \bigg|_{s = 1/2} \left ( j^{s + 1/2} \s_{-2s}(j) \right )
    \end{equation}
where the sum runs over degree $j$ isogenies $w \colon E \ra E_0$ of elliptic curves.
\end{corollary}

The notation $h_{\mrm{Fal}}(E)$ denotes the (stable) Faltings height of the elliptic curve $E$ after descent to any number field, and similarly for $E_0$. The quantity $j^{s+1/2} \s_{-2s}(j)$ is the product of the normalized non-Archimedean local Whittaker functions in the $j$-th Fourier coefficient $E^*_j(y,s)^{\circ}_2$ (with $m = 1$) as in \cref{equation:intro:arith_Siegel-Weil:weight_2_Eisenstein}. The derivative of the Archimedean local Whittaker function $W^*_{j,\infty}(y,s)^{\circ}_2$ at $s = 1/2$ was calculated explicitly and compared with its geometric counterpart (integral of Green function wedge Chern form on upper half-plane) in our companion paper \crefext{II:ssec:Archimedean_identity:case_n_is_2}.

In fact, Corollary \ref{corollary:part_IV:intro:results:n=2} holds if $E_0$ is replaced by any elliptic curve over $\overline{\Q}$.
Indeed, it is known that the quantity on the left-hand side of Corollary \ref{corollary:part_IV:intro:results:n=2}
is constant as $E_0$ varies over all elliptic curves over $\overline{\Q}$.
The case of isogenies of square-free degree appears in \cite[{Proposition 8.1}]{SU99},
and the case of cyclic isogenies appears in \cite[{Corollaire 3.3}]{Autissier03}.
Our Corollary \ref{corollary:part_IV:intro:results:n=2} may be extracted from either of those cited
results via a combinatorial calculation involving the usual Hecke algebra recurrence relations.
The actual proof that we give (as a byproduct of the proof of our main theorem) is different,
passing through local calculations on $p$-divisible groups and formal groups (e.g. as appearing in the
calculation of Nakkajima--Taguchi).

Our purely Archimedean result (for arbitrary $n$ and arbitrary rank $m^{\flat} \geq 1$) is the following.

\begin{theoremLetter}[Archimedean arithmetic Siegel--Weil, nonsingular]\label{theorem:part_IV:intro:results:Archimedean}
Consider any integer $m^{\flat}$ with $1 \leq m^{\flat} \leq n$, and consider any $T^{\flat} \in \mrm{Herm}_{m^{\flat}}(\Q)$ which is nonsingular and not positive definite.
    \begin{enumerate}[(1)]
        \addtocounter{enumi}{2}
        \item For any $y^{\flat} \in \mrm{Herm}_{m^{\flat}}(\R)_{>0}$, we have an equality of real numbers
            \begin{equation}\label{equation:part_IV:intro:results:Archimedean}
            \widehat{\deg}([\widehat{\mc{Z}}(T^{\flat})] \cdot \widehat{c}_1(\widehat{\mc{E}}^{\vee})^{n - m^{\flat}}) \coloneqq \int_{\mc{M}_{\C}} g_{T^{\flat},y^{\flat}} \wedge c_1(\widehat{\mc{E}}^{\vee}_{\C})^{n - m^{\flat}} = (-1)^{n - m^{\flat}} C \cdot \frac{h_F}{w_F} \frac{d}{ds} \bigg|_{s = s_0^{\flat}} E^*_{T^{\flat}}(y^{\flat}, s)^{\circ}_n
            \end{equation}
        where $s_0^{\flat} \coloneqq (n - m^{\flat}) / 2$. Here $C \in \Q_{>0}$ is the volume constant from \crefext{lemma:local_Siegel-Weil:uniformization_degree:main:1}, for the Hermitian space $V$ and $v_0 = \infty$ in the notation of loc. cit.. The constant $C$ may depend on $n$ and $m^{\flat}$ (and $F$), but does not otherwise depend on $T^{\flat}$.
    \end{enumerate}
\end{theoremLetter}

In Theorem \ref{theorem:part_IV:intro:results:Archimedean},
the assumption that $T^{\flat}$ is not positive definite ensures that the cycle $\mc{Z}(T^{\flat})$
is empty, so that the only contribution to the arithmetic degree $\widehat{\deg}([\widehat{\mc{Z}}(T^{\flat})] \cdot \widehat{c}_1(\widehat{\mc{E}}^{\vee})^{n - m^{\flat}})$
should be from the displayed Green current integral.
Theorem \ref{theorem:part_IV:intro:results:Archimedean} appears below (in stronger form)
as Theorem \ref{theorem:arithmetic_Siegel-Weil:main_results:Archimedean}. That version applies for all $n$ (even or not) and arbitrary level, as it is a statement about the complex Shimura variety. We gave the weaker version here to avoid more notation in the introduction. Due to non-properness of $\mc{M}_{\C} \ra \Spec \C$ for $n > 2$, the corresponding Archimedean result of \cite{GS19} does not apply here if $n > 2$.

When $m^{\flat} = n$, the preceding Archimedean theorem follows from Liu's result \cite[{Theorem 4.17}]{Liu11}. We do not have a new proof of this case. Instead, we deduce our general result from his by a certain limiting argument. This is also our method at non-Archimedean places (replacing Liu's Archimedean results with the non-Archimedean results of Li--Zhang \cite{LZ22unitary} and Li--Liu \cite{LL22II}). These local theorems were the main results of our companion papers \cite{corank1_ASW_I.pdf,corank1_ASW_II.pdf}. In this paper, we explain how to combine our local theorems at all places to prove the (global) Theorems \ref{theorem:part_IV:intro:results:main} and \ref{theorem:part_IV:intro:results:Archimedean}.

\begin{remark}
The ``global-to-local'' reduction step for heights 
from \cite{corank1_ASW_III.pdf} is the main reason why we
only prove an arithmetic Siegel--Weil result for imaginary quadratic $F / \Q$,
rather than a general CM field $F / F_0$.
Our local methods from \cite{corank1_ASW_I.pdf} carry through formally for 
arbitrary \'etale degree $2$ extensions $F / F_0$ for $F_0$ a $p$-adic field,
disallowing only the case $p = 2$ with $F / F_0$ nonsplit.
However, it remains to show that the analogous ``local geometric quantities''
in fact have global geometric meaning for general $F / F_0$.

In the general case, an ``absolute'' Rapoport--Zink space appears naturally
in $p$-adic uniformization of Shimura varieties, while our local methods seem more directly applicable
on certain ``relative'' Rapoport--Zink spaces.
There are known comparison isomorphisms between these Rapoport--Zink spaces,
as in \cite[{Theorem 1.4}]{Mihatsch22} (inert places)
and \cite[{\S 2.8}]{LL22II} (ramified places, unramified over $\Q_p$).
It remains to be understood how the height decomposition procedure
of this paper behaves with respect to these comparison isomorphisms.
To complete the global result,
presumably one needs to also develop a similar theory at split primes,
along with the related further study of ordinary Rapoport--Zink spaces,
generalizing the treatment in
\cite[{Sections \ref{I:sec:moduli_pDiv} to \ref{I:sec:can_and_qcan}}]{corank1_ASW_I.pdf}
and the Rapoport--Zink uniformization in ordinary situations as treated in
\crefext{III:sec:non-Arch_uniformization}.
\end{remark}

            \subsection{Sketch}
            \label{ssec:part_IV:intro:sketch}
                We next outline what remains for the proofs of Theorems \ref{theorem:part_IV:intro:results:main} and \ref{theorem:part_IV:intro:results:Archimedean},
after the previous three papers \cite{corank1_ASW_I.pdf,corank1_ASW_II.pdf,corank1_ASW_III.pdf} in this sequence.

We illustrate (a special case of)
the local-to-global reduction process on the analytic side of the arithmetic Siegel--Weil formula,
involving Eisenstein series and local Whittaker functions.
In the special case of $T = \mrm{diag}(0,T^{\flat})$ for $T^{\flat}$ nonsingular of rank $n - 1$
and $y = \mrm{diag}(1,y^{\flat})$, we give a local decomposition
    \begin{align}
    \frac{1}{2} \frac{d}{ds}\bigg|_{s = 0} E^*_T(y,s)^{\circ}_n & = \frac{d}{d s} \bigg|_{s = 0} \left ( \frac{\Lambda_n(s)_n^{\circ}}{\Lambda_{n - 1}(s + 1/2)^{\circ}_n} E^*_{T^{\flat}}(y^{\flat}, s + 1/2)^{\circ}_n \right ) \label{equation:part_II:intro:local-to-global:pre-unfold}
    \\
    & = \left ( \frac{d}{d s} \bigg|_{s = 0} \frac{\Lambda_n(s)_n^{\circ}}{\Lambda_{n - 1}(s + 1/2)^{\circ}_n} \right ) E^*_{T^{\flat}}(y^{\flat}, 1/2)^{\circ}_n \label{equation:part_II:intro:local-to-global:height_constant_term}
    \\
    & \mathrel{\hphantom{=}} - \frac{h_F}{w_F} \left ( \frac{d}{ds} \bigg|_{s = 1/2} W^*_{T^{\flat}, \infty}(y^{\flat},1/2)^{\circ}_n \right ) \prod_{p} W^*_{T^{\flat},p}(1/2)^{\circ}_n \label{equation:part_II:intro:local-to-global:Archimedean}
    \\
    & \mathrel{\hphantom{=}} - \frac{h_F}{w_F} \sum_p \left( \frac{d}{ds} \bigg|_{s = 1/2} W^*_{T,p}(s)^{\circ}_n \right) W^*_{T^{\flat},\infty}(y,1/2)^{\circ}_n \prod_{\ell \neq p} W^*_{T^{\flat},\ell}(1/2)^{\circ}_n
    \label{equation:part_II:intro:local-to-global:non-Arch}
    \end{align}
coming from Leibniz rule and an Euler product
    \begin{equation}
    E^*_{T^{\flat}}(y^{\flat},s)^{\circ}_n = W^*_{T^{\flat}, \infty}(y^{\flat},s)^{\circ}_n \prod_{p} W^*_{T^{\flat}, p}(s)^{\circ}_n
    \end{equation}
of local (normalized) Whittaker functions over all places (\cref{ssec:Eisenstein:global_normalized_Fourier:global_normalization}).
More generally, one can relate the $U(n,n)$ Eisenstein series Fourier coefficient $E^*_{T}(y,s)^{\circ}_n$
with the $U(n - 1, n - 1)$ Eisenstein series Fourier coefficients
$E^*_{T^{\flat}}(y^{\flat}, s + 1/2 )^{\circ}_n$ and $E^*_{T^{\flat}}(y^{\flat}, s - 1/2)^{\circ}_n$
(\cref{corollary:Eisenstein:singular_Fourier:corank_1}).
Here, the work is to explicitly pin down normalizing factors, constants, intertwining operators, etc.,
so that the final form of \cref{theorem:part_IV:intro:results:main} matches our explicit prediction.
These tasks are carried out in \cref{part:part_IV:Eisenstein}.

If $T$ is not in block diagonal form but ${}^t \overline{\gamma} T \gamma$
is block diagonal for some $\gamma \in \GL_n(\mc{O}_F)$,
a certain ``linear invariance'' property for Eisenstein series (\cref{ssec:Eisenstein:setup:Fourier_and_Whittaker})
gives a similar local decomposition.
If no such $\gamma$ exists, we instead take $\gamma \in \GL_n(\mc{O}_F \otimes_{\Z} \Z_{(p)})$
(``local diagonalizability''); this introduces a discrepancy of $\sum_{\ell \neq p} \Q \cdot \log \ell$,
but varying $p$ removes this discrepancy by $\R$-linear independence of $\log p$ over all primes $p$
(see \cref{ssec:arithmetic_Siegel-Weil:main_results}).
This ``$p$-local diagonalization'' argument is a new feature on the analytic side
for our local-to-global reduction of arithmetic Siegel--Weil.
This possible non-diagonalizability of $T$ is the main reason why the proof of
\cref{theorem:arithmetic_Siegel-Weil:main_results:main} is more than twice as long as this sketch.

To reduce \cref{theorem:part_IV:intro:results:main} (for $T$ singular of co-rank $1$) to our local main theorems from \crefext{I:sec:non-Arch_identity} and \crefext{II:sec:Archimedean_identity},
we recall the geometric quantities
    \begin{align}
        \mrm{Int}_{\infty, \mrm{global}}(T, y) & \coloneqq \int_{\mc{M}^{\mrm{an}}_{\C}} g_{T,y} \wedge c_1(\widehat{\mc{E}}_{\C}^{\vee})^{n - m}
        \\
        \sum_{p} \mrm{Int}_{\ms{H},p,\mrm{global}}(T) & \coloneqq \widehat{\deg}((\widehat{\mc{E}}^{\vee})^{n - \rank(T)}|_{\mc{Z}(T)_{\ms{H}}}) - (\deg_{\Z} \mc{Z}(T)_{\ms{H}}) \cdot h_{\widehat{\mc{E}}^{\vee}}^{\mrm{CM}}
        \\
        \mrm{Int}_{\ms{V},p,\mrm{global}}(T) & \coloneqq \deg_{\F_p} ({}^{\mbb{L}} \mc{Z}(T)_{\ms{V},p} \cdot (\mc{E}^{\vee})^{n - m}) \log p
    \end{align}
from \cite[{\cref{III:equation:Arch_uniformization:Archimedean}, \cref{III:equation:non-Arch_uniformization:horizontal:local_global_degrees}, \cref{III:equation:non-Arch_uniformization:vertical:local_global_degrees}}]{corank1_ASW_III.pdf},
as well as the ``total non-Archimedean intersection number''
$\opn{Int}_{p, \mrm{global}}(T) \coloneqq \opn{Int}_{\ms{H}, p, \mrm{global}}(T) + \opn{Int}_{\ms{V}, p, \mrm{global}}(T)$.
We recall that $h^{\mrm{CM}}_{\widehat{\mc{E}}^{\vee}}$ is a certain height constant
as in \crefext{III:equation:part_II:arith_cycle_classes:Hodge_bundles:taut_height_constants},
coming from Faltings height of elliptic curves with CM by $\mc{O}_F$.

We then proceed by matching
\cref{equation:part_II:intro:local-to-global:height_constant_term}, \cref{equation:part_II:intro:local-to-global:Archimedean}, and \cref{equation:part_II:intro:local-to-global:non-Arch}
to the geometric quantities
\begin{align}
    (\deg_{\Z} \mc{Z}(T)_{\ms{H}}) \cdot h^{\mrm{CM}}_{\widehat{\mc{E}}^{\vee}}
    \quad \quad
    \mrm{Int}_{\infty, \mrm{global}}(T, y)
    \quad \quad
    \sum_{p} \opn{Int}_{p, \mrm{global}}(T)
\end{align}
respectively, and the sums over primes $p$ match term by term.

The height constant $h^{\mrm{CM}}_{\widehat{\mc{E}}^{\vee}}$ is essentially the derivative appearing in \cref{equation:part_II:intro:local-to-global:height_constant_term} (see \cref{equation:arithmetic_Siegel-Weil:main_results:height_constant}), and the degree $\deg_{\Z} \mc{Z}(T)_{\ms{H}}$ is essentially the special value $E^*_{T^{\flat}}(y^{\flat},s)^{\circ}_n$ appearing in loc. cit. (by geometric Siegel--Weil as in \cref{remark:geometric_Siegel-Weil:degrees:self-dual} or \cref{equation:intro:arith_Siegel_Weil:geometric_Siegel-Weil}).
The comparison with \crefext{equation:part_II:intro:local-to-global:height_constant_term} emerges from this.

For the remaining terms, we recall the local-to-global relations
    \begin{align}
        \mrm{Int}_{\infty, \mrm{global}}(T, y) & = \mrm{Int}_{\infty}(T,y) \frac{h_F}{w_F} \cdot \deg \Biggl [ U(V)(\Q) \backslash \coprod_{\substack{\underline{x} \in V^m \\ (\underline{x}, \underline{x}) = T}} \mc{D}(\underline{x}_f) \Biggr ]
        \\
        \mrm{Int}_{p,\mrm{global}}(T) & = \mrm{Int}_{p}(T) \frac{h_F}{w_F} \cdot \deg \Biggl [ I_1(\Q) \backslash \Biggl ( \coprod_{\substack{\underline{\mbf{x}} \in \mbf{W}^m \\ (\underline{\mbf{x}}, \underline{\mbf{x}}) = T}} U(\mbf{W}^{\perp}_p) / K_{1,\mbf{L}^{\perp}_{p}} \times \mc{Z}(\underline{\mbf{x}}^p) \Biggr ) \Biggr ]
    \end{align}
from our companion paper,
as given in \crefext{III:equation:Arch_uniformization:Archimedean:local_global_degrees}
and (the combination of) \cite[{\eqref{III:equation:non-Arch_uniformization:horizontal:local_global_degrees},
\eqref{III:equation:non-Arch_uniformization:vertical:local_global_degrees}, \eqref{III:equation:non-Arch_uniformization:horizontal:total_intersection}}]{corank1_ASW_III.pdf},
with $[K_{L_0,f} : K_{0,f}] = 1$ in the notation of loc. cit.
Here $\opn{Int}_{\infty}(T,y)$ and $\opn{Int}_p(T)$ are ``local intersection numbers''
(Archimedean and non-Archimedean),
formulated on the Hermitian symmetric domain and Rapoport--Zink spaces respectively.
Each $\deg[ \cdots ]$ comes from some ``uniformization degrees'' (for complex or Rapoport--Zink uniformization)
for special cycles, and are essentially theta integrals.
For the reader's convenience, we recall that $\mc{D}(\underline{x}_f) \subseteq U(V)(\A_f) / K_{L,f}$
and $\mc{Z}(\underline{\mbf{x}}^p) \subseteq U(V)(\A_f^p) / K_{L,f}^p$ are certain (discrete) sets.
The symbol $\mbf{W}$ is a certain positive definite Hermitian space depending on $p$
(rank $n$ if $p$ is nonsplit, and rank $n - 1$ if $p$ is split),
with a ``perpendicular complement'' $\mbf{W}^{\perp}$
(rank $0$ if $p$ is nonsplit, and rank $1$ if $p$ is split)
with $I_1 \coloneqq U(\mbf{W}) \times U(\mbf{W}^{\perp})$.
The notation $K_{1,\mbf{L}^{\perp}_p} \subseteq U(\mbf{W}_p^{\perp})$
means the unique maximal compact subgroup.
For further discussion of these objects, we refer to the introduction \crefext{III:sec:part_III:intro}.

Our main Archimedean local theorem \crefext{II:theorem:local_identities:Archimedean_identity:statement:main_Archimedean} gives
    \begin{equation}
    \mrm{Int}_{\infty}(T,y) = - \frac{d}{ds} \bigg|_{s = 1/2} W^*_{T^{\flat},\infty}(y,s)^{\circ}_n
    \end{equation}
By ``local Siegel--Weil'' (\cref{ssec:local_Siegel-Weil:uniformization_degree}), we have
    \begin{equation}
    \deg \Biggl [ U(V)(\Q) \backslash \coprod_{\substack{\underline{x} \in V^m \\ (\underline{x}, \underline{x}) = T}} \mc{D}(\underline{x}_f) \Biggr ] = \frac{2 h_F}{w_F} \prod_p W^*_{T^{\flat}, p}(1/2)^{\circ}_n
    \end{equation}
for any $y^{\flat} \in \mrm{Herm}_{n - 1}(\R)_{>0}$, and the comparison with \crefext{equation:part_II:intro:local-to-global:Archimedean} now emerges.

Our local main non-Archimedean theorem from \crefext{I:sec:non-Arch_identity} gives
    \begin{equation}
    \mrm{Int}_{p}(T) = - e_p \frac{d}{ds} \bigg|_{s = 1/2} W^*_{T^{\flat},p}(s)^{\circ}_n
    \end{equation}
where $e_p = 1$ (resp. $e_p = 2$) if $p$ is unramified (resp. if $p$ is ramified).
By ``local Siegel--Weil'' (\cref{ssec:local_Siegel-Weil:uniformization_degree}), we have
    \begin{equation}
        \deg \Biggl [ I_1(\Q) \backslash \Biggl ( \coprod_{\substack{\underline{x} \in \mbf{W}^m \\ (\underline{x}, \underline{x}) = T}} U(\mbf{W}^{\perp}_p) / K_{1,\mbf{L}^{\perp}_{p}} \times \mc{Z}(\underline{\mbf{x}}^p) \Biggr ) \Biggr ] = \frac{2 h_F}{e_p w_F} W^*_{T^{\flat},\infty}(y^{\flat},1/2)^{\circ}_n \prod_{\ell \neq p} W^*_{T^{\flat}, \ell}(1/2)^{\circ}_n
    \end{equation}
and the comparison with \crefext{equation:part_II:intro:local-to-global:non-Arch} now emerges.

            \subsection{Outline}
            \label{ssec:part_IV:intro:outline}
                We briefly summarize the remaining content in this paper, and discuss the relation with our companion papers \cite{corank1_ASW_I.pdf,corank1_ASW_II.pdf,corank1_ASW_III.pdf}.
Further explanations may be found at the beginning of some sections.

The remaining sections are divided into \cref{part:part_IV:Eisenstein,part:part_IV:Siegel--Weil}.

In \cref{part:part_IV:Eisenstein} ``Eisenstein series'', we study $U(m,m)$ Siegel--Weil Eisenstein series. To formulate and prove our main results, it is extremely important that we normalize the Eisenstein series and local Whittaker functions (e.g. by certain $L$-factors). As in our main local theorems from \cite{corank1_ASW_I.pdf,corank1_ASW_II.pdf}, it seems that these normalized versions correspond more naturally to geometric objects (e.g. global and local special cycles). We pin down explicit precise normalizations, guided by special value formulas and symmetric functional equations. We also study (normalized) Fourier coefficients for singular $T$ (focusing on rank $m - 1$ and size $m \times m$), and give formulas needed for our main results.

In \cref{part:part_IV:Siegel--Weil} ``Siegel--Weil'', we first give some special value formulas (local and geometric Siegel--Weil, Sections \ref{sec:local_Siegel-Weil} and \ref{sec:geometric_Siegel-Weil}) which are needed to prove our arithmetic Siegel--Weil theorems. Our treatment of our geometric Siegel--Weil result is inspired by \cite[{Remark 4.6.2}]{LZ22unitary}, and the argument may be carried out using complex (Archimedean) uniformization or Rapoport--Zink (non-Archimedean) uniformization in similar fashions.
We also compare complex volumes of unitary Shimura varieties with the special value of our Eisenstein series normalizing factor (``geometric Siegel mass formula''), which may be of independent interest (\cref{ssec:geometric_Siegel-Weil:complex_volumes}).
The finale
occurs in Section \ref{ssec:arithmetic_Siegel-Weil:main_results}, where we collect our local main theorems to prove our (global) arithmetic Siegel--Weil theorems. This proof relies on results from almost all preceding sections, including our three other companion papers. The key inputs from those papers are our non-Archimedean local arithmetic Siegel--Weil theorems from \cite{corank1_ASW_I.pdf}, our Archimedean local arithmetic Siegel--Weil theorem from \cite{corank1_ASW_II.pdf}, and the geometric local-to-global reduction process (via Archimedean and non-Archimedean uniformization) from \cite{corank1_ASW_III.pdf}.
Section \ref{ssec:arithmetic_Siegel-Weil:Serre_tensor} contains a reformulation of our arithmetic Siegel--Weil results in the special case $n = 2$, via an exceptional comparison with Hecke translates of CM elliptic curves.
    
            \subsection{Acknowledgements}
            \label{part_IV:acknowledgements}
                I thank my advisor Wei Zhang for suggesting this topic, for his dedicated support and constant enthusiasm, for insightful discussions throughout the entire course of this project, and for helpful comments on earlier drafts. I thank Tony Feng, Qiao He, Benjamin Howard, Ishan Levy, Chao Li, Keerthi Madapusi, Andreas Mihatsch, Bjorn Poonen, Siddarth Sankaran, Ananth Shankar, Yousheng Shi, Tonghai Yang, Shou-Wu Zhang, and Zhiyu Zhang for helpful comments or discussions.
I also thank the anonymous referees for their helpful comments, line-by-line close reading, and many corrections.

This work was partly supported by the National Science Foundation Graduate Research Fellowship under Grant Nos. DGE-1745302 and DGE-2141064.
Parts of this work were completed at the Mathematical Sciences Research Institute (MSRI), now becoming the Simons Laufer Mathematical Sciences Institute (SLMath), and the Hausdorff Institute for Mathematics. I thank these institutes for their support and hospitality. The former is supported by the National Science Foundation (Grant No. DMS-1928930), and the latter is funded by the Deutsche Forschungsgemeinschaft (DFG, German Research Foundation) under Germany's Excellence Strategy -- EXC-2047/1 -- 390685813.
This manuscript was prepared partially during the period the author served as a Clay Research Fellow.

    \clearpage


    \part{Eisenstein series}
    \label{part:part_IV:Eisenstein}

        \section{Setup}
        \label{sec:setup}
            
            \subsection{The group \texorpdfstring{$U(m,m)$}{U(m,m)}} 
            \label{ssec:Eisenstein:setup:group}
                We fix notation for the unitary group $U(m,m)$.

Let $A \ra B$ be a finite locally free morphism of (commutative) rings, and suppose $B$ is given an involution $b \mapsto \overline{b}$ (``conjugation'') over $A$. We are mostly interested in the case where $F/F^+$ is a CM extension of number fields (with $F^+$ the index $2$ totally real subfield) and $B/A = \mc{O}_F/\mc{O}_{F^+}$ for the corresponding rings of integers (also the local analogues) etc..

Fix an integer $m \geq 0$. Write $1_m$ for the $m \times m$ identity matrix (sometimes we drop the subscript $m$), and let $H = U(m,m)$ be the unitary group
    \begin{equation}
    H = U(m,m) \coloneqq \left \{ h \in \Res_{B/A} \GL_{2m} : {}^t \overline{h} \begin{pmatrix} 0 & 1_m \\ -1_m & 0 \end{pmatrix} h = \begin{pmatrix} 0 & 1_m \\ -1_m & 0 \end{pmatrix} \right \}
    \end{equation}
where ${}^t \overline{h}$ denotes conjugate transpose (with $H$ the trivial group if $m = 0$, by convention). Equivalently, $H$ consists of block matrices
    \begin{equation}
    \begin{pmatrix}
    a & b \\
    c & d
    \end{pmatrix} \quad \text{satisfying} \quad {}^t \overline{a} c = {}^t \overline{c} a \quad \quad {}^t \overline{a} d - {}^t \overline{c} b = 1_m \quad \quad {}^t \overline{b} d = {}^t \overline{d} b
    \end{equation}
with $a, b, c, d \in \mrm{Res}_{B/A} M_{m \times m}$.
We refer to $H$ as the group $U(m,m)$ (for signature reasons when $B/A$ is $\C / \R$).

Given an integer $j$ with $0 \leq j \leq m$, we consider the injection
    \begin{equation}
    \mu_j^m \colon U(j,j) \ra U(m,m)
    \quad \quad
    \begin{pmatrix}
    a & b \\ c & d
    \end{pmatrix}
    \mapsto 
    \begin{pmatrix}
    1_{m-j} & 0 & 0 & 0 \\
    0 & a & 0 & b \\
    0 & 0 & 1_{m-j} & 0 \\
    0 & c & 0 & d
    \end{pmatrix}.
    \end{equation}
Consider the subgroups
    \begin{align}
    P & \coloneqq  \left \{ h = \begin{pmatrix} * & * \\ 0 & * \end{pmatrix} \in H \right \} \\
    M & \coloneqq  \left \{ m(a) = \begin{pmatrix} a & 0 \\ 0 & {}^t \overline{a}^{-1} \end{pmatrix} : a \in \Res_{B/A} \GL_m \right \} \\
    N & \coloneqq  \left \{ n(b) = \begin{pmatrix} 1_m & b \\ 0 & 1_m \end{pmatrix} : b \in \mrm{Herm}_m  \right \}
    \end{align}
of $H$. We have $P(R) = M(R) N(R)$ for all $A$-algebras $R$. We occasionally write $P_m, M_m, N_m$ to emphasize dependence on $m$.

Set
    \begin{equation}
    w_j = 
    \begin{pmatrix}
    1_{m - j} & 0 & 0 & 0 \\
    0 & 0 & 0 & 1_j \\
    0 & 0 & 1_{m-j} & 0 \\
    0 & -1_j & 0 & 0
    \end{pmatrix}
    \end{equation}
for $j$ with $0 \leq j \leq m$. We also write $w = w_m$ when $j = m$ and $m$ is understood.

Let $F_v$ be a finite \'etale algebra of degree $2$ over a local field $F^+_v$. Consider $B / A = \mc{O}_{F_v} / \mc{O}_{F^+_v}$ for the respective rings of integers (with $\mc{O}_{F^+_v} \coloneqq F^+_v$ and $\mc{O}_{F_v} \coloneqq F_v$ if $F^+_v$ is Archimedean).

If $F_v/F^+_v = \C/\R$, we consider the standard maximal compact subgroup
    \begin{equation}
    K_v^{\circ} \coloneqq \left \{ \begin{pmatrix} a & b \\ - b & a \end{pmatrix} \in H(\R) : a {}^t \overline{a} + b {}^t \overline{b} = 1_m \text{ and } a {}^t \overline{b} = b {}^t \overline{a} \right \} \subseteq H(\R)
    \end{equation}
We write $U(m) \subseteq \GL_m(\C)$ for (the real points of) the unitary group for the usual positive definite rank $m$ complex Hermitian space (specified by the Gram matrix $1_m$). There is an isomorphism $K_{v}^{\circ} \ra U(m) \times U(m)$ sending the displayed matrix to $(a + i b, a - i b) \in U(m) \times U(m)$ (see e.g. \cite[{\S 2.5.1}]{GS19}).

If $F^+_v$ is non-Archimedean, we consider the standard open compact subgroup
    \begin{equation}
    K_v^{\circ} \coloneqq H(\mc{O}_{F^+_v}) \subseteq H(F^+_v).
    \end{equation}
If $F_v/F^+_v = \C / \R$ or if $F^+_v$ is non-Archimedean, we have $H(F^+_v) = P(F^+_v) K_v^{\circ}$. 
If $F^+_v$ is non-Archimedean, and given any $n(b) \in N(F^+_v)$
with any decomposition
    \begin{equation}\label{equation:Eisenstein:element_decomp}
    w^{-1} n(b) w = m(a) n(b') k
    \end{equation}
with $m(a) \in M(F^+_v)$ and $n(b') \in N(F^+_v)$ and $k \in K_v^{\circ}$,
we may assume $\det a \in F^+_v$ (it a priori lies in $F_v$)
after possibly changing $k$,
as follows from \cite[{(13.4.1)}]{Shimura97} and the explanation
at the end of \cite[{\S 3.6}]{Shimura97}.

If $F/F^+$ is a CM extension of number fields and $B/A = \mc{O}_F/\mc{O}_{F^+}$, we write
    \begin{equation}
    K^{\circ} = \prod_v K_v^{\circ} \quad \quad K_{\infty}^{\circ} = \prod_{v \mid \infty} K_v^{\circ} \quad \quad K_f^{\circ} = \prod_{v < \infty} K_v^{\circ}
    \end{equation}
where the products run over places $v$ of $F^+$. 

For places $v$ of $F^+$, we use the notation $F_v \coloneqq \prod_{w \mid v} F_w$ where $w$ runs over places of $F$, similarly $\mc{O}_{F_v} \coloneqq \prod_{w \mid v} \mc{O}_{F_w}$, as well as $F^+_{\infty} = \prod_{v \mid \infty} F^+_v$ and $F_{\infty} = \prod_{w \mid \infty} F_w$, etc..
            
            \subsection{Ad\`elic and classical Eisenstein series} 
            \label{ssec:Eisenstein:setup:adelic_v_classical}
                Characters are assumed continuous and unitary unless specified otherwise.
Let $F_v$ be a degree $2$ \'etale algebra over a local field $F^+_v$, and form the corresponding unitary group $H = U(m,m)$ as in Section \ref{ssec:Eisenstein:setup:group}. If $F^+_v$ is Archimedean, we assume in Section \ref{ssec:Eisenstein:setup:adelic_v_classical} that $F_v/F^+_v$ is $\C/\R$.

Given a character $\chi_v \colon F_v^{\times} \ra \C^{\times}$ and $s \in \C$, we may form the local \emph{degenerate principal series}
    \begin{equation}
    I(s, \chi_v) \coloneqq \mrm{Ind}_{P(F^+_v)}^{H(F^+_v)}(\chi_v | - |_{F_v}^{s + m/2}).
    \end{equation}
This is an unnormalized induction, consisting of smooth and $K_v^{\circ}$-finite functions $\Phi_v \colon H(F^+_v) \ra \C$ satisfying
    \begin{equation}\label{equation:Eisenstein:setup:adelic_v_classical:local_degen_principal_series}
    \Phi_v(m(a) n(b) h, s) = \chi_v(a) |\det a|_{F_v}^{s + m/2} \Phi_v(h, s)
    \end{equation}
for all $m(a) \in M(F^+_v)$ and $n(b) \in N(F^+_v)$ and $h \in H(F^+_v)$. Here we wrote $\chi_v(a) \coloneqq \chi_v(\det a)$ for short. A section $\Phi_v(h,s)$ of $I(s, \chi_v)$ is \emph{standard} if $\Phi(k,s)$ is independent of $s$ for any fixed $k \in K_v^{\circ}$. We say $\Phi_v$ is \emph{spherical} if $\Phi_v(h k,s) = \Phi_v(h, s)$ for any $k \in K_v^{\circ}$.
If $F_v / F^+_v$ is non-Archimedean, then $I(s,\chi_v)$ contains nonzero spherical $\Phi_v$
if and only if $\chi_v$ is unramified.
In this case,
we write $\Phi_v^{\circ}$ for the unique spherical standard section satisfying $\Phi_v^{\circ}(1,s) = 1$ for all $s$, and call $\Phi^{\circ}_v$ the \emph{normalized spherical section}.

Next, suppose $F/F^+$ is a CM extension of number fields. We write $\A_F$ for the ad\`ele ring of $F$ and $\A$ for the ad\`ele ring of $F^+$.
Given a character $\chi \colon F^\times \backslash \A_F^{\times} \ra \C^{\times}$ and $s \in \C$, we similarly form the global \emph{degenerate principal series}
    \begin{equation}
    I(s, \chi) \coloneqq \mrm{Ind}_{P(\A)}^{H(\A)}(\chi | - |_F^{s + m/2})
    \end{equation}
which is an unnormalized induction, consisting of smooth and $K^{\circ}$-finite functions $\Phi \colon H(\A) \ra \C$ satisfying
    \begin{equation}
    \Phi(m(a) n(b) h, s) = \chi(a) |\det a|_F^{s + m/2} \Phi(h, s)
    \end{equation}
for all $m(a) \in M(\A)$ and $n(b) \in N(\A)$ and $h \in H(\A)$. Given characters $\chi_f \colon \A_{F,f}^{\times} \ra \C^{\times}$ and $\chi_{\infty} \colon \A_{F,\infty}^{\times} \ra \C^{\times}$, we similarly form $I(s, \chi_f)$ and $I(s, \chi_{\infty})$. We also speak of \emph{spherical sections} and the \emph{spherical standard section}, as above. We sometimes write $I_m(s, \chi)$ etc. to indicate dependence on $m$.

Given a standard section $\Phi(h,s)$ of the global degenerate principal series $I(s, \chi)$, we form the \emph{Siegel Eisenstein series}
    \begin{equation}
    E(h, s, \Phi) = \sum_{\g \in P(F^+) \backslash H(F^+)} \Phi(\g h, s)
    \end{equation}
which is absolutely convergent for $\Re(s) > m/2$. 
We also form $E(h,\Phi,s)$ when $\Phi$ is a finite meromorphic linear combination of standard sections by extending linearly.

Define another character $\check{\chi} \colon F^{\times} \backslash \A^{\times}_F \ra \C^{\times}$ as $\check{\chi}(a) \coloneqq \chi(\overline{a})^{-1}$. There is a functional equation
    \begin{equation}\label{equation:Eisenstein:setup:adelic_v_classical:global_functional_equation}
    E(h, -s,  M(\chi, s) \Phi) = E(h, s, \Phi)
    \end{equation}
where $M(\chi,s) \colon I(s, \chi) \ra I(-s, \check{\chi})$ is the intertwining operator
    \begin{equation}
    (M(s, \chi) \Phi)(h) = \int_{N(\A)} \Phi(w^{-1} n(b) h, s) ~ d n(b)
    \end{equation}
for $\Re(s) > m/2$ (see e.g. \cite{Tan99}).
Here the Haar measure on $N(\A)$ should be the one which is self-dual with respect
to any nontrivial additive character $\psi \colon F^+ \backslash \A \ra \C^{\times}$,
and the pairing $b,b' \mapsto \psi(\operatorname{tr}(b b'))$.
We occasionally write $M_m(s, \chi)$ to emphasize the understood $m$ (in $U(m,m)$).

Fix an identification of $F^+_v$-algebras $F_v \cong \C$ for each Archimedean place $v$ of $F^+$. We consider classical Eisenstein series on the \emph{Hermitian upper half-space}
    \begin{align}
    \mc{H}_m & \coloneqq \{ z \in M_{m, m}(F_{\infty}) : (2i)^{-1} (z - {}^t \overline{z}) > 0 \}
    \\
    & = \{ z = x + i y : x,y \in \mrm{Herm}_m(F^+_{\infty}) \text{ with } y > 0 \}, \label{equation:Eisenstein:setup:adelic_v_classical:Hermitian_upper-half_space}
    \end{align}
where the latter expression means that $x$ and $y$ are $m \times m$ Hermitian matrices with $y$ positive definite (at every place $v \mid \infty$ of $F^+_v$).
Given $z = x + i y \in \mc{H}_m$, we write $h_z \in H(F^+_{\infty}) \subseteq H(\A)$ for any element $h_z = n(x) m(a)$ where $a \in \mrm{GL}_m(F_{\infty})$ satisfies $a {}^t \overline{a} = y$. Note $h_z \cdot i 1_m = z$.

We restrict to $\Phi = \Phi_{\infty} \otimes \Phi_f$ for standard sections $\Phi_{\infty} \in I(s, \chi_{\infty})$ and $\Phi_f \in I(s, \chi_f)$. Fix an integer $n_v$ for each place $v \mid \infty$ of $F^+_v$, and assume $\chi_v|_{F^{+ \times}_v} = \mrm{sgn}(-)^{n_v}$ for every $v \mid \infty$. We also let $k(\chi_v) \in \Z$ be the integer satisfying
    \begin{equation}
    \chi_{v}(z) = (z / |z|^{1/2}_{F_v} )^{k(\chi_v)} \quad \quad \text{where $z \in F_v$,}
    \end{equation}
for each place $v \mid \infty$ of $F^+_v$. For such $v$, we let $\Phi_v = \Phi_v^{(n_v)}$ be the unique standard section of $I(s, \chi_v)$ of scalar weight
    \begin{equation}\label{equation:adelic_v_classical:scalar_weight}
    (n_1, n_2) \quad \quad \text{where} \quad \quad n_1 = \frac{n_v + k(\chi)}{2} \quad \text{and} \quad n_2 = \frac{-n_v + k(\chi)}{2}
    \end{equation}
such that $\Phi_{v}^{(n_v)}(1,s) = 1$ (as in \cite[{\S 3.2, \S 3.3}]{GS19}). The scalar weight condition means that $\Phi_{v}^{(n_v)}(h k,s) = \det(k_1)^{n_1} \det(k_2)^{n_2} \Phi_{v}^{(n_v)}(h,s)$ for all $h \in H(F^+_v)$ and $k \in K_{v}^{\circ}$ where
    \begin{equation}
    k = 
    \begin{pmatrix}
    a & b \\
    - b & a 
    \end{pmatrix} \in K_{v}^{\circ}
    \quad \quad
    k_1 = a + i b
    \quad \quad
    k_2 = a - i b.
    \end{equation}
Note that $\Phi_v^{(n_v)}$ does not depend on the choice of identification $F_v \cong \C$.

If $y = a {}^t \overline{a}$ for some $a \in \GL_m(F_v)$, a computation (omitted) shows
    \begin{equation}\label{equation:Eisenstein:setup:adelic_v_classical:explicit_Archimedean_scalar_weight_vector}
    \chi_v(a)^{-1} (\det y)^{- n_v/2} \Phi^{(n_v)}_v(w^{-1}n(b)m(a)) = (\det y)^{s - s_0} \det(-i y + b)^{-(s - s_0)} \det(iy + b)^{-(s - s_0) - n_v}
    \end{equation}
for any $b \in \mrm{Herm}_m(F^+_v)$, where $s_0 = (n_v - m)/2$ (reduce to the case $a = 1_m$ and write $w^{-1} n(b) = n(-b(1_m + b^2)^{-1}) m(b + i 1_m)^{-1} k$ for $k \in K_v^{\circ}$).
Equation \eqref{equation:Eisenstein:setup:adelic_v_classical:explicit_Archimedean_scalar_weight_vector} may be used to translate various statements from \cite{Shimura82} to statements about Archimedean Whittaker functions, etc. (see \crefext{II:ssec:Archimedean_identity:more_Whittaker} for more on this).

\begin{remark}\label{remark:Eisenstein:setup:adelic_v_classical:log_det_definition}
Given $g = x_g + i y_g \in M_{m,m}(\C)$ with $x_g, y_g$ Hermitian and $x_g$ positive definite, we define $\log \det(g)$ by the ``principal branch'' (such that $g \mapsto \log \det g$ is holomorphic, and $\log \det g \in \R$ if $y_g = 0$) as in \cite[{(1.11)}]{Shimura82} and the surrounding discussion of loc. cit.. If $y_g$ is positive definite and $x_g$ is only assumed Hermitian, we also take
    \begin{equation}
    \log \det g = \log \det (-i g) + m \log i \quad \quad \log \det \overline{g} = \log \det (i g) - m \log i
    \end{equation}
where $\log i \coloneqq \pi i/2$ (as in \cite[{(1.11)}]{Shimura82}).
This convention is implicit in \eqref{equation:Eisenstein:setup:adelic_v_classical:explicit_Archimedean_scalar_weight_vector}.
\end{remark}

We take $\Phi_{\infty} = \otimes_{v \mid \infty} \Phi_v^{(n_v)}$. We write $n = (n_v)_{v \mid \infty}$ for the collection of Archimedean weights (and will eventually focus on the case where all $n_v$ are equal to some fixed integer $n$). In the above situation, we write $E(h,s,\Phi)_n \coloneqq E(h,s,\Phi)$ and  consider an associated \emph{classical Eisenstein series}
    \begin{equation}\label{equation:Eisenstein:setup:adelic_v_classical:associated_classical_series}
    E(z, s, \Phi)_n \coloneqq E(z, s, \Phi) \coloneqq \chi_{\infty}(a)^{-1} \det(y)^{-n/2} E(h_z, s, \Phi)_n
    \end{equation}
where $z = x + i y$ and $h_z = n(x) m(a)$ with $a {}^t \overline{a} = y$ as above, and where $\det(y)^{-n/2}$ stands for $\prod_{v \mid \infty} \det(y_v)^{-n_v/2}$. This does not depend on the choice of $h_z$, i.e. $E(h_z k_{\infty}, s, \Phi)_n = E(h_z, s, \Phi)_n$ for any $k_{\infty} \in K_{\infty}^{\circ}$.

When $F^+ = \Q$ and $s_0 \coloneqq (n - m)/2$ (setting $n = n_{\infty}$ and $k(\chi) = k(\chi_{\infty})$), a computation (omitted) gives the more classical form
    \begin{align}
    E(z, s, \Phi)_n & = \sum_{\g \in P(\Z) \backslash U(m,m)(\Z)} \frac{\det(y)^{s - s_0} \det(\gamma)^{(n + k(\chi))/2}}{\det (c z + d)^n |\det ( c z + d )|^{2(s - s_0)}} \Phi_f(\g, s)
    \\
    & = \sum_{\g \in P_1(\Z) \backslash SU(m,m)(\Z)} \frac{\det(y)^{s - s_0}}{\det (c z + d)^n |\det ( c z + d )|^{2(s - s_0)}} \Phi_f(\g, s)  
    \end{align}
where 
    \begin{equation}
    \g = \begin{pmatrix} a & b \\ c & d \end{pmatrix},
    \end{equation}
where $SU(m,m) \subseteq U(m,m)$ is the determinant $1$ subgroup, and $P_1 \coloneqq SU(m,m) \cap P$. We have $P(\Q) \backslash H(\Q) = P(\Z) \backslash H(\Z) = P_1(\Q) \backslash H_1(\Q) = P_1(\Z) \backslash H_1(\Z)$ (e.g. \cite[{Proposition 12.6}]{Ikeda08}). 
When $m = 1$, the exceptional isomorphism $\SL_2 \ra SU(1,1)$ (over $\Spec \Q$) implies that the above expression is a classical Eisenstein series for $\SL_2$ on the upper half-plane.

Our main theorem (Theorem \ref{theorem:arithmetic_Siegel-Weil:main_results:main}) concerns Fourier coefficients of $E(z,s,\Phi)_n$ (normalized as in Section \ref{ssec:Eisenstein:global_normalized_Fourier:global_normalization}), but the variant
    \begin{equation}\label{equation:setup:adelic_v_classical:variant}
    \tilde{E}(a, s, \Phi)_n \coloneqq \chi(a)^{-1} |\det(a)|_F^{- n / 2} E(m(a), s, \Phi)_n \quad \text{for $a \in \GL_m(\A_F)$}
    \end{equation}
will be useful for studying Fourier coefficients of $E(z,s,\Phi)_n$ for singular $T$ (see below).
If $a \in \GL_m(F_{\infty})$ is any element satisfying $a {}^t \overline{a} = y$ for $y \in \mrm{Herm}_m(F^+_{\infty})_{>0}$, we have
    \begin{equation}
    E(i y, s, \Phi)_n = \tilde{E}(a, s, \Phi)_n.
    \end{equation}

            \subsection{Fourier expansion and local Whittaker functions} 
            \label{ssec:Eisenstein:setup:Fourier_and_Whittaker}
                Take notation as in Section \ref{ssec:Eisenstein:setup:adelic_v_classical}, e.g. $F/F^+$ is a CM extension of number fields. Choose a nontrivial additive character $\psi \colon F^+ \backslash \A \ra \C^\times$. We have a Fourier expansion
    \begin{equation}
    E(h,s,\Phi) = \sum_{T \in \mrm{Herm}_m(F^+)} E_T(h, s, \Phi)
    \end{equation}
where
    \begin{equation}
    E_T(h, s, \Phi) = \int_{N(F^+) \backslash N(\A)} E( n(b) h, s, \Phi) \psi(-\mrm{tr}(Tb)) ~ d n(b)
    \end{equation}
for $\mrm{Re}(s) > m/2$, and  where $d n(b)$ is the Haar measure on $N(\A)$ which is self-dual with respect to the pairing $(b, b') \mapsto \psi(\mrm{tr}(b b'))$. We refer to $E_T(h,s,\Phi)$ as the \emph{$T$-th Fourier term}.

For any $a \in \GL_m(F)$, a change of variables gives
    \begin{equation} \label{equation:Eisenstein:setup:Fourier_and_Whittaker:invariance}
    E_T(m(a) h, s, \Phi) = E_{{}^t \overline{a} T a}(h, s, \Phi).
    \end{equation}
We also have 
    \begin{equation} \label{equation:eisenstein:setup:Fourier_and_Whittaker:invariance_y}
    E_T(m(a) h, s, \Phi) = E_T(h, s, \Phi) \quad \text{for any} \quad \begin{pmatrix}
    1_{m - m^{\flat}} & * \\
    0 & 1_{m^{\flat}}
    \end{pmatrix}
    \in \GL_m(\A_F)
    \quad
    \text{if \quad}
    T = \begin{pmatrix}
    0 & 0 \\
    0 & T^{\flat}
    \end{pmatrix}
    \end{equation}
with the block matrix $T^{\flat} \in \mrm{Herm}_{m^{\flat}}(F^+)$ having $\det T^{\flat} \neq 0$ (here $m^{\flat}$ is arbitrary) (follows from \cite[{Lemma 5.4, (5.56)}]{GS19}).

Allowing arbitrary $T$ again, assume there is a factorization $\Phi = (\otimes_{v \mid \infty} \Phi_v) \otimes \Phi_f$. For each $v \mid \infty$, assume $\Phi_v = \Phi_v^{(n_v)}$ is the scalar weight standard section as in Section \ref{ssec:Eisenstein:setup:adelic_v_classical}, for some $n_v \in \Z$. Write $n = (n_v)_{v \mid \infty}$ for the resulting tuple of integers.

Consider $a = a_{\infty} a_f \in \GL_m(\A_F)$, with $a_{\infty} \in \GL_m(F_{\infty})$ and $a_f \in \GL_m(\A_{F,f})$. Set $y = a_{\infty} {}^t \overline{a}_{\infty}$ (temporary). We then have $T$-th \emph{Fourier coefficients} $E_T(y,s,\Phi)_n$ and $\tilde{E}_T(a, s, \Phi)_n$
characterized by the relations
    \begin{align} \label{equation:setup:Fourier_and_Whittaker:y_Fourier}
    E_T(y,s,\Phi)_n q^T & =  \chi_{\infty}(a_{\infty})^{-1} \det(y)^{-n / 2} E_T(n(x) m(a_{\infty}), s, \Phi)
    \\
    \tilde{E}_T(a,s,\Phi)_n \psi_f(\mrm{tr}(T b)) q^T & =  \chi(a)^{-1} |\det a|_F^{-n / 2} E_T(n(x + b) m(a), s, \Phi) \label{equation:setup:Fourier_and_Whittaker:a_Fourier}
    \end{align}
for any $x \in \mrm{Herm}_m(F_{\infty})$ and $b \in \mrm{Herm}_m(\A_f)$, with $z \coloneqq x + i y$, and with $q^T \coloneqq \psi_{\infty}(\mrm{tr}(Tz))$. These correspond to the classical Eisenstein series and its variants in \eqref{equation:Eisenstein:setup:adelic_v_classical:associated_classical_series} and \eqref{equation:setup:adelic_v_classical:variant}.

When $\det T \neq 0$ and $\Phi = \otimes_v \Phi_v$ is factorizable over all places, we have a factorization
    \begin{equation}
    E_T(h, s, \Phi) = \prod_v W_{T,\psi_v}(h_v, s, \Phi_v)
    \end{equation}
into \emph{local Whittaker functions} defined below \eqref{equation:local_Whittaker_function}.

We switch to local notation: let $F_v$ be a degree $2$ \'etale algebra over a local field $F^+_v$, with nontrivial involution $a \mapsto \overline{a}$. We assume $F^+_v$ has characteristic $0$ (because Karel assumes this \cite{Karel79}). If $F^+_v$ is Archimedean, we also assume $F_v / F^+_v = \C / \R$.

Let $\chi_v \colon F_v^{\times} \ra \C^{\times}$ and $\psi_v \colon F^+_v \ra \C^{\times}$ be characters with $\psi_v$ nontrivial, and suppose $\Phi_v \in I(s, \chi_v)$ is a standard section. Given $T \in \mrm{Herm}_m(F^+_v)$ with $\det T \neq 0$, there is a \emph{local Whittaker function} defined by the absolutely convergent integral
    \begin{equation}\label{equation:local_Whittaker_function}
    W_{T,\psi_v}(h, s, \Phi_v) \coloneqq \int_{N(F^+_v)} \Phi_v(w^{-1} n(b) h, s) \psi_v(- \mrm{tr}(Tb)) ~dn(b)
    \end{equation}
for $h \in H(F^+_v)$ and $s \in \C$ with $\mrm{Re}(s) > m/2$, where $d n(b)$ is the Haar measure which is self-dual with respect to the pairing $(b, b') \mapsto \psi_v(\mrm{tr}(b b'))$ on $\mrm{Herm}_m(F^+_v) \cong N(F^+_v)$. For each fixed $h$, the function $W_{T,\psi_v}(h, s, \Phi_v)$ admits holomorphic continuation to $s \in \C$ \cite[{Corollary 3.6.1}]{Karel79}\cite{KS97}\cite[{\S 6}]{Ichino04}. Extending linearly defines $W_{T,\psi_v}(h, s, \Phi_v)$ whenever $\Phi_v$ is a finite meromorphic linear combination of standard sections. For any $a \in \GL_m(F_v)$, a change of variables shows
    \begin{equation}\label{equation:setup:Fourier_and_Whittaker:Whittaker_linear_change}
    W_{T,\psi_v}(m(a) h, s, \Phi_v) = \check{\chi}_v(a) |\det a|_{F_v}^{-s + m / 2} W_{{}^t \overline{a} T a, \psi_v}(h, s, \Phi_v)
    \end{equation}
for $\check{\chi}_v(a) \coloneqq \chi_v(\overline{a})^{-1}$ as above. We use the shorthand $W_{T,\psi_v}(s, \Phi_v) \coloneqq W_{T,\psi_v}(1, s, \Phi_v)$.

If $F_v / F^+_v$ is $\C / \R$, we write $W_{T,v}$ instead of $W_{T,\psi_v}$
when $\psi_v(x) = e^{2 \pi i x}$ is the standard additive character.
If $F^+_v$ is non-Archimedean, we write $W_{T,v}$ instead of $W_{T,\psi_v}$
if $\psi_v$ is an understood unramified additive character (suppressed from notation).

\begin{lemma}\label{lemma:non-Arch_Whittaker_polynomial}
With notation as above, assume that $F^+_v$ is non-Archimedean with residue field of cardinality $q_v$. Suppose $\Phi_v \in I(s, \chi_v)$ is a standard section and $h \in H(F^+_v)$ is a fixed element.
    \begin{enumerate}[(1)]
        \item We have $W_{T,\psi_v}(h, s, \Phi_v) \in \C[q_v^{-s}, q_v^s]$.
        \item If $h \in K_v^{\circ}$, we have $W_{T,\psi_v}(h, s, \Phi_v) \in \C[q_v^{-2s}]$.
            \begin{enumerate}[(i)]
                \item If $\chi_v|_{F^{+ \times}_v}$ is a finite order character
                and $\Phi_v(k,s) \in \overline{\Q}$ for every $k \in K^{\circ}_v$,
                then we have $W_{T,\psi_v}(h, s, \Phi_v) \in \overline{\Q}[q_v^{-2s}]$ for any $h \in K^{\circ}_v$.
                \item If $\chi_v|_{F^{+ \times}_v}$ is a finite order character
                and $\Phi_v(k,s) \in \overline{\Q}$ for every $k \in K^{\circ}_v$,
                such that moreover every $\Phi_v(k,s)$
                is $\ell$-integral for some prime $\ell$ distinct from the residue
                characteristic of $F^+_v$, then $W_{T,\psi_v}(h, s, \Phi_v) \in \overline{\Q}[q_v^{-2s}]$
                also has coefficients which are $\ell$-integral.
                \item Suppose $\chi_v$ is unramified and $\Phi_v = \Phi^{\circ}_v$ is the spherical standard section.
                    Assume $\psi_v(\operatorname{tr}(Tb)) = 1$ for every $b \in \mrm{Herm}_m(\mc{O}_{F^+_v})$.
                    If $f(X) \in \C[X]$ is the polynomial satisfying $W_{T,\psi_v}(1,s,\Phi_v) = f(q_v^{-2s})$,
                    then $f(0) = \operatorname{vol}(\mrm{Herm}_m(\mc{O}_{F^+_v}))$.
            \end{enumerate}
        \item Suppose $\chi'_v \colon F_v^{\times} \ra \C^{\times}$ is another character satisfying $\chi'_v|_{F^{+ \times}_v} = \xi_v \chi_v|_{F^{+ \times}_v}$ for an unramified character $\xi_v \colon F^{+ \times}_v \ra \C^{\times}$.
        Assume $h \in K_v^{\circ}$, and suppose $\Psi_v \in I(s, \chi'_v)$ is a standard section satisfying $\Psi_v(k, s) = \Phi_v(k, s)$
        for any $k \in K^{\circ}_v$ (for all $s$).
        If $f(X) \in \C[X]$ is the polynomial satisfying $f(q_v^{-2s}) = W_{T,\psi_v}(h, s, \Phi_v)$, then we have $f(\xi_v(\varpi_0) q_v^{-2s}) = W_{T,\psi_v}(h, s, \Psi_v)$, where $\varpi_0 \in F^+_v$ is any uniformizer.
    \end{enumerate}
\end{lemma}
\begin{proof}
Claim (1) follows from a general result of Karel \cite[{Corollary 3.6.1}]{Karel79},
which implies that $W_{T,\psi_v}(h, s, \Phi_v) \in \C[q_v^{-s}, q_v^s]$,
and that $W_{T,\psi_v}(h, s, \Phi_v)$ may be computed for all $s$ as
the integral over a sufficiently large open compact subgroup of $N(F^+_v)$.

Next, given any $b \in N(F^+_v)$, write $w^{-1} n(b) w = m(a) n(b') k$
for $m(a) \in M(F^+_v)$ and $n(b') \in N(F^+_v)$ and $k \in K^{\circ}_v$,
as in \eqref{equation:Eisenstein:element_decomp}.
As explained in loc. cit., we can moreover assume $\det a \in F^+_v$.
We then have $\Phi_v(w^{-1} n(b) h, s) = \Phi_v(m(a) n(b') k w^{-1} h, s)
= \chi_v(\det a) |\det a|_{F_v}^{s + m/2} \Phi_v(k w^{-1} h, s)$.
Since $\det a \in F^+_v$, 
we must have $|\det a|_{F_v}^s$ being a nonnegative power of $q_v^{-2s}$.
Substituting into \eqref{equation:local_Whittaker_function}
(and replacing the integral over $N(F^+_v)$ with an integral
over a sufficiently large open compact subgroup)
now proves claim (2).
The claims about algebraicity and integrality
in (2)(i) and (2)(ii) follow because
$\psi_v$ and $\chi|_{F^{+ \times}_v}$ are both valued in roots of unity,
under the given hypotheses.
The claim in (2)(iii) follows because,
in the above notation,
we have $|\det a|_{F_v} = 1$ if and only if $b \in \mrm{Herm}_m(\mc{O}_{F^+_v})$
by inspection, or e.g. see the discussion at the end of \cite[{\S 3.6}]{Shimura97}
or \cite[{\S 13.5}]{Shimura97}.
The volume $\operatorname{vol}(\mrm{Herm}_m(\mc{O}_{F^+_v}))$
should be taken with respect to the self-dual Haar measure associated to $\psi_v$,
as above.

For (3), a similar computation gives
$\Psi_v(w^{-1} n(b) h, s) = \chi'_v(\det a) |\det a|^{s + m / 2}_{F_v} \Psi_v(k w^{-1} h, s)$,
for the same $a \in \GL_m(F_v)$ as above.
If $v_0$ is the valuation on $F^+_v$ normalized with $v_0(\varpi_0) = 1$,
we thus have
    \begin{align}
    \Phi_v(w^{-1} n(b) h, s) & = \chi_v(\det a) (q_v^{-2s})^{v_0(\det a)} \cdot q_v^{-m \cdot v_0(\det a)} \Phi_v(k w^{-1} h, s)
    \\
    \Psi_v(w^{-1} n(b) h, s) & = \chi_v(\det a) (\xi_v(\varpi_0) q_v^{-2s})^{v_0(\det a)} \cdot q_v^{-m \cdot v_0(\det a)} \Psi_v(k w^{-1} h, s).
    \end{align}
Claim (3) of the lemma now follows upon substitution into \eqref{equation:local_Whittaker_function}
(and again replacing the integral over $N(F^+_v)$ with an integral
over a sufficiently large open compact subgroup).
\end{proof}

In the case where $F^+_v$ is non-Archimedean, consider the case where $\chi_v$ is unramified and $\chi_v|_{F^{+ \times}_v} = \eta_v^n$ for some integer $n$, where $\eta_v \colon F^{+ \times}_v \ra \{ \pm 1\}$ is the quadratic character associated to $F_v / F^+_v$. Consider the normalized spherical standard section $\Phi^{\circ}_v \in I(s, \chi_v)$. We temporarily write $W_{T,\psi_v}(h, s, \Phi_v^{\circ})_n$ for the associated local Whittaker function, emphasizing the possible dependence on $n$.
By Lemma \ref{lemma:non-Arch_Whittaker_polynomial}(3), the implicit $\chi_v$-dependence of $W_{T,\psi_v}(h, s, \Phi_v^{\circ})_n$ is only on the restriction $\chi_v|_{F^{+ \times}_v}$. If $F_v / F^+_v$ is not inert, then $W_{T,\psi_v}(h, s, \Phi_v^{\circ})_n$ does not depend on $n$ (note $n$ must be even if $F_v / F^+_v$ is ramified). If $F_v / F^+_v$ is inert, then $W_{T,\psi_v}(h, s, \Phi_v^{\circ})_n$ depends only on the parity of $n$. The $\C$-algebra endomorphism of $\C[q_v^{-2s}]$ sending $q_v^{-2s} \mapsto -q_v^{-2s}$ swaps $W_{T,\psi_v}(h, s, \Phi_v^{\circ})_n$ and $W_{T,\psi_v}(h, s, \Phi_v^{\circ})_{n+1}$, by Lemma \ref{lemma:non-Arch_Whittaker_polynomial}(3).

\begin{remark}[Change of additive character]\label{remark:Eisenstein:setup:Fourier_and_Whittaker:change_of_additive_character}
We sometimes fix convenient choices of additive character $\psi$ or $\psi_v$.
This recovers the general case, as follows.
Given $\psi$, we write $E_{T, \psi}$ as needed
to emphasize dependence on the choice of additive character,
and similarly for normalized variants like $\tilde{E}_{T,\psi}$, etc..

Given $u \in (F^+)^{\times}$ or $u \in (F^+_v)^{\times}$,
consider $\psi'(x) \coloneqq \psi(u \cdot x)$ and $\psi_{v}'(x) \coloneqq \psi_v(u \cdot x)$.
We have
    \begin{align}
    E_{u T,\psi}(h,s,\Phi) &= E_{T, \psi'}(h,s,\Phi) 
    && E_{u T,\psi}(y,s,\Phi)_n = E_{T, \psi'}(y,s,\Phi)_n
    \\
    \tilde{E}_{u T,\psi}(a,s,\Phi)_n &= \tilde{E}_{T, \psi'}(a,s,\Phi)_n
    && |u|_{F^+_v}^{m^2/2} W_{u T,\psi_v}(h,s,\Phi_v) = W_{T, \psi_{v}'}(h,s,\Phi_v)
    \end{align}
with notation as above.
The additional factor $|u|^{m^2/2}_{F^+_v}$
arises from the change of Haar measure
on $N(F^+_v)$, from the requirements of self-duality
with respect to $\psi_v$ or $\psi'_v$.
\end{remark}

            \subsection{Singular Fourier coefficients} 
            \label{ssec:Eisenstein:setup:singular_Fourier}
                Retain notation from Section \ref{ssec:Eisenstein:setup:Fourier_and_Whittaker} (switching back to global notation). The Fourier terms $E_{T}(h,s,\Phi)$ for singular $T \in \mrm{Herm}_m(F^+)$ are known to be closely related with Fourier terms of Eisenstein series on smaller groups (e.g. \cite[{\S 5.2}]{GS19}). We focus on the case where $\rank T = m - 1$ (assume this throughout Section \ref{ssec:Eisenstein:setup:singular_Fourier}). On account of \eqref{equation:Eisenstein:setup:Fourier_and_Whittaker:invariance}, it will be enough to describe the case where $T$ is block diagonal of the form
    \begin{equation}\label{equation:Eisenstein:setup:singular_Fourier:block_diagonal}
    T = 
    \begin{pmatrix}
    0 & 0 \\
    0 & T^{\flat}
    \end{pmatrix}
    \end{equation}
with $\det T^{\flat} \neq 0$.

Assume $m \geq 1$, and fix an integer $n \in \Z$.
Let $\chi \colon F^\times \backslash \A_F^{\times} \ra \C^{\times}$ be a character satisfying $\chi|_{\A^{\times}} = \eta^n$, where $\eta$ is the quadratic character associated with $F / F^+$. Note $\check{\chi} = \chi$ in this case.

Take a factorizable standard section $\Phi = \otimes_v \Phi_v \in I(s, \chi)$, and assume $\Phi_v = \Phi_v^{(n)}$ is the normalized scalar weight standard section (Section \ref{ssec:Eisenstein:setup:adelic_v_classical}) for every Archimedean place $v$, with $n$ the fixed integer from above (same for every $v \mid \infty$). 

Take $T$ as in \eqref{equation:Eisenstein:setup:singular_Fourier:block_diagonal}. Given $a \in \GL_m(\A_F)$, we study the Fourier coefficient $\tilde{E}_T(a, s, \Phi)_n$. By the Iwasawa decomposition, every $a \in \GL_m(\A_F)$ admits a decomposition
    \begin{equation}
    a = 
    \begin{pmatrix}
    1_1 & * \\
    0 & 1_{m - 1} 
    \end{pmatrix}
    \begin{pmatrix}
    a^{\#} & 0 \\
    0 & a^{\flat}
    \end{pmatrix}
    k
    \end{equation}
with $a^{\#} \in \GL_1(\A_F)$, with $a^{\flat} \in \GL_{m - 1}(\A_F)$, and with $k \in \prod_{v \mid \infty} U(m) \times \prod_{v < \infty} \GL_m(\mc{O}_{F_v})$. We will be eventually interested in the case when $\Phi_f$ is spherical, which implies $\tilde{E}_T(a k, s, \Phi)_n = \tilde{E}_T(a, s, \Phi)$ for any $k \in \prod_{v \mid \infty} U(m) \times \prod_{v < \infty} \GL_m(\mc{O}_{F_v})$ (also using the fact that $\Phi_v$ is a scalar weight standard section for each $v \mid \infty$). In light of the invariance property in \eqref{equation:eisenstein:setup:Fourier_and_Whittaker:invariance_y}, it is thus harmless to restrict to the case of block diagonal $a = \mrm{diag}(a^{\#}, a^{\flat})$. Assume this for the rest of Section \ref{ssec:Eisenstein:setup:singular_Fourier} (but we do not assume $\Phi_f$ is spherical for now).

Set $m^{\flat} \coloneqq m - 1$ and $s_0 \coloneqq (n - m)/2$.
Consider the operators
\begin{equation}
    \begin{tikzcd}[row sep = tiny]
    I_m(s, \chi) \arrow{r}{\mu^{m *}_{m^{\flat}}(s, \chi)} & I_{m^{\flat}}(s+1/2, \chi) \\
    \Psi \arrow[mapsto]{r} & \Psi \circ \mu^m_{m^{\flat}}
    \end{tikzcd}
    \end{equation}
(with $\mu^m_{m^{\flat}} \colon U(m^{\flat}, m^{\flat}) \ra U(m,m)$ as in Section \ref{ssec:Eisenstein:setup:group}), and where
    \begin{equation}
    \begin{tikzcd}[row sep = tiny]
    I_m(s, \chi) \arrow{r}{U^m_{m^{\flat}}(s, \chi)} & I_{m^{\flat}}(s-1/2, \chi) \\
    \Psi \arrow[mapsto]{r} & \left ( h \mapsto \scalebox{2}{$\int$}_{\substack{b_1 \in \mrm{Herm}_{m - m^{\flat}}(\A) \\ b_{12} \in M_{m - m^{\flat}, m^{\flat}}(\A_F)}} \Psi\left ( w_m^{-1} \cdot n \begin{pmatrix} b_1 & b_{12} \\ {}^t \overline{b}_{12} & 0 \end{pmatrix} w_{m^{\flat}} \mu_{m^{\flat}}^m(h), s \right ) ~ db_1 ~d b_{12} \right )
    \end{tikzcd}
    \end{equation}
for $\mrm{Re}(s) > m/2$ (with meromorphic continuation to $s \in \C$).
A calculation shows
    \begin{equation}\label{equation:Eisenstein:setup:singular_Fourier}
    M_{m^{\flat}}(s - 1/2, \chi) \circ U^m_{m^{\flat}}(s, \chi) = \mu_{m^{\flat}}^{m *}(-s, \chi) \circ M_m(s, \chi);
    \end{equation}
compare \cite[{Lemma 5.5(iii)}]{GS19}
(this will not be used until Section \ref{ssec:Eisenstein:global_normalized_Fourier:singular_Fourier}).

\begin{lemma}
We have
    \begin{align}
    \label{equation:Eisenstein:setup:singular_Fourier:corank_1}
    \tilde{E}_{T}(a, s, \Phi)_n & = |\det a^{\#}|_F^{s - s_0} \tilde{E}_{T^{\flat}}(a^{\flat}, s + 1/2, \mu_{m^{\flat}}^{m*}(s, \chi) \Phi)_n 
    \\
    & \hphantom{=} + |\det a^{\#}|_F^{-s - s_0} \tilde{E}_{T^{\flat}}(a^{\flat}, s - 1/2, U_{m^{\flat}}^m(s, \chi) \Phi)_n \notag.
    \end{align}
\end{lemma}
\begin{proof}
We claim that arguing as in the proof of \cite[{Lemma 2.4}]{KR88} (see also \cite[{Lemma 5.4}]{GS19} and \cite[{Theorem 2.2}]{HSY21}) gives
    \begin{align}\label{equation:Eisenstein:setup:singular_Fourier:corank_1:unnormalized}
    E_{T}(a, s, \Phi) & = \chi(\det a^{\#}) |\det a^{\#}|_F^{s + m/2} E_{T^{\flat}}(a^{\flat}, s + 1/2, \mu_{m^{\flat}}^{m*}(s, \chi) \Phi) 
    \\
    \notag
    & \hphantom{=} + \chi(\det \overline{a}^{\#})^{-1} |\det a^{\#}|_F^{-s + m/2} E_{T^{\flat}}(a^{\flat}, s - 1/2, U_{m^{\flat}}^m(s, \chi) \Phi).
    \end{align}
The formula in the lemma statement follows from this,
upon using our normalization convention in
\eqref{equation:setup:adelic_v_classical:variant}
and using the equality $\check{\chi} = \chi$.

The special case of \eqref{equation:Eisenstein:setup:singular_Fourier:corank_1:unnormalized}
with $a = 1$ and $F^+ = \Q$ is explained in \cite[{Proposition 2.4}]{HSY21}
(taking their $a$ to equal $1$, with their $w_n$ being our $w^{-1}$,
with their $g$ being our $m(a)$, and noting that their $T_a$ need only be nonsingular),
and the case of general $F / F^+$ is similar.
In the general case (for arbitrary block diagonal $a = \diag(a^{\#}, a^{\flat})$
and arbitrary $F / F^+$),
we claim that the first line (resp. second line)
of our \eqref{equation:Eisenstein:setup:singular_Fourier:corank_1:unnormalized}
corresponds to $W^{(n-1)}_{T^{\flat}}(m(a), s, \Phi)$
and $W^{(n)}_{T}(m(a), s, \Phi)$, in the notation of \cite[{Theorem 2.2}]{HSY21}.
Reformulated in our notation, we have
$W^{(r)}_T(m(a), s, \Phi) = \int_{N_{m-1}(\A)} \Phi(w^{-1}_{m-1} n(\operatorname{diag}(0,b)) m(a), s) \psi(-\operatorname{tr}(T^{\flat} b)) ~dn(b)$.
We have $w^{-1}_{m-1} n(\operatorname{diag}(0,b)) m(a) = m(\diag(a^{\#}, 1)) w^{-1}_{m-1} n(\operatorname{diag}(0,b)) m(\diag(1, a^{\flat}))$,
so the first line of \eqref{equation:Eisenstein:setup:singular_Fourier:corank_1:unnormalized}
now follows upon applying \eqref{equation:Eisenstein:setup:adelic_v_classical:local_degen_principal_series},
with $m(a)$ of loc. cit. replaced by $m(\diag(a^{\#}, 1))$.
For the second line of \eqref{equation:Eisenstein:setup:singular_Fourier:corank_1:unnormalized},
consider the third-to-last line of the last displayed equation in the proof of
\cite[{Proposition 2.4}]{HSY21}. Taking $a = 1$ in loc. cit., with their $w_n$ being our $w^{-1}$,
and taking their $g$ to be our $m(a)$,
the second line of \eqref{equation:Eisenstein:setup:singular_Fourier:corank_1:unnormalized}
now follows from the formula
    \begin{align}
    & w^{-1} n(b) m(a) = m(\diag({}^t\overline{a}^{\# -1},1)) w^{-1} n(b') m(\diag(1, a^{\flat}))
    \\
    & \text{where $n(b') \coloneqq m(\diag(a^{\#}, 1)) \cdot n(b) \cdot m(\diag(a^{\#}, 1))^{-1}$}
    \end{align}
for any $b \in N(\A)$,
noting that $b$ and $b'$ have the same lower-right $(m-1) \times (m-1)$ block,
and again applying \eqref{equation:Eisenstein:setup:adelic_v_classical:local_degen_principal_series},
with $m(a)$ of loc. cit. replaced by $m(\diag(\overline{a}^{\# -1}, 1))$.
To obtain the second line of \eqref{equation:Eisenstein:setup:singular_Fourier:corank_1:unnormalized},
one should additionally note that the change of variables
$b' \mapsto b$ introduces a factor of $|\det a^{\#}|_F^{m}$,
so that we obtain $|\det a^{\#}|_F^{-s + m/2}$ rather than $|\det a^{\#}|_F^{-s - m/2}$.
\end{proof}

In Corollary \ref{corollary:Eisenstein:singular_Fourier:corank_1}, we rewrite \eqref{equation:Eisenstein:setup:singular_Fourier:corank_1} more explicitly when $\Phi_v$ is the normalized spherical standard section for every non-Archimedean $v$.

        \section{Weil representation}
        \label{sec:Eisenstein:Weil_rep}

            \subsection{Weil index} 
            \label{ssec:Eisenstein:Weil_rep:Weil_index}
                We recall \emph{Weil indices}, which are certain constants appearing in the Weil representation and other calculations below. We compute the instances which we need.

Suppose $F^+_v$ is a local field (arbitrary characteristic) with nontrivial additive character $\psi_v \colon F^+_v \ra \C^\times$, and suppose $V_v$ is a (finite dimensional) $F^+_v$ vector space equipped with a non-degenerate quadratic form $Q(-)$. 
The map $V_v \ra \C^{\times}$ given by $x \mapsto \psi_v(Q(x))$
is a ``non-degenerate character of the second degree'' in the sense of \cite{Weil64} \cite[{Appendix}]{Rao93}, so there is an associated \emph{Weil index} $\g_{\psi_v}(V_v) \in \C^\times$ (which is an eighth root of unity). The quantity $\gamma_{\psi_v}(V_v)$ depends only on $\psi_v$ and the isomorphism class of $V_v$, and we have
    \begin{equation}
    \gamma_{\overline{\psi}_v}(V_v) = \overline{\gamma_{\psi_v}(V_v)} \quad \quad \gamma_{\psi_v}(V_v \oplus V'_v) = \gamma_{\psi_v}(V_v) \gamma_{\psi_v}(V'_v)
    \end{equation}
for orthogonal direct sums $V_v \oplus V'_v$ (follows from the definition, see \cite[{Theorem A.2}]{Rao93}).
The Weil index also satisfies a global product formula \cite[{Proposition 5}]{Weil64}.

When $F^+_v$ has characteristic $\neq 2$ and $V_v$ has a bilinear pairing $(-,-)$, our convention is that $x \mapsto (x,x)$ is the associated quadratic form (and vice-versa).

\begin{lemma}\label{lemma:Weil_index:quadratic}
Let $F^+_v$ be a local field of characteristic $\neq 2$, let $\psi_v \colon F^+_v \ra \C^\times$ be a nontrivial additive character, and let $V_v$ be a finite dimensional $F^+_v$ vector space with non-degenerate bilinear pairing. Assume any of the following situations hold.
\begin{enumerate}[(1)]
    \item The bilinear pairing on $V_v$ is given by 
    \begin{equation}
    \begin{pmatrix}
    0 & 1_d \\
    1_d & 0
    \end{pmatrix}.
    \end{equation}
    \item The field $F^+_v$ is non-Archimedean with residue characteristic $\neq 2$, there exists a self-dual lattice in $V_v$, and $\psi_v$ is unramified.
\end{enumerate}
Then the Weil index is $\gamma_{\psi_v}(V_v) = 1$.
\end{lemma}
\begin{proof}
(1) By compatibility with orthogonal direct sums, we reduce to the case $d = 1$. Given a nonzero element $a \in F_v^{+ \times}$, we use the temporary notation $\gamma_{\psi_v}(a)$ for the Weil index of the one-dimensional quadratic space containing an element $x$ with $(x,x) = a$. We have $\gamma_{\psi_v}(V_v) = \gamma_{\psi_v}(a) \gamma_{\psi_v}(-a^{-1})$ for some $a \in F_v^{+ \times}$. We have $\gamma_{\psi_v}(a) \gamma_{\psi_v}(-a^{-1}) = 1$ (follows from \cite[{Theorem A.4}]{Rao93}, which relates Weil indices and the Hilbert symbol).

(2) By compatibility with orthogonal direct sums, it is enough to show $\gamma_{\psi_v}(a) = 1$ for $a \in \mc{O}_{F^+_v}^{\times}$. This follows from \cite[{Proposition A.11}]{Rao93}.
\end{proof}

\begin{remark}
The explicit formula of \cite[{Proposition A.12}]{Rao93} shows that Lemma \ref{lemma:Weil_index:quadratic}(2) is false if $F^+_v = \Q_2$ (e.g. if $V_v$ has rank one).
\end{remark}

Next, let $F_v$ be an \'etale algebra of degree $2$ over a local field $F^+_v$ of characteristic $\neq 2$ (residue characteristic $2$ allowed). Write $\eta_v \colon F^{+ \times}_v \ra \{ \pm 1 \}$ for the quadratic character associated to $F_v / F^+_v$ (trivial if $F_v/F^+_v$ is split), and write $a \mapsto \overline{a}$ for the nontrivial involution of $F_v$ over $F^+_v$. If $F^+_v$ is non-Archimedean, we write $\mf{d}$ (resp. $\Delta$) for the different (resp. discriminant) ideal for the extension $F_v/F^+_v$ (where $\mf{d} = \mc{O}_{F_v}$ and $\Delta = \mc{O}_{F^+_v}$ in the split case). We sometimes abuse notation and write $\mf{d}$ and $\Delta$ for understood/chosen generators of these ideals. We write $q_v$ for the residue cardinality of $F^+_v$ if $F^+_v$ is non-Archimedean.

Any non-degenerate $F_v/F^+_v$ Hermitian space $V_v$ with pairing $(-,-)$
has an associated $F^+_v$-bilinear pairing $\frac{1}{2} \mrm{tr}_{F_v/F^+_v}(-,-)$ and quadratic form $x \mapsto \frac{1}{2} \mrm{tr}_{F_v/F^+_v}(x,x)$. (Elsewhere, we typically normalize the trace bilinear pairing without the factor of $\frac{1}{2}$.) We write $\g_{\psi_v}(V_v)$ for the Weil index of this quadratic space with respect to a nontrivial additive character $\psi_v \colon F^+_v \ra \C^{\times}$. We know $\gamma_{\psi_v}(V_v)^4 = 1$ (see e.g. \cite[{Corollary A.5(4)}]{Rao93} and \cite[{Theorem 3.1}]{Kudla94}).

We write $\gamma_{\psi_v}(F_v)$ for the Weil index associated to the one-dimensional Hermitian space $F_v$ with pairing $(x,y) = \overline{x} y$. We write $\epsilon_v(s, \xi_v, \psi_v)$ for the local epsilon factor associated to a quasi-character $\xi_v \colon F^{+ \times}_v \ra \C^{\times}$ (as in \cite[{\S 3}]{Tate79}\cite{Tate67Thesis}, for the quasi-character $\xi_v |-|^s$ and the self-dual Haar measure for $\psi_v$).

If $F^+_v$ is non-Archimedean with uniformizer $\varpi_0$, we have
    \begin{equation}\label{equation:Eisenstein:Weil_rep:Weil_index:epsilon_factors_and_Weil_index}
    \epsilon_v(s, \eta_v, \psi_v) = |\varpi_0^{c(\psi_v)} \Delta|_{F^+_v}^{s - 1/2} \gamma_{\psi_v}(F_v)
    \end{equation}
where
    \begin{equation}\label{equation:Eisenstein:Weil_rep:Weil_index:additive_character_conductor}
    c(\psi_v) = \max \{j \in \Z : \psi_v|_{\varpi_0^{-j} \mc{O}_{F^+_v}} \text{ is trivial} \}.
    \end{equation}
If $F^+_v$ is Archimedean, we have
    \begin{equation}
    \epsilon_v(s, \eta_v, \psi_v) = |u|_{F^+_v}^{s - 1/2} \gamma_{\psi_v}(F_v)
    \end{equation}
where $u \in F^{+ \times}_v$ is such that
    \begin{equation}\label{equation:Eisenstein:Weil_rep:Weil_index:Archimedean}
    \psi_v(x) = e^{2 \pi i u x} \quad \text{if $F^+_v = \R$} \quad \text{and} \quad \psi_v(z) = e^{2 \pi i \mrm{tr}_{\C/\R} (u z)} \quad \text{if $F^+_v = \C$}.
    \end{equation}
In both non-Archimedean and Archimedean cases, if $\psi'_v(x) = \psi_v(ux)$ for some $u \in F^{+ \times}_v$,
we have
    \begin{equation}\label{equation:Eisenstein:Weil_rep:Weil_index:change_of_additive_character}
    \epsilon_v(s, \eta_v, \psi'_v) = \eta_v(u) |u|^{s - 1/2} \epsilon(s, \eta_v, \psi_v)
    \quad\quad
    \gamma_{\psi'_v}(F_v) = \eta_v(u) \gamma_{\psi_v}(F_v).
    \end{equation}
These identities follow from \cite[{Lemma 1.2(iii),(iv)}]{JL70} and properties of epsilon factors
as in \cite[{\S 3.2}]{Tate79}. For the reader's convenience, we recall $\gamma_{\psi_v}(\C) = i$ if $F^+_v = \R$ and $\psi_v(x) = e^{2 \pi i x}$.

In all cases, we have
    \begin{equation}
    \gamma_{\psi_v}(F_v)^2 = \epsilon_v(1/2, \eta_v, \psi_v)^2 = \eta_v(-1).
    \end{equation}
If $F^+_v$ is non-Archimedean, recall that $\epsilon_v(s,\xi_v, \psi_v) = 1$ if $\xi_v$ and $\psi_v$ are unramified. If $F^+_v = \R$ and $\psi_v(x) = e^{2 \pi i x}$, recall $\epsilon_v(s, \mrm{sgn}^j, \overline{\psi}_v) = 1$ (resp. $= -i$) if $j$ is even (resp. odd) where $\mrm{sgn} \colon \R^{\times} \ra \{ \pm 1\}$ is the sign character (these formulas will be used implicitly in Section \ref{ssec:Eisenstein:local_functional_equations:Arch}).
Recall our convention that self-duality for Hermitian lattices is understood with respect to the trace pairing (unless otherwise specified), see \crefext{I:ssec:Hermitian_conventions:lattices}.
If $F_v/F^+_v$ is ramified and $L$ is a self-dual Hermitian lattice, then $L$ must have even rank
(explained in loc. cit.).

By the \emph{standard additive character} $\psi_{\mrm{tr},v} \colon F^+_v \ra \C^{\times}$
for Archimedean $F^+_v$, we mean the characters of
\eqref{equation:Eisenstein:Weil_rep:Weil_index:Archimedean}
with $u = 1$.
By the \emph{standard additive character}
for non-Archimedean $F^+_v$, we mean $\psi_{\mrm{tr}}(x) \coloneqq \psi_p(\operatorname{tr}_{F^+_v / \Q_p}(x))$,
where the residue characteristic is $p$,
and where $\psi_p \colon \Q_p^{\times} \ra \C^{\times}$ is the unramified
character satisfying $\psi_p(x) = e^{-2 \pi i x}$ for any $x \in \Z[1/p]$.
In the non-Archimedean case, we write $\mf{d}_{F^+_v / \Q_p}$ for the different ideal
of $F^+_v / \Q_p$. An element $u \in F^+_v$
is a generator for $\mf{d}_{F^+_v / \Q_p}$
if and only if the additive character
$x \mapsto \psi_{\mrm{tr},v}(u^{-1} x)$ is unramified.

\begin{lemma}\label{lemma:Weil_index:Hermitian}
Let $F^+_v$ be a local field of characteristic $\neq 2$, let $\psi_v \colon F^+_v \ra \C^\times$ be a nontrivial additive character, and let $F_v/F^+_v$ be a degree $2$ \'etale algebra. Let $V_v$ be a finite dimensional non-degenerate $F_v/F^+_v$ Hermitian space. Assume any of the following situations hold.
\begin{enumerate}[(1)]
    \item The Hermitian space $V_v$ admits a basis with Gram matrix
    \begin{equation}
    \begin{pmatrix}
    0 & 1_d \\
    1_d & 0
    \end{pmatrix}.
    \end{equation}
    \item We have $F_v \cong F^+_v \times F^+_v$.
    \item The extension $F_v/F^+_v$ is unramified or $F^+_v$ has residue characteristic $\neq 2$. Moreover, the field $F^+_v$ is non-Archimedean, there exists a full-rank self-dual $\mc{O}_{F_v}$-lattice in $V_v$, and $V_v$ has even rank if $F_v / F^+_v$ is unramified.
    \item The field $F^+_v$ is non-Archimedean, the extension $F_v/F^+_v$ is unramified, there exists a full-rank self-dual lattice in $V_v$, and $\psi_v$ is unramified.
    \item The field $F^+_v$ is non-Archimedean, we have $\psi_v = \psi_{\mrm{tr},v}$
        being the standard additive character, and $V_v$ contains a full-rank $\mc{O}_{F_v}$-lattice
        $M$ which is self-dual for the $\Z_p$-bilinear
        pairing $x,y \mapsto \operatorname{tr}_{F_v / \Q_p}(x,y)$.
        Also assume that at least one of the following three conditions is satisfied:
        (i) The extension $F_v / F^+_v$ is unramified (ii) The residue characteristic of $F^+_v$ is 
        $\neq 2$ (iii) For some generator $u$ of $\mf{d}_{F^+_v / \Q_p}$,
        the Hermitian space $V_v$ with the modified pairing $u \cdot (-,-)$
        admits a basis with Gram matrix as in situation (1).
\end{enumerate}
Then the Weil index is $\gamma_{\psi_v}(V_v) = 1$.
\end{lemma}
\begin{proof}
We have (3) $\implies$ (1) (see \cite[{Lemma 2.12}]{LL22II} for the ramified situation, in which case $V_v$ automatically has even rank). This implication is false if $F_v/F^+_v$ is ramified with $F^+_v$ of residue characteristic $2$.

In situations (1) and (2) we may pick a basis $\{1, \a \}$ for $F_v$ as an $F^+_v$ vector space where $\mrm{tr}_{F_v/F^+_v}(\a) = 0$. Applying Lemma \ref{lemma:Weil_index:quadratic} proves the claims.

In situation (4), we may diagonalize the given self-dual lattice, hence reducing to the case where $V_v$ has rank one. In this case, we have $\gamma_{\psi_v}(V_v) = \gamma_{\psi_v}(F_v) = \epsilon(1/2, \eta_v, \psi_v) = 1$.

In situation (5), let $u \in \mf{d}_{F^+_v / \Q_p}$ be any generator
(in situation (5)(iii), take $u$ to be the given one). Take $\psi'_v(x) \coloneqq \psi_v(u^{-1} x)$,
and take $V'_v \coloneqq V_v$ but give $V'_v$ the modified pairing $(-,-)' \coloneqq u \cdot (-,-)$.
We find $\gamma_{\psi_v}(V_v) = \gamma_{\psi'_v}(V'_v)$,
so this reduces to situations (4), (3), and (1) above
(note that $M$ viewed as a lattice in $V'_v$ is self-dual for $(-,-)'$.)
\end{proof}
    
            \subsection{Weil representation} 
            \label{ssec:Eisenstein:Weil_rep:Weil_rep}
                Let $F_v/F^+_v$ and accompanying notation be as in Section \ref{ssec:Eisenstein:Weil_rep:Weil_index}.
Assume $F_v/F^+_v = \C / \R$ if $F^+_v$ is Archimedean. We also assume $F^+_v$ has characteristic $0$ (because \cite{Kudla94} assumes this).

Let $V_v$ be a non-degenerate $F_v/F^+_v$ Hermitian space of dimension $n \geq 0$. Choose a nontrivial additive character $\psi_v \colon F^+_v \ra \C^\times$, and let $\chi_v \colon F_v^{\times} \ra \C^{\times}$ be a character such that $\chi_v|_{F^{+ \times}_v} = \eta_v^{n}$. There is a local \emph{Weil representation} $\omega_v = \omega_{\chi_v,\psi_v}$ of $U(m,m)(F^+_v) \times U(V_v)(F^+_v)$ on the space of Schwartz functions $\mc{S}(V_v^m)$ (the Schr\"odinger model \cite[{\S 5}]{Kudla94}),\footnote{
    We explain how to extract our formulas from Kudla's.
    Our $F_v / F^+_v$ is Kudla's $E / F$ (with $D = E$),
    and we are in Kudla's ``Case $3_{+}$''
    in the notation of \cite[{\S 2}]{Kudla94}, so $\epsilon = +1$ in Kudla's notation.
    Kudla's $V$ is our $V_v$, and our $U(m,m)$ has a natural right action on Kudla's $W$
    (in the basis specified on \cite[{Page 364}]{Kudla94})
    via right multiplication with $W$ viewed as row vectors in loc. cit.
    Our $m$ is Kudla's $n$.
    Our $\chi_v$ is Kudla's $\xi$.
    We have $\kappa = 1/2$ in Kudla's notation.
    Our $\psi$ is Kudla's $\eta$,
    which is also Kudla's $u \mapsto \psi(\frac{1}{2} u)$.
    Our convention on $\psi$ is compatible with that
    of \cite[{\S 10}]{KR14}.
    (Alternatively, we could take our $\psi$
    to agree with Kudla's $\psi$,
    and instead identify the Hermitian pairing $(-,-)$ on our $V_v$
    with the pairing $\frac{1}{2}(-,-)$ on Kudla's $V$.)

    Our formula for $\omega_v(h) \varphi_v$ corresponds
    to the formula at the bottom of \cite[{Page 400}]{Kudla94}.
    The remaining formulas we give are all specializations of the integral formula
    for $\omega_V(g) \varphi(x)$ in \cite[{\S 5}]{Kudla94},
    as we explain below (our $\underline{x}$ is Kudla's $x$).
    To be historically accurate, we remark that Kudla's integral formula arises
    from combining Kudla's splitting of a metaplectic cover
    (also studied by Gelbart and Rogawski in the unitary case
    \cite{GR91}) with the Schr\"odinger model of the Weil representation
    for metaplectic groups, as studied by Weil \cite{Weil64},
    though we follow Kudla in using the integral formulas from Rao \cite{Rao93}.

    \vspace{\maxdimen}

    First consider the action of $\omega_v(m(a))$, in our notation.
    To calculate Kudla's $\beta_V(m(a))$ from \cite[{Theorem 3.1}]{Kudla94},
    note that $j = 0$ in Kudla's notation, so $\beta_V(m(a)) = \xi(x(m(a))) = \chi_v(\det a)$
    (first two expressions use Kudla's notation), because $x(m(a)) = \det a$.
    In Kudla's formula for $\omega_V(g) \varphi(x)$, we take $g = m(a)$,
    in which case $b = c = 0$ (in Kudla's notation) and Kudla's $a$ is our $a$.
    Since $c = 0$, we have $V_g^n = 0$ in Kudla's notation.
    In Kudla's integral expression,
    note that only the term $\varphi(x a)$ is nonzero in this case.
    To compute Kudla's integral expression, we must determine the measure
    $\mu_{\sigma}$ in his notation, for which Kudla directs the reader to Rao's 
    paper \cite{Rao93}.
    In the notation of \cite[{Theorem 3.5}]{Rao93} and its proof,
    Rao specifies that $\mu_{1}$ is the counting measure $\delta_0$ on $V_g^n$,
    and thus \cite[{(3.11)}]{Rao93} implies that $\mu_{m(a)} = |\det a|^{n / 2}_{F_v} \cdot \delta_0$,
    upon taking $p_1 = 1$ and $p_2 = m(a)$ in Rao's notation:
    note that Rao's $V^*$ corresponds to Kudla's $\mathbb{Y}$,
    and is hence isomorphic to $(F_v^{\oplus m}) \otimes_{F_v} V_v$
    as an $M(F^+_v)$-module, where $m(a)$ acts by the
    right multiplication $(F_v^{\oplus m}) \actedby {}^t \overline{a}^{-1}$.
    Putting this together verifies our claimed formula for $\omega_v(m(a))$.

    Next, consider the action of $\omega_v(n(b))$, in our notation.
    To calculate $\beta_V(n(b))$, we again note that $j = 0$ in this case,
    and we find $\beta_V(n(b)) = \xi(x(n(b))) = \chi_v(1) = 1$.
    In Kudla's notation, we again have $c = 0$, and so again $V^n_g = 0$
    in his notation. We also find that Rao's measure is $\mu_{n(b)} = \mu_1 = \delta_0$
    in this case.
    To finish the comparison with Kudla's formula, we
    note that $\opn{tr}(b \cdot (\underline{x}, \underline{x}))
    = \frac{1}{2} \opn{tr}_{F_v / F^+_v} (\underline{x}, \underline{x} b)$,
    where we use the fact that $b$ is Hermitian.

    Lastly, we consider the action of $\omega_v(w)$ in our notation.
    We first describe the action of $\omega_v(w^{-1})$.
    Our $w^{-1}$ is Kudla's $\tau_n$.
    For $g = \tau_n$, we thus have $\beta_V(g) = \gamma_F(\eta \circ RV)^{-n}$
    in the notation of \cite[{Theorem 3.1}]{Kudla94}.
    Kudla's $\gamma_F(\eta \circ RV)$ is our $\gamma_{\psi_v}(V_v)$,
    since we defined $\gamma_{\psi_v}(V_v)$
    with respect to the quadratic form 
    $x \mapsto \frac{1}{2} \operatorname{tr}_{F_v / F^+_v}(x, x)$
    (Section \ref{ssec:Eisenstein:Weil_rep:Weil_index})
    which is compatible with the convention in \cite[{Lemma 2.4}]{Kudla94}.
    In the notation of Kudla's integral formula,
    we have $V_g^n = V^n$.
    Inspecting Kudla's integral formula for $\omega_V(g) \phi(x)$,
    we find (in our notation)
    $(\omega_v(w^{-1}) \varphi_v)(\underline{x}) = \gamma_{\psi_v}(V_v)^{-m} C \cdot \widehat{\varphi}_v(-\underline{x})$
    for some positive constant $C$,
    arising from the choice of measure in Kudla's integral formula.
    Since $\omega_v$ defines an action,
    we must have $\omega_v(w^{-1}) \circ \omega_v(w^{-1}) = \omega_v(-1)$,
    and the previous computation implies $(\omega_v(-1) \varphi_v)(\underline{x})
    = \chi_v(-1)^{m} \varphi_v(- \underline{x})$.
    Since $\chi_v(-1)$ and $\gamma_{\psi_v}(V_v)$ both have norm $1$,
    this determines the constant $C$: by Fourier inversion,
    we must have $C = 1$, i.e. the measure in Kudla's integral
    must be the self-dual Haar measure from \eqref{equation:Eisenstein:Weil_rep:Weil_rep:Fourier_transform}.
    This reasoning also implies that $\gamma_{\psi_v}(V_v)^2 = \chi_v(-1)$ (consider the special case $m = 1$).
    On the other hand, since $\omega_v(w) = \omega_v(w^{-1}) \circ \omega_v(-1)$,
    we now obtain the claimed formula
    $(\omega_v(w) \varphi_v)(\underline{x}) = \gamma_{\psi_v}(V_v)^m \widehat{\varphi}_v(\underline{x})$.
}
which we normalize as
    \begin{align*}
    (\omega_v(m(a)) \varphi_v)(\underline{x}) & = \chi_v(\det a) |\det a|^{n/2}_{F_v} \varphi_v(\underline{x} \cdot a) && m(a) \in M(F^+_v) \\
    (\omega_v(n(b)) \varphi_v)(\underline{x}) & = \psi_v(\operatorname{tr} (b \cdot (\underline{x}, \underline{x}))) \varphi_v(\underline{x}) && n(b) \in N(F^+_v) \\
    (\omega_v(w) \varphi_v)(\underline{x}) & = \g_{\psi_v}(V_v)^m \widehat{\varphi}_v(\underline{x}) && \\
    (\omega_v(h) \varphi_v)(\underline{x}) & = \varphi_v(h^{-1} \cdot \underline{x}) && h \in U(V_v)(F^+_v)
    \end{align*}
for $\varphi_v \in \mc{S}(V_v^m)$ and $\underline{x} \in V_v^m$ (viewed as $n \times m$ matrices), where
    \begin{equation}\label{equation:Eisenstein:Weil_rep:Weil_rep:Fourier_transform}
    \widehat{\varphi}_v(\underline{x}) = \int_{V_v^m} \varphi_v(\underline{\smash{y}}) \psi_v(\operatorname{tr}_{F_v/F^+_v} \operatorname{tr}(\underline{x}, \underline{\smash{y}})) ~d\underline{\smash{y}}
    \end{equation}
is Fourier transform for the corresponding self-dual Haar measure on $V_v^m$,
and where we should view $\underline{x}$ as a length $m$ row vector with elements in $V_v$
in the above formulas.
The constant $\gamma_{\psi_v}(V_v)$ is the Weil index from Section \ref{ssec:Eisenstein:Weil_rep:Weil_index}.
The symbol $(\underline{x}, \underline{\smash{y}})$ denotes the matrix of inner products
with $i,j$-th entry $(x_i, y_j)$.

\begin{remark}[Change of additive character]\label{remark:Eisenstein:Weil_rep:Weil_rep:rescaling_trick}
We sometimes restrict our attention to Weil representations
$\omega_{\chi_v, \psi_v}$ for convenient choices of $\psi_v$
(e.g. unramified in the non-Archimedean case, or $\psi_v(x) = e^{2\pi i x}$
in the Archimedean case).
For some purposes, this loses no generality by the following rescaling trick.
If $\psi'_v(x) = \psi_v(u \cdot x)$ for some $u \in (F^+_v)^{\times}$,
and if $V'_v$ is $V_v$ but with the Hermitian pairing $u^{-1} (-,-)$,
then the Weil representations of
$\omega_{\chi_v, \psi_v}$ on $\mc{S}(V_v^m)$
and $\omega_{\chi_v, \psi'_v}$ on $\mc{S}(V_v^{\prime m})$
are compatible for the canonical identification $V_v \cong V_v^{\prime}$
as $F_v$-modules (not necessarily preserving Hermitian pairings).
\end{remark}

With $s_0 \coloneqq (n - m)/2$, there is a map $\mc{S}(V_v^m) \ra I(\chi_v, s_0)$ sending $\varphi_v \in \mc{S}(V_v^m)$ to the function $h \mapsto (\omega_v(h) \varphi_v)(0)$. The associated standard section $\Phi_{\varphi_v} \in I(\chi_v,s)$ is the \emph{Siegel--Weil section} for $\varphi_v$ \cite[{\S 5.1}]{GS19}.

If $F^+_v$ is non-Archimedean, choose a generator $\mf{d}$ of the different ideal of $F_v / F^+_v$, and let $M^{\circ}_2$ be the rank $2$ Hermitian $\mc{O}_{F_v}$-lattice admitting a basis with Gram matrix
    \begin{equation}\label{equation:Eisenstein:Weil_rep:Weil_rep:standard_self-dual_lattice}
    \begin{pmatrix}
    0 & \mf{d}^{-1} \\
    \overline{\mf{d}}^{-1} & 0
    \end{pmatrix}.
    \end{equation}
Note that $M_2^{\circ} = M_2^{\circ *}$ is self-dual (with respect to the $F^+_v$-bilinear pairing $\operatorname{tr}_{F_v/F^+_v}(-,-)$).

\begin{lemma}\label{lemma:spherical_Siegel-Weil_section}
In the situation above, assume moreover that $\chi_v$ is unramified, and that $F^+_v$ is non-Archimedean. Suppose $\varphi_v = \pmb{1}_{M}^{m}$ where $\pmb{1}_{M}$ is the characteristic function of a full rank $\mc{O}_{F_v}$-lattice $M \subseteq V_v$ in any of the following situations.
    \begin{enumerate}[(1)]
    \item The additive character $\psi_v$ is unramified, and the lattice $M$ is self-dual. Moreover, the extension $F_v/F^+_v$ is unramified, or $F^+_v$ has residue characteristic $\neq 2$.
    
    \item The additive character $\psi_v$ is unramified, and we have $M \cong (M^{\circ}_2)^{\oplus d}$ (orthogonal direct sum) for some $d \geq 0$.
    
    \item We have $\psi_v = \psi_{\mrm{tr}, v}$ being the standard additive character, and $M$ is self-dual for the $\Z_p$-bilinear
        pairing $x,y \mapsto \operatorname{tr}_{F_v / \Q_p}(x,y)$.
        Also assume that at least one of the following three conditions is satisfied:
        (i) The extension $F_v / F^+_v$ is unramified (ii) The residue characteristic of $F^+_v$ is 
        $\neq 2$ (iii) For some generator $u$ of $\mf{d}_{F^+_v / \Q_p}$,
        the lattice $M$ with the modified pairing $u \cdot (-,-)$
        is isomorphic to $(M^{\circ}_2)^{\oplus d}$ for some $d \geq 0$.
    \end{enumerate} 
Then the associated Siegel--Weil section $\Phi_{\varphi_v}$ is the normalized spherical section $\Phi_v^{\circ}$, i.e. $K_v^{\circ}$-fixed with $\Phi_{\varphi_v}(1) = 1$.
\end{lemma}
\begin{proof}
If either (1) or (2) holds, the lemma follows from the explicit formulas above, since $w$ and $P(\mc{O}_{F^+_v})$ generate $K_v^{\circ} = U(m,m)(\mc{O}_{F^+_v})$ and since the Weil index $\gamma_{\psi_v}(V_v)$ is $1$ (Lemma \ref{lemma:Weil_index:Hermitian}),
where we use the fact that any self-dual $M$ must have even rank if $F_v / F^+_v$ is ramified
with residue characteristic $2$.

In case (3), pick any generator $u$ of $\mf{d}_{F^+_v / \Q_p}$
(in situation (3)(iii), pick the $u$ that is given).
Consider the modified pairing $(-,-)' \coloneqq u \cdot (-,-)$
on $M$. Then $\psi'_v(x) \coloneqq \psi_{\mrm{tr},v}(u^{-1} x)$ is unramified,
and the lattice $M$ with pairing $(-,-)'$ satisfies
either condition (1) or (2).
The claim now follows from
Remark \ref{remark:Eisenstein:Weil_rep:Weil_rep:rescaling_trick}.
\end{proof}

In the statement of Lemma \ref{lemma:spherical_Siegel-Weil_section},
condition (1) implies condition (2) 
if $M$ has even rank (the ramified case is \cite[{Lemma 2.12}]{LL22II}).

Next, consider the case where $F_v/F^+_v = \C/\R$. Suppose the $n$-dimensional Hermitian space $V_v$ is positive definite, with Hermitian pairing $(-,-)$.
If $\psi_v(x) = e^{2 \pi i u x}$ for some $u > 0$, the Gaussian
    \begin{equation}
    \varphi_v(\underline{x}) = e^{-2 \pi u \cdot \mrm{tr}(\underline{x}, \underline{x})} \in \mc{S}(V_v^m)
    \end{equation}
for $\underline{x} = (x_1, \ldots, x_m) \in V_v^m$ has associated Siegel--Weil section
    \begin{equation}
    \Phi_{\varphi_v} = \Phi_v^{(n)}
    \end{equation}
where $\Phi_v^{(n)}$ is the scalar weight standard section described surrounding \eqref{equation:adelic_v_classical:scalar_weight}, see \cite[{(2.68)}]{GS19}
(there in the case $u = 1$, and the general case follows from Remark \ref{remark:Eisenstein:Weil_rep:Weil_rep:rescaling_trick}).

\begin{remark}\label{remark:Weil_rep:Weil_rep:chi_independence}
Suppose $F / F^+$ is a CM extension of number fields with associated quadratic character $\eta$ and accompanying notation as in Section \ref{ssec:Eisenstein:setup:adelic_v_classical}. With $m$ and $n$ as above, choose any character $\chi \colon F^{\times} \backslash \A_F^{\times} \ra \C^{\times}$ satisfying $\chi|_{\A^{\times}} = \eta^n$. Choose nontrivial additive characters $\psi_v \colon F^+_v \ra \C^{\times}$ for each place $v$ (the $\psi_v$ need not come from a global character). Suppose we are given a collection of local Weil representations $\omega_{\chi_v, \psi_v}$ on some $\mc{S}(V_v^m)$ for each place $v$ of $F^+_v$ (where the collection $(V_v)_v$ of local Hermitian spaces need not come from a global Hermitian space). Choose $\varphi_v \in \mc{S}(V_v^m)$ for each place $v$, and assume $\varphi_v = \pmb{1}_{L_v}^m$ for some full-rank self-dual lattice $L_v \subseteq V_v$ for all but finitely many $v$. Set $\Phi \coloneqq \bigotimes_v \Phi_{\varphi_v}$.

In this situation, the Eisenstein series variant $\tilde{E}(a, s, \Phi)_n$ \eqref{equation:setup:adelic_v_classical:variant} does not depend on the choice of $\chi$. This follows upon inspecting the Weil representation, particularly the action of $m(a)$.

This remark also has a local version, i.e. the Whittaker function variants $\tilde{W}^*_{T,\psi_v}(a, s)^{\circ}_n$ and $\tilde{W}^*_{T,\psi_v}(a, s, \Phi_{\varphi_v})_n$ (Sections \ref{ssec:Eisenstein:local_Whittaker:Archimedean} and \ref{ssec:Eisenstein:local_Whittaker:non-Arch}) do not depend on the choice of $\chi_v$.
\end{remark}

        \section{Local Whittaker functions} 
        \label{sec:local_Whittaker}
            Let $F_v/F^+_v$ and accompanying notation be as in Section \ref{ssec:Eisenstein:Weil_rep:Weil_rep}. If $F^+_v$ has residue characteristic $2$, we also assume $F_v / F^+_v$ is unramified.
Let $\chi_v \colon F_v^{\times} \ra \C^{\times}$ be a character satisfying $\chi_v|_{F^{+ \times}_v} = \eta_v^n$ for some integer $n \in \Z$, with $n$ even if $F_v / F^+_v$ is ramified in the non-Archimedean case.\footnote{
    In the non-Archimedean case, for $F_v / F^+_v$ ramified and $n$ odd, there is no spherical standard section $\Phi^{\circ}_v$ in $I(\chi_v, s)$,
    so the setup below would need modification.
}
Assume $\chi_v$ is unramified if $F^+_v$ is non-Archimedean.
Let $\psi_v \colon F^+_v \ra \C^{\times}$ be an unramified nontrivial additive character. Assume $\psi_v(x) = e^{2 \pi i x}$ if $F^+_v = \R$. These are our default hypotheses, but weaker hypotheses often suffice (as will be indicated below).

Let $\Phi^{\circ}_v \in I(s, \chi_v)$ be the normalized spherical standard section if $F^+_v$ is non-Archimedean. Let $\Phi^{(n)}_v \in I(s, \chi_v)$ be the normalized scalar weight standard section from Section \ref{ssec:Eisenstein:setup:adelic_v_classical} if $F_v / F^+_v = \C / \R$.

Given an integer $m \geq 0$ (we do not assume $m \leq n$, unless otherwise specified) and given $T \in \mrm{Herm}_m(F^+_v)$ with $\det T \neq 0$, we define \emph{normalized local Whittaker functions}
    \begin{align}
    W_{T,v}^*(h, s)^{\circ}_n & \coloneqq \Lambda_{T,v}(s)^{\circ}_n W_{T,v}(h, s, \Phi^{\circ}_v) && \text{for $F^+_v$ non-Archimedean}
    \\
    W_{T,v}^*(h, s)^{\circ}_n & \coloneqq \Lambda_{T,v}(s)^{\circ}_n W_{T,v}(h, s, \Phi^{(n)}_v) && \text{for $F^+_v$ Archimedean}
    \end{align}
for certain normalizing factors $\Lambda_{T,v}(s)^{\circ}_n$ (see \eqref{equation:Eisenstein:local_Whittaker:non-Arch:normalizing_factor} and \eqref{equation:Eisenstein:local_Whittaker:Archimedean:normalizing_factor} below). 

The preceding normalization gives $W_{T,v}^*(h, s)^{\circ}_n$ a clean functional equation (Section \ref{sec:Eisenstein:local_functional_equations}).
Moreover, the normalized function $W^*_{T,v}(h, s)^{\circ}_n$ (as opposed to the unnormalized versions) seems to correspond more naturally to local information about special cycles (e.g. local contributions to arithmetic degrees) in arithmetic (and non-arithmetic) Siegel--Weil formulas. For example, our main local theorems \crefext{I:sec:non-Arch_identity} \and \crefext{II:sec:Archimedean_identity} are proved in terms of the derivative of $W^*_{T,v}(1, s)^{\circ}_n$ and not $W_{T,v}(1, s, \Phi^{\circ}_v)$ or $W_{T,v}(1, s, \Phi^{(n)}_v)$.

The normalizing factors $\Lambda_{T,v}(s)^{\circ}_n$ also carry geometric information. For example, consider an imaginary quadratic field $F / \Q$ of odd discriminant, suppose $m = n$ is even, and form the product $2 \prod_v \Lambda_{T,v}(s)^{\circ}_n$ over all places $v$ of $\Q$. If $n \equiv 0 \pmod{4}$, evaluation at $s = 0$ returns the degree of a certain $0$-dimensional unitary complex Shimura variety (stack), giving a case of a unitary analogue of the Siegel mass formula. If $n \equiv 2 \pmod{4}$, evaluation at $s = 0$ returns the volume of a certain $(n - 1)$-dimensional unitary complex Shimura variety (stack). These volume identities will be discussed in Section \ref{ssec:geometric_Siegel-Weil:complex_volumes} (but are not needed for our main theorems on arithmetic Siegel--Weil).
    
            \subsection{Local \texorpdfstring{$L$}{L}-factors} 
            \label{ssec:local_Whittaker:local_L-factors}
                We use the following (standard) local factors as in \cite[{\S 3}]{Tate79}.

If $F^+_v$ is a local field (allowing arbitrary characteristic in Section \ref{ssec:local_Whittaker:local_L-factors}) and $\xi_v \colon F^{+ \times}_v \ra \C^{\times}$ is a quasi-character, we write $L_v(s, \xi_v)$ for the corresponding local $L$-factor (for the quasi-character $\xi_v |-|_{F^+_v}^s$). Given any nontrivial additive character $\psi_v \colon F^+_v \ra \C^{\times}$, we write $\epsilon_v(s, \xi_v, \psi_v)$ for the corresponding local epsilon factor (as appeared in Section \ref{ssec:Eisenstein:Weil_rep:Weil_index}) and $\rho_v(s, \xi_v, \psi_v)$ for the local factor from Tate's thesis \cite[{Theorem 2.4.1}]{Tate67Thesis}, which satisfies
    \begin{equation}
    \rho_{v}(s, \xi_v, \psi_v) = \epsilon_v(s, \xi_v, \psi_v)^{-1} L_v(1-s, \xi_v^{-1})^{-1} L_v(s, \xi_v).
    \end{equation}

If $F^+$ is a global field with a quasi-character $\xi \colon F^{+ \times} \backslash \A_{F^+} \ra \C^{\times}$ and nontrivial additive character $\psi \colon F^+ \backslash \A_{F^+} \ra \C^{\times}$, we write
    \begin{equation}
    \Lambda(s, \xi) = \prod_v L_v(s, \xi_v) \quad \quad L(s, \xi) = \prod_{v < \infty} L_v(s, \xi) \quad \quad \epsilon(s, \xi) = \prod_v \epsilon_v(s, \xi_v, \psi_v)
    \end{equation}
and have $\Lambda(s, \xi) = \epsilon(s, \xi) \Lambda(1-s, \xi^{-1})$.
For the reader's convenience, we recall the formulas
    \begin{align*}
    L_v(s, \xi_v) & = 
    \begin{cases}    
    (1 - \xi_v(\varpi_0) |\varpi_0|^s_{F^+_v})^{-1} & \text{if $\xi_v$ is unramified} \\
    1 & \text{if $\xi_v$ is ramified}
    \end{cases}
    \quad \quad 
    \parbox[t]{0.3\textwidth}{if $F^+_v$ is non-Archimedean with uniformizer $\varpi_0 \in F^+_v$}
    \\
    L_v(s, \mrm{sgn}^j) & = 
    \begin{cases}
    \pi^{-s / 2} \Gamma(s/2) & \text{if $j$ is even} \\
    \pi^{-(s+1)/2} \Gamma((s+1)/2) & \text{if $j$ is odd}
    \end{cases}
    \quad \quad 
    \parbox[t]{0.3\textwidth}{if $F^+_v = \R$ and $\mrm{sgn}$ denotes the sign character.}
    \end{align*}

            \subsection{Normalized Archimedean Whittaker functions} 
            \label{ssec:Eisenstein:local_Whittaker:Archimedean}
                With notation as above, assume $F_v/F^+_v$ is $\C / \R$. The symbol $h$ will denote an element of $U(m,m)(F^+_v)$. Fix integers $n,m$ with $m \geq 0$.

Let $\psi_v \colon F^+_v \ra \C^{\times}$
be an arbitrary nontrivial additive character,
which is necessarily of the form $\psi_v(x) = e^{2\pi i u x}$
for some $u \in \R^{\times}$.
Consider $T \in \mrm{Herm}_m(F^+_v)$ with $\det T \neq 0$. With $s_0 \coloneqq (n-m)/2$ as above, we define the normalizing factor
    \begin{equation}\label{equation:Eisenstein:local_Whittaker:Archimedean:normalizing_factor}
    \Lambda_{T,\psi_v}(s)^{\circ}_n \coloneqq \frac{(2 \pi)^{m(m-1)/2}}{(-2 \pi i)^{nm}} \pi^{m(- s + s_0)} \left ( \prod_{j = 0}^{m - 1} \Gamma(s - s_0 + n - j) \right ) |\det T|^{- s - s_0}_{F^+_v} |u|_{F^+_v}^{-m (s + s_0) - m^2/2}
    \end{equation}
(compare \cite[{(3.3.14)}]{GS19}, also Shimura \cite{Shimura82}) where $\Gamma$ is the usual gamma function.
If $\psi_v(x) = e^{2 \pi i x}$ is the standard additive character, we
write $\Lambda_{T,v}(s)^{\circ}_n \coloneqq \Lambda_{T,\psi_v}(s)^{\circ}_n$.
For general $\psi_v(x) = e^{2\pi i u x}$, we have
$\Lambda_{T,\psi_v}(s)^{\circ}_n = |u|_{F^+_v}^{-m^2/2} \Lambda_{u T,v}(s)^{\circ}_n$.

For arbitrary $\psi_v$, we define a \emph{normalized Archimedean Whittaker function}
    \begin{equation}
    W^*_{T,\psi_v}(h, s)^{\circ}_n \coloneqq \Lambda_{T,\psi_v}(s)^{\circ}_n W_{T,\psi_v}(h, s, \Phi^{(n)}_v).
    \end{equation}
For $a \in \GL_m(F_v)$, we also consider the variant
    \begin{equation}
    \tilde{W}^*_{T,\psi_v}(a, s)^{\circ}_n \coloneqq \chi_v(a)^{-1} |\det a|_{F_v}^{-n / 2} W^*_{T,\psi_v}(m(a), s)^{\circ}_n \cdot q^{-T} \quad \quad q^{-T} \coloneqq \psi_v(-\operatorname{tr}(i T y))
    \end{equation}
with $y \coloneqq a {}^t \overline{a}$ (temporary notation). This is a (normalized) local analogue of \eqref{equation:setup:Fourier_and_Whittaker:a_Fourier}.
When $\psi_v(x) = e^{2 \pi i x}$ is the standard additive character,
we write $W^*_{T,v}$ and $\tilde{W}^*_{T,v}$ instead of $W^*_{T,\psi_v}$ and $\tilde{W}^*_{T,\psi_v}$.
For general $\psi_v(x) = e^{2\pi i u x}$, 
we have $W^*_{T,\psi_v} = W^*_{u T, v}$ and $\tilde{W}^*_{T, \psi_v} = \tilde{W}^*_{u T, v}$
(consider Remark \ref{remark:Eisenstein:setup:Fourier_and_Whittaker:change_of_additive_character}).

For any $a \in \GL_m(F_v)$ and $k \in U(m)$, we have the ``linear invariance'' properties
    \begin{equation}\label{equation:local_Whittaker:Archimedean:linear_invariance}
    \tilde{W}^*_{T,\psi_v}(a, s)^{\circ}_n = \tilde{W}_{{}^t \overline{a} T a, \psi_v}^*(1, s)^{\circ}_n \quad \quad 
    \tilde{W}^*_{T,\psi_v}(1,s)^{\circ}_n = \tilde{W}^*_{T,\psi_v}(k,s)^{\circ}_n.
    \end{equation}
The first expression follows from \eqref{equation:setup:Fourier_and_Whittaker:Whittaker_linear_change}, and the second expression follows from the scalar weight property of $\Phi^{(n)}_v$.
Given $y \in \mrm{Herm}_m(\R)_{>0}$, we also set $W^*_{T,\psi_v}(y,s)^{\circ}_n \coloneqq \tilde{W}^*_{T,\psi_v}(a, s)^{\circ}_n$
for any $a \in \GL_m(F_v)$ satisfying $a {}^t \overline{a} = y$ (does not depend on the choice of $a$).\footnote{With this notation, there is possible ambiguity for the meaning of $W^*_{T,\psi_v}(1,s)^{\circ}_n$, which could refer to either $W^*_{T,\psi_v}(h,s)^{\circ}_n$ or $W^*_{T,\psi_v}(y,s)^{\circ}_n$ evaluated at $h = 1_{2m}$ or $y = 1_m$. To avoid confusion, we will avoid the symbol $W^*_{T,\psi_v}(1,s)^{\circ}_n$ when $v$ is Archimedean.}
We use the shorthand $W_{T,\psi_v}^*(s)^{\circ}_n \coloneqq \tilde{W}^*_{T,\psi_v}(1,s)^{\circ}_n$.
We continue to use our usual convention
of using the subscript ``$v$'' instead of ``$\psi_v$''
if $\psi_v$ is the standard character.

To simplify exposition,
assume $\psi_v(x) = e^{2 \pi i u x}$ with $u > 0$ for the moment.
For all $n \in \Z$, we have the functional equation
    \begin{equation}\label{equation:Eisenstein:local_Whittaker:Archimedean:functional_equation}
    W_{T,\psi_v}^*(h, s)^{\circ}_n = \eta_v(\det T)^{n - m - 1} W_{T,\psi_v}^*(h,-s)^{\circ}_n.
    \end{equation}
The case when $T$ is positive definite follows from \cite[{Theorem 3.1}]{Shimura82} (via \eqref{equation:Eisenstein:setup:adelic_v_classical:explicit_Archimedean_scalar_weight_vector}, see also \cite[{(3.54)}]{GS19}). The case of general $T$ (still with $\det T \neq 0$) should follow from \cite[{Theorem 4.2, (4.34.K)}]{Shimura82}, though we will give an alternative proof (Lemma \ref{lemma:Eisenstein:local_functional_equations:Archimedean:scalar_weight}). Here $\eta_v$ is the sign character $\operatorname{sgn}(-)$.

Write $(r_1,r_2)$ for the signature of $T$ (temporary notation). If either $n \geq r_1$ or $r_2 = 0$, then the function $W_{T,\psi_v}^*(h,s)^{\circ}_n$ is holomorphic for all $s \in \C$, for fixed $h$ (follows from \cite[{Theorem 4.2, (4.34.K)}]{Shimura82}).
For any $a \in \GL_m(F_v)$, we also have
    \begin{equation}\label{equation:Eisenstein:local_Whittaker:Archimedean:special_value}
    \tilde{W}^*_{T,\psi_v}(a, s_0)^{\circ}_n
    =
    \begin{cases}
    1 & \text{if $T$ is positive definite} \\
    0 & \text{if $m \leq n$ and $T$ is not positive definite.}
    \end{cases}
    \end{equation}
For the case when $T$ is positive definite, see \cite[{(3.15)}]{Shimura97} (also the proof of \cite[{Proposition 3.2}]{GS19}). The case of $T$ not positive definite with $m \leq n$ follows from \cite[{Theorem 4.2, (4.34.K)}]{Shimura82} (see also \cite[{Proposition 3.3(i)}]{GS19}).

For the general case of $\psi_v(x) = e^{2 \pi i u x}$
with possibly $u < 0$, the formula
in \eqref{equation:Eisenstein:local_Whittaker:Archimedean:functional_equation}
remains valid if ``$\eta_v(\det T)$'' is replaced with ``$\eta_v(\det uT)$'',
and \eqref{equation:Eisenstein:local_Whittaker:Archimedean:special_value}
remains valid if both instances of ``$T$'' on the right are replaced
by ``$uT$''.

            \subsection{Normalized non-Archimedean Whittaker functions} 
            \label{ssec:Eisenstein:local_Whittaker:non-Arch}
                With $n$, $\chi_v$, $\psi_v$, $\eta_v$, etc. as at the beginning of Section \ref{sec:local_Whittaker}, assume $F^+_v$ is non-Archimedean.
For the moment, we only assume $F^+_v$ has characteristic $\neq 2$, and allow both $\chi_v$ and $\psi_v$
to be possibly ramified. We can also allow $F_v/F^+_v$ to be ramified with $F^+_v$ of residue characteristic $2$ in Section \ref{ssec:Eisenstein:local_Whittaker:non-Arch}. The symbol $h$ will denote an element of $U(m,m)(F^+_v)$.

Let $\varpi_0$ be any uniformizer of $F^+_v$, and let $q_v$ be the residue cardinality of $F^+_v$. Consider $T \in \mrm{Herm}_m(F^+_v)$ with $\det T \neq 0$.
We define the \emph{local normalizing factor}
    \begin{equation}\label{equation:Eisenstein:local_Whittaker:non-Arch:normalizing_factor}
    \Lambda_{T,\psi_v}(s)^{\circ}_n \coloneqq |\Delta|_{F^+_v}^{- m(m-1)/4} \left ( \prod_{j = 0}^{m - 1} L_v(2s + m - j, \eta_v^{j + n})  \right ) |(\det T)\Delta^{\lfloor m/2 \rfloor}|_{F^+_v}^{-s - s_0} (q_v^{c(\psi_v)})^{m(s + s_0) + m^2/2}
    \end{equation}
where $c(\psi_v)$ is the conductor of $\psi_v$, as in \eqref{equation:Eisenstein:Weil_rep:Weil_index:additive_character_conductor}.
The local $L$-factors appearing in $\Lambda_{T,\psi_v}(s)^{\circ}_n$ should be compared with e.g. \cite[{\S 6}]{HKS96}.
If $\psi_v$ is any unramified nontrivial additive character (i.e. $c(\psi_v) = 0$), we
write $\Lambda_{T,v}(s)^{\circ}_n \coloneqq \Lambda_{T,\psi_v}(s)^{\circ}_n$.
For general $\psi'_v(x) = \psi_v(u x)$
for some unramified $\psi_v$ and $u \in F^{+ \times}_v$, we have
$\Lambda_{T,\psi'_v}(s)^{\circ}_n = |u|_{F^+_v}^{-m^2 / 2}\Lambda_{u T,v}(s)^{\circ}_n$.

Suppose $n \geq 0$, and suppose $V_v$ is an $n$-dimensional non-degenerate $F_v / F^+_v$ Hermitian space. Consider a full-rank lattice $L_v \subseteq V_v$, and take the Schwartz function $\varphi_v = \pmb{1}_{L_v}^m \in \mc{S}(V_v^m)$.
Form the associated Siegel--Weil standard section $\Phi_{\varphi_v} \in I(\chi_v, s)$,
which could depend on $\psi_v$.
Let $S$ be the Gram matrix of any basis for $L_v$. 

We consider the \emph{normalized local Whittaker function} $W^*_{T,\psi_v}$ and the variant $\tilde{W}^*_{T,\psi_v}$ 
    \begin{align}\label{equation:local_Whittaer:non-Arch:normalized_any_lattice}
    W^*_{T,\psi_v}(h, s, \Phi_{\varphi_v})_n & \coloneqq \gamma_{\psi_v}(V_v)^m \mrm{vol}(L_v)^{-m} \Lambda_{T,\psi_v}(s)^{\circ}_n W_{T,\psi_v}(h, s, \Phi_{\varphi_v})
    \\
    \label{equation:local_Whittaer:non-Arch:normalized_any_lattice_tilde}
    \tilde{W}^*_{T,\psi_v}(a, s, \Phi_{\varphi_v})_n & \coloneqq \chi_v(a)^{-1} |\det a|_{F_v}^{-n / 2} W^*_{T,\psi_v}(m(a), s, \Phi_{\varphi_v})
    \end{align}
for $a \in \GL_m(F_v)$. The volume $\mrm{vol}(L_v)$ is taken with respect to the self-dual Haar measure with respect to the pairing $x,y \mapsto \psi_v(\mrm{tr}(x,y))$ on $V_v$ (compare Lemma \ref{lemma:local_density_Whittaker_interpolation}). The variant $\tilde{W}^*_{T,\psi_v}$ is a local analogue of \eqref{equation:setup:adelic_v_classical:variant}. These will depend on $n$ in general. For any $a \in \GL_m(F_v)$ and $k \in \GL_m(\mc{O}_{F^+_v})$, we have the ``linear invariance'' property
    \begin{equation}\label{equation:local_Whittaker:non-Arch:linear_invariance}
    \tilde{W}^*_{T,\psi_v}(a, s, \Phi_{\varphi_v})_n = \tilde{W}_{{}^t \overline{a} T a, \psi_v}^*(1, s, \Phi_{\varphi_v})_n \quad \quad W^*_{T,\psi_v}(1, s, \Phi_{\varphi_v})_n = \tilde{W}^*_{T,\psi_v}(1,s, \Phi_{\varphi_v})_n = \tilde{W}^*_{T,\psi_v}(k, s, \Phi_{\varphi_v})_n. 
    \end{equation}
The left expression follows from \eqref{equation:setup:Fourier_and_Whittaker:Whittaker_linear_change}. The right expression follows as $\chi_v(k)^{-1} \omega_v(m(k)) \varphi_v = \varphi_v$ holds for all $k \in \GL_m(\mc{O}_{F^+_v})$, where $\omega_v$ is the local Weil representation (Section \ref{ssec:Eisenstein:Weil_rep:Weil_rep}).
If $\psi'_v$ is another additive character satisfying $\psi'_v(x) = \psi(u x)$ for some $u \in F^{+ \times}_v$,
if $V'_v$ is the Hermitian space obtained from $V_v$ after rescaling by $u^{-1}$
(Remark \ref{remark:Eisenstein:Weil_rep:Weil_rep:rescaling_trick}),
and if $\varphi'_v \in \mc{S}(V_v^{\prime m})$ is identified with $\varphi_v$
using the canonical isomorphism $V_v \cong V'_v$ (generally not preserving Hermitian pairings),
then we have $W^*_{uT, \psi_v}(h, s, \Phi_{\varphi_v})_n = W^*_{T, \psi'_v}(h, s, \Phi_{\varphi'_v})_n$
and $\tilde{W}^*_{uT, \psi_v}(h, s, \Phi_{\varphi_v})_n = \tilde{W}^*_{T, \psi'_v}(h, s, \Phi_{\varphi'_v})_n$
(consider Remark \ref{remark:Eisenstein:setup:Fourier_and_Whittaker:change_of_additive_character}).

Drop the assumption $n \geq 0$, assume $\chi_v$ is unramified, and recall the normalized spherical standard section $\Phi^{\circ}_v \in I(\chi_v, s)$.
We set
    \begin{equation}
    W^*_{T,\psi_v}(h, s)^{\circ}_n \coloneqq \Lambda_{T,\psi_v}(s)^{\circ}_n W_{T,\psi_v}(h, s, \Phi^{\circ}_v)_n 
    \quad \quad
    \tilde{W}^*_{T,\psi_v}(a, s)^{\circ}_n \coloneqq \chi_v(a)^{-1} |\det a|_{F_v}^{-n / 2} W^*_{T,\psi_v}(m(a), s)^{\circ}_n
    \end{equation}
for $h \in H(F^+_v)$ and $a \in \GL_m(F_v)$.
If $\psi'_v$ is another additive character satisfying $\psi'_v(x) = \psi(u x)$ for some $u \in F^{+ \times}_v$,
then we have $W^*_{uT, \psi_v}(h,s)^{\circ}_n = W^*_{T,\psi'_v}(h,s)^{\circ}_n$
and $\tilde{W}^*_{uT, \psi_v}(a, s)^{\circ}_n = \tilde{W}^*_{T, \psi'_v}(a, s)^{\circ}_n$
(again consider Remark \ref{remark:Eisenstein:setup:Fourier_and_Whittaker:change_of_additive_character}).

The alternative normalization
    \begin{equation}\label{equation:Eisenstein:local_Whittaker:non-Arch:parentheses-star-vs-star}
    W^{(*)}_{T,\psi_v}(h,s)^{\circ}_n \coloneqq (q_v^{c(\psi_v)})^{-m(s + s_0) - m^2 / 2} |(\det T)\Delta^{\lfloor m/2 \rfloor}|_{F^+_v}^{s + s_0} W^{*}_{T,\psi_v}(h, s)^{\circ}_n
    \end{equation}
will also be useful.
We use the shorthand $W^*_{T,\psi_v}(s)^{\circ}_n \coloneqq W^*_{T,\psi_v}(1,s)^{\circ}_n$ and $W^{(*)}_{T,\psi_v}(s)^{\circ}_n \coloneqq W^{(*)}_{T,\psi_v}(1,s)^{\circ}_n$.
We continue to use our usual convention
of using the subscript ``$v$'' instead of ``$\psi_v$''
if $\psi_v$ is an understood unramified character, suppressed from notation.

Now assume $\psi_v$ is unramified and $n \geq 0$.
If $L_v$ is self-dual, we have $\Phi_{\varphi_v} = \Phi^{\circ}_v$ (Section \ref{ssec:Eisenstein:Weil_rep:Weil_rep}), at least outside the case of $F_v / F^+_v$ ramified with residue characteristic $2$. If $F_v / F^+_v$ is ramified of residue characteristic $2$, this still holds if $L_v = (M_2^{\circ})^{\oplus d}$ for some $d \geq 0$ (with $M_2^{\circ}$ the ``standard'' self-dual lattice from \eqref{equation:Eisenstein:Weil_rep:Weil_rep:standard_self-dual_lattice}). Note that $\gamma_{\psi_v}(V_v) = 1$ in these cases.
In these situations, we have
    \begin{equation*}
    W^*_{T,\psi_v}(h, s)^{\circ}_n = W^*_{T,\psi_v}(h, s, \Phi_{\varphi_v})_n 
    \quad \quad
    \tilde{W}^*_{T,\psi_v}(a, s)^{\circ}_n = \tilde{W}^*_{T,\psi_v}(a, s, \Phi_{\varphi_v})_n
    \end{equation*}
for $h \in H(F^+_v)$ and $a \in \GL_m(F_v)$.

We further describe these functions in the following sections (e.g. special values and functional equations). We are mostly interested in the spherical local Whittaker function $W^*_{T,\psi_v}(h, s)^{\circ}_n$, and the case of general $\varphi_v$ plays a very limited role in the present work.
             
            \subsection{Local densities} 
            \label{ssec:Eisenstein:local_Whittaker:local_densities}
                We relate non-Archimedean Whittaker functions with local densities. This should be essentially known, but we restate the result for clarity (Lemma \ref{lemma:local_density_Whittaker_interpolation}).\footnote{\label{footnote:Whittaker_interpolation}The proof is essentially as in \cite[{Proposition 10.1}]{KR14}, with a few modifications. 
In the ramified situation, we should use $M^{\circ}_2$ (from Section \ref{ssec:Eisenstein:Weil_rep:Weil_rep}) instead of $L_{1,1}$ (in the proof of loc. cit.); the proposition statement changes correspondingly, see \cite[{Proposition 9.7}]{Shi22}. 
Moreover, the quantity $\g_p(V)^n$ appearing before \cite[{(10.3)}]{KR14} should be $\g_p(V)^{-n}$ for consistency with the Schr\"odinger model of the Weil representation from \cite[{Theorem 3.1 \S 3, \S 5}]{Kudla94} (and the same applies to \cite[{Proposition 9.7}]{Shi22}). 
The interpolation of $W_{T,v}(s, \Phi_{\varphi_v})$ in the two cited references should also be shifted by $s_0 = (n - m)/2$ in the $s$-variable. The cited results also restrict to the case $F^+_v = \Q_p$, but the result and (modified) proof hold more generally.
} In Section \ref{ssec:Eisenstein:local_Whittaker:local_densities}, we do not need to assume $\chi_v$ is unramified (but still require $\chi_v|_{F^{+ \times}_v} = \eta_v^n$).
We do, however, assume that $\psi_v$ is unramified.

Retain notation and assumptions from Section \ref{ssec:Eisenstein:local_Whittaker:non-Arch}. In Section \ref{ssec:Eisenstein:local_Whittaker:local_densities}, we now require $F^+_v$ to have characteristic $0$, exclude the case where $F_v/F^+_v$ is ramified with $F^+_v$ of residue characteristic $2$, and take $n \geq 0$.
We write
    \begin{align}
    \mrm{Herm}_m(\mc{O}_{F^+_v})^{*} \coloneqq {} & \{b \in \mrm{Herm}_m(F^+_v) : \mrm{tr}(bc) \in \mc{O}_{F^+_v} \text{ for all $c \in \mrm{Herm}_m(\mc{O}_{F^+_v})$}  \}
    \\
    = {} & \{b \in \mrm{Herm}_m(F^+_v) : b_{i,j} \in \mc{O}_{F^+_v} \text{ if $i = j$} \text{ and } b_{i,j} \in \mf{d}^{-1} \mc{O}_{F_v} \text{ if $i \neq j$} \}. \notag
    \end{align}
Given nonsingular Hermitian matrices $S \in \mrm{Herm}_n(F^+_v)$ and $T \in \mrm{Herm}_m(F^+_v)$, we consider the \emph{local representation density}
(or just \emph{local density})
    \begin{equation}\label{equation:Eisenstein:local_Whittaker:local_densities:def}
    \mrm{Den}(S,T) \coloneqq \lim_{k \ra \infty} \frac{\operatorname{vol}(\{x \in M_{n,m}(\mc{O}_{F_v}) : {}^t \overline{x} S x - T \in \varpi_0^k \mrm{Herm}_m(\mc{O}_{F^+_v})^{*}\})}{q_v^{-k m^2}}
    \end{equation}
where $M_{n,m}(\mc{O}_{F_v})$ is given the Haar measure of total volume $1$.
The limit argument stabilizes for $k \gg 0$ (follows from the proof of Lemma \ref{lemma:local_density_Whittaker_interpolation}).\footnote{
    We have the following argument in mind. The general result of \cite{Karel79} guarantees
    that the argument of the limit $\lim_{k \rightarrow \infty}$ from
    the proof of Lemma \ref{lemma:local_density_Whittaker_interpolation} stabilizes
    for $k \gg 0$,
    i.e. that the Whittaker integral may be computed
    as an integral with $N(F^+_v)$ replaced by any sufficiently large
    open compact subgroup.
    In the case of Lemma \ref{lemma:local_density_Whittaker_interpolation}
    with $r = 0$, the argument of the last limit which appears in the proof of
    Lemma \ref{lemma:local_density_Whittaker_interpolation}
    is, after dividing by the factor $\opn{vol}(\opn{Herm}_m(\mc{O}_{F^+_v})) \opn{vol}(L_{v,r,r}^m)$,
    exactly the argument in the limit in \eqref{equation:Eisenstein:local_Whittaker:local_densities:def}.
}
The local density $\mrm{Den}(S,T)$ depends only on the isomorphism classes of the Hermitian lattices defined by $S$ and $T$. If $n < m$, then $\mrm{Den}(S,T) = 0$.

If $S \in \mrm{Herm}_n(\mc{O}_{F^+_v})^{*}$, we have
    \begin{equation}
    \mrm{Den}(S,T) = \lim_{k \ra \infty} \frac{\# \{x \in M_{n,m}(\mc{O}_{F_v}/\varpi_0^k \mc{O}_{F_v}) : {}^t \overline{x} S x - T \in \varpi_0^k \mrm{Herm}_m(\mc{O}_{F^+_v})^{*}\}}{q_v^{k \cdot m(2n - m)}}.
    \end{equation}
If $S \in \mrm{Herm}_n(\mc{O}_{F^+_v})^{*}$ and $T \not \in \mrm{Herm}_m(\mc{O}_{F^+_v})^{*}$, we have $\mrm{Den}(S,T) = 0$.

\begin{remark}
If $S \in \mrm{Herm}_n(\mc{O}_{F^+_v})^{*}$ and $T \in \mrm{Herm}_m(\mc{O}_{F^+_v})^{*}$ with $m \leq n$, the local density $\mrm{Den}(S,T)$ admits the following equivalent formulation. Suppose $M$ (resp. $L$) is a Hermitian $\mc{O}_{F_v}$-lattice which admits a basis with Gram matrix $S$ (resp. $T$). Write $\mf{d}$ for any trace-zero generator of the different ideal $\mf{d}$ of $F_v/F^+_v$, and let $M'$ (resp. $L'$) be the skew-Hermitian lattice with pairing $\mf{d} S$ (resp. $\mf{d} T$). If $\mrm{Herm}(M', L')$ is the \emph{scheme of skew-Hermitian module homomorphisms} given by 
    \begin{equation}
    \mrm{Herm}(M', L')(R) \coloneqq \mrm{Herm}(M' \otimes R, L' \otimes R)
    \end{equation}
for $\mc{O}_{F^+_v}$-algebras $R$ (where the right-hand side means $\mc{O}_{F_v}$-linear homomorphisms preserving the skew-Hermitian pairing), we have
    \begin{equation}
    \# \mrm{Herm}(M',L')(\mc{O}_{F^+_v} / \varpi_0^k \mc{O}_{F^+_v}) = \# \{x \in M_{n,m}(\mc{O}_{F_v}/\varpi_0^k \mc{O}_{F_v}) : {}^t \overline{x} S x - T \in \varpi_0^k \mrm{Herm}_m(\mc{O}_{F^+_v})^{*}\}
    \end{equation}
and also $m (2n - m) = \dim (\mrm{Herm}(M',L') \times \Spec F^+_v)$ (and the right-hand side is nonempty).
This recovers the formulations in \cite[{\S 3.1}]{LZ22unitary} (inert), \cite[{\S 2.3}]{FYZ22SW} (inert and split), and \cite[{\S 5.1}]{HLSY22} (ramified).
\end{remark}

Return to the situation of general nonsingular $S$ and $T$ (and possibly $m > n$). Fix characters $\chi_v \colon F_v^{\times} \ra \C^\times$ and $\psi_v \colon F^+_v \ra \C^{\times}$ as above (with $\psi_v$ unramified). Let $M$ be a Hermitian $\mc{O}_{F_v}$-lattice admitting a basis whose Gram matrix is $S$. Write $V_v = M \otimes_{\mc{O}_{F_v}} F_v$ for the associated $F_v/F^+_v$ Hermitian space, and let $\varphi_v \in \mc{S}(V_v^m)$ be the function $\varphi_v = \pmb{1}_{M}^{m}$, where $\pmb{1}_{M}$ is the characteristic function of $M$. Let $\Phi_{\varphi_v} \in I(s, \chi_v)$ be the associated Siegel--Weil section, and form the local Whittaker function $W_{T,v}(h, s, \Phi_{\varphi_v})$ as in Section \ref{ssec:Eisenstein:setup:Fourier_and_Whittaker}. Set $W_{T,v}(s, \Phi_{\varphi_v}) \coloneqq W_{T,v}(1, s, \Phi_{\varphi_v})$.

With $M^{\circ}_2$ being the rank $2$ self-dual Hermitian lattice from \eqref{equation:Eisenstein:Weil_rep:Weil_rep:standard_self-dual_lattice}, let $S_{r,r}$ be the Gram matrix of a basis for $L_{v,r,r} \coloneqq M \oplus (M^{\circ}_2)^{\oplus r}$ (orthogonal direct sum).
When $F_v/F^+_v$ is not ramified, we also let $S_r$ be the Gram matrix of a basis for $L_{v,r} \coloneqq M \oplus \langle 1 \rangle^{\oplus r}$ (orthogonal direct sum), where $\langle 1 \rangle$ is a rank one self-dual lattice. The notations $L_{v,r,r}$ and $L_{v,r}$ will only be used in the proof of the next lemma.

\begin{lemma}\label{lemma:local_density_Whittaker_interpolation}
With notation as above, there exists $\mrm{Den}(S,T,X) \in \Q[X]$ (necessarily unique) such that
    \begin{align}\label{equation:local_density_Whittaker_interpolation}
    W_{T,v}(s_0 + s, \Phi_{\varphi_v}) & = \gamma_{\psi_v}(V_v)^{-m} |\det S|^{m}_{F^+_v} |\Delta|^e_{F^+_v} \mrm{Den}(S,T,q_v^{-2s}) && \text{for all $s \in \C$} \\
    \mrm{Den}(S_{r,r}, T) & = \mrm{Den}(S,T,q_v^{-2r}) && \text{for all $r \in \Z_{\geq 0}$}
    \end{align}
where $\gamma_{\psi_v}(V_v)$ is the Weil index, $s_0 = (n - m)/2$, and $e = nm/2 + m(m-1)/4$. For all $r \in \Z_{\geq 0}$, we also have
    \begin{align}
    \mrm{Den}(S_r, T) & = \mrm{Den}(S, T, (-q_v)^{-r}) && \text{if $F_v/F^+_v$ is inert} \\
    \mrm{Den}(S_r, T) & = \mrm{Den}(S, T, q_v^{-r}) && \text{if $F_v/F^+_v$ is split}.
    \end{align}
\end{lemma}
\begin{proof}
As mentioned above (Footnote \ref{footnote:Whittaker_interpolation}), this is a restatement of a result which should be essentially known \cite[{Proposition 10.1}]{KR14} \cite[{Proposition 9.7}]{Shi22}, up to a few modifications. The modified version stated here may be proved by a similar interpolation argument, as explained below. For any $r \in \Z_{\geq 0}$, set $V_{v,r,r} \coloneqq L_{v,r,r} \otimes_{\mc{O}_{F_v}} F_v$, and let $\varphi_{v,r,r} = \pmb{1}_{L_{v,r,r}}^m$. Equip $\mrm{Herm}_m(\mc{O}_{F^+_v})$ and $V_{v,r,r}$ with the self-dual Haar measures with respect to $(b,c) \mapsto \psi_v(\mrm{tr}(bc))$ and $\psi_v(\mrm{tr}_{F_v / F^+_v}(\mrm{tr}(-,-)))$ respectively. Using the Weil representation, we compute
    \begin{align}
    & W_{T,v}(s_0 + r, \Phi_{\varphi_v})
    \\
    & = \lim_{k \ra \infty} \int_{\varpi_0^{-k} \mrm{Herm}_m(\mc{O}_{F^+_v})} \psi_v(- \mrm{tr}(Tb)) \Phi_{\varphi_v}(w^{-1} n(b), s_0 + r) ~ dn(b) \notag 
    \\
    & = \gamma_{\psi_v}(V_v)^{-m} \lim_{k \ra \infty} \int_{\varpi_0^{-k} \mrm{Herm}_m(\mc{O}_{F^+_v})} \psi_v(- \mrm{tr}(Tb)) \int_{V_{v,r,r}^m} \psi_v(\mrm{tr}(b(\underline{x},\underline{x}))) \varphi_{v,r,r}(\underline{x}) ~ d\underline{x} ~ dn(b) \notag
    \\
    & = \gamma_{\psi_v}(V_v)^{-m} \lim_{k \ra \infty} \mrm{vol}(\varpi_0^{-k} \mrm{Herm}_m(\mc{O}_{F^+_v})) \int_{\substack{\underline{x} \in V_{v,r,r}^m \\ (\underline{x},\underline{x}) - T \in \varpi_0^k \mrm{Herm}_m(\mc{O}_{F^+_v})^{*}}} \varphi_{v,r,r}(\underline{x}) ~d\underline{x} \notag
    \\
    & = \gamma_{\psi_v}(V_v)^{-m} \opn{vol}(\opn{Herm}_m(\mc{O}_{F^+_v})) \opn{vol}(L_{v,r,r}^m) \opn{Den}(S_{r,r},T). \notag
    \end{align}
We have the volume identities
    \begin{equation}
    \mrm{vol}(\mrm{Herm}_m(\mc{O}_{F^+_v})) = |\Delta|_{F^+_v}^{m(m-1)/4} \quad \quad \mrm{vol}(L_{v,r,r}) = |\det S|_{F^+_v} |\Delta|_{F^+_v}^{n/2}
    \end{equation}
with respect to the self-dual Haar measures described above. We already know $W_{T,v}(s, \Phi_{\varphi_v}) \in \C[q_v^{-2s}]$ by Lemma \ref{lemma:non-Arch_Whittaker_polynomial}. Since $\mrm{Den}(S_{r,r},T) \in \Q$ for all $r \geq 0$, we conclude $W_{T,v}(s, \Phi_{\varphi_v}) \in \Q[q_v^{-2s}]$.
The additional claims involving $\mrm{Den}(S_r,T)$ in the unramified case may be proved similarly, using $L_{v,r}$ instead of $L_{v,r,r}$.
\end{proof}

            \subsection{Local densities and spherical non-Archimedean Whittaker functions} 
            \label{ssec:Eisenstein:local_Whittaker:local_densities_spherical}
                Take $F_v/F^+_v$, $\psi_v$, and $\chi_v$ as in Section \ref{ssec:Eisenstein:local_Whittaker:local_densities}, and continue to assume $n \geq 0$ for the moment. Set $s_0 = (n - m) / 2$.
We assume $\chi_v$ is unramified (the character $\psi_v$ was already assumed unramified).

Let $M^{\circ}$ be a self-dual Hermitian $\mc{O}_{F_v}$-lattice of rank $n$. 
This characterizes $M^{\circ}$ uniquely up to isomorphism, and forces $n$ to be even if $F_v / F^+_v$ is ramified.
If $F_v / F^+_v$ is ramified, then $M^{\circ}$ is isomorphic to an orthogonal
direct sum of $n / 2$ copies of the rank $2$ lattice
from \eqref{equation:Eisenstein:Weil_rep:Weil_rep:standard_self-dual_lattice}.
We also have $\g_{\psi_v}(V_v) = 1$ (Lemma \ref{lemma:Weil_index:Hermitian}).

Set $V_v = M^{\circ} \otimes_{\mc{O}_{F_v}} F_v$, and let $\varphi_v \in \mc{S}(V_v^m)$ be the characteristic function of $M^{\circ m}$. Then the associated Siegel--Weil section $\Phi_{\varphi_v} \in I(s, \chi_v)$ coincides with the normalized spherical standard section $\Phi_v^{\circ}$ (Lemma \ref{lemma:spherical_Siegel-Weil_section}).

\begin{remark}\label{remark:Eisenstein:local_Whittaker:character_could_ramify}
Even if $\chi_v$ is possibly ramified, we still have $W_{T,v}(s, \Phi_{\varphi_v}) = W_{T,v}(s, \Phi^{\circ}_v)$ for any $T \in \mrm{Herm}_m(F^+_v)$ with $\det T \neq 0$ (by Lemma \ref{lemma:non-Arch_Whittaker_polynomial}(3) or Lemma \ref{lemma:local_density_Whittaker_interpolation}), where $\Phi^{\circ}_v \in I(s, \chi'_v)$ is the standard normalized spherical section for an unramified $\chi'_v$.
\end{remark}

Suppose $T \in \mrm{Herm}_m(F^+_v)$ with $\det T \neq 0$. If $S$ is the Gram matrix of any basis for $M^{\circ}$, Lemma \ref{lemma:local_density_Whittaker_interpolation} gives
    \begin{equation}\label{equation:local_Whittaker:local_densities_spherical:no-normalization}
    W_{T,v}(s_0 + s, \Phi_v^{\circ}) = |\Delta|_{F^+_v}^{m(m-1)/4} \mrm{Den}(S, T, q_v^{-2s})
    \end{equation}
for all $s \in \C$.
Here, we are also using the equality $|\det S|_{F^+_v} = |\Delta|_{F^+_v}^{- n / 2}$.
Applying Lemma \ref{lemma:non-Arch_Whittaker_polynomial}(2)(ii)
shows that $\mrm{Den}(S,T,X) \in \Z[1/q_v][X]$.
Note that the right-hand side of \eqref{equation:local_Whittaker:local_densities_spherical:no-normalization}
manifestly has no dependence on $\chi_v$ or $\psi_v$
(our convention at the moment is that $\psi_v$ is unramified).

Suppose $M^{\circ \prime}$ is a rank $m$ Hermitian $\mc{O}_{F_v}$-lattice such that
    \begin{equation}
    \begin{cases}
    M^{\circ \prime} \text{ is self-dual} & \text{if $F_v/F^+_v$ is unramified or if $m$ is even} \\
    M^{\circ \prime} \text{ is almost self-dual} & \text{if $F_v/F^+_v$ is ramified and $m$ is odd.}
    \end{cases}
    \end{equation}
Let $S' \in \mrm{Herm}(F^+_v)$ be the Gram matrix of a basis for $M^{\circ \prime}$. We have
    \begin{align}\label{equation:local_Whittaker:spherical:local_den_and_L-functions}
    \left ( \prod_{j = 0}^{m - 1} L_v(2(s + s_0) + m - j, \eta_v^{j + n})  \right )^{-1} 
    & = \mrm{Den}(S,S',X)|_{X = q_v^{-2s}}.
    \end{align}
See \cite[{(3.2.0.1)}]{LZ22unitary} (inert), \cite[{Theorem 2.2}]{FYZ22SW} (split and inert), \cite[{Lemma 2.15}]{LL22II} (ramified).

We normalize $\mrm{Den}(S,T,X)$ and define a
(normalized) \emph{local density polynomial} $\mrm{Den}(X,T)_n \in \Z[1/q_v][X]$,
characterized by the equality
    \begin{equation}\label{equation:Eisenstein:local_Whittaker:local_densities_spherical:normalized_local_density_polynomial}
    W^{(*)}_{T,v}(s + s_0)^{\circ}_n = \mrm{Den}(q_v^{-2s}, T)_n
    \end{equation}
for all $s \in \C$ (with $W^{(*)}_{T,v}$ as in Section \ref{ssec:Eisenstein:local_Whittaker:non-Arch}).
These local density polynomials have more explicit descriptions
via ``Cho--Yamauchi formulas'' proved
in \cite[{Theorem 3.5.1}]{LZ22unitary} (inert), \cite[{Theorem 2.2}]{FYZ22SW} (split and inert), and \cite[{Lemma 2.15}]{LL22II} (ramified). Note that our convention differs slightly from \cite{LL22II} in the ramified case, where they consider polynomials in $q_v^{-s}$ instead.

The polynomial $\mrm{Den}(X,T)_n$ is nonzero if and only if $T \in \mrm{Herm}(\mc{O}_{F^+_v})^{*}$, in which case $\mrm{Den}(X,T)_n$ has constant term $1$.
When $m = n$, we have $\mrm{Den}(X,T)_n \in \Z[X]$ for any $T$. More classically, see \cite[{Theorem 13.6}]{Shimura97}, which implies that $\mrm{Den}(q_v^{n} X, T)_n \in \Z[X]$ with constant term $1$.

We have
    \begin{align}
    \mrm{Den}(X, T)_{n+1} &= \mrm{Den}(q_v^{-1} X, T)_n && \text{if $F_v/F^+_v$ is split} \notag \\
    \mrm{Den}(X, T)_{n+1} &= \mrm{Den}(-q_v^{-1} X, T)_n && \text{if $F_v/F^+_v$ is inert} \label{equation:Eisenstein:local_Whittaker:local_densities_spherical:change_of_n} \\
    \mrm{Den}(X, T)_{n+2} &= \mrm{Den}(q_v^{-2} X, T)_n && \text{if $F_v/F^+_v$ is ramified}. \notag
    \end{align}
For $n < 0$, we define $\mrm{Den}(X,T)_n$ using \eqref{equation:Eisenstein:local_Whittaker:local_densities_spherical:change_of_n}. Note that \eqref{equation:Eisenstein:local_Whittaker:local_densities_spherical:normalized_local_density_polynomial} continues to hold. For the rest of Section \ref{ssec:Eisenstein:local_Whittaker:local_densities_spherical}, we allow general $n \in \Z$ (assumed even if $F_v/F^+_v$ is ramified).

Similarly, there is a (normalized) \emph{local density (Laurent) polynomial} $\mrm{Den}^*(X,T)_n \in \Z[1/q_v][X, X^{-1/2}]$ such that
    \begin{equation}\label{equation:local_Whittaker:local_densities_spherical:Den_star}
    W^*_{T,v}(s + s_0)^{\circ}_n = \mrm{Den}^*(q_v^{-2s}, T)_n
    \end{equation} 
for all $s \in \C$ (with $W^*_{T,v}$ as in Section \ref{ssec:Eisenstein:local_Whittaker:non-Arch}). 

\begin{remark}\label{remark:local_Whittaker:local_densities_spherical:square_root_notation}
On the right-hand side of \eqref{equation:local_Whittaker:local_densities_spherical:Den_star}, we mean evaluating $\mrm{Den}^*(X,T)_n$ at $X^{1/2} = q_v^{-s}$.
We similarly abuse notation elsewhere. For example, $\mrm{Den}^*(q_v X, T)_n \in \Z[1 / q_v^{1/2}][X, X^{-1/2}]$ is obtained from $\mrm{Den}^*(X, T)_n$ by replacing $X^{1/2}$ with $q_v^{1/2} X^{1/2}$. The notation $\frac{d}{d X} \colon \Q[X, X^{-1/2}] \ra \Q[X, X^{-1/2}]$ means the $\Q$-linear map $X^{j/2} \mapsto (j/2) X^{j/2 - 1}$.
\end{remark}

If $T$ defines a self-dual Hermitian lattice when $m$ is even or $F_v/F^+_v$ is unramified (resp. almost self-dual Hermitian lattice when $m$ is odd and $F_v/F^+_v$ is ramified), we have 
    \begin{equation}
    \label{equation:Eisenstein:local_Whittaker:local_densities_spherical:often_1}
    W^*_{T,v}(s)^{\circ}_n = W^{(*)}_{T,v}(s)^{\circ}_n = 1 \quad \quad \mrm{Den}^*(X, T)_n = \mrm{Den}(X, T)_n = 1
    \end{equation}
(follows from \eqref{equation:local_Whittaker:spherical:local_den_and_L-functions}). For such $T$, an application of Lemma \ref{lemma:non-Arch_Whittaker_polynomial}(3) also shows that
    \begin{equation}\label{equation:Eisenstein:local_Whittaker:local_densities_spherical:typical_T}
    W_{T,v}(s, \Phi_v^{\circ}) = |\Delta|^{m(m-1)/4}_{F^+_v} \prod_{j = 0}^{m - 1} L_v(2s + m - j, \eta_v^{j} \chi'_v|_{F^{+ \times}_v})^{-1}
    \end{equation}
if $\Phi_v^{\circ} \in I(s, \chi'_v)$ is the normalized spherical section for any unramified character $\chi'_v \colon F_v^{\times} \ra \C^{\times}$ (not assuming $\chi'_v|_{F^{+ \times}_v} = \eta_v^n$).

If $L$ is a $\mc{O}_{F_v}$ Hermitian lattice of rank $m$, and if $L$ admits a basis with Gram matrix $T$ (allowing arbitrary $T \in \mrm{Herm}_m(F^+_v)$ with $\det T \neq 0$ again), we write $\mrm{Den}(X,L)_n \coloneqq \mrm{Den}(X, T)_n$ and similarly $\mrm{Den}^*(X, L)_n \coloneqq \mrm{Den}^*(X, T)_n$.
We claim that we have
    \begin{align}\label{equation:local_Whittaker:local_densities_spherical:comparison}
    \mrm{Den}^*(X,L)_n
    & = 
    (q_v^{2 s_0} X^{- 1/2})^{\mrm{val}'(L)} \mrm{Den}(X,L)_n
    \\
    \label{equation:local_Whittaker:local_densities_spherical:val_prime}
    \mrm{val}'(L) \coloneqq \lfloor \mrm{val}(L) \rfloor & =
    \begin{cases}
    \mrm{val}(L) - 1/2 & \text{if $F_v / F^+_v$ is ramified and $m$ is odd} \\
    \mrm{val}(L) & \text{else}
    \end{cases}
    \end{align}
with $\mathrm{val}(L)$ as defined in \crefext{I:ssec:Hermitian_conventions:lattices}.\footnote{
    For the reader's convenience: the quantity $\mrm{val}(L) \in \frac{1}{2} \Z$
    can be characterized by the formula $q_v^{2 \mrm{val}(L)} = |L^* / L|$
    for integral lattices $L$, where $L^*$ is the dual for the trace pairing as in loc. cit.,
    and by the additional formula
    $\mrm{val}(\varpi_0 L) = 2 m + \mrm{val}(L)$
    for arbitrary $L$, where $\varpi_0$ is any uniformizer for $\mc{O}_{F^+_v}$.
}
Indeed, upon recalling \eqref{equation:Eisenstein:local_Whittaker:non-Arch:parentheses-star-vs-star}
(and the definitions
\eqref{equation:Eisenstein:local_Whittaker:local_densities_spherical:normalized_local_density_polynomial}
and \eqref{equation:local_Whittaker:local_densities_spherical:Den_star}),
the formula in \eqref{equation:local_Whittaker:local_densities_spherical:comparison}
follows from the computation that $|(\det T) \Delta^{\lfloor m / 2 \rfloor}|_{F^+_v} = q_v^{-\mrm{val}'(L)}$.

The local densities satisfy the following cancellation property (which we will use).

\begin{lemma}
If $L^{\circ}$ is a self-dual Hermitian lattice of rank $d$, then for any non-degenerate Hermitian lattice $L$ and every integer $r \in \Z$ (assume $r$ is even if $F_v/F^+_v$ is ramified), we have
    \begin{equation}\label{equation:local_Whittaker:local_densities_spherical:cancellation}
    \mrm{Den}(X,L \oplus L^{\circ})_{r + d} = \mrm{Den}(X,L)_{r} \quad \quad \mrm{Den}^*(X,L \oplus L^{\circ})_{r + d} = \mrm{Den}^*(X,L)_{r}
    \end{equation}
where $L \oplus L^{\circ}$ is the orthogonal direct sum.
\end{lemma}
\begin{proof}
This can be seen from the Cho--Yamauchi type formulas
mentioned after \eqref{equation:Eisenstein:local_Whittaker:local_densities_spherical:normalized_local_density_polynomial},
as we explain below.

First assume $\rank(L) = r$
(if $F_v / F^+_v$ is ramified, we will deal with the 
odd rank case separately).
Then the cited Cho--Yamauchi formulas may be formulated
in the following uniform way:
we have
    \begin{align}\label{equation:Eisenstein:non-Arch_identity:horizontal_identity:Cho--Yamauchi}
    \mrm{Den}(X,L)_r & = \sum_{L \subseteq M \subseteq M^*} X^{\ell(M / L)} \mrm{Den}(X,M)_r^{\circ}
    \\
    \mrm{Den}(X,M)_r^{\circ} & \coloneqq \prod_{i = 0}^{t(M) - 1} (1 - \eta^{*i}(\varpi_0) q_v^i X) \label{equation:non-Arch_identity:horizontal_identity:Cho--Yamauchi:"primitive"}
    \end{align}
where $\eta^*(\varpi_0) \coloneqq \eta^{*i}(\varpi_0) \coloneqq -1, 0, 1$ if $i$ is odd (resp. $\eta^{*i}(\varpi_0) \coloneqq 1$ if $i$ is even) in the inert, ramified, split cases respectively,
and $\mrm{Den}(X,L)_r \in \Z[X]$ is the local density polynomial normalized as
in \eqref{equation:Eisenstein:local_Whittaker:local_densities_spherical:normalized_local_density_polynomial}.
The displayed sum runs over lattices $M \subseteq L_{F_v} \coloneqq L \otimes_{\mc{O}_{F_v}} F_v$.
The notations $\ell(M / L)$ and $t(M)$ are as defined in
\crefext{I:ssec:Hermitian_conventions:lattices};
for the reader's convenience, we recall
$q^{\ell(M / L)} = |M/L|$ (right-hand side means cardinality)
and recall that $t(M)$ is the number of nonzero fundamental invariants of $M$
(and depends only on $M^* / M$).
The above formula \eqref{equation:Eisenstein:non-Arch_identity:horizontal_identity:Cho--Yamauchi}
also appeared in \crefext{I:sec:non-Arch_identity},
where it was used in a crucial way for our main local results.

In the case where $F_v / F^+_v$ is ramified and $\rank(L) = r - 1$,
a similar formula for $\mrm{Den}(X,L)_r$ is given
in \crefext{I:ssec:non-Arch_identity:horizontal_identity},
upon taking $L^{\flat}$ in loc. cit. to be our current $L$
(the formula in loc. cit. is extracted from \cite[{Lemma 2.15}]{LL22II}).
For our present purposes, we only need to know that the formula is
also a sum over $M$ satisfying $L \subseteq M \subseteq M^*$ as above,
where the polynomials being summed depend only on $t(M)$ and $\ell(M / L)$.

In the case $r = \rank_{\mc{O}_{F_v}}(L)$,
the left equality in \eqref{equation:local_Whittaker:local_densities_spherical:cancellation}
now follows from \eqref{equation:Eisenstein:non-Arch_identity:horizontal_identity:Cho--Yamauchi}
(replacing $L$ with $L \oplus L^{\circ}$)
via the following linear algebra fact: every lattice
$M \subseteq (L \oplus L^{\circ}) \otimes_{\mc{O}_{F_v}} F_v$ satisfying
$L^{\circ} \subseteq M \subseteq M^{*}$
admits an orthogonal direct sum decomposition
$M = L^{\circ} \oplus M'$ for some sublattice $M'$.
The case of $r = \rank_{\mc{O}_{F_v}}(L) + 1$ and $F_v / F^+_v$ ramified
follows by similar reasoning.
The case of general $r$ now follows by shifting
as in \eqref{equation:Eisenstein:local_Whittaker:local_densities_spherical:change_of_n}.

The right equality in \eqref{equation:local_Whittaker:local_densities_spherical:cancellation}
now follows, using \eqref{equation:local_Whittaker:local_densities_spherical:comparison}
and the fact that $\mrm{val}(L \oplus L^{\circ}) = \mrm{val}(L)$.
\end{proof}

        \section{Local functional equations} 
        \label{sec:Eisenstein:local_functional_equations}
            Let $F_v$ be a degree $2$ \'etale algebra over a local field $F^+_v$ of characteristic $\neq 2$, with notation $\mf{d}$, $\Delta$, $\eta_v$, and $a \mapsto \overline{a}$ as above. If $F^+_v$ is Archimedean, we also assume $F_v / F^+_v$ is $\C / \R$. Fix an integer $m \geq 0$.

Consider a character $\chi_v \colon F^{\times}_v \ra \C^{\times}$ and a nontrivial additive character $\psi_v \colon F^+ \ra \C^{\times}$ (for the moment, we do not require $\chi_v|_{F^{+ \times}_v} = \eta_v^n$, and allow $\chi_v$ and $\psi_v$ to be ramified).

Set $\check{\chi}_v(a) \coloneqq \chi_v(\overline{a})^{-1}$. There is a local intertwining operator
    \begin{equation}
    M(s,\chi_v) \colon I(s, \chi_v) \ra I(-s, \check{\chi}_v)
    \end{equation}
(where $I(s, \chi_v)$ and $I(-s, \check{\chi}_v)$ are degenerate local principal series for $U(m,m)$) defined by the integral
    \begin{equation}
    M(s,\chi_v)\Phi_v(h) = \int_{N(F^+_v)} \Phi_v(w^{-1} n(b) h, s) ~dn(b)
    \end{equation}
for $\mrm{Re}(s) > m/2$, with meromorphic continuation to $\C$ (e.g. see \cite{KS97} in the non-Archimedean case).

Given $T \in \mrm{Herm}_m(F^+_v)$, we define the quantity
    \begin{align}
    \kappa_T(s,\chi_v,\psi_v) & \coloneqq \chi_v(-1)^m \chi_v(\det T)^{-1} |\det T|_{F^+_v}^{-2s} \g_{\psi_v}(F_v)^{m(m-1)/2} \eta_v(\det T)^{m-1}  \notag
    \\
    & \mathrel{\phantom{=}} \cdot \prod_{j = 0}^{m-1} \rho_{v}(2s + j - m + 1, \eta_v^{j} \cdot \chi_v|_{F_v^{+ \times}}, \overline{\psi}_v)
     \end{align}
where $\gamma_{\psi_v}(F_v)$ is a Weil index (Section \ref{ssec:Eisenstein:Weil_rep:Weil_index},
with conventions as explained before \eqref{equation:Eisenstein:Weil_rep:Weil_index:epsilon_factors_and_Weil_index})
and $\rho_v$ is a local factor as in Tate's thesis (Section \ref{ssec:local_Whittaker:local_L-factors}). This factor is taken from \cite[{\S 3}]{KS97}\footnote{
    The factor $\kappa_T(s, \chi_v, \psi_v)$ is given there in the non-Archimedean case, but we will use the same formula in the Archimedean case. For comparing our formulas with theirs, note the different convention used to define $W_{T,v}$ and $M(s, \chi_v)$ ($\psi_v$ versus $\overline{\psi}_v$ and $w$ vs $w^{-1}$).
    Up to these changes in notation/conventions,
    our $\kappa_T$ should be compared with the $\kappa_T$ of \cite[{Proposition 6.3}]{HKS96},
    or the $\gamma(s, \chi, \beta)$ from \cite[{(3.5)}]{KS97},
    which unfolds via \cite[{(3.25),(3.9)}]{KS97} (also \cite[{Page 303}]{KS97} for the split case).
    For comparison with \cite[{Proposition 6.3}]{HKS96},
    we remark that our $\gamma_{\psi_v}(F_v)$
    equals $\gamma(\Delta, \overline{\psi}_v)$
    from loc. cit. (upon taking their $E / F$ to be our $F_v / F^+_v$),
    as follows e.g. from the formula at the bottom of \cite[{Page 387}]{Kudla94}
    (apply that formula with $m = 1$ and with their $V$ equal to our Hermitian
    space $F_v$, with conventions as explained before \eqref{equation:Eisenstein:Weil_rep:Weil_index:epsilon_factors_and_Weil_index}).
    The additional $\chi_v(-1)^m$ in our formula is accounted for
    because we define Whittaker functions and intertwining operators using $w^{-1}$
    instead of $w$ as in loc. cit.; note that $w^{-1} = -w$,
    and then recall that $\Phi_v(m(-1) h, s) = \chi_v(-1)^m \Phi_v(h, s)$
    for all $h \in U(m,m)(F^+_v)$
    and any $\Phi_v$, by \eqref{equation:Eisenstein:setup:adelic_v_classical:local_degen_principal_series}.
    Also note $\check{\chi}(-1) = \chi(-1)$.
}
(see also \cite[{Proposition 6.3}]{HKS96}).
    
            \subsection{Non-Archimedean} 
            \label{ssec:Eisenstein:local_functional_equations:non-Arch}
                Suppose $F^+_v$ is non-Archimedean (with notation as above). For any $T \in \mrm{Herm}_m(F^+_v)$ with $\det T \neq 0$ and any standard section $\Phi_v$ of $I(s, \chi_v)$, there is a functional equation
    \begin{equation}\label{equation:Eisenstein:local_functional_equations:non-Arch:Kudla-Sweet}
    W_{T,\psi_v}(h, -s, M(s, \chi_v) \Phi_v) = \kappa_T(s,\chi_v,\psi_v) W_{T,\psi_v}(h, s, \Phi_v)
    \end{equation}
as in \cite[{\S 3, \S 7}]{KS97}.

We next consider spherical Whittaker functions. Assume $\chi_v$ is unramified.
We require $F^+_v$ to be characteristic $0$ (because \cite[{\S 13}]{Shimura97} assumes this).
With $\Phi^{\circ}_v$ denoting the normalized spherical sections of $I(s, \chi_v)$ and $I(s, \check{\chi}_v)$, we have
    \begin{equation}\label{equation:Eisenstein:local_functional_equations:non-Arch:spherical_intertwining}
    M(s,\chi_v) \Phi^{\circ}_v(s) = q_v^{-c(\psi_v) \cdot m^2 / 2} |\Delta|_{F^+_v}^{m(m-1)/4} \prod_{j=0}^{m-1} \frac{L_v(2s + j - m + 1, \eta_v^{j} \chi_v|_{F^{+ \times}_v})}{L_v(2s + m - j, \eta_v^{j} \chi_v|_{F^{+ \times}_v})} \Phi^{\circ}_v(-s),
    \end{equation}
see \cite[{Theorem 13.6}]{Shimura97}.\footnote{Take $\zeta = 0$ in the notation of loc. cit.. Strictly speaking, the statement there is only for $\chi_v|_{F^{+ \times}_v}$ trivial, but the general case follows from this; see \eqref{equation:Eisenstein:element_decomp} and the proof of Lemma \ref{lemma:non-Arch_Whittaker_polynomial}(3).
The statement in loc. cit. also restricts to unramified $\psi_v$, but the general case follows immediately from this,
as the factor $q_v^{-c(\psi_v) \cdot m^2 / 2}$ accounts for the $\psi_v$-dependence of the self-dual Haar measure
on $N(F^+_v)$.}

Now, we further restrict to the situation where $\chi_v|_{F^{+ \times}_v} = \eta_v^n$ for some $n \in \Z$, with $n$ assumed even if $F_v/F^+_v$ is ramified. Note $\check{\chi}_v = \chi_v$.
We also assume $\psi_v$ is unramified, to simplify exposition.
Combining \eqref{equation:Eisenstein:local_functional_equations:non-Arch:spherical_intertwining} with the 
identities stated above (including the relation between Weil indices and epsilon factors in \eqref{equation:Eisenstein:Weil_rep:Weil_index:epsilon_factors_and_Weil_index}), a straightforward computation (omitted)
yields the functional equations
    \begin{align}
    W^{(*)}_{T,v}(h, -s)^{\circ}_n &= |(\det T) \Delta^{\lfloor m/2 \rfloor}|_{F^+_v}^{-2s} \eta_v((-1)^{m(m-1)/2} \det T)^{n - m - 1} W^{(*)}_{T,v}(h, s)^{\circ}_n \\
    W^*_{T,v}(h, -s)^{\circ}_n &= \eta_v((-1)^{m(m-1)/2} \det T)^{n - m - 1} W^*_{T,v}(h, s)^{\circ}_n \label{equation:local_functional_equations:non-Arch:Whittaker_star}
    \end{align}
with $W^{(*)}_{T,v}(h,s)^{\circ}_n$ and $W^*_{T,v}(h,s)^{\circ}_n$ as in Section \ref{ssec:Eisenstein:local_Whittaker:non-Arch}.

\begin{remark}
We also obtain functional equations in the case when $\psi_v$
is not necessarily unramified.
Let $\psi_v$ be an arbitrary nontrivial additive character,
and pick any $u \in F^{+ \times}_v$
with the property that $x \mapsto \psi_v(u^{-1} x)$
is unramified, i.e. $u \in \varpi_0^{c(\psi_v)} \mc{O}_{F^+_v}^{\times}$
for $\varpi_0$ a uniformizer of $\mc{O}_{F^+_v}$.
Using the relation $W^*_{T,\psi_v}(h,s)^{\circ}_n = W^*_{u T, v}(h,s)^{\circ}_n$
(see Section \ref{ssec:Eisenstein:local_Whittaker:non-Arch}),
the functional equation in \eqref{equation:local_functional_equations:non-Arch:Whittaker_star}
immediately gives
    \begin{equation}\label{equation:local_functional_equations:non-Arch:arbitrary_psi:Whittaker_star}
    W^*_{T,\psi_v}(h,-s)^{\circ}_n
    =
    \eta_v((-1)^{m(m-1)/2} \det T)^{n - m - 1} \eta_v(\varpi_0)^{c(\psi_v) \cdot m(n - m - 1)}
        W^*_{T,\psi_v}(h,s)^{\circ}_n.
    \end{equation}
Note that the sign $\eta_v(\varpi_0)^{c(\psi_v) \cdot m(n - m - 1)}$
always equals $1$ if $n$ is even.
In particular, the sign $\eta_v(\varpi_0)^{c(\psi_v) \cdot m(n - m - 1)}$
does not depend on the choice of uniformizer $\varpi_0$, because
we are assuming $n$ is even if $F_v / F^+_v$ is ramified.
\end{remark}

Next, assume that $F_v/F^+_v$ is unramified or that $F^+_v$ has residue characteristic $\neq 2$. 
If $L$ is a Hermitian $\mc{O}_{F_v}$-lattice, we thus have
    \begin{align}\label{equation:Eisenstein:local_functional_equations:non-Arch:non-star_density}
    \mrm{Den}(q_v^{2 s_0} X^{-1},L)_n & = \varepsilon(L)^{n - m - 1} X^{- \mrm{val}'(L)} \mrm{Den}(q_v^{2 s_0} X,L)_n \\
    \mrm{Den}^*(q_v^{2 s_0} X^{-1},L)_n & = \varepsilon(L)^{n - m - 1} \mrm{Den}^*(q_v^{2 s_0} X,L)_n
    \label{equation:local_functional_equations:non-Arch:Den_star}
    \end{align}
with $\mrm{val}'(L) \coloneqq \lfloor \mrm{val}(L) \rfloor$ as in \eqref{equation:local_Whittaker:local_densities_spherical:val_prime} (both $\varepsilon(L)$ and $\mrm{val}(L)$ were defined in \crefext{I:ssec:Hermitian_conventions:lattices}).

In the case where $\chi_v|_{F^{+\times}_v}$ is trivial, these functional equations are essentially \cite[{Corollary 3.2}]{Ikeda08}.

            \subsection{Archimedean} 
            \label{ssec:Eisenstein:local_functional_equations:Arch}
                Suppose $F_v / F^+_v$ is $\C / \R$ (with notation as above). For any $T \in \mrm{Herm}_m(F^+_v)$ with $\det T \neq 0$ and any standard section $\Phi_v$ of $I(s, \chi_v)$, we have            
    \begin{equation}\label{equation:Eisenstein:local_functional_equations:Archimedean:Kudla-Sweet}
    W_{T,v}(h, -s, M(s, \chi_v) \Phi_v) = \kappa_T(s,\chi_v,\psi_v) W_{T, v}(h, s, \Phi_v).
    \end{equation}
This may be deduced, e.g. by combining the non-Archimedean analogue \eqref{equation:Eisenstein:local_functional_equations:non-Arch:Kudla-Sweet} with the global functional equation \eqref{equation:Eisenstein:setup:adelic_v_classical:global_functional_equation}.

In the rest of Section \ref{ssec:Eisenstein:local_functional_equations:Arch}, we require $\chi_v|_{F^{+ \times}_v} = \eta_v^n$ for some $n \in \Z$, and let $\psi_v(x) = e^{2\pi i x}$.
Recall that we have defined a normalized Archimedean Whittaker function $W^*_{T,v}(h,s)^{\circ}_n$ (Section \ref{ssec:Eisenstein:local_Whittaker:Archimedean}).

\begin{lemma}\label{lemma:Eisenstein:local_functional_equations:Archimedean:scalar_weight}
For any $T \in \mrm{Herm}_m(F^+_v)$ with $\det T \neq 0$, we have the functional equation
    \begin{equation}\label{equation:Eisenstein:local_functional_equations:Archimedean:scalar_weight}
    W^*_{T,v}(h, -s)^{\circ}_n = \eta_v(\det T)^{n - m -1} W^*_{T,v}(h, s)^{\circ}_n.
    \end{equation}
\end{lemma}
\begin{proof}
By \eqref{equation:Eisenstein:local_functional_equations:Archimedean:Kudla-Sweet}, we must have $W^*_{T,v}(h,-s)^{\circ}_n = \eta_v(\det T)^{n - m - 1} f(s) W^*_{T,v}(h,s)^{\circ}_n$
for some meromorphic factor $f(s)$ which is independent of $T$. When $T$ is positive definite, we have $f(s) = 1$ (see Section \ref{ssec:Eisenstein:local_Whittaker:Archimedean}), so we obtain the claimed functional equation for all $T \in \mrm{Herm}_m(F^+_v)$ with $\det T \neq 0$. Note that $\eta_v$ is simply the sign character $\operatorname{sgn}(-)$.
\end{proof}

\begin{remark}
We also obtain functional equations in the case when $\psi_v$
is not necessarily the standard additive character,
as stated in \eqref{equation:Eisenstein:local_Whittaker:Archimedean:functional_equation}
and the discussion at the end of Section \ref{ssec:Eisenstein:local_Whittaker:Archimedean}.
Indeed, if $\psi_v(x) = e^{2 \pi i u x}$
is an arbitrary nontrivial additive character
for some $u \in \R^{\times}$,
the claimed functional equation for $W^*_{T,\psi_v}(h,s)^{\circ}_n$
is immediately implied by \eqref{equation:Eisenstein:local_functional_equations:Archimedean:scalar_weight}
and the relation $W^*_{T,\psi_v}(h,s)^{\circ}_n = W^*_{u T, v}(h,s)^{\circ}_n$
(Section \ref{ssec:Eisenstein:local_Whittaker:Archimedean}).
\end{remark}

Recall that $\Phi_v^{(n)} \in I(s, \chi_v)$ is our notation for a certain scalar weight standard section, as in Section \ref{ssec:Eisenstein:setup:adelic_v_classical}. For verifying the next lemma, it may be helpful to recall the relation between local epsilon factors $\epsilon_v(-)$ and Weil indices $\gamma_v(-)$ (Section \ref{ssec:Eisenstein:Weil_rep:Weil_index}).

\begin{lemma}\label{lemma:Eisenstein:local_functional_equations:Archimedean:intertwining_on_scalar_weight}
We have
    \begin{align} \label{equation:Eisenstein:local_functional_equations:Archimedean:intertwining_on_scalar_weight}
    & M(s, \chi_v) \Phi_v^{(n)}(s)
    \\
    & = \left ( \prod_{j = 0}^{m-1} \frac{L_v(2s + j - m + 1, \eta_v^{n+j}) \Gamma(-s - s_0 + n - j)}{\epsilon_v(2s + j - m + 1, \eta_v^{n+j}, \overline{\psi}_v) L_v(-2s - j + m, \eta_v^{n+j}) \Gamma(s - s_0 + n - j)} \right ) \notag
    \\
    & \mathrel{\phantom{=}} \cdot (-1)^{nm} i^{m(m-1)/2} \pi^{2 m s} \Phi_v^{(n)}(-s) \notag
    \end{align}
with $s_0 = (n - m)/2$ as above.
\end{lemma}
\begin{proof}
A priori, the displayed identity holds up to some meromorphic scale factor. We may compute this scale factor by combining \eqref{equation:Eisenstein:local_functional_equations:Archimedean:Kudla-Sweet} and Lemma \ref{lemma:Eisenstein:local_functional_equations:Archimedean:scalar_weight} (take $T = 1_m$).
\end{proof}

\begin{remark} 
Lemma \ref{lemma:Eisenstein:local_functional_equations:Archimedean:intertwining_on_scalar_weight} should be a reformulation (with alternative proof) of a case of \cite[{(1.31)}]{Shimura82} (translating into Shimura's setup via \eqref{equation:Eisenstein:setup:adelic_v_classical:explicit_Archimedean_scalar_weight_vector}). Shimura's computation in loc. cit. implies
    \begin{equation}
    M(s, \chi_v) \Phi_v^{(n)}(s) = \left ( \frac{i^{-mn} (2 \pi)^{m^2} \pi^{-m(m-1)/2}}{2^{m(m-1)/2 + 2ms}} \prod_{j = 0}^{m-1} \frac{\Gamma(2s - j)}{\Gamma(s - s_0 + n - j) \Gamma(s - s_0 - j)} \right ) \Phi^{(n)}_v(-s).
    \end{equation}
Similarly, the functional equation in Lemma \ref{lemma:Eisenstein:local_functional_equations:Archimedean:scalar_weight} should follow from \cite[{Theorem 4.2, (4.34.K)}]{Shimura82} (alternative proof) after some rearranging.
    
For our later calculations, we prefer to use these results as stated in Lemmas \ref{lemma:Eisenstein:local_functional_equations:Archimedean:scalar_weight} and \ref{lemma:Eisenstein:local_functional_equations:Archimedean:intertwining_on_scalar_weight}.
\end{remark}

\begin{remark}
If instead $\psi_v(x) = e^{2\pi i u x}$ for some $u \in \R^{\times}$,
then the formula in \eqref{equation:Eisenstein:local_functional_equations:Archimedean:intertwining_on_scalar_weight}
remains valid after multiplying by $|u|_{F^+_v}^{m^2 / 2}$ on the right,
to account for the self-dual Haar measure on $\mrm{Herm}_m(F^+_v)$,
which depends on $\psi_v$.
\end{remark}
        
        \section{Normalized Fourier coefficients} 
        \label{sec:Eisenstein:global_normalized_Fourier}
    
            \subsection{Global normalization} 
            \label{ssec:Eisenstein:global_normalized_Fourier:global_normalization}
                With notation as in Section \ref{ssec:Eisenstein:setup:adelic_v_classical} and Section \ref{ssec:Eisenstein:setup:Fourier_and_Whittaker}, let $F / F^+$ be a CM extension of number fields. We allow $2$-adic places of $F^+$ to ramify in $F$. Write $\mf{d}$ (resp. $\Delta$) for the different ideal (resp. discriminant ideal) of $F / F^+$.
Write $\mf{d}_{F^+ / \Q}$ (resp. $\Delta_{F^+ / \Q}$) for the different ideal (resp. discriminant ideal) of $F^+ / \Q$.
Let $\eta \colon F^{+ \times} \backslash \A^{\times} \ra \{ \pm 1\}$ be the quadratic character associated with $F / F^+$.

Fix integers $m$ and $n$ with $m \geq 0$, with $s_0 \coloneqq (n - m)/2$ as above.
Let $\chi \colon F^{\times} \backslash \A_F^{\times} \ra \C^{\times}$ be a character satisfying $\chi|_{\A^{\times}} = \eta^n$. To simplify, we assume that $\chi$ is unramified at every non-Archimedean place,
as the spherical standard section $\Phi^{\circ}_v$ for non-Archimedean $v$ will not exist if $\chi_v$ is ramified
(but see also Remark \ref{remark:Eisenstein:local_Whittaker:character_could_ramify}).
Thus, if any non-Archimedean place of $F^+$ ramifies in $F$,
these assumptions force $n$ to be even.
For arbitrary $F / F^+$, there certainly exists a character
$\chi$ as above if $n$ is even, in which case $\chi$ could be the trivial character.

Take the standard section
    \begin{equation}
    \Phi^{(n) \circ} \coloneqq \bigg (\bigotimes_{v \mid \infty} \Phi_v^{(n)} \bigg ) \otimes \bigg (\bigotimes_{v < \infty} \Phi^{\circ}_v \bigg) \in I(s, \chi)
    \end{equation}
(scalar weight at Archimedean places and spherical at non-Archimedean places). Form the associated Eisenstein series $E(h,s,\Phi^{(n) \circ})$ and its variants $E(z,s,\Phi^{(n) \circ})_n$ and $\tilde{E}(a, s, \Phi^{(n) \circ})_n$ as in Section \ref{ssec:Eisenstein:setup:adelic_v_classical}. The Eisenstein series variant $\tilde{E}(a,s, \Phi^{(n) \circ})_n$ does not depend on the choice of $\chi$ (Remark \ref{remark:Weil_rep:Weil_rep:chi_independence}).

Define the \emph{global normalizing factor}
    \begin{align}\label{equation:global_normalized_Fourier:global_normalization:normalizing_factor}
    \Lambda_m(s)^{\circ}_n &\coloneqq \left (\frac{(2 \pi)^{m(m-1)/2}}{(-2 \pi i)^{nm}} \pi^{m(- s + s_0)} \right ) ^{[F^+ : \Q]} |N_{F^+/\Q}(\Delta)|^{m(m-1)/4} |N_{F^+/\Q}(\Delta^{\lfloor m/2 \rfloor})|^{s + s_0} \notag
    \\
    & \mathrel{\phantom{\coloneqq}} \cdot \left ( \prod_{j = 0}^{m - 1} \Gamma(s - s_0 + n - j)^{[F^+ : \Q]} \cdot L(2s + m - j, \eta^{j + n})  \right ) \cdot |\Delta_{F^+ / \Q}|^{m(s + s_0) + m^2 / 2}.
    \end{align}
We define the \emph{normalized Eisenstein series} and its variants
    \begin{align}
    & E^*(h,s)^{\circ}_n \coloneqq \Lambda_m(s)^{\circ}_n E(h, s, \Phi^{(n) \circ}) 
    \\
    & E^*(z, s)_n^{\circ} \coloneqq \Lambda_m(s)^{\circ}_n E(z,s,\Phi^{(n) \circ})_n \quad \quad \tilde{E}^*(a, s)_n^{\circ} \coloneqq \Lambda_m(s)^{\circ}_n \tilde{E}(a,s,\Phi^{(n) \circ})_n
    \end{align}
where $h \in U(m,m)(\A)$ and $z \in \mc{H}_m$ and $a \in \GL_m(\A_F)$.
For $T \in \mrm{Herm}_m(F^+)$,
and given any nontrivial additive character $\psi \colon F^+ \backslash \A \ra \C^{\times}$,
we similarly define 
    \begin{align}
    & E^*_{T,\psi}(h,s)^{\circ}_n \coloneqq \Lambda_m(s)^{\circ}_n E_{T,\psi}(h, s, \Phi^{(n) \circ})
    \\
    & E^*_{T,\psi}(y, s)^{\circ}_n \coloneqq \Lambda_m(s)^{\circ}_n E_{T,\psi}(y, s, \Phi^{(n) \circ})_n \quad \quad \tilde{E}^*_{T,\psi}(a,s)^{\circ}_n \coloneqq \Lambda_m(s)^{\circ}_n \tilde{E}_{T,\psi}(a, s, \Phi^{(n) \circ})_n.
    \end{align}
The latter two are \emph{normalized Fourier coefficients}.
We often suppress $\psi$ from notation.

Given any $T \in \mrm{Herm}_m(F^+)$ with $\det T \neq 0$,
and given any nontrivial additive character $\psi \colon F^+ \backslash \A \ra \C^{\times}$,
the local normalizing factors from Sections \ref{ssec:Eisenstein:local_Whittaker:Archimedean} and \ref{ssec:Eisenstein:local_Whittaker:non-Arch} satisfy
    \begin{equation}
    \Lambda_m(s)^{\circ}_n \coloneqq \prod_{v} \Lambda_{T,\psi_v}(s)^{\circ}_n
    \end{equation}
where the product (over all places $v$ of $F^+$) is convergent for $\mrm{Re}(s) > 0$.
Indeed, it is enough to check the case where $\psi = \psi_{\mrm{tr}}$ is the standard character
satisfying $\psi_v(x) = e^{2 \pi i x}$ for all $v \mid \infty$,
and this case follows because $q_v^{-c(\psi_v)} = |N_{F^+_v / \Q_p}(\mf{d}_{F^+_v / \Q_p})|_p$
for all $p$-adic places $v$,
where we also use the relation $N_{F^+ / \Q}(\mf{d}_{F^+ / \Q}) = \Delta_{F^+ / \Q}$.
The case of general $\psi$ follows by product formula.
For nonsingular $T$ as above, and for arbitrary $\psi$, 
we have factorizations into (normalized) local Whittaker functions
    \begin{equation}\label{equation:Eisenstein:global_normalized_Fourier:global_normalization:factorization}
    E^*_{T,\psi}(h, s)^{\circ}_n = \prod_v W^*_{T,\psi_v}(h_v, s)^{\circ}_n \quad \quad \tilde{E}^*_{T,\psi}(a, s)^{\circ}_n = \prod_{v} \tilde{W}^*_{T,\psi_v}(a_v ,s)^{\circ}_n
    \end{equation}
where all but finitely many factors are identically equal to $1$ (as functions of $s$) for fixed $T$, $h$, and $n$.

The next lemma involves a sign $\eta(\mf{d}_{F^+/\Q})^{m(n - m - 1)} \in \{ \pm 1 \}$
which is abuse of notation, as we now explain.
We define
    \begin{equation}
    \eta(\mf{d}_{F^+/\Q})^{m(n - m - 1)}
    \coloneqq
    \prod_{v < \infty} \eta_v(\mf{d}_{F^+_v / \Q_p})^{m(n - m - 1)}
    \end{equation}
where the product runs over all non-Archimedean places of $F^+$,
and $p$ denotes the residue characteristic of $F^+_v$.
The notation $\eta_v(\mf{d}_{F^+_v / \Q_p})^{m(n - m - 1)}$
is defined to mean $\eta_v(u)^{m(n - m - 1)}$ for any generator $u$ of $\mf{d}_{F^+_v / \Q_p}$;
note that the choice of generator does not matter, since $n$ is even if $\eta_v$ is ramified.
We find
    \begin{equation}\label{equation:Eisenstein:global_normalized_Fourier:global_normalization:different_ideal_sign}
    \eta(\mf{d}_{F^+/\Q})^{m(n - m - 1)}
    =
    \begin{cases}
    +1 & \text{if $mn$ is even}
    \\    
    \displaystyle\prod_{v < \infty} \gamma_{\psi_{\mrm{tr},v}}(F_v)
    =
    \displaystyle\prod_{v < \infty} \epsilon_v(1/2, \eta_v, \psi_{\mrm{tr},v})
    &
    \text{if $mn$ is odd}
    \end{cases}
    \end{equation}
using \eqref{equation:Eisenstein:Weil_rep:Weil_index:epsilon_factors_and_Weil_index} and
\eqref{equation:Eisenstein:Weil_rep:Weil_index:change_of_additive_character}
and Lemma \ref{lemma:Weil_index:Hermitian}(4),
since we are assuming that $F / F^+$ is unramified at all non-Archimedean places if $n$ is odd.
Note also $\eta(\mf{d}_{F^+/\Q})^{m(n - m - 1)} = + 1$ if $F^+ = \Q$.

\begin{lemma}\label{lemma:Eisenstein:global_normalized_Fourier:global_normalization:global_functional_equation}
We have
    \begin{align}
    E^*(h,-s)^{\circ}_n & = \pm E^*(h,s)^{\circ}_n
    \end{align}
where the sign ``$\pm$'' is
    \begin{equation}
    (-1)^{m(m-1)(n-m-1)[F^+ : \Q]/2} \eta(\mf{d}_{F^+/\Q})^{m(n - m - 1)}.
    \end{equation}
\end{lemma}
\begin{proof}
Let $\psi \colon F^+ \backslash \A \ra \C^{\times}$
be any nontrivial additive character.
Given $T \in \mrm{Herm}_m(F^+)$ with $\det T \neq 0$, the local functional equations (Section \ref{sec:Eisenstein:local_functional_equations}) and the factorization from \eqref{equation:Eisenstein:global_normalized_Fourier:global_normalization:factorization} imply
    \begin{equation}\label{equation:Eisenstein:global_normalization:Fourier_functional_equation}
    E^*_{T,\psi}(h,-s)^{\circ}_n = (-1)^{m(m-1)(n - m - 1)[F^+ : \Q]/2} \eta(\mf{d}_{F^+/\Q})^{m(n - m - 1)} E^*_{T,\psi}(h,s)^{\circ}_n.
    \end{equation}
The global functional equation \eqref{equation:Eisenstein:setup:adelic_v_classical:global_functional_equation} implies that $E^*(h,-s)^{\circ}_n = f(s) E^*(h,s)^{\circ}_n$ for some meromorphic function $f(s)$ (temporary notation) independent of $T$.
There exists $T$ with $\det T \neq 0$ and $E^*_{T,\psi}(h,s)^{\circ}_n$ not identically zero,
e.g. we can take $T = 1_m$ if $\psi$ is the standard character,
as follows from \eqref{equation:Eisenstein:local_Whittaker:Archimedean:special_value}
and Lemma \ref{lemma:non-Arch_Whittaker_polynomial}(2)(iii).
So $f(s)$ is identically $1$ and \eqref{equation:Eisenstein:global_normalization:Fourier_functional_equation} holds for all $T \in \mrm{Herm}_m(F^+)$.
\end{proof}
    
            \subsection{Singular Fourier coefficients} 
            \label{ssec:Eisenstein:global_normalized_Fourier:singular_Fourier}
                Retain notation and assumptions from Section \ref{ssec:Eisenstein:global_normalized_Fourier:global_normalization}. In this section, the main result is Corollary \ref{corollary:Eisenstein:singular_Fourier:corank_1} on singular Fourier terms of co-rank $1$.

We use various subscripts to emphasize $m$-dependence (in the implicit $U(m,m)$). For example, we write $\Phi_{m,v}^{\circ}$ rather than just $\Phi_v^{\circ}$ for non-Archimedean $v$ (resp. $\Phi_{m,v}^{(n)}$ instead of $\Phi_v^{(n)}$ for Archimedean $v$), similarly $\Phi_m^{(n) \circ}$ instead of $\Phi^{(n) \circ}$ for the global standard section from Section \ref{ssec:Eisenstein:global_normalized_Fourier:global_normalization}, also $M_m(s, \chi)$ instead of $M(s, \chi)$ for the intertwining operator, etc..

Suppose $m \geq 1$ and set $m^{\flat} = m - 1$. Recall the operators $\mu^{m *}_{m^{\flat}}(s, \chi)$, $M_m(s, \chi)$, $M_{m^{\flat}}(s, \chi)$ and $U^m_{m^{\flat}}(s, \chi)$ as in Section \ref{ssec:Eisenstein:setup:singular_Fourier}.

\begin{lemma}\label{lemma:Eisenstein:singular_Fourier:operators}
We have
    \begin{align}
    \mu_{m^{\flat}}^{m *}(s, \chi) \Phi_m^{(n) \circ}(s) & = \Phi^{(n) \circ}_{m^{\flat}}(s + 1/2)
    \\
    \label{equation:Eisenstein:singular_Fourier:operators:U-operator}
    U^m_{m^{\flat}}(s, \chi) \Phi_m^{(n) \circ}(s) & = \pm \frac{\Lambda_{m^{\flat}}(s - 1/2)^{\circ}_n \Lambda_m(-s)^{\circ}_n}{\Lambda_m(s)^{\circ}_n \Lambda_{m^{\flat}}(- s + 1/2)^{\circ}_n} \Phi_{m^{\flat}}^{(n) \circ}(s - 1/2)
    \\
    M_{m}(s, \chi) \Phi_{m}^{(n) \circ}(s) & = |\Delta_{F^+ / \Q}|^{- m^2 / 2} |N_{F^+/\Q}(\Delta)|^{-m(m-1)/4} ( (-1)^{nm} i^{m(m-1)/2} \pi^{2 m s} )^{[F^+ : \Q]}
    \\
    & \mathrel{\phantom{=}} \cdot \left ( \prod_{j=0}^{m-1} \frac{L(2s + j - m + 1, \eta^{n + j})}{L(2s + m - j, \eta^{n + j})} \right ) \notag
    \\
    & \mathrel{\phantom{=}} \cdot \left ( \prod_{v \mid \infty} \prod_{j = 0}^{m-1} \frac{L_v(2s + j - m + 1, \eta_v^{n+j}) }{\epsilon_v(2s + j - m + 1, \eta_v^{n+j}, \overline{\psi}_v) L_v(-2s - j + m, \eta_v^{n+j}) } \right ) \notag
    \\
    & \mathrel{\phantom{=}} \cdot \left ( \prod_{j = 0}^{m - 1} \frac{\Gamma(-s - s_0 + n - j)}{\Gamma(s - s_0 + n - j)}\right )^{[F^+ : \Q]} \notag
    \\
    & \mathrel{\phantom{=}} \cdot \Phi^{(n) \circ}_m(-s), \notag
    \end{align}
allowing $m = 0$ for $M_m(s, \chi)$,  where $\psi_v(x) \coloneqq e^{2 \pi i x}$ for $v \mid \infty$, and
where the sign ``$\pm$'' in \eqref{equation:Eisenstein:singular_Fourier:operators:U-operator} is
    \begin{equation}
    (-1)^{(m(m-1)(n-m-1)/2 - m^{\flat} (m^{\flat} - 1)(n - m^{\flat} - 1)/2)[F^+ : \Q]}
    \eta(\mf{d}_{F^+ / \Q})^{m(n - m - 1)} \eta(\mf{d}_{F^+ / \Q})^{m^{\flat}(n - m^{\flat} - 1)}.
    \end{equation}
\end{lemma}
\begin{proof}
Each identity holds a priori up to a meromorphic scale factor. We may compute this scale factor by evaluating both sides at $1 \in U(m^{\flat}, m^{\flat})$ or $1 \in U(m,m)$ as appropriate.

The identity for $\mu^{m*}_{m^{\flat}}(s, \chi)$ is then clear. For $M_m(s, \chi)$, the identity follows directly upon combining \eqref{equation:Eisenstein:local_functional_equations:non-Arch:spherical_intertwining} and \eqref{equation:Eisenstein:local_functional_equations:Archimedean:intertwining_on_scalar_weight},
e.g. considering the standard character $\psi_{\mrm{tr}}$ (note that the Lemma statement
does not need to mention a choice of global additive character).

Define the temporary notation $\a_m(s)_n$ for the meromorphic function (in the lemma statement) satisfying $M_m(s, \chi) \Phi^{(n) \circ}_m(s) = \a_m(s)_n \Phi^{(n) \circ}_m(-s)$. 
By \eqref{equation:Eisenstein:setup:singular_Fourier}, proving the claimed identity for $U^m_{m^{\flat}}(s, \chi)$ is equivalent to showing
    \begin{equation}
    \frac{\a_{m}(s)_n}{\a_{m^{\flat}}(s - 1/2)_n} = \pm \frac{\Lambda_{m^{\flat}}(s - 1/2)^{\circ}_n \Lambda_m(-s)^{\circ}_n}{\Lambda_m(s)^{\circ}_n \Lambda_{m^{\flat}}(- s + 1/2)^{\circ}_n}
    \end{equation}
with sign $\pm$ as in the lemma statement. This may be computed explicitly as follows. Some rearranging yields
    \begin{align*}
    \frac{\Lambda_{m^{\flat}}(s - 1/2)^{\circ}_n \Lambda_m(-s)^{\circ}_n}{\Lambda_m(s)^{\circ}_n \Lambda_{m^{\flat}}(- s + 1/2)^{\circ}_n} 
    & =
    (\pi^{2 m s} \pi^{(-2s + 1) m^{\flat}})^{[F^+ : \Q]} |N_{F^+ / \Q}(\Delta^{\lfloor m/2 \rfloor})|^{-2s} |N_{F^+ / \Q}(\Delta^{\lfloor m^{\flat}/2 \rfloor})|^{2s-1}
    \\
    & \mathrel{\phantom{=}} \cdot \Gamma(s - s_0 + n)^{-[F^+ : \Q]} \Gamma(- s - s_0 + n - m + 1)^{[F^+ : \Q]}
    \\
    & \mathrel{\phantom{=}} \cdot L(2s, \eta^{m + n}) L(2s + m, \eta^n)^{-1} L(2s + m -1, \eta^{n+1})^{-1}
    \\
    & \mathrel{\phantom{=}} \cdot L(-2s + 1, \eta^{m + n + 1}) \cdot |\Delta_{F^+ / \Q}|^{-2 s - m + 1}
    \end{align*}
and
    \begin{align*}
    \frac{\a_{m}(s)_n}{\a_{m^{\flat}}(s - 1/2)_n} & = |\Delta_{F^+ / \Q}|^{-m + 1/2} |N_{F^+/\Q}(\Delta)|^{-(m - 1)/2} ( (-1)^{n} i^{m-1} \pi^{2 m s - 2 m^{\flat}(s - 1/2)} )^{[F^+ : \Q]}
    \\
    & \mathrel{\phantom{=}} \cdot L(2s, \eta^{n + m + 1}) \cdot L(2s + m, \eta^{n})^{-1} L(2s + m - 1, \eta^{n + 1})^{-1} L(2s, \eta^{n + m})
    \\
    & \mathrel{\phantom{=}} \cdot \left ( \prod_{v \mid \infty} \frac{L_v(2s, \eta_v^{n + m + 1})}{\epsilon_v(2s, \eta_v^{n + m + 1}, \overline{\psi}_v) L_v(-2s + 1, \eta_v^{n + m + 1})} \right )
    \\
    & \mathrel{\phantom{=}} \cdot \Gamma(- s - s_0 + n - m + 1)^{[F^+ : \Q]} \Gamma(s - s_0 + n)^{-[F^+ : \Q]}.
    \end{align*}
Let $\psi = \psi_{\mrm{tr}} \colon F^+ \backslash \A \ra \C^{\times}$.
Using the global functional equation $\Lambda(s, \eta^{n + m + 1}) = \epsilon(s, \eta^{n + m + 1}) \Lambda(1 - s, \eta^{n + m + 1})$ (notation as in Section \ref{ssec:local_Whittaker:local_L-factors}),
we find
    \begin{align*}
    \frac{\a_{m}(s)_n}{\a_{m^{\flat}}(s - 1/2)_n} & = |\Delta_{F^+ / \Q}|^{-m + 1/2} |N_{F^+/\Q}(\Delta)|^{-(m - 1)/2} ( (-1)^{n} i^{m-1} \pi^{2 m s - 2 m^{\flat}(s - 1/2)} )^{[F^+ : \Q]}
    \\
    & \mathrel{\phantom{=}} \cdot L(-2s +1, \eta^{n + m + 1}) \cdot L(2s + m, \eta^{n})^{-1} L(2s + m - 1, \eta^{n + 1})^{-1} L(2s, \eta^{n + m})
    \\
    & \mathrel{\phantom{=}} \cdot \left ( \prod_{v < \infty} \epsilon_v(2s, \eta_v^{n + m + 1}, \overline{\psi}_v) \right )
    \\
    & \mathrel{\phantom{=}} \cdot \Gamma(- s - s_0 + n - m + 1)^{[F^+ : \Q]} \Gamma(s - s_0 + n)^{-[F^+ : \Q]}.
    \end{align*}
It remains to show
    \begin{align}\label{equation:Eisenstein:singular_Fourier:operators:numerology}
    \begin{split}
    & |N_{F^+ / \Q}(\Delta)|^{-(m-1)/2} ((-1)^n i^{m - 1})^{[F^+ : \Q]} \prod_{v < \infty} \epsilon_v(2s, \eta_v^{n + m + 1}, \overline{\psi}_v)
    \\
    & = \pm |N_{F^+ / \Q}(\Delta^{\lfloor m/2 \rfloor})|^{-2s} |N_{F^+ / \Q}(\Delta^{\lfloor m^{\flat}/2 \rfloor})|^{2s-1} \cdot |\Delta_{F^+ / \Q}|^{-2 s + 1/2}
    \end{split}
    \end{align}
with sign $\pm$ as in the lemma statement.
For studying this sign, we find it convenient to rewrite
    \begin{align}\label{equation:Eisenstein:singular_Fourier:operators:numerology_sign_partial}
    m (m - 1)(n - m - 1) - m^{\flat} (m^{\flat} - 1) (n - m^{\flat} - 1) = (m - 1)(- m - (- 2n + 4)).
    \end{align}
The rest of the proof consists of casework on
the parity of $m$ and $n$ (four cases).

First consider the case when $n + m + 1$ is odd.
We compute
    \begin{align}
    \prod_{v < \infty} \epsilon_v(2s, \eta_v^{n + m + 1}, \overline{\psi}_v)
    &= |\Delta_{F^+ / \Q}|^{1/2 - 2s} |N_{F^+ / \Q}(\Delta)|^{-2s + 1/2} \prod_{v < \infty} \gamma_{\overline{\psi}_v}(F_v)
    \\
    & = |\Delta_{F^+ / \Q}|^{1/2 - 2s} |N_{F^+ / \Q}(\Delta)|^{-2s + 1/2} \prod_{v \mid \infty} \gamma_{\overline{\psi}_v}(F_v)^{-1}
    \\
    & = |\Delta_{F^+ / \Q}|^{1/2 - 2s} |N_{F^+ / \Q}(\Delta)|^{-2s + 1/2} i^{[F^+ : \Q]}
    \end{align}
where the first equality follows from
\eqref{equation:Eisenstein:Weil_rep:Weil_index:epsilon_factors_and_Weil_index},
the second equality follows from the global product formula for Weil indices,
and the third equality follows from the formula for Archimedean Weil indices
given after \eqref{equation:Eisenstein:Weil_rep:Weil_index:Archimedean},
i.e. $\gamma_{\overline{\psi}_v}(\C) = -i$.

Suppose $n$ and $m$ are both odd. Since $n$ is odd,
we must have $|N_{F^+ / \Q}(\Delta)| = 1$ (by convention).
Then \eqref{equation:Eisenstein:global_normalized_Fourier:global_normalization:different_ideal_sign}
gives
    \begin{equation}\label{equation:Eisenstein:singular_Fourier:operators:mn_odd:numerology_sign_different}
    \prod_{v < \infty} \gamma_{\overline{\psi}_v}(F_v)
    =
    \eta(\mf{d}_{F^+ / \Q})^{m(n - m - 1)} \eta(\mf{d}_{F^+ / \Q})^{m^{\flat}(n - m^{\flat} - 1)}.
    \end{equation}
Since the above quantity (lying in $\pm 1$) is also $i^{[F^+ : \Q]}$, we conclude
that $[F^+ : \Q]$ is even,
so $((-1)^n i^{m - 1})^{[F^+ : \Q]} = +1$.
Since $[F^+ : \Q]$ is even, we also find that the sign $\pm$
of the lemma statement is given by
\eqref{equation:Eisenstein:singular_Fourier:operators:mn_odd:numerology_sign_different},
so \eqref{equation:Eisenstein:singular_Fourier:operators:numerology} now follows.

Suppose $n$ and $m$ are both even.
The previous computation shows that the left-hand side
of \eqref{equation:Eisenstein:singular_Fourier:operators:numerology}
equals $|N_{F^+ / \Q}(\Delta)|^{-(m-1)/2 - 2s + 1/2} (-1)^{m \cdot [F^+ : \Q] / 2} |\Delta_{F^+ / \Q}|^{1/2 - 2s}$.
The factors of $|N_{F^+ / \Q}(\Delta)|$
on both sides of
\eqref{equation:Eisenstein:singular_Fourier:operators:numerology}
agree because
$-(m-1)/2 - 2s + 1/2 = \lfloor m / 2 \rfloor(- 2 s) + \lfloor m^{\flat} / 2 \rfloor (2s - 1)$.
The signs also agree because
$m / 2$ has the same parity as
$(m - 1)(- m - (- 2n + 4)) / 2$,
where the latter comes from \eqref{equation:Eisenstein:singular_Fourier:operators:numerology_sign_partial}.

Next consider the case where $n + m + 1$ is even.
We compute
    \begin{align}
    \prod_{v < \infty} \epsilon_v(2s, \eta_v^{n + m + 1}, \overline{\psi}_v)
    & = |\Delta_{F^+ / \Q}|^{1/2 - 2s}
    \end{align}
using \eqref{equation:Eisenstein:Weil_rep:Weil_index:epsilon_factors_and_Weil_index}
and Lemma \ref{lemma:Weil_index:Hermitian}(2).

Suppose $n$ is odd and $m$ is even.
Since $n$ is odd, we must have $|N_{F^+ / \Q}(\Delta)| = 1$.
As in the case where $n$ and $m$ were both odd,
we find
    \begin{equation}
    i^{[F^+ : \Q]}
    =
    \eta(\mf{d}_{F^+ / \Q})^{m(n - m - 1)} \eta(\mf{d}_{F^+ / \Q})^{m^{\flat}(n - m^{\flat} - 1)}
    \end{equation}
and that $[F^+ : \Q]$ is even.
To see that the signs in \eqref{equation:Eisenstein:singular_Fourier:operators:numerology}
agree, it thus suffices to note
$((-1)^n i^m)^{[F^+ : \Q]} = +1$.

Lastly, suppose $n$ is even and $m$ is odd.
The factors of $|N_{F^+ / \Q}(\Delta)|$ on both
sides of \eqref{equation:Eisenstein:singular_Fourier:operators:numerology}
agree since $(m - 1)/2 = \lfloor m / 2 \rfloor = \lfloor m^{\flat} / 2 \rfloor$
because $m$ is odd.
The signs also agree because $(m - 1)/2$ and
$(m - 1)(- m - (- 2n + 4)) / 2$ have the same parity,
where the latter comes from \eqref{equation:Eisenstein:singular_Fourier:operators:numerology_sign_partial}.
\end{proof}

\begin{corollary}\label{corollary:Eisenstein:singular_Fourier:corank_1}
Consider any $a = \mrm{diag}(a^{\#}, a^{\flat}) \in \GL_m(\A_F)$ with $a^{\#} \in \GL_1(\A_F)$ and $a^{\flat} \in \GL_{m^{\flat}}(\A_F)$.
Take any nontrivial additive character $\psi \colon F^+ \backslash \A \ra \C^{\times}$.
For any $T \in \mrm{Herm}_m(F^+)$ with $\rank T = m - 1$ and $T = \mrm{diag}(0, T^{\flat})$ being block diagonal with $\det T^{\flat} \neq 0$, we have
    \begin{align*}
    \tilde{E}_T^*(a, s)^{\circ}_n & = |\det a^{\#}|_{F}^{s - s_0} \frac{\Lambda_m(s)^{\circ}_n}{\Lambda_{m^{\flat}}(s + 1/2)^{\circ}_n} \tilde{E}^*_{T^{\flat}}(a^{\flat}, s + 1/2)^{\circ}_n 
    \\
    & \mathrel{\phantom{=}} \pm | \det a^{\#} |_F^{- s - s_0} \frac{\Lambda_m(- s)^{\circ}_n}{\Lambda_{m^{\flat}}(- s + 1/2)^{\circ}_n} \tilde{E}^*_{T^{\flat}}(a^{\flat}, s - 1/2)^{\circ}_n
    \end{align*}
with sign ``$\pm$'' as in Lemma \ref{lemma:Eisenstein:singular_Fourier:operators}.
\end{corollary}
\begin{proof}
This follows immediately from Lemma \ref{lemma:Eisenstein:singular_Fourier:operators}, \eqref{equation:Eisenstein:setup:singular_Fourier:corank_1}, and the definition of the normalized Fourier coefficients $\tilde{E}^*_T(a,s)^{\circ}_n$ and $\tilde{E}^*_{T^{\flat}}(a^{\flat}, s)^{\circ}_n$ (Section \ref{ssec:Eisenstein:global_normalized_Fourier:global_normalization}).
\end{proof}

\begin{remark}
In the situation of Corollary \ref{corollary:Eisenstein:singular_Fourier:corank_1}, the functional equation 
    \begin{equation}
    \tilde{E}^*_T(a, s)^{\circ}_n = (-1)^{m(m-1)(n - m - 1) [F^+ : \Q]/2} \eta(\mf{d}_{F^+ / \Q})^{m(n - m - 1)}\tilde{E}^*_T(a, -s)^{\circ}_n
    \end{equation}
for the singular matrix $T$ is a visible consequence of the functional equation
    \begin{equation}
    \tilde{E}^*_{T^{\flat}}(a^{\flat}, s)^{\circ}_n = (-1)^{m^{\flat}(m^{\flat} - 1)(n - m^{\flat} - 1) [F^+ : \Q]/2}  \eta(\mf{d}_{F^+ / \Q})^{m^{\flat}(n - m^{\flat} - 1)} \tilde{E}^*_{T^{\flat}}(a^{\flat}, -s)^{\circ}_n
    \end{equation}
for the nonsingular matrix $T^{\flat}$.
\end{remark}

    \clearpage


    \part{Siegel--Weil}
    \label{part:part_IV:Siegel--Weil}
        Our main results (arithmetic Siegel--Weil) are in Section \ref{sec:arithmetic_Siegel-Weil}. We also give some explicit formulas for special values (local Siegel--Weil and geometric Siegel--Weil) in Sections \ref{sec:local_Siegel-Weil} and \ref{sec:geometric_Siegel-Weil}. These special value formulas will be needed as ingredients in the proofs of our arithmetic Siegel--Weil results.

        \section{Local Siegel--Weil} 
        \label{sec:local_Siegel-Weil}
    
            \subsection{Volume forms} 
            \label{ssec:local_Siegel-Weil:volume_forms}
                Given a scheme $X$ which is smooth and equidimensional over a field $A$, a \emph{volume form} (or \emph{gauge form}) on $X$ will mean a nowhere vanishing (algebraic) differential form of top degree on $X$. When $X$ is also affine and $A$ is a local field, the set $X(A)$ has the natural structure of an $A$-analytic manifold (in the sense of \cite[{Part II, Chapter III}]{Serre64}). In this case, a volume form on $X$ defines a Borel measure on $X(A)$ in a standard way (see \cite[{\S 2.2}]{Weil82}).

We use volume forms to normalize various Haar measures. Let $B$ be a degree $2$ \'etale algebra over a field $A$ of characteristic $\neq 2$, and write $b \mapsto \overline{b}$ for the nontrivial involution on $B$.
Let $V$ be a $B/A$ Hermitian space which is free of rank $n$, and set $G = U(V)$. Fix a nonnegative integer $m \leq n$, and choose translation-invariant volume forms $\a$ and $\b$ on $V^m$ and $\mrm{Herm}_m$ respectively (viewed as group schemes over $A$). The forms $\a$ and $\b$ have degrees $2nm$ and $m^2$ respectively. 

Consider the moment map
    \begin{equation}
    \begin{tikzcd}[row sep = tiny]
    V^m \arrow{r}{\mc{T}} & \mrm{Herm}_m \\
    \underline{x} \arrow[mapsto]{r} & (\underline{x},\underline{x}).
    \end{tikzcd}
    \end{equation}
We assume $n \geq m$, and write $V^m_{\mrm{reg}} \subseteq V^m$ for the open subscheme where $\det \mc{T}$ is invertible. A tangent space calculation shows that $\mc{T}$ is smooth when restricted to $V^m_{\mrm{reg}}$.

Given $T \in \mrm{Herm}_m(A)$, we write $\Omega_T \subseteq V^m$ for the fiber of the moment map over $T$. If $\underline{x} \in V^m(A)$ has Gram matrix $T = (\underline{x}, \underline{x})$, then $g \mapsto g^{-1} \cdot \underline{x}$ defines a morphism $\iota_{\underline{x}} \colon G \ra \Omega_T$. If $\det T$ is invertible, then a dimension count and tangent space calculation show that $\iota_{\underline{x}}$ is smooth. If $\det T$ is invertible, if $A$ is a local field, and if $G_{\underline{x}} \subseteq G$ denotes the stabilizer of $\underline{x}$, then the induced map $G_{\underline{x}}(A) \backslash G(A) \ra \Omega_T(A)$ is a homeomorphism (surjectivity is from Witt's theorem, and openness is from the submersivity of $G(A) \ra \Omega_T(A)$, which in turn comes from smoothness of $\iota_{\underline{x}}$).

 In part (3) of the lemma below, $x \in V^m_{\mrm{reg}}$ means $x \in V^m_{\mrm{reg}}(S)$
 for some suppressed $A$-scheme $S$, and we similarly abuse notation in part (2).
 In part (3), we also use the usual notation $\G_a$ for the additive group scheme. 

\begin{lemma}\label{lemma:local_Siegel-Weil:volume_forms:exists_nu}
There exists an algebraic differential form $\nu$ on $V^m_{\mrm{reg}}$ of degree $m(2n - m)$ satisfying the following conditions.
    \begin{enumerate}[(1)]
        \item We have $\a = \mc{T}^*(\b) \wedge \nu$.
        \item For the $G \times \Res_{B/A} \GL_m$ action on $V^m_{\mrm{reg}}$ given by $x \mapsto g x h^{-1}$ for $(g,h) \in G \times \Res_{F/F^+} \GL_m$, we have $(g, h)^* \nu = \det({}^t \overline{h} h)^{m-n} \nu$.
        \item For each $x \in V^m_{\mrm{reg}}$, the restriction of $\nu \colon T_x(V^m) \ra \G_a$ to $\ker d \mc{T}_x$ is nonzero.
        \item For any fixed non-degenerate subspace $V^{\flat} \subseteq V$ which is free of rank $m$, and for $\underline{x} \in V^{\flat m}_{\mrm{reg}}(A)$, the differential form
            \begin{equation}
            \det(\underline{x}, \underline{x})^{m - n} \iota_{\underline{x}}^* \nu
            \end{equation}
        on $G$ is independent of the choice of $\underline{x}$. This form is right $G$-invariant.
    \end{enumerate}
\end{lemma}
\begin{proof}
The case $m = n$ is stated in \cite[{\S 10}]{KR14}. The analogue of that case for orthogonal groups is discussed in \cite[{Lemmas 5.3.1, 5.3.2}]{KRY06} (there stated and proved for three-dimensional quadratic spaces). The present lemma may be proved by a similar computation.

Part (4) follows from part (2) (where ``non-degenerate subspace'' means that the restriction of the Hermitian pairing is non-degenerate).
\end{proof}
                
            \subsection{Special value formula} 
            \label{ssec:local_Siegel-Weil:local_Siegel-Weil}
                We retain notation from Section \ref{ssec:local_Siegel-Weil:volume_forms}, and specialize to the case where $B/A$ is the extension $F_v/F^+_v$ where $F^+_v$ is a local field of characteristic $\neq 2$. If $F^+_v$ is Archimedean, we assume $F_v/F^+_v$ is $\C / \R$. We often use subscripts $v$ to emphasize $F^+_v$ being a local field, e.g. we write $\underline{x}_v$ for elements of $V^m_{\mrm{reg}}(F^+_v)$.

Fix a nontrivial additive character $\psi_v \colon F^+_v \ra \C^{\times}$. We write $d b_v$ for the self-dual Haar measure on $\mrm{Herm}_m(F^+_v)$ with respect to the trace pairing $(b,c) \mapsto \psi_v(\mrm{tr}(bc))$. We also write $d \underline{x}_v$ for the self-dual Haar measure on $V^m(F^+_v)$ with respect to the pairing $\psi_v(\mrm{tr}_{F_v/F^+_v}(\mrm{tr}(-,-)))$.

Fix translation-invariant volume forms $\a$ and $\b$ as in Section \ref{ssec:local_Siegel-Weil:volume_forms}. These determine Haar measures $d_{\b} b_v$ and $d_{\a} \underline{x}_v$ on $\mrm{Herm}_m(F^+_v)$ and $V^m(F^+_v)$ respectively. Define positive real constants $c_v(\a, \psi_v)$ and $c_v(\b, \psi_v)$ such that
    \begin{equation}
    d_{\a} \underline{x}_v = c_v(\a, \psi_v) d \underline{x}_v \quad \quad d_{\b} b_v = c_v(\b, \psi_v) d b_v.
    \end{equation}
    
Suppose $T \in \mrm{Herm}_m(F^+_v)$ is a matrix with $\det T \neq 0$. For the rest of Section \ref{ssec:local_Siegel-Weil:local_Siegel-Weil}, fix a differential form $\nu$ as in Lemma \ref{lemma:local_Siegel-Weil:volume_forms:exists_nu}. The restriction of $\det (T)^{m - n} \nu$ to $\Omega_T$ is a $G$-invariant volume form on $\Omega_T$, and we write $d_{T, \nu} \underline{x}_v$ for the resulting measure on $\Omega_T(F^+_v)$.

Take any character $\chi_v \colon F^{\times}_v \ra \C^{\times}$
with $\chi_v |_{F^{+ \times}_v} = \eta_v^n$,
where $\eta_v \colon F^{+ \times}_v \ra \{ \pm 1\}$
is our usual notation for the quadratic character associated to $F_v / F^+_v$.
It is known that there exists a constant $c_{T,v}$ (depending on $T$, the measure $d_{T, \nu} \underline{x}_v$, and the character $\psi_v$) such that
    \begin{equation}\label{equation:local_Siegel-Weil:local_Siegel-Weil:up-to-constant}
    W_{T,v}(s_0, \Phi_{\varphi_v}) = c_{T,v} \int_{\Omega_T(F^+_v)} \varphi_v(\underline{x}_v) ~ d_{T, \nu} \underline{x}_v
    \end{equation}
holds for any Schwartz function $\varphi_v \in \mc{S}(V^m(F^+_v))$ (see \cite[{Lemma 5.1, Lemma 5.2}]{Ichino04}). Here we set $s_0 \coloneqq (n - m)/2$ as usual, and $\Phi_{\varphi_v}$ is the Siegel--Weil section associated with $\varphi_v$ (Section \ref{ssec:Eisenstein:Weil_rep:Weil_rep}).
If $\Omega_T(F^+_v) = \emptyset$, we thus have $W_{T,v}(s_0, \Phi_{\varphi_v}) = 0$ for all $\varphi_v$.
The function $W_{T,v}(s, \Phi_{\varphi_v})$ indeed does not depend on the choice of $\chi_v$
by Lemma \ref{lemma:local_density_Whittaker_interpolation}, as the right-hand side
of \eqref{equation:local_density_Whittaker_interpolation} has no $\chi_v$ dependence.\footnote{
    Strictly speaking, we required $\psi_v$ to be unramified in loc. cit..
    However, since $V_v$ in loc. cit. was arbitrary,
    we can reduce the case of general $\psi_v$ to the case of
    unramified $\psi_v$ by rescaling as in
    Remark \ref{remark:Eisenstein:Weil_rep:Weil_rep:rescaling_trick}.
}

We may compute the constant $c_{T,v}$ by evaluating \eqref{equation:local_Siegel-Weil:local_Siegel-Weil:up-to-constant} on any nonzero nonnegative Schwartz function $\varphi_v$. We may take $\varphi_v$ to have support which is compact and contained in $V^m_{\mrm{reg}}(F^+_v)$. The relation $\a = \mc{T}^*(\b) \wedge \nu$ and an ``integrate along the fibers of $\mc{T}$'' computation (similar to the proof of \cite[{Proposition 5.3.3}]{KRY06}) gives
    \begin{equation}\label{equation:local_Siegel-Weil:local_Siegel-Weil:c_Tv}
    c_{T,v} = \frac{\gamma_{\psi_v}(V)^{-m} c_v(\b, \psi_v)}{c_v(\a, \psi_v)} |\det T|_{F^+_v}^{n - m}.
    \end{equation}
Here $\gamma_{\psi_v}(V)$ is the Weil index, as appearing in the Weil representation (Section \ref{ssec:Eisenstein:Weil_rep:Weil_rep}).

\begin{lemma}[Local Siegel--Weil]\label{lemma:local_Siegel-Weil:local_Siegel-Weil:unique_Haar_measure}
Let $V$ be an $F_v/F^+_v$ Hermitian space of rank $n$, and let $\psi_v \colon F^+_v \ra \C^{\times}$ be a nontrivial additive character. Fix a non-degenerate subspace $V^{\flat} \subseteq V$ which is free of rank $m$, and fix a Haar measure on $U(V^{\flat \perp})(F^+_v)$. 

There exists a unique Haar measure $d g_v$ on $G(F^+_v)$ such that, for any basis $\underline{x}_v \in V^{\flat m}$ of $V^{\flat}$ and any Schwartz function $\varphi_v \in \mc{S}(V^m(F^+_v))$, we have
	\begin{equation}
        W_{T,v}(s_0, \Phi_{\varphi_v}) = \gamma_{\psi_v}(V)^{-m} |\det T|_{F^+_v}^{n-m}  \int_{G_{\underline{x}_v}(F^+_v) \backslash G(F^+_v)} \varphi_v (g_v^{-1} \underline{x}_v) ~ d g_v
	\end{equation}
for the corresponding quotient measure, where $T = (\underline{x}_v, \underline{x}_v)$ is the Gram matrix of $\underline{x}_v$ (and where the Haar measure on $G_{\underline{x}_v}(F^+_v)$ is induced by the canonical identification $G_{\underline{x}_v} \cong U(V^{\flat \perp})$).
\end{lemma}
\begin{proof}
Select any basis $\underline{x}_v$ of $V^{\flat}$. Set $\omega_1 = \det(\underline{x}_v, \underline{x}_v)^{m - n} \iota_{\underline{x}}^* \nu$ (temporary notation). 
We know $\omega_1$ does not depend on the choice of $\underline{x}_v$, by Lemma \ref{lemma:local_Siegel-Weil:volume_forms:exists_nu}(4). Let $\omega_2$ be a right $G$-invariant differential form of degree $(n - m)^2$ on $G$ such that $\omega_1 \wedge \omega_2$ is a nowhere vanishing differential form of top degree $n^2$ (also right $G$-invariant).
The volume form $\omega_1 \wedge \omega_2$ on $G$ defines a Haar measure on $G(F^+_v)$. The restriction $\omega_2|_{G_{\underline{x}}}$ is a volume form on $G_{\underline{x}}$ (by smoothness of $\iota_{\underline{x}}$), and defines a Haar measure on $G_{\underline{x}_v}(F^+_v)$. The resulting quotient measure on $G_{\underline{x}}(F^+_v) \backslash G(F^+_v) \cong \Omega_T(F^+_v)$ is precisely the measure for the volume form $(\det T)^{m - n} \nu|_{\Omega_T}$ on $\Omega_T$.

The lemma then follows from \eqref{equation:local_Siegel-Weil:local_Siegel-Weil:up-to-constant} and the constant calculated in \eqref{equation:local_Siegel-Weil:local_Siegel-Weil:c_Tv}.
\end{proof}

\begin{remark}\label{remark:local_Siegel-Weil:unique_Haar_measure:invariance}
Consider the situation of Lemma \ref{lemma:local_Siegel-Weil:local_Siegel-Weil:unique_Haar_measure}, and suppose $V^{\flat \prime} \subseteq V$ is a subspace which is isomorphic to $V^{\flat}$ as a Hermitian space. Suppose $f_v \in U(V)(F^+_v)$ satisfies $f_v(V^{\flat}) = V^{\flat \prime}$, and equip $U(V^{\flat \prime \perp})(F^+_v)$ with the Haar measure induced from $U(V^{\flat \perp})(F^+_v)$ via $f_v$. If $d g_v$ and $d g'_v$ are the induced Haar measures on $G(F^+_v)$ corresponding to $V^{\flat}$ and $V^{\flat \prime}$ respectively (Lemma \ref{lemma:local_Siegel-Weil:local_Siegel-Weil:unique_Haar_measure}), a change of variables shows $d g_v = d g'_v$.
\end{remark}

            \subsection{Explicit Haar measures} 
            \label{ssec:local_Siegel-Weil:non-Arch}
                For our application to uniformization of special cycles (Section \ref{ssec:local_Siegel-Weil:uniformization_degree}), we need to explicitly compute the Haar measures from Lemma \ref{lemma:local_Siegel-Weil:local_Siegel-Weil:unique_Haar_measure} in a few cases. The main result of this subsection is Lemma \ref{lemma:local_Siegel-Weil:non-Arch:spherical_and_measures}, and the other lemmas are auxiliary.

We retain notation from Section \ref{ssec:local_Siegel-Weil:local_Siegel-Weil}. In addition, we assume that $F^+_v$ is non-Archimedean and that $\psi_v$ is unramified. 
Let $\varpi_0$ be a uniformizer of $F^+_v$.
For $F_v/F^+_v$ ramified in Section \ref{ssec:local_Siegel-Weil:non-Arch},
we assume that $F^+_v$ has residue characteristic $\neq 2$
and write $\varpi \in F_v$ for a uniformizer satisfying $\varpi^{\s} = - \varpi$,
where $\varpi^{\s}$ denotes the Galois conjugate of $\varpi$ under the Galois involution $\s$
of $F_v$ over $F^+_v$. 

Let $M^{\circ}_2$ be the rank $2$ self-dual lattice described in Section \ref{ssec:Eisenstein:Weil_rep:Weil_rep}, and write $U(M^{\circ}_2)$ for the group of (unitary) automorphisms of $M^{\circ}_2$. Let $q_v$ be the residue cardinality of $F^+_v$.

The next lemma should be compared with Witt's theorem for lattices with quadratic forms, as in \cite{Morin-Strom79}.

\begin{lemma}\label{lemma:local_Siegel-Weil:non-Arch:lattice_Witt_theorem}
For any given $c \in \mc{O}_{F^+_v}^{\times}$, the group $U(M^{\circ}_2)$ acts transitively on the set
    \begin{equation}
    \{ x \in M^{\circ}_2 : (x, x) = c \}.
    \end{equation}
If $F_v/F^+_v$ is inert, the same holds for any $c \in \varpi_0 \mc{O}_{F^+_v}^{\times}$.
\end{lemma}
\begin{proof}
Given $y \in M^{\circ}_2$, we write $\langle y \rangle \subseteq M^{\circ}_2$ for the submodule generated by $y$. If $F_v/F^+_v$ is ramified, we view $\varpi$ as a generator of the different ideal $\mf{d}$, and we otherwise view $1$ as a generator of $\mf{d}$. Choose a basis $e_1, e_2$ of $M^{\circ}_2$ with Gram matrix given by \eqref{equation:Eisenstein:Weil_rep:Weil_rep:standard_self-dual_lattice}. In this basis, we also consider the elements
    \begin{equation}
    w' = \begin{pmatrix} 0 & 1 \\ \epsilon & 0 \end{pmatrix} 
    \quad \quad m(a) = \begin{pmatrix} a & 0 \\ 0 & \overline{a}^{-1} \end{pmatrix} 
    \quad \quad n(b) = \begin{pmatrix} 1 & b \\ 0 & 1 \end{pmatrix}
    \quad \quad \epsilon = \begin{cases} -1 & \text{if $F_v/F^+_v$ is ramified} \\ 1 & \text{else} \end{cases}
    \end{equation}
of $U(M^{\circ}_2)$ (acting on column vectors), where $a \in \mc{O}_{F_v}^{\times}$ and $b \in \mc{O}_{F_v}$ satisfies $\overline{b} = - \epsilon b$.

\textit{Case 1.} Assume $F_v / F^+_v$ is unramified and $c \in \mc{O}_{F^+_v}^{\times}$. Given any $x \in M^{\circ}_2$ with $(x, x) = c$, there exists an orthogonal direct sum decomposition $M^{\circ}_2 = \langle x \rangle \oplus \langle y \rangle$ for some $y \in M^{\circ}_2$ with $(y,y) = 1$ (by self-dualness). 
Via this decomposition, the lemma is clear in this case.

\textit{Case 2.} Assume $F_v / F^+_v$ is ramified and $c \in \mc{O}_{F^+_v}^{\times}$. Suppose $x = a_1 e_1 + a_2 e_2 \in M^{\circ}_2$ with $(x,x) = c$. Without loss of generality, we may assume $a_2 \in \mc{O}_{F_v}^{\times}$ (replace $x$ with $w' x$ if necessary), and we may further assume $a_2 = 1$ (replace $x$ with $m(\overline{a}_2) x$). We then have $\mrm{tr}_{F_v/F^+_v}(\varpi^{-1} a_1) = - c$. Given another $x' = a'_1 e_1 + e_2 \in M^{\circ}_2$ with $(x', x') = c$, we take $b = a'_1 - a_1$ and have $n(b) x = x'$.

\textit{Case 3.} Assume $F_v / F^+_v$ is inert and $c \in \varpi_0 \mc{O}_{F^+_v}^{\times}$. Suppose $x = a_1 e_1 + a_2 e_2 \in M^{\circ}_2$ with $(x,x) = c$. Without loss of generality, we may assume $a_2 = 1$ and $\mrm{tr}_{F_v/F^+_v}(a_1) = c$ (argue as in Case 2). Given another $x' = a'_1 e_1 + e_2 \in M^{\circ}_2$ with $(x', x') = c$, we take $b = a'_1 - a_1$ and have $n(b) x = x'$.
\end{proof}

\begin{remark}
If $F_v / F_v^+$ is nonsplit, the proof of Lemma \ref{lemma:local_Siegel-Weil:non-Arch:lattice_Witt_theorem}
shows that the $U(M^{\circ}_2)$ orbits in $(M^{\circ}_2 \otimes_{\mc{O}_{F_v}} F_v) \setminus \{ 0 \}$
are classified by pairs $(c, r) \in F^+_v \times \Z$, in the sense that the map
    \begin{equation}
        \begin{tikzcd}[row sep = tiny]
        U(M^{\circ}_2) \backslash ((M^{\circ}_2 \otimes_{\mc{O}_{F_v}} F_v) \setminus \{ 0 \} )
        \arrow{r} & 
        F^+_v \times \Z
        \\
        x \arrow[mapsto]{r} & (x,x), \min \{r \in \Z : \varpi^r x \in M^{\circ}_2 \}
        \end{tikzcd}
    \end{equation}
is an injection (say, setting $\varpi = \varpi_0$ in the inert case).
\end{remark}

\begin{lemma}\label{lemma:local_Siegel-Weil:non-Arch:transitive_dual_self_dual}
Let $L$ be a self-dual Hermitian $\mc{O}_{F_v}$-lattice of rank $n$. Any isomorphism between self-dual sublattices of $L$ extends to a (unitary) automorphism of $L$.
The same holds for rank $n - 1$ almost self-dual sublattices of $L$.
\end{lemma}
\begin{proof}
Any self-dual lattice $L^{\flat} \subseteq L$ admits a (unique) orthogonal direct sum decomposition $L = L^{\flat} \oplus L^{\#}$ where $L^{\#}$ is also self-dual. This immediately implies the claim for self-dual sublattices of $L$, as self-dual lattices are unique up to isomorphism (for a fixed rank).

Next, assume that $F_v / F^+_v$ is nonsplit and that $L^{\flat} \subseteq L$ is almost self-dual of rank $n - 1$. There is an orthogonal direct sum decomposition $L^{\flat} = L^{\flatflat} \oplus L^{\flat \kern-1.4pt \#}$, where $L^{\flatflat}$ is self-dual of rank $n - 2$ and $L^{\flat \kern-1.4pt \#}$ is almost self-dual of rank $1$. We also have an orthogonal direct sum decomposition $L = L^{\flatflat} \oplus L^{\#}$ where $L^{\#}$ is self-dual of rank $2$.

Suppose $L^{\flat \prime} \subseteq L$ is another almost self-dual lattice of rank $n - 1$, equipped with an isomorphism $L^{\flat} \ra L^{\flat \prime}$. Applying the result just proved above (in the case of rank $n - 2$ self-dual sublattices), we may assume there is an orthogonal decomposition $L^{\flat \prime} = L^{\flatflat} \oplus L^{\flat \prime  \#}$ where $L^{\flat \kern-1.4pt \#} \cong L^{\flat \prime \#}$. We thus reduce to the case $n = 2$ (the claim for $L^{\#}$), which was proved in Lemma \ref{lemma:local_Siegel-Weil:non-Arch:lattice_Witt_theorem}.
\end{proof}

\begin{lemma}\label{lemma:local_Siegel-Weil:non-Arch:almost_self_dual_in_self_dual}
Assume $F_v/ F^+_v$ is nonsplit, and let $V$ be an $F_v/F^+_v$ Hermitian space of rank $n$, and assume that $V$ contains a full-rank self-dual lattice. Suppose $L^{\flat} \subseteq V$ is a non-degenerate lattice of rank $n - 1$ satisfying $L^{\flat} \subseteq L^{\flat *}$ and $t(L^{\flat}) \leq 1$. Then $L^{\flat}$ is contained in a self-dual lattice of rank $n$.
\end{lemma}
\begin{proof}
Recall that $t(L^{\flat}) \coloneqq \dim_{k} ((L^{\flat *}/L^{\flat}) \otimes k)$ where $k$ is the residue field of $\mc{O}_{F_v}$. 

Let $L^{\flat} \subseteq V$ be as in the lemma statement. The existence of such $L^{\flat}$ implies $n \geq 2$. There exists an orthogonal decomposition $L^{\flat} = L^{\flatflat} \oplus L^{\flat \kern-1.4 pt \#}$ where $L^{\flatflat}$ is self-dual of rank $n - 2$. Replacing $V$ with the orthogonal complement of $L^{\flatflat}$, we reduce immediately to the case $n = 2$, which we now assume.

Let $\varpi$ be a uniformizer for $F_v$ (take $\varpi = \varpi_0$ if $F_v/F^+_v$ is inert). We may take $V = M^{\circ}_2 \otimes F_v$, where $M^{\circ}_2$ is as in Lemma \ref{lemma:local_Siegel-Weil:non-Arch:lattice_Witt_theorem}. We also choose a standard basis $e_1, e_2$ for $M^{\circ}_2$ and consider the elements $w', m(a), n(b) \in U(V)$ as in the proof of that lemma (now allowing $a \in F_v^{\times}$ and allowing $b \in F_v$ satisfying $\overline{b} = - \epsilon b$). 

The rank one lattice $L^{\flat}$ is generated by an element $x = a_1 e_1 + a_2 e_2$ for some $a_1, a_2 \in F_v$ (such that $(x,x)$ is nonzero and lies in $\mc{O}_{F^+_v}$).  It is enough to check that the orbit $U(V) \cdot x$ intersects $M^{\circ}_2$. Acting on $x$ by $m(a) \in U(V)$ for suitable $a$, we see that it is enough to check the case where $a_2 = 1$ and $a_1 \in F_v^{\times}$. 

If $F_v/F^+_v$ is inert, there exists $a' \in \mc{O}_{F_v}$ such that $\mrm{tr}_{F_v/F^+_v}(a') = (x,x)$ since $\mc{O}_{F_v}$ is self-dual with respect to the trace pairing. If $F_v/F^+_v$ is ramified, there exists $a' \in \mc{O}_{F_v}$ such that $\mrm{tr}_{F_v/F^+_v}(-\varpi^{-1} a') = (x,x)$ since $\mc{O}_{F_v}$ and $\varpi^{-1} \mc{O}_{F_v}$ are dual. In either case, we can take $b = a' - a_1$, and have $n(b) x \in M^{\circ}_2$.
\end{proof}

\begin{lemma}\label{lemma:local_Siegel-Weil:non-Arch:compact_open_intersection}
In the situations of Lemma \ref{lemma:local_Siegel-Weil:non-Arch:lattice_Witt_theorem}, choose $x \in M^{\circ}_2$ with $(x,x) = c$ and form the orthogonal complement lattice $x^{\perp} \subseteq M^{\circ}_2$ (of rank one). We view both $U(M^{\circ}_2)$ and $U(x^{\perp})$ as subgroups of $U(M^{\circ}_2 \otimes F_v)$.

Viewing $U(x^{\perp})$ as the norm-one subgroup of $\mc{O}_{F_v}^{\times}$, we have
    \begin{equation}
    U(M^{\circ}_2) \cap U(x^{\perp}) = \{ \a \in \mc{O}_{F_v}^{\times} : \a \overline{\a} = 1, \text{ and } \a \equiv 1 \pmod {c \mf{d} \mc{O}_{F_v}} \} \subseteq U(x^{\perp}).
    \end{equation}
The subgroup $U(M^{\circ}_2) \cap U(x^{\perp}) \subseteq U(x^{\perp})$ has index
    \begin{equation}
    \begin{cases}
    1 & \text{if $c \in \mc{O}_{F^+_v}^{\times}$ and $F_v / F^+_v$ is unramified}
    \\
    2 & \text{if $c \in \mc{O}_{F^+_v}^{\times}$ and $F_v / F^+_v$ is ramified}
    \\
    q_v + 1 & \text{if $c \in \varpi_0 \mc{O}_{F^+_v}^{\times}$ and $F_v / F^+_v$ is inert}.
    \end{cases}
    \end{equation}
\end{lemma}
\begin{proof}
We express elements of $U(M^{\circ}_2 \otimes F_v)$ in a standard basis $e_1, e_2$ of $M^{\circ}_2$, as in the proof of Lemma \ref{lemma:local_Siegel-Weil:non-Arch:lattice_Witt_theorem}.

\textit{Case 1.} Assume $F_v / F^+_v$ is unramified and $c \in \mc{O}_{F^+_v}^{\times}$. We then have $U(M^{\circ}_2) \cap U(x^{\perp}) = U(x^{\perp})$, as follows immediately from an orthogonal direct sum decomposition $M^{\circ}_2 = \langle x \rangle \oplus \langle y \rangle$ as in the proof of Lemma \ref{lemma:local_Siegel-Weil:non-Arch:lattice_Witt_theorem} Case 1.

\textit{Case 2.} Assume $F_v / F^+_v$ is ramified and $c \in \mc{O}_{F^+_v}^{\times}$. By the proof of Lemma \ref{lemma:local_Siegel-Weil:non-Arch:lattice_Witt_theorem} Case 2, we may assume (after conjugating $U(M^{\circ}_2 \otimes F_v)$ by an appropriate element of $U(M^{\circ}_2)$) that $x = a_1 e_1 + e_2$ for some $a_1 \in \mc{O}_{F_v}$, where $a_1 - \overline{a}_1 = - \varpi c$. Then $\overline{a}_1 e_1 + e_2$ is orthogonal to $x$. For every $\a \in \mc{O}_{F_v}^{\times}$, the matrix
    \begin{equation}
    \begin{pmatrix}
    a_1 & \overline{a}_1 \\
    1 & 1
    \end{pmatrix}
    \begin{pmatrix}
    1 & 0 \\
    0 & \a
    \end{pmatrix}
    \begin{pmatrix}
    a_1 & \overline{a}_1 \\
    1 & 1
    \end{pmatrix}^{-1}
    =
    (- \varpi c)^{-1} 
    \begin{pmatrix}
    a_1 - \overline{a}_1 \a & (- 1 + \a) a_1 \overline{a}_1 \\
    1 - \a & - \overline{a}_1 + a_1 \a
    \end{pmatrix}
    \end{equation}
lies in $U(M^{\circ}_2)$ if and only if $\a \equiv 1 \pmod{\varpi \mc{O}_{F_v}}$. The claim about index follows from surjectivity of the reduction modulo $\varpi$ map 
    \begin{equation}
    \{ \a \in \mc{O}_{F_v}^{\times} : \a \overline{\a} = 1 \} \ra \{ \a \in \F_{q_v}^{\times} : \a^2 = 1 \}
    \end{equation}
(surjectivity is by smoothness of the corresponding unitary group over $\Spec \mc{O}_{F^+_v}$).

\textit{Case 3.} Assume $F_v / F^+_v$ is inert and $c \in \varpi_0 \mc{O}_{F^+_v}^{\times}$. By the proof of Lemma \ref{lemma:local_Siegel-Weil:non-Arch:lattice_Witt_theorem} Case 3, we may assume (after conjugating $U(M^{\circ}_2 \otimes F_v)$ by an appropriate element of $U(M^{\circ}_2)$) that $x = a_1 e_1 + e_2$ for some $a_1 \in \mc{O}_{F_v}$, where $a_1 + \overline{a}_1 = c$. Then $-\overline{a}_1 e_1 + e_2$ is orthogonal to $x$. For every $\a \in \mc{O}_{F_v}^{\times}$, the matrix
    \begin{equation}
    \begin{pmatrix}
    a_1 & -\overline{a}_1 \\
    1 & 1
    \end{pmatrix}
    \begin{pmatrix}
    1 & 0 \\
    0 & \a
    \end{pmatrix}
    \begin{pmatrix}
    a_1 & -\overline{a}_1 \\
    1 & 1
    \end{pmatrix}^{-1}
    =
    c^{-1} 
    \begin{pmatrix}
    a_1 + \overline{a}_1 \a & (1 - \a) a_1 \overline{a}_1 \\
    1 - \a & \overline{a}_1 + a_1 \a
    \end{pmatrix}
    \end{equation}
lies in $U(M^{\circ}_2)$ if and only if $\a \equiv 1 \pmod{\varpi_0 \mc{O}_{F_v}}$. The claim about index follows from surjectivity of the reduction modulo $\varpi_0$ map 
    \begin{equation}
    \{ \a \in \mc{O}_{F_v}^{\times} : \a \overline{\a} = 1 \} \ra \{ \a \in \F_{q_v^2}^{\times} : \a \overline{\a} = 1 \}
    \end{equation}
(surjectivity is by smoothness of the corresponding unitary group over $\Spec \mc{O}_{F^+_v}$).
\end{proof}

The following lemma can be interpreted
as giving an explicit description
of the Haar measure from Lemma \ref{lemma:local_Siegel-Weil:local_Siegel-Weil:unique_Haar_measure}
in some special cases.
Indeed, since $W^*_{T,v}(s)^{\circ}_n = \Lambda_{T,v}(s)^{\circ}_n W_{T,v}(s, \Phi_{\varphi_v})$
in the situation of Lemma \ref{lemma:local_Siegel-Weil:non-Arch:spherical_and_measures}
(using the normalizations defined in Section \ref{ssec:Eisenstein:local_Whittaker:non-Arch},
and using Lemma \ref{lemma:spherical_Siegel-Weil_section}),
we find that the Haar measure in Lemma \ref{lemma:local_Siegel-Weil:local_Siegel-Weil:unique_Haar_measure}
must be
    \begin{equation}
    e^{-1} \cdot |\det T|_{F^+_v}^{m - n} \cdot (\Lambda_{T,v}(s_0)^{\circ}_n)^{-1} \cdot dg_v
    \end{equation}
with $T$ and $e$ and $d g_v$ and $s_0$ as in Lemma \ref{lemma:local_Siegel-Weil:non-Arch:spherical_and_measures}
(using also Lemma \ref{lemma:Weil_index:Hermitian}
to compute the Weil index).
The above expression does not depend on $T$ (as expected),
as can be verified upon inspecting \eqref{equation:Eisenstein:local_Whittaker:non-Arch:normalizing_factor}.

\begin{lemma}\label{lemma:local_Siegel-Weil:non-Arch:spherical_and_measures}
Take $m = n - 1$ or $m = n$ and $s_0 \coloneqq (n - m)/2$.
Assume the rank $n$ Hermitian space $V$ contains a full-rank self-dual lattice $L$. 
Let $K_v \subseteq G = U(V)$ be the stabilizer of such a lattice $L$. 

Consider any $\underline{x}_v \in V^m(F^+_v)$ with nonsingular Gram matrix $T = (\underline{x}_v, \underline{x}_v) \in \mrm{Herm}_m(F^+_v)$. Let $\pmb{1}_{L}$ be the characteristic function of $L$, and set $\varphi_v = \pmb{1}_{L}^{m} \in \mc{S}(V^m(F^+_v))$. 

Give $G(F^+_v)$ the Haar measure which assigns volume $1$ to $K_v$. Give $G_{\underline{x}_v}(F^+_v)$ the Haar measure which assigns volume $1$ to the (unique) maximal open compact subgroup. We have
    \begin{equation}\label{equation:local_Siegel-Weil:non-Arch:integral_compare_star_Whittaker}
    W_{T,v}^*(s_0)^{\circ}_n = \frac{1}{e} \int_{G_{\underline{x}_v}(F^+_v) \backslash G(F^+_v)} \varphi_v(g_v^{-1} \underline{x}_v) ~ d g_v
    \quad \quad 
    e \coloneqq
    \begin{cases}
    2 & \text{if $F_v/F^+_v$ is ramified and $m = n - 1$} \\
    1 & \text{else}
    \end{cases}
    \end{equation}
with respect to the associated quotient measure.
\end{lemma}
\begin{proof}
Recall that $W_{T,v}^*(s)^{\circ}_n$ is our notation for a certain normalized spherical Whittaker function (Section \ref{ssec:Eisenstein:local_Whittaker:non-Arch}), which is a rescaled version of $W_{T,v}(s, \Phi_{\varphi_v})$. 

In the lemma statement, the stabilizer in $G(F^+_v)$ of any full-rank self-dual lattice in $V$ has volume $1$ (because any such stabilizer is conjugate to $K_v$),
i.e. the measure on $G(F^+_v)$ does not depend on the choice of $L$.
To verify \eqref{equation:local_Siegel-Weil:non-Arch:integral_compare_star_Whittaker}, we can (and will) replace $L$ by another full-rank self-dual lattice in $V$.
Lemma \ref{lemma:local_Siegel-Weil:local_Siegel-Weil:unique_Haar_measure} implies that
if \eqref{equation:local_Siegel-Weil:non-Arch:integral_compare_star_Whittaker}
holds for $\varphi_v = \pmb{1}_L^{\otimes m}$
for some full-rank self-dual $L$, then it must hold for all full-rank self-dual $L$.

Let $V^{\flat} \subseteq V$ be the rank $m$ subspace spanned by $\underline{x}_v$. Then $V^{\flat}$ is free of rank $m$. By Lemma \ref{lemma:local_Siegel-Weil:local_Siegel-Weil:unique_Haar_measure}, it is enough to show \eqref{equation:local_Siegel-Weil:non-Arch:integral_compare_star_Whittaker} holds for one choice of basis 
$\underline{x}_v$ for $V^{\flat}$. We choose $\underline{x}_v$ to be a basis for a full-rank lattice $L^{\flat} \subseteq V^{\flat}$ which is
    \begin{equation}
    \begin{cases}
    \text{self-dual} & \text{if $V^{\flat}$ contains a full-rank self-dual lattice}
    \\
    \text{almost self-dual} & \text{else}.    
    \end{cases}
    \end{equation}
Note that $V^{\flat}$ always contains a full-rank self-dual lattice if $F_v / F^+_v$ is split.

\textit{Case 1.} Assume $L^{\flat}$ is self-dual. There exists a rank $n-m$ self-dual lattice $L^{\#} \subseteq V$ which is orthogonal to $L^{\flat}$. Form the rank $n$ self-dual lattice $L = L^{\flat} \oplus L^{\#}$. Any isomorphism between self-dual sublattices of $L$ lifts to an element of $K_v = U(L)$ (Lemma \ref{lemma:local_Siegel-Weil:non-Arch:transitive_dual_self_dual}). This implies that $g_v \mapsto \varphi_v(g_v^{-1} \underline{x}_v)$ is the characteristic function of $G_{\underline{x}_v}(F^+_v) \backslash ( G_{\underline{x}_v}(F^+_v) K_v)$.

We know that $K_v \cap G_{\underline{x}_v}(F^+_v)$ is the unique maximal open compact subgroup in $G_{\underline{x}_v}(F^+_v)$ (i.e. $U(L^{\#})$). We compute
    \begin{equation}
    \int_{G_{\underline{x}_v}(F^+_v) \backslash G(F^+_v)} \varphi_v(g_v^{-1} \underline{x}_v) ~ d g_v = \mrm{vol}(G_{\underline{x}_v}(F^+_v) \backslash ( G_{\underline{x}_v}(F^+_v) K_v)) = \frac{\mrm{vol}(K_v)}{\mrm{vol}(K_v \cap G_{\underline{x}_v}(F^+_v))} = 1.
    \end{equation}
Since $T = (\underline{x}_v, \underline{x}_v)$ and $\underline{x}_v$ is a basis for the self-dual lattice $L^{\flat}$, we also know $W_{T,v}^*(s_0)_n = 1$ (see \eqref{equation:Eisenstein:local_Whittaker:local_densities_spherical:often_1}; note that $V^{\flat}$ containing a self-dual lattice means that $F_v / F^+_v$ is unramified if $m$ is odd).

\textit{Case 2.} Assume that $L^{\flat}$ is almost self-dual and that $F_v / F^+_v$ is ramified. Then $n \geq 2$ and $m = n - 1$. There is an orthogonal direct sum decomposition $L^{\flat} = L^{\flatflat} \oplus L^{\flat \kern-1.4pt  \#}$, where $L^{\flatflat}$ is self-dual of rank $m - 1$ and $L^{\flat \kern-1.4pt \#}$ is almost self-dual of rank $1$. There exists a rank $2$ self-dual lattice $L^{\#} \subseteq V$ which is orthogonal to $L^{\flatflat}$. We can assume $L^{\flat \kern-1.4pt \#} \subseteq L^{\#}$ (Lemma \ref{lemma:local_Siegel-Weil:non-Arch:almost_self_dual_in_self_dual}). Form the rank $n$ self-dual lattice $L = L^{\flatflat} \oplus L^{\#}$. Any isomorphism between rank $n - 1$ almost self-dual sublattices in $L$ lifts to an element of $K_v = U(L)$ (Lemma \ref{lemma:local_Siegel-Weil:non-Arch:transitive_dual_self_dual}). This implies that $g_v \mapsto \varphi_v(g_v^{-1} \underline{x}_v)$ is the characteristic function of $G_{\underline{x}_v}(F^+_v) \backslash ( G_{\underline{x}_v}(F^+_v) K_v)$.

We know that $K_v \cap G_{\underline{x}_v}(F^+_v) = U(L^{\#}) \cap G_{\underline{x}_v}(F^+_v)$ has index $2$ inside the unique maximal open compact subgroup of $G_{\underline{x}_v}(F^+_v)$ (reduces immediately to the case $n = 2$, which is Lemma \ref{lemma:local_Siegel-Weil:non-Arch:compact_open_intersection}). 
We compute
    \begin{equation}
    \int_{G_{\underline{x}_v}(F^+_v) \backslash G(F^+_v)} \varphi_v(g_v^{-1} \underline{x}_v) ~ d g_v = \mrm{vol}(G_{\underline{x}_v}(F^+_v) \backslash ( G_{\underline{x}_v}(F^+_v) K_v)) = \frac{\mrm{vol}(K_v)}{\mrm{vol}(K_v \cap G_{\underline{x}_v}(F^+_v))} = 2.
    \end{equation}
Since $T = (\underline{x}_v, \underline{x}_v)$ and since $\underline{x}_v$ is a basis for the almost self-dual lattice $L^{\flat}$, we also know $W_{T,v}^*(s_0)^{\circ}_n = 1$ \eqref{equation:Eisenstein:local_Whittaker:local_densities_spherical:often_1}.

\textit{Case 3.} Assume that $L^{\flat}$ is almost self-dual and that $F_v / F^+_v$ is inert. This implies $n \geq 2$ and $m = n - 1$. Arguing as in Case 2 (use the same notation; the first paragraph applies verbatim), again apply Lemma \ref{lemma:local_Siegel-Weil:non-Arch:transitive_dual_self_dual} and Lemma \ref{lemma:local_Siegel-Weil:non-Arch:compact_open_intersection} to compute
    \begin{equation}
    \int_{G_{\underline{x}_v}(F^+_v) \backslash G(F^+_v)} \varphi_v(g_v^{-1} \underline{x}_v) ~ d g_v = \mrm{vol}(G_{\underline{x}_v}(F^+_v) \backslash ( G_{\underline{x}_v}(F^+_v) K_v)) = \frac{\mrm{vol}(K_v)}{\mrm{vol}(K_v \cap G_{\underline{x}_v}(F^+_v))} = q_v + 1.
    \end{equation}
When $n = 2$, we have $\mrm{Den}^*(X, L^{\flat})_n = q_v X^{-1/2} + X^{1/2}$ (follows from the relevant Cho--Yamauchi type formula; see \cite[{Example 3.5.2}]{LZ22unitary} \cite[{Theorem 2.2}]{FYZ22SW}). The ``cancellation'' property for local densities and self-dual lattices \eqref{equation:local_Whittaker:local_densities_spherical:cancellation} implies $\mrm{Den}^*(X, L^{\flat})_n = q_v X^{-1/2} + X^{1/2}$ for $n \geq 2$. We thus have $W^*_{T,v}(s_0)^{\circ}_n = \mrm{Den}^*(1,L^{\flat})_n = q_v + 1$.
\end{proof}

            \subsection{Uniformization degrees for special cycles} 
            \label{ssec:local_Siegel-Weil:uniformization_degree}
                The purpose of this section is to express the groupoid cardinality of \eqref{equation:local_Siegel-Weil:uniformization_degree:main} in terms of special values of local Whittaker functions, with explicit constants (Lemma \ref{lemma:local_Siegel-Weil:uniformization_degree:main}). This groupoid has already appeared as a ``uniformization degree'' for special cycles (see \crefext{III:equation:non-Arch_uniformization:quotient:unitary_stacky_degree}, also \crefext{III:ssec:non-Arch_uniformization:vertical,III:ssec:non-Arch_uniformization:horizontal} and \crefext{III:ssec:Arch_uniformization:Archimedean}). This calculation will be needed for our main arithmetic Siegel--Weil results (Section \ref{ssec:arithmetic_Siegel-Weil:main_results}).

Let $F/F^+$ be a CM extension of number fields, with respective ad\`ele rings $\A_F$ and $\A$ and finite ad\`ele rings $\A_{F,f}$ and $\A_f$, etc.. As in \cref{sec:setup,sec:Eisenstein:Weil_rep,sec:local_Whittaker,sec:Eisenstein:local_functional_equations,sec:Eisenstein:global_normalized_Fourier}, we write $v$ for places of $F^+$ with completions $F^+_v$, and set $F_v \coloneqq F \otimes_{F^+} F^+_v$.

Let $T \in \mrm{Herm}_m(F^+)$ be a Hermitian matrix (with $F$-coefficients) for any integer $m \geq 0$. Set $m^{\flat} \coloneqq \rank(T)$. For each place $v$, select any $a_v \in \GL_m(F_v)$ such that ${}^t \overline{a}_v^{-1} T a_v^{-1} = \mrm{diag}(0,T_v^{\flat})$ for some $T_v^{\flat} \in \mrm{Herm}_{m^{\flat}}(F^+_v)$ with $\det T^{\flat}_v \neq 0$. For each $v$, choose any decomposition (Iwasawa decomposition)
    \begin{equation}\label{equation:local_Siegel-Weil:uniformization_degree:a_v_Iwasawa_decomp}
    a_v = 
    \begin{pmatrix}
    1_{m - m^{\flat}} & * \\
    0 & 1_{m^{\flat}}
    \end{pmatrix}
    \begin{pmatrix}
    a_v^{\#} & 0 \\
    0 & a_v^{\flat}
    \end{pmatrix}
    k_v
    \quad \quad
    k_v \in
    \begin{cases}
    \GL_m(\mc{O}_{F_v}) & \text{if $v$ is non-Archimedean} \\
    U(m) & \text{if $v$ is Archimedean,}
    \end{cases}
    \end{equation}
where $a_v^{\#} \in \GL_{m - m^{\flat}}(F_v)$, $a_v^{\flat} \in \GL_{m^{\flat}}(F_v)$, and $U(m) \subseteq \GL_m(\C)$ is the unitary group for the standard diagonal positive definite Hermitian pairing.

Let $L$ be a non-degenerate Hermitian $\mc{O}_F$-lattice of any rank $n$, set $V \coloneqq L \otimes_{\mc{O}_F} F$, and let $G = U(V)$ be the associated unitary group. Set $s_0^{\flat} \coloneqq (n - m^{\flat}) / 2$. 
For any place $v$ of $F^+$, we set $V_v \coloneqq V \otimes_{F^+} F^+_v$. Let $K_{L,f} = \prod K_{L,v} \subseteq U(V)(\A_f)$ be the ad\`elic stabilizer of $L$ (i.e. $K_{L,v}$ is the stabilizer of $L_v \coloneqq L \otimes_{\mc{O}_{F^+}} \mc{O}_{F^+_v}$ for every place $v < \infty$ of $F^+_v$). 
Fix a place $v_0$ of $F^+$ (Archimedean or non-Archimedean). Assume $V_v$ is positive definite for every Archimedean $v \neq v_0$. 

Given $\underline{x}_f^{v_0} \in (V \otimes_{F^+} \A_f^{v_0})^m$, we define the ``away from $v_0$ special cycle'' (compare Sections \crefext{III:ssec:non-Arch_uniformization:away_p_special_cycle} and \crefext{III:ssec:Arch_uniformization:away_infty_special_cycle})
    \begin{equation}
    \mc{Z}(\underline{x}^{v_0}_f) \coloneqq \{g_f \in G(\A_f^{v_0}) / K_{L,f}^{v_0} : g_{f,v}^{-1} \underline{x}_v \in L_v \text{ for all non-Archimedean $v \neq v_0$}\}
    \end{equation}
where $\underline{x}_v \in V_v^m$ is the $v$-component of $\underline{x}_f^{v_0}$.

Fix a nontrivial additive character $\psi_v$ for each place $v$. Assume $\psi_v$ is unramified if $v < \infty$, and assume $\psi_v(x) = e^{2 \pi i x}$ if $F^+_v = \R$.
For every non-Archimedean place $v \neq v_0$, set $\varphi_v \coloneqq \pmb{1}_{L_v}^{m}$ (characteristic function of $L_v^{m} \subseteq V_v^{m}$) and set 
    \begin{equation}
    \varphi_f^{v_0} = \otimes_{\substack{v < \infty \\ v \neq v_0}} \varphi_v \in \mc{S}(V(\A_f^{v_0})^{m}).        
    \end{equation}
Similarly set $\varphi_v^{\flat} \coloneqq \pmb{1}_{L_v}^{m^{\flat}} \in \mc{S}(V(F^+_v)^{m^{\flat}})$ for such $v$.

For every place $v$ of $F^+_v$, let $\eta_v \colon F^{+ \times}_v \ra \{ \pm 1 \}$ be the quadratic character associated to $F_v / F^+_v$. Let $\chi_v \colon F^{\times}_v \ra \C^{\times}$ be any character satisfying $\chi_v |_{F^{+ \times}_v} = \eta_v^n$. Form the associated Siegel--Weil standard section $\Phi_{\varphi_v} \in I(s, \chi_v)$ (Section \ref{ssec:Eisenstein:Weil_rep:Weil_rep}) for every place $v < \infty$ with $v \neq v_0$. To simplify slightly, we assume that $2$-adic places of $F^+$ are unramified in $F$ for the rest of Section \ref{ssec:local_Siegel-Weil:uniformization_degree}.

For $v < \infty$ with $v \neq v_0$, the local Whittaker function variant $\tilde{W}^*_{T^{\flat}_v, v}(a^{\flat}_v, s, \Phi_{\varphi_v})_n$ does not depend on the choice of $a_v$ or $a_v^{\flat}$.
Indeed, the $\GL_{m^{\flat}}(\mc{O}_{F_v})$-equivalence class
of the Hermitian matrix ${}^t \overline{a}^{\flat}_v T^{\flat}_v a^{\flat}_v$
does not depend on the choice of $a_v$ or $a^{\flat}_v$
(then apply \eqref{equation:local_Whittaker:non-Arch:linear_invariance}).
For $v \mid \infty$ with $v \neq v_0$,
the local Whittaker function variant $\tilde{W}^*_{T^{\flat}_v, v}(a^{\flat}_v, s)^{\circ}_n$
similarly does not depend on the choice of $a_v$ or $a^{\flat}_v$,
as the $U(m^{\flat})$-equivalence class of ${}^t \overline{a}^{\flat}_v T^{\flat}_v a^{\flat}_v$
is similarly well-defined (then apply \eqref{equation:local_Whittaker:Archimedean:linear_invariance}).

Given any tuple $\underline{x} \in V^m$ which spans a non-degenerate Hermitian space, we write $G_{\underline{x}} \subseteq G$ for the stabilizer of $\underline{x}$ (i.e. the unitary group of the orthogonal complement $\mrm{span}(\underline{x})^{\perp} \subseteq V$). We write $\underline{x}_f^{v_0}$
for the image of $\underline{x}$ in $(V \otimes_{F^+} \A_f^{v_0})^m$.

Suppose there exists $\underline{x} \in V^m$ with Gram matrix $(\underline{x}, \underline{x}) = T$.
Fix such an $\underline{x}$, and assume $\mrm{span}(\underline{x})^{\perp}$ is definite at every Archimedean place. Let $K_{\underline{x},v_0} \subseteq G_{\underline{x}}(F^+_{v_0})$ be any open compact subgroup, and assume $K_{\underline{x},v_0}(F^+_{v_0}) = G_{\underline{x}}(F^+_{v_0})$ if $v_0$ is Archimedean. 

We are mostly interested in applying Lemma \ref{lemma:local_Siegel-Weil:uniformization_degree:main} below when $m^{\flat} \geq n - 1$ and $L_v$ is self-dual for all $v < \infty$ with $v \neq v_0$. The result and proof is simpler in that case, and the lemma may not be optimal otherwise.

\begin{lemma}\label{lemma:local_Siegel-Weil:uniformization_degree:main}
Consider the groupoid quotient
    \begin{equation}\label{equation:local_Siegel-Weil:uniformization_degree:main}
    \bigg [G_{\underline{x}}(F^+) \backslash \bigg ( G_{\underline{x}}(F^+_{v_0})/K_{\underline{x},v_0} \times \mc{Z}(\underline{x}_f^{v_0}) \bigg ) \bigg ].
    \end{equation}
The displayed groupoid has finite automorphism groups and finitely many isomorphism classes. Its groupoid cardinality is
        \begin{equation}
        C \cdot 
        \prod_{\substack{v \mid \infty \\ v \neq v_0}} \tilde{W}^*_{T^{\flat}_v, v}(a^{\flat}_v, s_0^{\flat})^{\circ}_n
        \prod_{\substack{v < \infty \\ v \neq v_0}} \tilde{W}^*_{T^{\flat}_v, v}(a^{\flat}_v, s_0^{\flat}, \Phi_{\varphi_v^{\flat}})_n
        \end{equation}
for some volume constant $C \in \Q_{>0}$ which we describe in the following three situations.
\begin{enumerateThm}
    \item \plabel{lemma:local_Siegel-Weil:uniformization_degree:main:1}
    Suppose $v_0$ is Archimedean. Assume the local characters $(\psi_v)_v$ and $(\chi_v)_v$ arise from global characters $\psi \colon F^+ \backslash \A \ra \C^{\times}$ and $\chi \colon F^{\times} \backslash \A_F^{\times} \ra \C^{\times}$. The constant $C$ may depend on $V$, $n$, $m^{\flat}$, $F$, and the isomorphism classes of the local Hermitian lattices $\{L_v\}_{v < \infty}$. The constant $C$ does not otherwise depend on $T$ or $V^{\flat}$ or $\underline{x}$.
    \item Suppose $m^{\flat} = n$ (with $v_0$ not necessarily Archimedean). Then
        \begin{equation}
        C = \prod_{\substack{v < \infty \\ v \neq v_0}} c_v
        \end{equation}
    for some constants $c_v \in \Q_{>0}$, all but finitely many of which are $1$. For any given $v < \infty$ with $v \neq v_0$, the constant $c_v$ may depend on the local Hermitian lattice $L_v$ and the quadratic extension $F_v / F^+_v$, but otherwise does not depend on $T$ or $V$ or $\underline{x}$ or $v_0$ or $F / F^+$. 
    
    If $L_v$ is self-dual, then $c_v = 1$.

    \item Suppose $m^{\flat} = n - 1$ (with $v_0$ not necessarily Archimedean). Assume $K_{\underline{x},v_0} \subseteq G_{\underline{x}}(F^+_v)$ is the unique maximal open compact subgroup. Then there are constants $c'_v \in \Q_{>0}$ such that
        \begin{equation}
        C = \frac{2^{1 - o(\Delta)} h_F}{w_F h_{F^+} \cdot \# (\mc{O}_{F}^{\times} / (W \mc{O}_{F^+}^{\times}))} \prod_{\substack{v < \infty \\ v \neq v_0}} c'_v
        \end{equation}
    where $o(\Delta)$ is the number of prime ideals of $\mc{O}_{F^+}$ which ramify in $\mc{O}_{F}$, where $h_F$ (resp. $h_{F^+}$) is the class number of $F$ (resp. $F^+$), where $w_F$ (resp. $W$) is the number of (resp. group of) roots of unity in $F$. All but finitely many $c'_v$ are equal to $1$.
    
    For each $v < \infty$ with $v \neq v_0$, the constant $c'_v$ may depend on the local Hermitian lattice $L_v$, the quadratic extension $F_v / F^+_v$, and the local invariant $\varepsilon(V^{\flat}_v) \in \{ \pm 1\}$. The constant $c'_v$ does not otherwise depend on $T$ or $V$ or $V^{\flat}$ or $\underline{x}$ or $v_0$ or $F / F^+$. 
    
    If $L_v$ is self-dual, then $c'_v = 1$ if $F_v / F^+_v$ is unramified (resp. $c'_v = 2$ if $F_v / F^+_v$ is ramified).
\end{enumerateThm}
\end{lemma}
\begin{proof}
For the moment, we allow $v_0$ Archimedean or not. The groupoid in the lemma statement indeed has finite stabilizer groups, by discreteness of $G_{\underline{x}}(F^+)$. Take any factorizable open compact subgroup $K_{\underline{x}} = \prod_v K_{\underline{x}, v} \subseteq G_{\underline{x}}(\A)$. Assume $K_{\underline{x},v} = G_{\underline{x},v}(F^+_v)$ for every Archimedean $v$, and assume $K_{\underline{x},v} = K_{\underline{x},v_0}$ is the open compact subgroup fixed in the lemma statement when $v = v_0$. For each $v$, define $\underline{x}_v^{\flat} = [x_{1,v}^{\flat}, \ldots, x^{\flat}_{m^{\flat},v}] \in V_v^{m^{\flat}}$ to be the tuple satisfying $\underline{x} \cdot a_v^{-1} = [0,\ldots,0, x_{1,v}^{\flat}, \ldots, x_{m^{\flat},v}^{\flat}]$ (so $T^{\flat}_v = (\underline{x}^{\flat}_v, \underline{x}^{\flat}_v)$).

We have $\tilde{W}^*_{T^{\flat}_v,v}(a^{\flat}_v, s_0)^{\circ}_n = 1$ for all Archimedean $v \neq v_0$ by positive definite-ness of $T^{\flat}_v$ (Section \ref{ssec:Eisenstein:local_Whittaker:Archimedean}).
For all but finitely many $v$, the Hermitian matrix ${}^t \overline{a}^{\flat}_v T^{\flat}_v a^{\flat}_v$ defines a self-dual $\mc{O}_{F_v}$-lattice
(first check the case where the collection $(a_v)_v$ comes from a single element $a \in \GL_m(F)$;
then recall that the $\GL_{m^{\flat}}(\mc{O}_{F_v})$-equivalence class
of ${}^t \overline{a}^{\flat}_v T^{\flat}_v a^{\flat}_v$ does not depend on the choice of $a_v$ or $a^{\flat}_v$).
For such non-Archimedean $v \neq v_0$, we have $\tilde{W}^*_{T^{\flat}_v, v}(a^{\flat}_v, s, \Phi_{\varphi_v^{\flat}})_n = \tilde{W}^*_{T^{\flat}_v, v}(a^{\flat}_v, s)^{\circ}_n = 1$ if $L_v$ is self-dual (see \eqref{equation:Eisenstein:local_Whittaker:local_densities_spherical:often_1}, Remark \ref{remark:Eisenstein:local_Whittaker:character_could_ramify}, and the invariance property in \eqref{equation:local_Whittaker:non-Arch:linear_invariance}). Hence $\tilde{W}^*_{T^{\flat}_v, v}(a^{\flat}_v, s, \Phi_{\varphi_v^{\flat}})_n = 1$ for all but finitely many $v$.

Choose Haar measures $d g_{x,v}$ on $G_{\underline{x}}(F^+_v)$ for each $v$. Assume that $\mrm{vol}_{d g_{x,v}}(K_{\underline{x},v}) \in \Q$ for all $v$, that $\mrm{vol}_{d g_{x,v}}(K_{\underline{x},v}) = 1$ for all but finitely many $v$, and that $\mrm{vol}_{d g_x, v}(K_{\underline{x},v}) = 1$ if $v = v_0$ or if $v \mid \infty$.

For $v < \infty$ with $v \neq v_0$, we give $G(F^+_v)$ the unique Haar measure $d g_v$ such that
    \begin{equation}\label{equation:local_Siegel-Weil:uniformization_degree:Haar_measure_choice}
    W^*_{T^{\flat \prime}_v, v}(1, s_0^{\flat}, \Phi_{\varphi^{\flat}_v})_n = \int_{G_{\underline{x}_v}(F^+_v) \backslash G(F^+_v)} \varphi_v^{\flat} (g_v^{-1} \underline{x}_v') ~ d g_v
    \end{equation}
for any tuple $\underline{x}'_v \in V_v^m$ (temporary notation) with nonsingular Gram matrix $T^{\flat \prime}_v \coloneq (\underline{x}'_v, \underline{x}'_v)$ (Lemma \ref{lemma:local_Siegel-Weil:local_Siegel-Weil:unique_Haar_measure}). The integral is taken with respect to the quotient measure induced by $d g_{x,v}$. This measure $d g_v$ on $G(F^+_v)$ may depend on $n$, $m^{\flat}$, the isomorphism class of $L_v$ (as the normalization defining $\tilde{W}^*_{T^{\flat}, v}$ depended on $L_v$) as well as the local invariant $\varepsilon(V^{\flat}_v)$ (Remark \ref{remark:local_Siegel-Weil:unique_Haar_measure:invariance}). The measure $d g_v$ does not otherwise depend on $T^{\flat}_v$.
Note $\mrm{vol}_{d g_v}(K_{L,v}) \in \Q_{>0}$ for any $v < \infty$ with $v \neq v_0$, since the left-hand side of \eqref{equation:local_Siegel-Weil:uniformization_degree:Haar_measure_choice} lies in $\Q$ by Lemma \ref{lemma:local_density_Whittaker_interpolation}. We have $\mrm{vol}_{d g_v}(K_{L,v}) = 1$ for all but finitely many $v$ (cf. the proof of Lemma \ref{lemma:local_Siegel-Weil:non-Arch:spherical_and_measures}; we have $W^*_{T^{\flat}_v, v}(s_0^{\flat})^{\circ}_n = 1$ for all but finitely many $v$).
We equip $G(\A_f^{v_0})$ with the product measure $dg = \prod_{\substack{v < \infty \\ v \neq v_0}} d g_v$.

Using the Haar measures specified above, we may unfold the groupoid cardinality as
    \begin{align}
    & \deg \bigg [G_{\underline{x}}(F^+) \backslash \bigg ( G_{\underline{x}}(F^+_{v_0})/K_{\underline{x},v_0} \times \mc{Z}(\underline{x}^{v_0}_f) \bigg ) \bigg ]
    \\
    & = \mrm{vol}_{dg}(K_{L,f}^{v_0})^{-1} \int_{G_{\underline{x}}(F^+) \backslash ((\prod_{\substack{v = v_0 \\ \text{or } v \mid \infty}} G_{\underline{x}}(F^+_{v})) \times G(\A_f^{v_0}))} \varphi_f^{v_0}(g^{-1} \underline{x}) ~dg
    \\
    & = 
    \mrm{vol}(G_{\underline{x}}(F^+) \backslash G_{\underline{x}}(\A)) \mrm{vol}_{dg}(K_{L,f}^{v_0})^{-1} \left ( \int_{G_{\underline{x}}(\A_f^{v_0}) \backslash G(\A_f^{v_0})} \varphi_f^{v_0} (g^{-1} \underline{x}) ~dg \right )
    \\
    & = \mrm{vol}(G_{\underline{x}}(F^+) \backslash G_{\underline{x}}(\A)) \mrm{vol}_{dg}(K_{L,f}^{v_0})^{-1} \prod_{\substack{v < \infty \\ v \neq v_0}} \int_{G_{\underline{x}}(F^+_v) \backslash G(F^+_v)} \varphi_v^{\flat}(g_v^{-1} \underline{x}_v^{\flat} a^{\flat}_v) ~dg_v
    \\
    & = C \prod_{\substack{v < \infty \\ v \neq v_0}} \tilde{W}^*_{T^{\flat}_v, v}(a^{\flat}_v, s_0^{\flat}, \Phi_{\varphi_v^{\flat}})_n \label{equation:local_Siegel-Weil:uniformization_degree:main:unfold_local} 
    \end{align}
with
    \begin{equation}
    C \coloneqq \mrm{vol}(G_{\underline{x}}(F^+) \backslash G_{\underline{x}}(\A)) \prod_{\substack{v < \infty \\ v \neq v_0}} \mrm{vol}_{d g_v}(K_{L,v})^{-1}.
    \end{equation}
Note that the integrals are absolutely convergent, since the integrands are continuous and compactly supported. This unfolding also shows that the groupoid in \eqref{equation:local_Siegel-Weil:uniformization_degree:main} has finitely many isomorphism classes.

\begin{enumerate}[(1)]
    \item Suppose $v_0$ is Archimedean. Recall that the Tamagawa number of any nontrivial unitary group is $2$ \cite[{Section 4}]{Ichino04}. After scaling one of the non-Archimedean local measures $d g_{x,v}$ by an element of $\Q_{>0}$, we may assume $\prod_v d g_{x, v}$ is the Tamagawa measure on $G_{\underline{x}}(\A)$. If $v \mid \infty$, let $d g_v$ be the Haar measure on $G(F^+_v)$ given by Lemma \ref{lemma:local_Siegel-Weil:local_Siegel-Weil:unique_Haar_measure} (induced by $d g_{x,v}$). For $v \mid \infty$, the local invariant $\varepsilon(V^{\flat}_v)$ is already determined by $V$ and the requirement that $V^{\flat \perp}_v$ is definite. Hence the measures $d g_v$ for $v \mid \infty$ do not depend on $V^{\flat}$ (apply Remark \ref{remark:local_Siegel-Weil:unique_Haar_measure:invariance}).
    
    By construction of the measures in Lemma \ref{lemma:local_Siegel-Weil:local_Siegel-Weil:unique_Haar_measure} (via invariant differentials), we find that $\prod_v d g_v$ equals the Tamagawa measure on $G(\A)$ up to scaling by a constant which may depend on the lattices $\{L_v\}_{v < \infty}$ as well as $n$ and $m^{\flat}$ (coming from our normalization of local Whittaker functions $\tilde{W}^*_{T^{\flat},v}$, Section \ref{ssec:Eisenstein:local_Whittaker:non-Arch}). We conclude that the measure $dg$ on $G(\A_f)$ may depend on $V$, $n$, $m^{\flat}$, $F$, and the lattices $\{L_v\}_{v < \infty}$, but it does not otherwise depend on $T$ or $V^{\flat}$ or $\underline{x}$.
    
    \item Suppose $m^{\flat} = n$. Then $G_{\underline{x}}$ is the trivial group. Take $\mrm{vol}_{d g_{x,v}}(K_{\underline{x},v}) = 1$ for all $v$. Consider $v < \infty$ with $v \neq v_0$ and set $c_v = \mrm{vol}_{d g_v}(K_{L,v})^{-1}$. If $L_v$ is self-dual, then $c_v = 1$ by Lemma \ref{lemma:local_Siegel-Weil:non-Arch:spherical_and_measures}. In general, $d g_v$ may depend on $L_v$ (but not on $T$ or $T^{\flat}_v$).
    
    \item Suppose $m^{\flat} = n - 1$. Then $G_{\underline{x}}$ is isomorphic to the norm-one torus inside $\Res_{F/F^+} \G_m$. Assume $K_{\underline{x},v} \subseteq G_{\underline{x}}(F^+_v)$ is the unique maximal open compact subgroup for every $v$. Take $\mrm{vol}_{d g_{x,v}}(K_{\underline{x},v}) = 1$ for all $v$. Consider $v < \infty$ with $v \neq v_0$ and set $c'_v = \mrm{vol}_{d g_v}(K_{L,v})^{-1}$. If $L_v$ is self-dual, then $c'_v = 1$ if $F_v / F^+_v$ is unramified (resp. $c'_v = 2$ if $F_v / F^+_v$ is ramified) by Lemma \ref{lemma:local_Siegel-Weil:non-Arch:spherical_and_measures}. In general, $d g_v$ may depend on $L_v$, $m^{\flat}$ and the local invariant $\varepsilon(V^{\flat}_v)$ (but not on $T$ or $T^{\flat}_v$).

    We have
    \begin{equation}\notag
    \mrm{vol}(G_{\underline{x}}(F^+) \backslash G_{\underline{x}}(\A)) = \deg [G_{\underline{x}}(F^+) \backslash (G_{\underline{x}}(\A) / K_{\underline{x}})] = \frac{\deg (G_{\underline{x}}(F^+) \backslash G_{\underline{x}}(\A) / K_{\underline{x}}) }{w_F}
    \end{equation}
    where $\deg [ - ]$ denotes groupoid cardinality and $\deg ( - )$ denotes set cardinality. We have 
    \begin{equation}
    \deg (G_{\underline{x}}(F^+) \backslash G_{\underline{x}}(\A) / K_{\underline{x}}) = 2^{u - t} h_F h_{F^+}^{-1},
    \end{equation}
    where $t$ is the number of prime ideals of $F^+$ which ramify in $F$, and where $u \in \Z$ is such that $H^1(\Gal(F/F^+), \mc{O}_F^{\times}) \cong (\Z / 2 \Z)^u$ \cite[{(9)}]{Ono85}.
    A group cohomology computation (omitted) shows that $2^{-u} =  \# (\mc{O}_F^{\times} / (W \mc{O}_{F^+}^{\times}))/2$ (where $\#$ also means cardinality). \qedhere
\end{enumerate}
\end{proof}
    
        \section{Geometric Siegel--Weil}
        \label{sec:geometric_Siegel-Weil}
        
            \subsection{Degrees of \texorpdfstring{$0$}{0}-cycles}
            \label{ssec:geometric_Siegel-Weil:degrees}
                For the rest of the paper, we assume $F$ is an imaginary quadratic extension
of $\Q$, i.e. $F^+ = \Q$.

Let $L$ be any non-degenerate Hermitian $\mc{O}_F$-lattice of signature $(n - 1, 1)$ (not assuming $n$ is even). Let $\mc{M} \ra \Spec \mc{O}_F[1/d_L]$ be the associated moduli stack (\crefext{III:ssec:ab_var:integral_models}). Recall that $d_L \in \Z$ is a certain integer associated to $L$, with $d_L = 1$ if $L$ is self-dual when $2 \nmid \Delta$. Let $V \coloneqq L \otimes_{\mc{O}_F} F$ be the associated $F / \Q$ Hermitian space.

Consider an integer $m$ with $m = n$ or $m = n - 1$. Pick any embedding $F \ra \C$, and set $\mc{M}_{\C} \coloneqq \mc{M} \times_{\Spec \mc{O}_F} \Spec \C$, etc.. Given $T \in \mrm{Herm}_m(\Q)$ with $\rank T = n - 1$, recall that there is an associated Kudla--Rapoport special cycle $\mc{Z}(T) \ra \mc{M}$ \crefext{III:definition:part_II:global_special_cycles}. The base change $\mc{Z}(T)_{\C}$ is smooth, proper, and quasi-finite (and of dimension zero) over $\Spec \C$ \crefext{III:lemma:special_cycles_generically_smooth,III:lemma:horizontal_global_special_cycles_proper_quasi-finite}. 

For each place $v$ of $\Q$, select any $a_v \in \GL_m(F_v)$ such that ${}^t \overline{a}_v^{-1} T a_v^{-1} = \mrm{diag}(0,T^{\flat}_v)$ for some $T^{\flat}_v \in \mrm{Herm}_{n - 1}(F^+_v)$ with $\det T^{\flat}_v \neq 0$. Choose any $a_v^{\flat} \in \GL_{n - 1}(F_v)$ associated to $a_v$ via the Iwasawa decomposition, as in \crefext{equation:local_Siegel-Weil:uniformization_degree:a_v_Iwasawa_decomp} (if $m = n -1$, we can just take $a^{\flat}_v = a_v$).

For formation of local Whittaker functions, we use the standard additive character $\psi \colon \Q \backslash \A \ra \C^{\times}$ with $\psi_{\infty}(x) = e^{2 \pi i x}$. Suppose $\chi \colon F^{\times} \backslash \A_F^{\times} \ra \C^{\times}$ is a character satisfying $\chi|_{\A^{\times}} = \eta^n$, where $\eta$ is the quadratic character associated to $F / \Q$. 
For each prime $p$, we let $\varphi_p^{\flat} = \pmb{1}_{L_p}^{n - 1} \in \mc{S}(V(\Q_p)^{n - 1})$ where $\pmb{1}_{L_p}$ is the characteristic function of the lattice $L_p \subseteq V(\Q_p)$.

\begin{proposition}\label{proposition:geometric_Siegel-Weil:degrees}
Let $C \in \Q_{>0}$ be the volume constant from Lemma \ref{lemma:local_Siegel-Weil:uniformization_degree:main}(3), for the Hermitian space $V$ and with $v_0 = \infty$ in the notation of loc. cit.. In the situation above, we have
    \begin{equation}
    \deg \mc{Z}(T)_{\C} = \frac{h_F}{w_F} C \cdot \tilde{W}^*_{T^{\flat}_{\infty}, \infty}(a^{\flat}_{\infty}, 1/2)^{\circ}_n
        \prod_{p} \tilde{W}^*_{T^{\flat}_p, p}(a^{\flat}_p, 1/2, \Phi_{\varphi_v^{\flat}})_n.
    \end{equation}
\end{proposition}
\begin{proof}
As in Section \ref{ssec:local_Siegel-Weil:volume_forms}, we write $\Omega_T(R) \coloneqq \{\underline{x} \in (V \otimes_{\Q} R)^m : (\underline{x}, \underline{x}) = T\}$ for $\Q$-algebras $R$. Here $\deg \mc{Z}(T)_{\C}$ denotes the (stacky) degree of $\mc{Z}(T)_{\C}$ over $\Spec \C$, as explained at the end of \crefext{I:appendix:K0:K0_stacky}.

Suppose there is no tuple $\underline{x} \in V^m$ such that $(\underline{x}, \underline{x}) = T$. 
By the Hasse principle, we conclude $\Omega_T(\Q_{v_0}) = \emptyset$ for some place $v_0$ of $\Q$. Since $\rank(T) < n$, we must have $v_0 = \infty$ (i.e. for $v < \infty$, any non-degenerate Hermitian $F_v$ vector space of rank $n - 1$ embeds isometrically into any non-degenerate Hermitian $F_v$ vector space of rank $n$). We conclude that $T^{\flat}_{\infty}$ (and ${}^t \overline{a}^{\flat}_{\infty} T^{\flat}_{\infty} a^{\flat}_{\infty}$) has signature $(n - 1 - r, r)$ for some $r \geq 2$. The proposition holds in this case because $\tilde{W}^*_{T^{\flat}_{\infty}, \infty}(a^{\flat}_{\infty}, 1/2)^{\circ}_n = 0$ (by \crefext{equation:Eisenstein:local_Whittaker:Archimedean:special_value} or \crefext{equation:local_Siegel-Weil:local_Siegel-Weil:up-to-constant}).

Suppose there exists $\underline{x} \in V^m$ such that $(\underline{x}, \underline{x}) = T$. For such $\underline{x}$, write $\underline{x}_{\infty} \in V_{\R}^m$ and $\underline{x}_f \in (V \otimes_{\Q} \A_f)^m$ for the respective images.
By complex uniformization of special cycles \crefext{III:equation:Arch_uniformization:quotient_iso},\footnote{
    There is a typo in loc. cit., where ``framed'' should be removed from right-hand side.
} we have
    \begin{equation}\label{equation:geometric_Siegel-Weil:degrees:geometric}
    \deg \mc{Z}(T)_{\C} = \frac{h_F}{w_F} \cdot \deg \mc{D}(\underline{x}_{\infty}) \cdot \deg \Biggl [ U(V)(\Q) \backslash \coprod_{\substack{\underline{x} \in V^m \\ (\underline{x}, \underline{x}) = T}} \mc{D}(\underline{x}_f) \Biggr ],
    \end{equation}
i.e. \eqref{equation:geometric_Siegel-Weil:degrees:geometric} is the groupoid cardinality of
\cite[{(5.3.5)}]{corank1_ASW_III.pdf}.
Here $\deg \mc{D}(\underline{x}_{\infty})$ is the degree of the Archimedean local special cycle $\mc{D}(\underline{x}_{\infty}) \subseteq \mc{D}$ \crefext{II:ssec:Hermitian_domain:local_cycles} for any $\underline{x} \in V^m$ with $(\underline{x}, \underline{x}) = T$. We know $\mc{D}(\underline{x}_{\infty})$ is a single point if $T$ is positive semidefinite, and empty otherwise. Hence $\deg \mc{D}(\underline{x}_{\infty}) = \tilde{W}^*_{T^{\flat}_{\infty}, \infty}(a^{\flat}_{\infty}, 1/2)^{\circ}_n$ (by \cref{equation:Eisenstein:local_Whittaker:Archimedean:special_value}, the right-hand side is $1$ if $T^{\flat}_{\infty}$ is positive definite and $0$ otherwise).

We then use Lemma \ref{lemma:local_Siegel-Weil:uniformization_degree:main} to evaluate the groupoid cardinality in \cref{equation:geometric_Siegel-Weil:degrees:geometric}.
For the reader's convenience, we note that the present $\mathcal{D}(\underline{x}_f)$
is the $\mc{Z}(\underline{x}_f^{v_0})$ from loc. cit. (taking $v_0$ to be the Archimedean place),
and we are using transitivity of $U(V)(\Q)$
acting on the set of $\underline{x} \in V^m$ satisfying $(\underline{x}, \underline{x}) = T$
(Witt's theorem). The latter uses the fact that $\rank T = n - 1$,
which implies that any such $\underline{x}$ must span a non-degenerate subspace of $V$.
\end{proof}

\begin{remark}\label{remark:geometric_Siegel-Weil:degrees:self-dual}
Suppose $2 \nmid \Delta$ and that $L$ is self-dual (for the trace pairing, as is our running convention). We then have $C = 2 h_F / w_F$ in Proposition \ref{proposition:geometric_Siegel-Weil:degrees}. Take any $a \in \GL_m(F)$ such that ${}^t \overline{a}^{-1} T a^{-1} = \mrm{diag}(0, T^{\flat})$ where $T^{\flat} \in \mrm{Herm}_{n - 1}(\Q)$ with $\det T^{\flat} \neq 0$. For each place $v$ of $\Q$, let $a_v \coloneqq a \in \GL_m(F_v)$. Set $a^{\flat} = (a^{\flat}_v)_v \in \GL_m(\A_F)$ (running over places $v$ of $\Q$) in the notation above. The proposition then states
    \begin{equation}
    \deg \mc{Z}(T)_{\C} = 2 \frac{h_F^2}{w_F^2} \cdot \tilde{E}^*_{T^{\flat}}(a^{\flat}, 1/2)^{\circ}_n.
    \end{equation}
\end{remark}

\begin{remark}
As observed by Li and Zhang \cite[{Remark 4.6.2}]{LZ22unitary}, Proposition \ref{proposition:geometric_Siegel-Weil:degrees} may be proved using Rapoport--Zink non-Archimedean uniformization in essentially the same way. Indeed, the horizontal local special cycle $\mc{Z}(T)_{\ms{H}} \ra \Spec \mc{O}_F[1/d_L]$ is proper, quasi-finite, and flat \crefext{III:lemma:horizontal_global_special_cycles_proper_quasi-finite}, so we may calculate its degree in the fiber over any geometric point of $\Spec \mc{O}_F[1/d_L]$. 
Fix a geometric point in characteristic $p > 0$. Assume $p \neq 2$ if $2$ is nonsplit in $\mc{O}_F$, assume $L_p$ is self-dual, and assume either $p \nmid \Delta$ or that $L$ is self-dual and $2 \nmid \Delta$. Consider the $n$-dimensional positive definite non-degenerate Hermitian space $\mbf{V}$ with $\varepsilon(\mbf{V}_p) = -1$ and $\varepsilon(\mbf{V}_{\ell}) = \varepsilon(V_{\ell})$ for any $\ell \neq p$.

Using non-Archimedean uniformization, we may then argue as in the proof of Proposition \ref{proposition:geometric_Siegel-Weil:degrees} (see \crefext{III:equation:non-Arch_uniformization:horizontal:local_global_degrees:non-arithmetic}), using the special value formula for degrees of local special cycles \crefext{I:lemma:non-Arch_identity:statement:special_value}, and the formula for uniformization degrees (Lemma \ref{lemma:local_Siegel-Weil:uniformization_degree:main}) for $\mbf{V}$ and $v_0 = p$.
\end{remark}
            
            \subsection{Complex volumes}
            \label{ssec:geometric_Siegel-Weil:complex_volumes}
                Assume $2$ is unramified in $\mc{O}_F$. For even integers $n \in \Z_{>0}$, we show that the global normalizing factors $\Lambda_n(s)^{\circ}_n$ (Section \ref{ssec:Eisenstein:global_normalized_Fourier:global_normalization}) encode complex volumes of certain unitary Shimura varieties (Propositions \ref{proposition:geometric_Siegel-Weil:complex_volumes:positive_definite} and \ref{proposition:geometric_Siegel-Weil:complex_volumes:n-1_1}).

First consider $n \equiv 0 \pmod{4}$. Let $V$ be the unique $F / \Q$ Hermitian space of signature $(n,0)$ which satisfies $\varepsilon(V_p) = 1$ for all primes $p$ (with $\varepsilon$ as in \crefext{I:ssec:Hermitian_conventions:lattices}, i.e. $V$ contains a full-rank self-dual $\mc{O}_F$-lattice). Set $G \coloneqq U(V)$, let $L \subseteq V$ be a full-rank self-dual lattice, and write $K_{L,f} \subseteq G(\A_f)$ for the ad\`elic stabilizer of $L$. The following proposition should be a special case of a unitary analogue of the classical Siegel mass formula. It is included for comparison with the analogous volume identity for a signature $(n - 1, 1)$ unitary complex Shimura variety. The left-hand side counts self-dual positive definite $\mc{O}_F$-lattices of rank $n$, weighted by the inverses of the sizes of their automorphism groups.

\begin{proposition}\label{proposition:geometric_Siegel-Weil:complex_volumes:positive_definite}
We have
    \begin{equation}
    \# [G(\Q) \backslash ( G(\A_f) / K_{L,f})] = 2 \Lambda_n(0)^{\circ}_n
    \end{equation}
where the left-hand side denotes groupoid cardinality.
\end{proposition}
\begin{proof}
Let $\psi \colon \Q \backslash \A \ra \C^{\times}$ be the standard additive character with $\psi_{\infty}(x) = e^{2 \pi i x}$. Let $\chi \colon \A_F^{\times} \ra \C^{\times}$ be the trivial character. 

For $v = \infty$, let $\varphi_v(\underline{x}) = e^{2 \pi i \mrm{tr}(\underline{x},\underline{x})} \in \mc{S}(V(\R)^n)$ and let $T \in \mrm{Herm}_n(\R)$ be an arbitrary positive definite matrix. For $v < \infty$ corresponding to a prime $p$, let $\varphi_v = \pmb{1}_{L_v}^n \in \mc{S}(V(\Q_p)^n)$ and let $T$ be the Gram matrix for any basis of $L_v$.
For such $T$, we have $W^*_{T,v}(s_0)^{\circ}_n = 1$ for all $v$ (in both Archimedean and non-Archimedean cases),
by \eqref{equation:Eisenstein:local_Whittaker:Archimedean:special_value} and \eqref{equation:Eisenstein:local_Whittaker:local_densities_spherical:often_1}.
Recall $W^*_{T,v}(s)^{\circ}_n = \Lambda_{T,v}(s)^{\circ}_n W_{T,v}(s, \Phi_{\varphi_v})$ if $v < \infty$ (resp. $W^*_{T,v}(s)^{\circ}_n e^{- 2 \pi \mrm{tr}(T)} = \Lambda_{T,v}(s)^{\circ}_n W_{T,v}(s, \Phi_{\varphi_v})$ if $v = \infty$); see Section \ref{ssec:Eisenstein:Weil_rep:Weil_rep}. 

Using these data, the local Siegel--Weil formula and its proof (Lemma \ref{lemma:local_Siegel-Weil:local_Siegel-Weil:unique_Haar_measure}) for each place $v$ of $\Q$ shows that $\mrm{vol}(G(\R) \times K_{L,f})^{-1} = \Lambda_n(0)^{\circ}_n$ for the Tamagawa measure on $G(\A)$.
Since $G$ has Tamagawa number $2$ \cite[{\S 4}]{Ichino04}, the proposition follows.
\end{proof}

Next, consider $n \equiv 2 \pmod{4}$. Let $V$ be the unique $n$-dimensional $F / \Q$ Hermitian space of signature $(n - 1, 1)$ which satisfies $\varepsilon(V_p) = 1$ for all primes $p$. Again, set $G \coloneqq U(V)$, let $L \subseteq V$ be a full-rank self-dual lattice, and write $K_{L,f} \subseteq G(\A_f)$ for the ad\`elic stabilizer of $L$. For sufficiently small open compact $K_f \subseteq G(\A_f)$, there is a complex (analytic) Shimura variety
    \begin{equation}
    \mrm{Sh}_{K_f, \C} = G(\A) \backslash (\mc{D} \times G(\A_f) / K_f)
    \end{equation}
of dimension $n - 1$, where $\mc{D}$ is the Hermitian symmetric domain from \crefext{II:ssec:Hermitian_domain:setup} (parameterizing maximal negative definite subspaces of $V_{\R}$); the $V$ of loc. cit. is our $V_{\R}$, with $\C = F \otimes_{\Q} \R$-action. The metrized tautological bundle $\widehat{\mc{E}}^{\vee}$ of loc. cit. descends to $\mrm{Sh}_{K_f, \C}$. For any open compact $K'_f \subseteq G(\A_f)$ and any sufficiently small $K_f \subseteq K'_{f}$, we set
    \begin{equation}
    \mrm{vol}(\mrm{Sh}_{K_f, \C}) \coloneqq \int_{\mrm{Sh}_{K_f, \C}} c_1(\widehat{\mc{E}})^{n - 1} \quad \quad \mrm{vol}(\mrm{Sh}_{K_{f}',\C}) \coloneqq \frac{1}{[K_{f}' : K_f]} \mrm{vol}(\mrm{Sh}_{K_f, \C}).
    \end{equation}

If $K_{L',f} \subseteq G(\A_f)$ is the ad\`elic stabilizer of a full-rank lattice $L' \subseteq V$ which is self-dual for the Hermitian pairing, the quantity $\mrm{vol}(\mrm{Sh}_{K'_{L,f},\C})$ was computed explicitly in \cite[{Theorem A}]{BH21}. We show that the level $K_{L,f}$ (self-dual for the \emph{trace pairing}) removes the additional factors at ramified primes in loc. cit.,
and that the resulting complex volume agrees with $2 \Lambda_n(0)^{\circ}_n$ exactly. 

The volume identity should also follow from \cite[{Footnote 11}]{LL22I} (or possibly other geometric Siegel--Weil results). We instead compute $\mrm{vol}(\mrm{Sh}_{K_{L,f},\C})$ using \cite[{Theorem A}]{BH21} by calculating the ``change of level'' via the following lemma.

\begin{lemma}\label{lemma:Eisenstein:global_normalized_Fourier:volume_identity:comparison_volume_lattice_stabilizers}
Let $E^+_v$ be a non-Archimedean local field of odd residue cardinality $q_v$, and let $E_v/E^+_v$ be a ramified quadratic extension with involution $a \mapsto a^{\s}$.

Let $W$ be a rank $2d$ non-degenerate $E_v/E^+_v$ Hermitian space, and assume $W$ contains a full-rank lattice $M \subseteq W$ which is self-dual (for the trace pairing). Let $M' \subseteq W$ be any full-rank lattice which is self-dual for the Hermitian pairing.

If $K, K' \subseteq U(W)$ are the stabilizers of $M$ and $M'$ respectively, we have
    \begin{equation}
    \frac{\mrm{vol}(K)}{\mrm{vol}(K')} = 2^{-1} (1 + q_v^{d})
    \end{equation}
for any Haar measure on $U(W)$.
\end{lemma}
\begin{proof}
We know that any two full-rank lattices in $W$ which are self-dual (resp. self-dual for the Hermitian form) are isomorphic \cite[{Proposition 8.1}]{Jacobowitz62} (false if $E^+_v$ is allowed to have residue characteristic $2$). 
Hence $\mrm{vol}(K)/\mrm{vol}(K')$ does not depend on the choice of $M$ and $M'$ (nor the choice of Haar measure).

Let $\varpi$ be a uniformizer of $E_v$, and assume $\varpi^{\s} = - \varpi$. The lattices $M$ and $M'$ admit bases with Gram matrices
    \begin{equation}
    \begin{pmatrix}
    0 & \varpi^{-1} \\
    -\varpi^{-1} & 0
    \end{pmatrix}
    \quad \quad
    \begin{pmatrix}
    0 & 1 \\
    1 & 0
    \end{pmatrix}
    \end{equation}
respectively. Choose a basis $e_1, \ldots, e_{2d}$ for $M$ with Gram matrix as above. We may assume that $M'$ is the lattice with basis $e_1, \ldots, e_d, \varpi e_{d+1}, \ldots, \varpi e_{2d}$. Let $\overline{W}$ (resp. $\overline{W}'$) be the $2d$-dimensional vector space over $\F_{q_v}$ with symplectic pairing (resp. bilinear pairing) given by the block matrices
    \begin{equation}
    \begin{pmatrix}
    0 & 1 \\
    -1 & 0 
    \end{pmatrix}
    \quad \quad
    \text{resp.} \quad \begin{pmatrix}
    0 & 1 \\
    1 & 0
    \end{pmatrix}.
    \end{equation}
If $P_W \subseteq \mrm{Sp}(\overline{W})$ and $P_{\overline{W}'} \subseteq \mrm{O}(\overline{W})$ are the subgroups of upper triangular matrices (in $d \times d$ blocks), we have
    \begin{align}
    \# (K / (K \cap K')) & = \# (\mrm{Sp}(\overline{W})(\F_{q_v}) / P_{\overline{W}}(\F_{q_v})) \\
    \# (K' / (K \cap K')) & = \# (\mrm{O}(\overline{W}')(\F_{q_v}) / P_{\overline{W}'}(\F_{q_v})).
    \end{align}
The lemma now follows from the formulas
    \begin{align}
    & \# \mrm{Sp}(\overline{W})(\F_{q_v}) = q_v^{d^2} \prod_{i = 1}^d (q_v^{2i} - 1)
    && \# \mrm{O}(\overline{W}')(\F_{q_v}) = 2 q_v^{d(d-1)} (q_v^d + 1)^{-1} \prod_{i = 1}^d (q_v^{2i} - 1)  \notag \\
    & \# P_{\overline{W}}(\F_{q_v}) = q_v^{d(d+1)/2} \prod_{i = 1}^d (q_v^d - q_v^{i - 1})    
    && \# P_{\overline{W}'}(\F_{q_v}) = q_v^{d(d-1)/2} \prod_{i = 1}^d (q_v^d - q_v^{i - 1}). \notag \qedhere
    \end{align}
\end{proof}

We return to the global situation with $F / \Q$ as above and $L \subseteq V$ a self-dual lattice.

\begin{proposition}\label{proposition:geometric_Siegel-Weil:complex_volumes:n-1_1}
We have
    \begin{equation}
    \mrm{vol}(\mrm{Sh}_{K_{L,f}, \C}) = 2 \Lambda_n(0)^{\circ}_n.
    \end{equation}
\end{proposition}
\begin{proof}
If $K_{L',f} \subseteq G(\A_f)$ is the ad\`elic stabilizer of a full-rank lattice $L' \subseteq V$ which is self-dual for the Hermitian pairing, the result \cite[{Theorem A}]{BH21} (see also \cite[{Theorem 5.5.1}]{BH21} to compare $c_1(\widehat{\mc{E}})$ with the Chern form of the metrized Hodge bundle; note our $\widehat{\mc{E}}$ is $\widehat{\mc{L}}$ in loc. cit. (up to restricting)) gives
    \begin{equation}
    \mrm{vol}(\mrm{Sh}_{K_{L',f}, \C}) = \bigg [ 2^{1 - o(\Delta)} \prod_{\ell \mid \Delta} (1 + \varepsilon(V_{\ell}) \ell^{-n/2}) \prod_{j = 1}^{n} \frac{\Delta^{j/2} \Gamma(s+j) L(2s + j, \eta^j)}{2^j \pi^{s + j}} \bigg ]_{s = 0}
    \end{equation}
where $o(\Delta)$ is the number of primes dividing $\Delta$. We assumed $\varepsilon(V_{\ell}) = 1$ for all $\ell$, and a direct computation shows
    \begin{equation}
    \Lambda_n(s)^{\circ}_n = \Delta^{n/2(s - 1)} \prod_{j = 1}^{n} \frac{\Delta^{j/2} \Gamma(s+j) L(2s + j, \eta^j)}{2^j \pi^{s + j}}
    \end{equation}
(using even-ness of $n$). 
The claim now follows from the computation of $\mrm{vol}(K_{L,f})/\mrm{vol}(K_{L',f})$ (for any Haar measure on $G(\A_f)$) from Lemma \ref{lemma:Eisenstein:global_normalized_Fourier:volume_identity:comparison_volume_lattice_stabilizers}. Note that the only discrepancy between $\mrm{vol}(K_{L,f})$ and $\mrm{vol}(K_{L',f})$ is at ramified primes, since self-dual lattices for the Hermitian pairing are the same as self-dual lattices at unramified primes.
\end{proof}
    
        \section{Arithmetic Siegel--Weil}
        \label{sec:arithmetic_Siegel-Weil}
    
            \subsection{Main theorems}
            \label{ssec:arithmetic_Siegel-Weil:main_results}
                This section contains the statements and proofs of our main global results (Theorem \ref{theorem:arithmetic_Siegel-Weil:main_results:main} and the secondary Theorem \ref{theorem:arithmetic_Siegel-Weil:main_results:Archimedean}). Theorem \ref{theorem:arithmetic_Siegel-Weil:main_results:main} relies on essentially all preceding results in our four-part sequence of papers. We necessarily heavily cite our companion papers \cite{corank1_ASW_I.pdf,corank1_ASW_II.pdf,corank1_ASW_III.pdf}. In the proof, we explain how to combine our local main results (proved in \crefext{I:sec:non-Arch_identity} \crefext{II:sec:Archimedean_identity}) and a (new) ``local diagonalization'' argument to deal with singular $T$ (including those which are not necessarily $\GL_n(\mc{O}_F)$-conjugate to a block diagonal matrix with nonsingular diagonal blocks).

Assume $2 \nmid \Delta$, and let $L$ be any non-degenerate self-dual Hermitian $\mc{O}_F$-lattice of signature $(n - 1, 1)$. Set $n \coloneqq \rank L$, and note $n \equiv 2 \pmod{4}$ (by the global product formula for local invariants of Hermitian spaces; note $\varepsilon(L_p) = 1$ for all primes $p$).

Form the associated (smooth) moduli stack $\mc{M} \ra \Spec \mc{O}_F$ (\crefext{III:ssec:ab_var:integral_models} and \crefext{I:ssec:part_I:arith_intersections:integral_models}). We are imposing ``no level structure'' on $\mc{M}$ (i.e. $K_{0,f} \times K_f = K_{L_0,f} \times K_{L,f}$ in the notation of \crefext{III:ssec:ab_var:level_structure}).

For any $m$, given $T \in \mrm{Herm}_m(\Q)$ (with $F$-coefficients), and given $y \in \mrm{Herm}_m(\R)_{>0}$ (with $\C$-coefficients), recall that there is an arithmetic special cycle class $[\widehat{\mc{Z}}(T)] \in \arithCh^m(\mc{M})_{\Q}$ \crefext{III:sec:arith_cycle_classes} and a normalized $T$-th Fourier coefficient $E^*_{T}(y,s)^{\circ}_n$ (\cref{ssec:Eisenstein:global_normalized_Fourier:global_normalization}) of a $U(m,m)$ Eisenstein series.
Recall that the construction of $[\widehat{\mc{Z}}(T)]$ involves a choice of Green current $g_{T,y}$.
If $\rank(T) \geq n - 1$ or if $T$ is nonsingular and not positive definite, we are using the current $g_{T,y}$ from \crefext{III:ssec:Arch_uniformization:Archimedean} (constructed from the local analogue in \crefext{II:ssec:part_II:Hermitian_domain:corank1_modification}). The class $[\widehat{\mc{Z}}(T)]$ thus implicitly depends on $y$.

For special cycles $\mc{Z}(T)$ which are proper over $\Spec \mc{O}_F$, recall that we have defined certain arithmetic degrees without boundary contributions \cref{equation:part_IV:intro:results:if_proper}. These are the arithmetic degrees appearing in our main theorem below.

For use below, we record the expression
    \begin{equation}\label{equation:arithmetic_Siegel-Weil:main_results:corank_1_L-factor}
    \frac{\Lambda_n(s)_n^{\circ}}{\Lambda_{n - 1}(s + 1/2)^{\circ}_n} = - \frac{1}{2} L(2s + 1, \eta) \Gamma(s + 1) |\Delta|^{s + 1/2} \pi^{- s - 1}
    \end{equation}
which follows from our formula for the normalizing factor $\Lambda_m(s)^{\circ}_n$ \eqref{equation:global_normalized_Fourier:global_normalization:normalizing_factor}. We thus have
    \begin{equation}\label{equation:arithmetic_Siegel-Weil:main_results:height_constant}
    \frac{\Lambda_n(0)_n^{\circ}}{\Lambda_{n - 1}(1/2)^{\circ}_n} = - \frac{h_F}{w_F} \quad \quad \frac{d}{ds} \bigg|_{s = 0} \left ( \frac{\Lambda_n(s)_n^{\circ}}{\Lambda_{n - 1}(s + 1/2)^{\circ}_n} \right ) = 2 \frac{h_F}{w_F} h^{\mrm{CM}}_{\widehat{\mc{E}}^{\vee}}
    \end{equation}
where the left expression follows from the analytic class number formula, and $h^{\mrm{CM}}_{\widehat{\mc{E}}^{\vee}}$ is the height constant from \crefext{III:equation:part_II:arith_cycle_classes:Hodge_bundles:taut_height_constants}, arising from Faltings heights of elliptic curves with complex multiplication by $\mc{O}_F$ (for the purpose of the present paper, one could also take \cref{equation:arithmetic_Siegel-Weil:main_results:height_constant} as the definition of $h^{\mrm{CM}}_{\widehat{\mc{E}}^{\vee}}$).

\begin{theorem}[Co-rank $1$ arithmetic Siegel--Weil]\label{theorem:arithmetic_Siegel-Weil:main_results:main}
Assume the prime $2$ splits in $\mc{O}_F$.
\begin{enumerate}[(1)]
    \item For any $T \in \mrm{Herm}_n(\Q)$ with $\rank(T) = n - 1$ and any $y \in \mrm{Herm}_n(\R)_{>0}$, we have
        \begin{equation}\label{equation:arithmetic_Siegel-Weil:main_results:main:singular}
        \widehat{\deg}([\widehat{\mc{Z}}(T)]) = \frac{h_F}{w_F} \frac{d}{d s} \bigg |_{s = 0} E^*_{T}(y,s)^{\circ}_n.
        \end{equation}
    \item For any $T^{\flat} \in \mrm{Herm}_{n - 1}(\Q)$ with $\det T^{\flat} \neq 0$ and any $y^{\flat} \in \mrm{Herm}_{n-1}(\R)_{>0}$, we have
        \begin{equation}\label{equation:arithmetic_Siegel-Weil:main_results:main:nonsingular}
        \widehat{\deg}([\widehat{\mc{Z}}(T^{\flat})] \cdot \widehat{c}_1(\widehat{\mc{E}}^{\vee})) = 2 \frac{h_F}{w_F} \frac{d}{d s} \bigg|_{s = 0} \left ( \frac{\Lambda_n(s)_n^{\circ}}{\Lambda_{n - 1}(s + 1/2)^{\circ}_n} E^*_{T^{\flat}}(y^{\flat}, s + 1/2)^{\circ}_n \right ).
        \end{equation}
\end{enumerate}
\end{theorem}
\begin{proof}
In the theorem statement, $[\widehat{\mc{Z}}(T)]$ and $[\widehat{\mc{Z}}(T^{\flat})]$ are implicitly formed with respect to $y$ and $y^{\flat}$, respectively. Note that $E^*_T(y,s)^{\circ}_n$ is a normalized Fourier coefficient for a $U(n,n)$ Eisenstein series, while $E^*_{T^{\flat}}(y^{\flat}, s)^{\circ}_n$ is a normalized Fourier coefficient for a $U(n - 1, n - 1)$ Eisenstein series. In the theorem statement, note that $\mc{Z}(T) \ra \Spec \mc{O}_F$ and $\mc{Z}(T^{\flat}) \ra \Spec \mc{O}_F$ are both proper \crefext{III:lemma:global_special_cycles_proper}, so we may use \cref{equation:part_IV:intro:results:if_proper} to define arithmetic degrees without boundary contributions.

Note that \cref{theorem:arithmetic_Siegel-Weil:main_results:main}(2) is the special case of \cref{theorem:arithmetic_Siegel-Weil:main_results:main}(1) when $T = \mrm{diag}(0, T^{\flat})$ and $y = \mrm{diag}(1, y^{\flat})$. This follows from the unfolding of Fourier coefficients in Corollary \ref{corollary:Eisenstein:singular_Fourier:corank_1} (also the functional equation in Lemma \ref{lemma:Eisenstein:global_normalized_Fourier:global_normalization:global_functional_equation}) and from the definition of arithmetic degrees in \eqref{equation:part_IV:intro:results:if_proper}.

Fix $T$ and $y$ as in the statement of part (1) (not necessarily block diagonal). Fix any prime $p$. It is enough to show that \eqref{equation:arithmetic_Siegel-Weil:main_results:main:singular} holds modulo $\sum_{\ell \neq p} \Q \cdot \log \ell$ (i.e. as elements of the additive quotient $\R / (\sum_{\ell \neq p} \Q \cdot \log \ell)$), where the sum runs over primes $\ell \neq p$. Varying the prime $p$ removes this discrepancy (giving an equality as elements of $\R$) because the real numbers $\log \ell$ (ranging over all primes $\ell$ in $\Z$) form a $\Q$-linearly independent set.

(\textit{Step 1: Diagonalize}) For convenience, we fix an embedding $F \ra \C$. Pick any $b \in \GL_m(F)$ such that ${}^t \overline{b}^{-1} T b^{-1} = \mrm{diag}(0, T^{\flat})$ for some $T^{\flat} \in \mrm{Herm}_{n - 1}(\Q)$ with $\det T^{\flat} \neq 0$. We may (and do) assume $b \in \GL_n(\mc{O}_F \otimes_{\Z} \Z_{(p)})$ as well. The proof below will show that the theorem holds modulo $\Q \cdot \log \ell$ for primes $\ell$ such that $b \not \in \GL_n(\mc{O}_F \otimes_{\Z} \Z_{(\ell)})$.

For each place $v$ of $\Q$, select any $b_v^{\#} \in \GL_1(F_v)$ and $b_v^{\flat} \in \GL_{n - 1}(F_v)$ associated to an Iwasawa decomposition of $b_v \in \GL_n(F_v)$, as in \eqref{equation:local_Siegel-Weil:uniformization_degree:a_v_Iwasawa_decomp} (where $b_v$ denotes the image of $b$). Also consider the (unique) decomposition
    \begin{equation}
    b y {}^t \overline{b} = 
    \begin{pmatrix}
    1 & c \\
    0 & 1
    \end{pmatrix}
    \begin{pmatrix}
    y^{\#} & 0 \\
    0 & y^{\flat}
    \end{pmatrix}
    \begin{pmatrix}
    1 & 0 \\
    {}^t \overline{c} & 1
    \end{pmatrix}
    \end{equation}
as in \crefext{II:equation:Arch_uniformization:Archimedean:decomp_y_diagonal}, where $c \in M_{1,n-1}(\C)$, $y^{\#} \in \R_{>0}$, and $y^{\flat} \in \mrm{Herm}_{n - 1}(\R)_{>0}$. Pick any $a_{\infty}^{\#} \in \GL_1(\C)$ and $a_{\infty}^{\flat} \in \GL_{n - 1}(\C)$ such that $a_{\infty}^{\#} {}^t \overline{a}_{\infty}^{\#} = y^{\#}$ and $ a_{\infty}^{\flat} {}^t \overline{a}_{\infty}^{\flat} = y^{\flat}$. 

Let $a^{\#} \in \GL_1(\A_F)$ be the element with component $a^{\#}_v \coloneqq b^{\#}_v$ for places $v < \infty$ and $a^{\#}_v \coloneqq a^{\#}_{\infty}$ for the place $v = \infty$. Similarly define $a^{\flat} \in \GL_{n - 1}(\A_F)$, and set $a \coloneqq \mrm{diag}(a^{\#}, a^{\flat}) \in \GL_n(\A_F)$.

By unfolding for co-rank $1$ Fourier coefficients (Corollary \ref{corollary:Eisenstein:singular_Fourier:corank_1}) and Fourier coefficient invariance properties (see \eqref{equation:Eisenstein:setup:Fourier_and_Whittaker:invariance}, \eqref{equation:eisenstein:setup:Fourier_and_Whittaker:invariance_y},
\eqref{equation:local_Whittaker:Archimedean:linear_invariance}, and \eqref{equation:local_Whittaker:non-Arch:linear_invariance} for $U(m)$ invariance when $v \mid \infty$ and $\GL_m(\mc{O}_{F^+_v})$ invariance when $v < \infty$), we find
    \begin{align*}
    E^*_T(y,s)^{\circ}_n & = \chi_{\infty}(d)^{-1} \det(y)^{- n / 2} E^*_{T'}(m(a), s)^{\circ}_n = \tilde{E}^*_{T'}(a, s)^{\circ}_n
    \\
    & = |a^{\#}|_{F}^{s} \frac{\Lambda_n(s)^{\circ}_n}{\Lambda_{n - 1}(s + 1/2)^{\circ}_n} \tilde{E}^*_{T^{\flat}}(a^{\flat}, s + 1/2)^{\circ}_n 
    \\
    & \mathrel{\phantom{=}} - | a^{\#} |_F^{- s} \frac{\Lambda_n(- s)^{\circ}_n}{\Lambda_{n - 1}(- s + 1/2)^{\circ}_n} \tilde{E}^*_{T^{\flat}}(a^{\flat}, s - 1/2)^{\circ}_n
    \end{align*}
where $T' \coloneqq \mrm{diag}(0, T^{\flat})$.
The sign $\pm$ from Corollary \ref{corollary:Eisenstein:singular_Fourier:corank_1} is $-1$,
by the assumptions $F^+ = \Q$ and $n \equiv 2 \pmod{4}$.
We remind the reader that the notation $E^*_T(-,s)^{\circ}_n$ is overloaded (Section \ref{ssec:Eisenstein:global_normalized_Fourier:global_normalization}, also end of Section \ref{ssec:Eisenstein:setup:adelic_v_classical}) and has slightly different meaning when ``$-$'' is $y \in \mrm{Herm}_m(\R)_{>0}$ versus $h \in U(m,m)(\A)$ (e.g. $h = m(a)$).

(\textit{Step 2: Leibniz rule}) Since $n \equiv 2 \pmod{4}$, the functional equation for $\tilde{E}^*_{T^{\flat}}(a^{\flat},s)^{\circ}_n$ (Lemma \ref{lemma:Eisenstein:global_normalized_Fourier:global_normalization:global_functional_equation}, noting
that the sign is $+1$) implies
    \begin{align}
    & \frac{d}{ds} \bigg|_{s = 0} E^*_{T}(y,s)^{\circ}_n
    \\
    & = 2 \frac{d}{ds} \bigg|_{s = 0} \left ( |a^{\#}|_{F}^{s} \frac{\Lambda_n(s)^{\circ}_n}{\Lambda_{n - 1}(s + 1/2)^{\circ}_n} \tilde{E}^*_{T^{\flat}}(a^{\flat}, s + 1/2)^{\circ}_n \right ). \notag
    \end{align}
Since $\det T^{\flat} \neq 0$, we may factorize $\tilde{E}^*_{T^{\flat}}(a^{\flat}, s + 1/2)^{\circ}_n$ into a product of (variants of) normalized local Whittaker functions \eqref{equation:Eisenstein:global_normalized_Fourier:global_normalization:factorization}. Also recall the formulas in \eqref{equation:arithmetic_Siegel-Weil:main_results:height_constant}. We have $|a^{\#}_{\ell}|_{\ell} = 1$ ($\ell$-adic norm of $a^{\#}_{\ell}$) for any prime $\ell$ such that $b \in \GL_n(\mc{O}_F \otimes_{\Z} \Z_{(\ell)})$ (by construction, this includes $\ell = p$). By the Leibniz rule, we thus find
    \begin{align}
    & \left ( \frac{ 2 h_F}{w_F} \right )^{-1} \frac{d}{ds} \bigg|_{s = 0} E^*_{T}(y,s)^{\circ}_n 
    \label{equation:arithmetic_Siegel-Weil:main_results:main:Eisenstein_decomp:total}
    \\
    & = 2 h_{\widehat{\mc{E}}^{\vee}}^{\mrm{CM}} \tilde{E}^*_{T^{\flat}}(a^{\flat}, 1/2)^{\circ}_n 
    \label{equation:arithmetic_Siegel-Weil:main_results:main:Eisenstein_decomp:height_constant}
    \\
    & \mathrel{\phantom{=}} - \left (\frac{d}{ds} \bigg|_{s = 1/2} \left ( |a^{\#}_{\infty}|_{\infty}^s \tilde{W}^*_{T^{\flat}, \infty}(a^{\flat}_{\infty}, s)^{\circ}_n \right ) \right ) \prod_{\ell} \tilde{W}^*_{T^{\flat},\ell}(a^{\flat}_{\ell}, 1/2)^{\circ}_n
    \label{equation:arithmetic_Siegel-Weil:main_results:main:Eisenstein_decomp:Archimedean}
    \\
    & \mathrel{\phantom{=}} - \left ( \frac{d}{ds} \bigg|_{s = 1/2} \tilde{W}^*_{T^{\flat}, p}(a^{\flat}_{p}, s)^{\circ}_n \right ) \prod_{v \neq p} \tilde{W}^*_{T^{\flat},v}(a^{\flat}_v, 1/2)^{\circ}_n
    \label{equation:arithmetic_Siegel-Weil:main_results:main:Eisenstein_decomp:at_p}
    \\
    & \mathrel{\phantom{=}} - \sum_{\ell \neq p} \left ( \frac{d}{ds} \bigg|_{s = 1/2} |a^{\#}_{\ell}|^s_{\ell} \tilde{W}^*_{T^{\flat}, \ell}(a^{\flat}_{\ell}, s)^{\circ}_n \right ) \prod_{v \neq \ell} \tilde{W}^*_{T^{\flat},v}(a^{\flat}_v, 1/2)^{\circ}_n.
    \label{equation:arithmetic_Siegel-Weil:main_results:main:Eisenstein_decomp:away_p}
    \end{align}
The product in \eqref{equation:arithmetic_Siegel-Weil:main_results:main:Eisenstein_decomp:Archimedean} runs over all primes $\ell$ (not including the Archimedean place $\infty$). The products in \eqref{equation:arithmetic_Siegel-Weil:main_results:main:Eisenstein_decomp:at_p} and \eqref{equation:arithmetic_Siegel-Weil:main_results:main:Eisenstein_decomp:away_p} run over all places $v$ of $\Q$ (with $v \neq p$ or $v \neq \ell$ as indicated), including $v = \infty$. The sum in \eqref{equation:arithmetic_Siegel-Weil:main_results:main:Eisenstein_decomp:away_p} runs over all primes $\ell \neq p$. We remind the reader that $|a^{\#}_{\infty}|_{\infty} = \overline{a}^{\#}_{\infty} a^{\#}_{\infty} \in \R_{>0}$, by definition.

For all but finitely many primes $\ell$, the Hermitian matrix ${}^t \overline{a}^{\flat}_{\ell} T^{\flat} a^{\flat}_{\ell} \in \mrm{Herm}_{n - 1}(\Q_{\ell})$ defines a (non-degenerate) self-dual Hermitian $\mc{O}_F \otimes_{\Z} \Z_{\ell}$-lattice. For such $\ell$, we have $\tilde{W}^*_{T^{\flat},\ell}(a^{\flat}_{\ell}, s)^{\circ}_n$ identically equal to $1$ (as a function in the $s$-variable). This follows from \eqref{equation:Eisenstein:local_Whittaker:local_densities_spherical:often_1} and an invariance property for local Whittaker functions \eqref{equation:local_Whittaker:non-Arch:linear_invariance}. In particular, the sums and products are finite in the right-hand side of \eqref{equation:arithmetic_Siegel-Weil:main_results:main:Eisenstein_decomp:total}.

For every prime $\ell$, we have $\tilde{W}^*_{T^{\flat},\ell}(a^{\flat}_\ell, s + 1/2)^{\circ}_n \in \Z[\ell^{-1}, \ell^{-s}, \ell^{s}]$ (see \eqref{equation:local_Whittaker:local_densities_spherical:Den_star}, and again the invariance property in \eqref{equation:local_Whittaker:non-Arch:linear_invariance}). We also have $\tilde{W}^*_{T^{\flat}, v}(a^{\flat}_{v}, 1/2)^{\circ}_n \in \Q$ for all places $v$ of $\Q$ (if $v \mid \infty$, this quantity is $1$ if $T^{\flat}$ is positive definite and $0$ otherwise by \eqref{equation:Eisenstein:local_Whittaker:Archimedean:special_value}).
The quantity in \eqref{equation:arithmetic_Siegel-Weil:main_results:main:Eisenstein_decomp:at_p} thus lies in $\Q \cdot \log p$, and the quantity in \eqref{equation:arithmetic_Siegel-Weil:main_results:main:Eisenstein_decomp:away_p} thus lies in $\sum_{\ell \neq p} \Q \cdot \log \ell$.

As we explain below, every quantity on the right-hand side of \eqref{equation:arithmetic_Siegel-Weil:main_results:main:Eisenstein_decomp:total} has geometric meaning via our main local results, at least modulo $\Q \cdot \log \ell$ for primes $\ell$ such that $b \not \in \GL_n(\mc{O}_F \otimes_{\Z} \Z_{(\ell)})$.

(\textit{Step 3a: Local geometric interpretation: complex degree}) Set $\mc{Z}(T)_{\C} = (\mc{Z}(T) \times_{\Spec \mc{O}_F} \Spec \C)$ for the embedding $F \ra \C$ fixed above. We have $\deg \mc{Z}(T)_{\C} = (\deg_{F} \mc{Z}(T) \times_{\Spec \mc{O}_F} \Spec F) = 2 \deg_{\Q} (\mc{Z}(T) \times_{\Spec \Z} \Spec \Q) \eqqcolon \deg_{\Z} \mc{Z}(T)_{\ms{H}}$. Here $\deg_F$ and $\deg_{\Q}$ denote stacky degrees over $\Spec F$ and $\Spec \Q$, respectively, as defined at the end of \crefext{I:appendix:K0:K0_stacky}.

By the geometric Siegel--Weil formula for Kudla--Rapoport $0$-cycles over $\C$ (Proposition \ref{proposition:geometric_Siegel-Weil:degrees}, also Remark \ref{remark:geometric_Siegel-Weil:degrees:self-dual}), we conclude
    \begin{equation}\label{equation:arithmetic_Siegel-Weil:main_results:main:geometric_decomp:height_constant:global}
    \deg_{\Z} \mc{Z}(T)_{\ms{H}} = 2 \deg \mc{Z}(T)_{\C} = \frac{4 h_F^2}{w_F^2} \tilde{E}^*_{T^{\flat}}(a^{\flat}, 1/2)^{\circ}_n.    
    \end{equation}
This gives a geometric interpretation of \eqref{equation:arithmetic_Siegel-Weil:main_results:main:Eisenstein_decomp:height_constant}.

(\textit{Step 3b: Local geometric interpretation: at $\infty$}) 
We claim that
    \begin{equation}\label{equation:arithmetic_Siegel-Weil:main_results:main:geometric_decomp:Archimedean:local}
    \mrm{Int}_{\infty}(T,y) = - \frac{d}{ds} \bigg|_{s = 1/2} \left ( |a^{\#}_{\infty}|^s_{\infty} \tilde{W}^*_{T^{\flat}, \infty}(a^{\flat}_{\infty}, s)^{\circ}_n \right ) \mod \sum_{\substack{\ell \text{ such that} \\ b \not \in \GL_n(\mc{O}_F \otimes_{\Z} \Z_{(\ell)})}} \Q \cdot \log \ell
    \end{equation}
where $\mrm{Int}_{\infty}(T,y)$ is the geometric quantity defined in \crefext{III:equation:Arch_uniformization:Archimedean:Int_infty_local}. 

Indeed, \crefext{II:equation:arch_uniformization:Archimedean:singular_diagonalized_current} implies 
    \begin{align}
    & \mrm{Int}_{\infty}(T,y) = 
    \\
    & \begin{cases}
    \mrm{Int}_{\infty}(T^{\flat}, a^{\flat}_{\infty} {}^t \overline{a}^{\flat}_{\infty}) - \log(|a^{\#}_{\infty}|_{\infty}) \mod \sum_{\substack{\ell \text{ such that} \\ b \not \in \GL_n(\mc{O}_F \otimes_{\Z} \Z_{(\ell)})}} \Q \cdot \log \ell & \text{if $T^{\flat} > 0$}
    \\
    \mrm{Int}_{\infty}(T^{\flat}, a^{\flat}_{\infty} {}^t \overline{a}^{\flat}_{\infty}) \mod \sum_{\substack{\ell \text{ such that} \\ b \not \in \GL_n(\mc{O}_F \otimes_{\Z} \Z_{(\ell)})}} \Q \cdot \log \ell & \text{if $T^{\flat} \not > 0$.}
    \end{cases}
    \end{align}
The notation $T^{\flat} > 0$ (resp. $T^{\flat} \not > 0$) means that $T^{\flat}$ is positive definite (resp. not positive definite).
We have $\mrm{Int}_{\infty}(T^{\flat}, a^{\flat}_{\infty} {}^t \overline{a}^{\flat}_{\infty}) = \mrm{Int}_{\infty}({}^t \overline{a}^{\flat}_{\infty} T^{\flat} a^{\flat}_{\infty}, 1)$ \crefext{II:equation:Arch_uniformization:Archimedean_nonsingular_current}.
By our main Archimedean local identity \crefext{II:theorem:local_identities:Archimedean_identity:statement:main_Archimedean}, we have $\mrm{Int}_{\infty}({}^t \overline{a}^{\flat}_{\infty} T^{\flat} a^{\flat}_{\infty}, 1) = \frac{d}{ds} \big|_{s = -1/2} W^*_{{}^t \overline{a}^{\flat}_{\infty} T a^{\flat}_{\infty},\infty}(s)^{\circ}_n$. 

The Whittaker function invariance property \eqref{equation:local_Whittaker:Archimedean:linear_invariance} implies $W^*_{{}^t \overline{a}^{\flat}_{\infty} T^{\flat} a^{\flat}_{\infty}}(s)^{\circ}_n = \tilde{W}^*_{T^{\flat}}(a^{\flat}_{\infty}, s)^{\circ}_n$. By the Archimedean local functional equation \eqref{lemma:Eisenstein:local_functional_equations:Archimedean:scalar_weight} we have $\frac{d}{ds} \big|_{s = -1/2} \tilde{W}^*_{T^{\flat},\infty}(a^{\flat}_{\infty}, s)^{\circ}_n = - \frac{d}{ds} \big|_{s = 1/2} \tilde{W}^*_{T^{\flat},\infty}(a^{\flat}_{\infty}, s)^{\circ}_n$. This is still true when $T^{\flat}$ has signature $(n - 1 - r, r)$ for $r \geq 2$, as both sides are zero in this case (by definition for the geometric side, and by \crefext{II:equation:Archimedean_identity:statement:higher_vanishing_order} for the local Whittaker function). 
As already mentioned, recall that $\tilde{W}^*_{T^{\flat}}(a^{\flat}, 1/2)^{\circ}_n$ is $1$ if $T^{\flat}$ is positive definite, and is $0$ is $T^{\flat}$ is not positive definite \eqref{equation:Eisenstein:local_Whittaker:Archimedean:special_value}.
Now \eqref{equation:arithmetic_Siegel-Weil:main_results:main:geometric_decomp:Archimedean:local} follows from what we have just discussed.

Next, recall the global Archimedean intersection number $\mrm{Int}_{\infty, \mrm{global}}(T,y) = \int_{\mc{M}_{\C}} g_{T,y}$ (where $g_{T,y}$ is a current associated with $T$ and $y$) as in \crefext{III:equation:Arch_uniformization:Archimedean}. Recall the relation \crefext{III:equation:Arch_uniformization:Archimedean:local_global_degrees}
    \begin{equation}
    \mrm{Int}_{\infty, \mrm{global}}(T, y) = \frac{h_F}{w_F} \mrm{Int}_{\infty}(T,y) \cdot \deg \Biggl [ U(V)(\Q) \backslash \coprod_{\substack{\underline{x} \in V^n \\ (\underline{x}, \underline{x}) = T}} \mc{D}(\underline{x}_f) \Biggr ]
    \end{equation}
where $V \coloneqq L \otimes_{\mc{O}_F} F$ and $\mc{D}(\underline{x}_f)$ is a certain ``away-from-$\infty$'' local special cycle (it is a discrete set), defined in \crefext{III:ssec:Arch_uniformization:away_infty_special_cycle}. The displayed groupoid cardinality $\deg [ \cdots ]$ describes certain ``complex uniformization degrees'' \crefext{III:equation:Arch_uniformization:Archimedean}. If there exists $\underline{x} \in V^n$ with $(\underline{x}, \underline{x}) = T$, the groupoid cardinality is
    \begin{equation}
    \deg \Biggl [ U(V)(\Q) \backslash \coprod_{\substack{\underline{x} \in V^n \\ (\underline{x}, \underline{x}) = T}} \mc{D}(\underline{x}_f) \Biggr ] = \frac{2 h_F}{w_F} \prod_{\ell} \tilde{W}^*_{T^{\flat},\ell}(a^{\flat}_{\ell}, 1/2)^{\circ}_n
    \end{equation}
by local Siegel--Weil as in Lemma \ref{lemma:local_Siegel-Weil:uniformization_degree:main} (with $v_0 = \infty$ in the notation of loc. cit.,
and with $C$ as in Remark \ref{remark:geometric_Siegel-Weil:degrees:self-dual}). If there does not exist such $\underline{x}$, then the Hasse principle implies that $T^{\flat}$ has signature $(n - 1 - r, r)$ for some $r \geq 2$ (compare the proof of Proposition \ref{proposition:geometric_Siegel-Weil:degrees}). 
In this case, we have $\frac{d}{ds} \big|_{s = 1/2} \tilde{W}^*_{T^{\flat}, \infty}(a^{\flat}_{\infty}, s)^{\circ}_n = 0$ \crefext{II:equation:Archimedean_identity:statement:higher_vanishing_order}. In all cases, we thus have
    \begin{equation}\label{equation:arithmetic_Siegel-Weil:main_results:main:geometric_decomp:Archimedean:global}
    \mrm{Int}_{\infty, \mrm{global}}(T, y) = - \frac{2 h_F^2}{w_F^2} \left (\frac{d}{ds} \bigg|_{s = 1/2} \left ( |a^{\#}_{\infty}|_{\infty}^s \tilde{W}^*_{T^{\flat}, \infty}(a^{\flat}_{\infty}, s)^{\circ}_n \right ) \right ) \prod_{\ell} \tilde{W}^*_{T^{\flat},\ell}(a^{\flat}_{\ell}, 1/2)^{\circ}_n
    \end{equation}
modulo $\sum_{\ell} \Q \cdot \log \ell$ for primes $\ell$ such that $b \not \in \GL_n(\mc{O}_F \otimes_{\Z} \Z_{(\ell)})$. This gives a geometric interpretation of \eqref{equation:arithmetic_Siegel-Weil:main_results:main:Eisenstein_decomp:Archimedean}.

(\textit{Step 3c: Local geometric interpretation: at $p$}) Recall $\mrm{Int}_p(T) \coloneqq \mrm{Int}_{\ms{H},p}(T) + \mrm{Int}_{\ms{V},p}(T)$ \crefext{III:equation:non-Arch_uniformization:horizontal:total_intersection}, where $\mrm{Int}_{\ms{H},p}(T)$ is a ``horizontal local intersection number'' \crefext{III:equation:non-Arch_uniformization:horizontal:Int_horizontal_local} and $\mrm{Int}_{\ms{V},p}(T)$ is a ``vertical local intersection number'' \crefext{III:equation:non-Arch_uniformization:horizontal:Int_vertical_local} associated with $T$. The former describes ``local change of tautological (or Faltings) height'' and the latter describes degrees for ``components in positive characteristic'' in terms of local special cycles on Rapoport--Zink spaces.

We claim that
    \begin{equation}\label{equation:arithmetic_Siegel-Weil:main_results:main:geometric_decomp:at_p:local}
    \mrm{Int}_p(T) = - e_p \frac{d}{d s} \bigg|_{s = 1/2} \tilde{W}^*_{T^{\flat},p}(a^{\flat}_p, s)^{\circ}_n
    \end{equation}
where $e_p = 1$ if $p$ is unramified (resp. $e_p = 2$ if $p$ is ramified).

First note that the functional equation \eqref{equation:local_functional_equations:non-Arch:Whittaker_star} implies $- \frac{d}{d s} \big|_{s = 1/2} \tilde{W}^*_{T^{\flat},p}(a^{\flat}_p, s)^{\circ}_n = \frac{d}{d s} \big|_{s = -1/2} \tilde{W}^*_{T^{\flat},p}(a^{\flat}_p, s)^{\circ}_n$. The invariance property for Whittaker functions \eqref{equation:local_Whittaker:non-Arch:linear_invariance} implies $\tilde{W}^*_{T^{\flat},p}(a^{\flat}_p, s)^{\circ}_n = \tilde{W}^*_{{}^t \overline{a}^{\flat}_p T^{\flat} a^{\flat}_p, p}(s)^{\circ}_n$.

Form the positive definite $F / \Q$ Hermitian spaces $\mbf{W} \subseteq \mbf{V}$ as in \crefext{III:sec:non-Arch_uniformization} (recall $\varepsilon(\mbf{V}_p) = - 1$ and $\varepsilon(\mbf{V}_{\ell}) = \varepsilon(V_{\ell})$ for all $\ell \neq p$). Set $\mc{O}_{F,p} \coloneqq \mc{O}_F \otimes_{\Z} \Z_p$. For any $\underline{\mbf{x}}_p \in \mbf{W}_p^n$ with Gram matrix $T$ (such $\underline{\mbf{x}}_p$ exists because $\rank(T) \leq n - 1$; recall $\mbf{W}$ has rank $n$ if $p$ is nonsplit and rank $n - 1$ if $p$ is split), there exists a basis of $L^{\flat}_p \coloneqq \mrm{span}_{\mc{O}_{F,p}}(\underline{x}_p)$ with Gram matrix ${}^t \overline{a}^{\flat}_p T^{\flat} a^{\flat}_p$. Indeed, we have $a_p \in \GL_n(\mc{O}_{F_p})$ and $a^{\flat}_p \in \GL_{n - 1}(\mc{O}_{F_p})$ by construction (and recall ${}^t \overline{a}_p^{-1} T a_p^{-1} = \mrm{diag}(0, T^{\flat})$ by definition). 
We remind the reader that \eqref{equation:local_Whittaker:local_densities_spherical:Den_star} may be used to pass between (normalized) local densities and local Whittaker functions. We also pass between the notation $\mrm{Den}^*(X,L^{\flat}_p)_n = \mrm{Den}^*(X,{}^t \overline{a}^{\flat}_p T^{\flat} a^{\flat}_p)_n$ as explained in Section \ref{ssec:Eisenstein:local_Whittaker:local_densities_spherical}. 
Now \eqref{equation:arithmetic_Siegel-Weil:main_results:main:geometric_decomp:at_p:local} follows from our main non-Archimedean local identity \crefext{I:theorem:non-Arch_identity:statement}.

Next, recall the horizontal and vertical global intersection numbers $\mrm{Int}_{\ms{H},p,\mrm{global}}(T)$ and $\mrm{Int}_{\ms{V},p,\mrm{global}}(T)$ at $p$, associated with $T$ (see \crefext{III:equation:non-Arch_uniformization:horizontal:local_global_degrees} and \crefext{III:equation:non-Arch_uniformization:vertical:local_global_degrees}). These are elements of $\Q \cdot \log p$. Recall the $F / \Q$ Hermitian space $\mbf{W}^{\perp}$ defined in \crefext{III:ssec:non-Arch_uniformization:framing}, which satisfies $\mbf{V} = \mbf{W} \oplus \mbf{W}^{\perp}$ (orthogonal direct sum). In particular, $\mbf{W}^{\perp} = 0$ if $p$ is nonsplit and $\dim_F \mbf{W}^{\perp} = 1$ if $p$ is split.

By \crefext{III:equation:non-Arch_uniformization:horizontal:local_global_degrees} and \crefext{III:equation:non-Arch_uniformization:vertical:local_global_degrees} (and in the notation of loc. cit.), we have
    \begin{equation}
    \mrm{Int}_{p, \mrm{global}}(T) = \frac{h_F}{w_F} \mrm{Int}_p(T) \cdot \deg \Biggl [ I_1(\Q) \backslash \Biggl ( \coprod_{\substack{\underline{\mbf{x}} \in \mbf{W}^n \\ (\underline{\mbf{x}}, \underline{\mbf{x}}) = T}} U(\mbf{W}^{\perp}_p) / K_{1, \mbf{L}^{\perp}_p} \times \mc{Z}(\underline{\mbf{x}}^p) \Biggr ) \Biggr ].
    \end{equation}
The notation $\mc{Z}(\underline{\mbf{x}}^p)$ means a certain ``away-from-$p$'' local special cycle (a discrete set), defined in \crefext{III:ssec:non-Arch_uniformization:away_p_special_cycle}. Recall that $K_{1,\mbf{L}^{\perp}_p} \subseteq U(\mbf{W}^{\perp}_p)$ is the unique maximal open compact subgroup and $I_1 = U(\mbf{W}) \times U(\mbf{W}^{\perp})$ as algebraic groups over $\Q$ \crefext{III:ssec:non-Arch_uniformization:quotient}. The displayed groupoid cardinality $\deg [ \cdots ]$ encodes certain ``Rapoport--Zink non-Archimedean uniformization degrees''.

If there exists $\underline{\mbf{x}} \in \mbf{W}^n$ with Gram matrix $T$, then local Siegel--Weil (Lemma \ref{lemma:local_Siegel-Weil:uniformization_degree:main}) implies
    \begin{equation}\label{equation:arithmetic_Siegel-Weil:main_results:main:geometric_decomp:at_p:uniformization_degrees}
    \deg \Biggl [ I_1(\Q) \backslash \Biggl ( \coprod_{\substack{\underline{\mbf{x}} \in \mbf{W}^n \\ (\underline{\mbf{x}}, \underline{\mbf{x}}) = T}} U(\mbf{W}^{\perp}_p) / K_{1, \mbf{L}^{\perp}_p} \times \mc{Z}(\underline{\mbf{x}}^p) \Biggr ) \Biggr ] = \frac{2 h_F}{e_p w_F} \prod_{v \neq p} \tilde{W}^*_{T^{\flat},v}(a^{\flat}_v, 1/2)^{\circ}_n.
    \end{equation}
(in the notation of Lemma \ref{lemma:local_Siegel-Weil:uniformization_degree:main}, take $v_0 = p$ and use the Hermitian space $\mbf{V}$ for the $V$ in loc. cit.).

Set $\Omega_T(R) \coloneqq \{ \underline{\mbf{x}} \in (\mbf{W} \otimes_{\Q} R)^n : (\underline{\mbf{x}}, \underline{\mbf{x}}) = T \}$ for $\Q$-algebras $R$. If $\Omega_T(\Q) = \emptyset$, then the Hasse principle implies $\Omega_T(\Q_v) = \emptyset$ for some place $v$ of $\Q$. We have $\Omega_T(\Q_p) \neq \emptyset$ (either $p$ is nonsplit and $\mbf{W} = \mbf{V}$ and the claim follows because $\rank T < \rank \mbf{W}$ (compare the proof of Proposition \ref{proposition:geometric_Siegel-Weil:degrees}),
or $p$ is split and $\Omega_T(\Q_p) \neq \emptyset$ automatically). For all places $v$, we have $\Omega_T(\Q_v) = \emptyset$ if and only if $\Omega_{{}^t \overline{a}_v^{\flat} T^{\flat} a_v^{\flat}}(\Q_v) = \emptyset$ (where $\Omega_{{}^t \overline{a}_v^{\flat} T^{\flat} a_v^{\flat}}$ is defined like $\Omega_T$ but for $(n - 1)$-tuples); this follows from our diagonalization of $T$ (e.g. ${}^t \overline{a}_v^{-1} T a_v^{-1} = \mrm{diag}(0,T^{\flat})$ for all $v < \infty$).

If $\Omega_T(\Q_v) = \emptyset$, we thus conclude $\tilde{W}^*_{T^{\flat},v}(a^{\flat}_v, 1/2)^{\circ}_n = \tilde{W}^*_{{}^t \overline{a}^{\flat}_v T^{\flat} a^{\flat}_v, v}(1/2)^{\circ}_n = 0$ by the invariance property for local Whittaker functions (see \eqref{equation:local_Whittaker:Archimedean:linear_invariance} and \eqref{equation:local_Whittaker:non-Arch:linear_invariance}) and by local Siegel--Weil \eqref{equation:local_Siegel-Weil:local_Siegel-Weil:up-to-constant}.
Hence \eqref{equation:arithmetic_Siegel-Weil:main_results:main:geometric_decomp:at_p:uniformization_degrees} holds even if there is no $\underline{\mbf{x}} \in \mbf{W}^n$ such that $(\underline{\mbf{x}}, \underline{\mbf{x}}) = T$ (both sides are $0$ in this case).

We have shown
    \begin{equation}\label{equation:arithmetic_Siegel-Weil:main_results:main:geometric_decomp:at_p:global}
    \mrm{Int}_{p, \mrm{global}}(T) = - \frac{2 h^2_F}{w^2_F} \left ( \frac{d}{ds} \bigg|_{s = 1/2} \tilde{W}^*_{T^{\flat}, p}(a^{\flat}_{p}, s)^{\circ}_n \right ) \prod_{v \neq p} \tilde{W}^*_{T^{\flat},v}(a^{\flat}_v, 1/2)^{\circ}_n.
    \end{equation}
This gives a geometric interpretation for \eqref{equation:arithmetic_Siegel-Weil:main_results:main:Eisenstein_decomp:at_p}.

(\textit{Step 4: Finish}) Recall the definition of arithmetic degree without boundary contributions $\widehat{\deg}([\widehat{\mc{Z}}(T)])$ \eqref{equation:part_IV:intro:results:if_proper}. In our current situation, this is
    \begin{align}
    \widehat{\deg}([\widehat{\mc{Z}}(T)]) \coloneqq \left ( \int_{\mc{M}_{\C}} g_{T,y} \right )+ \widehat{\deg}(\widehat{\mc{E}}^{\vee}|_{\mc{Z}(T)_{\ms{H}}}) + \sum_{\ell} \deg_{\F_{\ell}}({}^{\mathbb{L}}\mc{Z}(T)_{\ms{V},\ell}) \log \ell \notag
    \end{align}
where the sum runs over all primes $\ell$. By definition, we have
    \begin{equation}\label{equation:arithmetic_Siegel-Weil:main_results:finish_decomp}
    \begin{gathered}
    \int_{\mc{M}_{\C}} g_{T,y} = \mrm{Int}_{\infty, \mrm{global}}(T, y) \quad \quad \deg_{\F_{\ell}}({}^{\mathbb{L}}\mc{Z}(T)_{\ms{V},\ell}) \log \ell = \mrm{Int}_{\ms{V}, \ell, \mrm{global}}(T)
    \\
    \widehat{\deg}(\widehat{\mc{E}}^{\vee}|_{\mc{Z}(T)_{\ms{H}}}) = (\deg_{\Z} \mc{Z}(T)_{\ms{H}}) \cdot h^{\mrm{CM}}_{\widehat{\mc{E}}^{\vee}} + \sum_{\ell} \mrm{Int}_{\ms{H}, \ell, \mrm{global}}(T)
    \end{gathered}
    \end{equation}
where $h^{\mrm{CM}}_{\widehat{\mc{E}}^{\vee}}$ is the height constant from \eqref{equation:arithmetic_Siegel-Weil:main_results:height_constant}. See \crefext{III:equation:Arch_uniformization:Archimedean} (Archimedean), \crefext{III:equation:non-Arch_uniformization:vertical:local_global_degrees} (vertical), and \crefext{III:equation:non-Arch_uniformization:horizontal:local_global_degrees:decomp} (horizontal). For all primes $\ell$, we have $\mrm{Int}_{\ms{V},\ell,\mrm{global}}(T) \in \Q \cdot \log \ell$ and $\mrm{Int}_{\ms{H}, \ell, \mrm{global}}(T) \in \Q \cdot \log \ell$. These quantities are $0$ for all but finitely many $\ell$.

After multiplying both sides of \eqref{equation:arithmetic_Siegel-Weil:main_results:main:Eisenstein_decomp:total} by $2 (h_F / w_F)^2$, we apply the results of Steps 3a, 3b, and 3c above (see \eqref{equation:arithmetic_Siegel-Weil:main_results:main:geometric_decomp:height_constant:global}, \eqref{equation:arithmetic_Siegel-Weil:main_results:main:geometric_decomp:Archimedean:global}, and \eqref{equation:arithmetic_Siegel-Weil:main_results:main:geometric_decomp:at_p:global}) to find
    \begin{equation}
    \frac{h_F}{w_F} \frac{d}{ds} \bigg|_{s = 0} E^*_T(y,s)^{\circ}_n = \widehat{\deg}([\widehat{\mc{Z}}(T)])
    \end{equation}
as elements of $\R/(\sum_{\ell \neq p} \Q \cdot \log \ell)$. As we already discussed, varying $p$ shows that this identity holds as an equality of real numbers.
\end{proof}

\begin{remark}[Nonsingular central-point arithmetic Siegel--Weil]\label{remark:arithmetic_Siegel-Weil:main_results:nonsingular}
In the setup above (in particular, $n \equiv 2 \pmod{4}$), consider any $T \in \mrm{Herm}_n(\Q)$ with $\det T \neq 0$ and any $y \in \mrm{Herm}_n(\R)_{>0}$. Assuming the prime $2$ is split in $\mc{O}_F$, we still have
    \begin{equation}\label{equation:arithmetic_Siegel-Weil:main_results:nonsingular}
    \widehat{\deg}([\widehat{\mc{Z}}(T)]) = \frac{h_F}{w_F} \frac{d}{d s} \bigg |_{s = 0} E^*_{T}(y,s)^{\circ}_n
    \end{equation}
where the Green current for $[\widehat{\mc{Z}}(T)]$ is formed with respect to $y$, and where $\widehat{\deg}([\widehat{\mc{Z}}(T)])$ again denotes the arithmetic degree without boundary contributions as in \eqref{equation:part_IV:intro:results:if_proper}. This should be compared with our preceding main theorem for singular $T$ of co-rank $1$ (Theorem \ref{theorem:arithmetic_Siegel-Weil:main_results:main}).

Using the local theorems of Liu, Li--Zhang, and Li--Liu (cited below), one can prove \eqref{equation:arithmetic_Siegel-Weil:main_results:nonsingular} by a local decomposition as in the proof of Theorem \ref{theorem:arithmetic_Siegel-Weil:main_results:main} (no diagonalization procedure is necessary here) using the volume constant calculated in Lemma \ref{proposition:geometric_Siegel-Weil:degrees}. This is possibly considered known to experts up to a volume constant by the cited local theorems. Nevertheless, the global statement is not available in the literature
(but see some closely related variants at the end of this remark), so we have stated it. A sketch is provided below.

Decomposing $E^*_T(y,s)^{\circ}_n$ into a product of local Whittaker functions (Section \ref{ssec:Eisenstein:global_normalized_Fourier:global_normalization}), we find
    \begin{align}
    \frac{d}{ds} \bigg|_{s = 0} E^*_{T}(y,s)^{\circ}_n & = \left ( \frac{d}{ds} \bigg|_{s = 0} W^*_{T, \infty}(y, s)^{\circ}_n \right ) \prod_{\ell} W^*_{T,\ell}(0)^{\circ}_n 
    \\
    & \mathrel{\phantom{=}} + \sum_p \left ( \frac{d}{ds} \bigg|_{s = 0} W^*_{T, p}(s)^{\circ}_n \right ) W^*_{T,\infty}(y,0)^{\circ}_n \prod_{\ell \neq p} W^*_{T,\ell}(0)^{\circ}_n
    \end{align}
At most one of the summands is nonzero (see below), and all but finitely many $W^*_{T,\ell}(s)^{\circ}_n$ are identically equal to $1$ as functions of $s$.

For the quantities $\mrm{Int}_{\infty, \mrm{global}}(T, y)$
from \crefext{III:equation:Arch_uniformization:Archimedean}
(a Green current integral)
and $\mrm{Int}_{p, \mrm{global}}(T)$
from the end of \crefext{III:ssec:non-Arch_uniformization:vertical}
(degrees of cycle classes supported in characteristic $p$),
we also have
    \begin{align}
    \widehat{\deg}([\widehat{\mc{Z}}(T)]) & = \mrm{Int}_{\infty, \mrm{global}}(T, y) + \sum_p \mrm{Int}_{p, \mrm{global}}(T),
    \end{align}
largely by definition.
In contrast with our main theorem, these intersection numbers $\mrm{Int}_{p, \mrm{global}}(T)$ are ``purely vertical'', without a mixed characteristic contribution.

In this setup, the local Archimedean theorem \cite[{Theorem 4.1.7}]{Liu11} (restated in our notation in \crefext{II:theorem:local_identities:Archimedean_identity:statement:main_Archimedean}) and the local Kudla--Rapoport theorems \cite[{Theorem 1.2.1}]{LZ22unitary} (inert) and \cite[{Theorem 2.7}]{LL22II} (ramified, exotic smooth, even $n$) take the place of our main local identities (which were for co-rank $1$ singular $T$).
These local theorems may be reformulated, in our notation, as
    \begin{equation}\label{equation:arithmetic_Siegel_Weil:main_results:central_point:ASW}
    \mrm{Int}_{\infty}(T, y) = \frac{d}{ds} \bigg|_{s = 0} W^*_{T, \infty}(y, s)^{\circ}_n
    \quad\quad
    \mrm{Int}_{p}(T) = \frac{d}{ds} \bigg|_{s = 0} W^*_{T, p}(s)^{\circ}_n
    \end{equation}
for nonsplit $p$, where the left-hand sides are certain ``local intersection numbers'',
with $\mrm{Int}_{\infty}(T, y)$ calculated on the symmetric space $\mc{D}$
\crefext{III:equation:Arch_uniformization:Archimedean:Int_infty_local},
and $\mrm{Int}_{p}(T)$ calculated on suitable Rapoport--Zink spaces
\crefext{III:equation:non-Arch_uniformization:horizontal:total_intersection}.

Let $V$ be the $F / \Q$ Hermitian space of signature $(n - 1, 1)$
which is split at every prime, meaning $\varepsilon_p(V_p) = +1$
for all primes $p$, with conventions as in \crefext{I:ssec:Hermitian_conventions:lattices}
(since $V$ has even rank, being split at every prime is equivalent to
$V$ containing a self-dual lattice of full-rank).
Let $L \subseteq V$ be a full-rank self-dual lattice.
Given a nonsplit prime $p$, let $\mbf{V}$ be the ``nearby Hermitian space"
with signature $(n,0)$, with $\varepsilon_p(V_p) = -1$,
and with $\varepsilon_{\ell}(V_{\ell}) = +1$ for all $\ell \neq p$.
The space $\mbf{V}$ certainly depends on $p$,
though we suppress this from notation.

By complex and Rapoport--Zink uniformization respectively,
the global intersection numbers unfold into the local quantities
    \begin{align}
    \mrm{Int}_{\infty, \mrm{global}}(T, y)
    & =
    \mrm{Int}_{\infty}(T, y)
    \cdot
    \frac{h_F}{w_F}
    \cdot
    \#
    \bigg [
    U(V)(\Q) \backslash \bigg ( \coprod_{\substack{\underline{x} \in V^n \\ (\underline{x},\underline{x}) = T}}
    \mc{Z}(\underline{x}^{\infty}) \bigg )
    \bigg ]
    \\
    \mrm{Int}_{p, \mrm{global}}(T)
    & =
    \mrm{Int}_{p}(T)
    \cdot
    \frac{h_F}{w_F}
    \cdot
    \#
    \bigg [
    U(\mbf{V})(\Q) \backslash \bigg ( \coprod_{\substack{\underline{x} \in \mbf{V}^n \\ (\underline{x},\underline{x}) = T}}
    \mc{Z}(\underline{x}^{p}) \bigg )
    \bigg ]
    \end{align}
for nonsplit $p$,
as explained in \cite[{\eqref{III:equation:Arch_uniformization:Archimedean:local_global_degrees}, (6.8.7)}]{corank1_ASW_III.pdf},
where the notation on the right is that of
Lemma \ref{lemma:local_Siegel-Weil:uniformization_degree:main},
with $\underline{x}_f^{v_0}$ of loc. cit. replaced by $\underline{x}^{\infty}$
and $\underline{x}^p$.

By Lemma \ref{lemma:local_Siegel-Weil:uniformization_degree:main}(2)
(we have $C = 1$ in the notation of loc. cit.)
we have
    \begin{align}
    \label{equation:arithmetic_Siegel_Weil:main_results:central_point:uniformization_degree:Arch}
    \#
    \bigg [
    U(V)(\Q) \backslash \bigg ( \coprod_{\substack{\underline{x} \in V^n \\ (\underline{x},\underline{x}) = T}}
    \mc{Z}(\underline{x}^{\infty}) \bigg )
    \bigg ]
    &= \prod_p W^*_{T,p}(0)^{\circ}_n
    \\
    \label{equation:arithmetic_Siegel_Weil:main_results:central_point:uniformization_degree:non-Arch}
    \#
    \bigg [
    U(\mbf{V})(\Q) \backslash \bigg ( \coprod_{\substack{\underline{x} \in \mbf{V}^n \\ (\underline{x},\underline{x}) = T}}
    \mc{Z}(\underline{x}^{p}) \bigg )
    \bigg ]
    &= W^*_{T,\infty}(y,0)^{\circ}_n \prod_{\ell \neq p} W^*_{T,\ell}(0)^{\circ}_n
    \end{align}
for all nonsplit $p$,
where the first formula is conditional on existence of $\underline{x} \in V^n$
with $(\underline{x}, \underline{x}) = T$,
and the second formula, for each $p$, is conditional on existence
of $\underline{x} \in \mbf{V}^n$ with $(\underline{x}, \underline{x}) = T$.
The left-hand sides vanish, by definition, if no such $\underline{x}$ exists.
To handle existence/non-existence of such $\underline{x}$,
we have in mind a (presumably routine) Hasse principle argument (compare \cite[{\S 9}]{KR14})
as we explain next.
For any prime $p$, set $\varepsilon_p(T) \coloneqq \eta_p((-1)^{n(n-1)/2} \det T)$ (the usual local invariant, by our conventions in \crefext{I:ssec:Hermitian_conventions:lattices}), where $\eta_p \colon \Q_p^{\times} \ra \{ \pm 1\}$ is the local quadratic character associated to $F / \Q$. 

We have $\mrm{Int}_{\infty, \mrm{global}}(T, y) = 0$ unless $T$ has signature $(n - 1, 1)$ and $\varepsilon_p(T) = 1$ for all $p$,
by construction (the Green current is zero otherwise).
In particular, we have $\mrm{Int}_{\infty, \mrm{global}}(T, y) = 0$ unless
there exists $\underline{x} \in V$ with $(\underline{x}, \underline{x}) = T$.
For such $T$, the special cycle $\mc{Z}(T)$ is empty (but may have a nontrivial Green current).

We have $\mrm{Int}_{p, \mrm{global}}(T) = 0$ unless $T$ is positive definite, $\varepsilon_p(T) = -1$, and $\varepsilon_{\ell}(T) = 1$ for all primes $\ell \neq p$. For such $T$, the special cycle $\mc{Z}(T)$ is supported in characteristic $p$ (or empty). For all other $T$, the special cycle $\mc{Z}(T)$ is empty.
These claims follow from e.g. uniformization of special cycles (e.g. Sections \crefext{III:ssec:Arch_uniformization:Archimedean} (Archimedean) and \crefext{III:ssec:non-Arch_uniformization:vertical} (non-Archimedean)) and the Hasse principle (e.g. applied to $\mbf{V}$ from loc. cit. in the non-Archimedean case).
In particular, $\mrm{Int}_{p, \mrm{global}}(T) = 0$ if $p$ is split in $\mc{O}_F$, and $\mc{Z}(T)$ is empty over any split $p$.
For nonsplit primes $p$, we also find that $\mrm{Int}_{p, \mrm{global}}(T) = 0$ unless
there exists $\underline{x} \in \mbf{V}$ with $(\underline{x}, \underline{x}) = T$
(recall that $\mbf{V}$ depends on $p$).

On the analytic side, we have $W^*_{T,p}(0)^{\circ}_n = 0$ if $\varepsilon_p(T) = - 1$ (by local Siegel--Weil \eqref{equation:local_Siegel-Weil:local_Siegel-Weil:up-to-constant}, or the functional equation \eqref{equation:local_functional_equations:non-Arch:Whittaker_star}).
If $T$ has signature $(n - r, r)$ for $r \geq 2$, we have $\frac{d}{ds} \big|_{s = 0} W^*_{T,\infty}(y,s)^{\circ}_n = 0$ \crefext{II:equation:Archimedean_identity:statement:higher_vanishing_order}.
In particular, if there is no $\underline{x} \in V^n$ with $(\underline{x}, \underline{x}) = T$,
then either the right-hand side of \eqref{equation:arithmetic_Siegel_Weil:main_results:central_point:uniformization_degree:Arch}
is $0$, or the derivative in the left equation of \eqref{equation:arithmetic_Siegel_Weil:main_results:central_point:ASW}
is $0$.

We also have $W^*_{T,\infty}(y,0)^{\circ}_n = 0$ if $T$ is not positive definite (local Siegel--Weil again, or \eqref{equation:Eisenstein:local_Whittaker:Archimedean:special_value}).
Thus, if there is no $\underline{x} \in \mbf{V}^n$ with $(\underline{x}, \underline{x}) = T$
(for some understood nonsplit $p$),
the right-hand side of \eqref{equation:arithmetic_Siegel_Weil:main_results:central_point:uniformization_degree:non-Arch}
is also $0$.

We then find
    \begin{align}
    \mrm{Int}_{\infty, \mrm{global}}(T, y) &= \frac{h_F}{w_F} \left ( \frac{d}{ds} \bigg|_{s = 0} W^*_{T, \infty}(y, s)^{\circ}_n \right ) \prod_p W^*_{T,p}(0)^{\circ}_n
    \\
    \mrm{Int}_{p, \mrm{global}}(T) &= \frac{h_F}{w_F} \left ( \frac{d}{ds} \bigg|_{s = 0} W^*_{T, p}(s)^{\circ}_n \right ) W^*_{T,\infty}(y,0)^{\circ}_n \prod_{\ell \neq p} W^*_{T,\ell}(0)^{\circ}_n
    \end{align}
for all nonsplit $p$.
The displayed formula also holds for split $p$:
in this case, we have $\mrm{Int}_{p,\mrm{global}}(T) = 0$
as explained above, and
we have $W^*_{T,\infty}(y,0)^{\circ}_n = 0$
if $T$ is not positive definite,
and $W^*_{T,\ell}(0)^{\circ}_n = 0$
if $\varepsilon_{\ell}(T) = -1$ (at least one of these must hold,
i.e. we cannot have $T$ positive definite
and $\varepsilon_{\ell}(T) = +1$ for all $\ell$,
by the global product formula for local invariants).
This completes the proof of \eqref{equation:arithmetic_Siegel-Weil:main_results:nonsingular}.

For the analogous global result (still $\det T \neq 0$ and $T \in \mrm{Herm}_n$, central derivative) for an unramified CM extension of number fields $F / F^+$ where all $2$-adic places are split (forcing $F^+ \neq \Q$) and a lattice $L$ which is self-dual for the Hermitian pairing, see \cite[{Theorem 15.5.1}]{LZ22unitary} (at least up to a volume constant). For the analogous global result (still $\det T \neq 0$ and $T \in \mrm{Herm}_n$, central derivative) for possibly ramified $F / F^+$ where all $2$-adic places are split, on Kr\"amer integral models (semistable reduction at ramified primes), and again $L$ self-dual for the Hermitian pairing, see \cite[{Theorem 10.1}]{HLSY22} (at least up to a volume constant). For the result on Kr\"amer models, one needs to correct the Eisenstein series derivative by special values of other Eisenstein series.
\end{remark}

\begin{remark}\label{remark:arithmetic_Siegel-Weil:main_results:0_mod_4}
When $n \equiv 0 \pmod{4}$, there is no non-degenerate self-dual signature $(n - 1, 1)$ Hermitian $\mc{O}_F$-lattice. In this case, Theorem \ref{theorem:arithmetic_Siegel-Weil:main_results:main}(1) still holds in the sense that $\frac{d}{d s} \big |_{s = 0} E^*_{T}(y,s)^{\circ}_n = 0$ (by the functional equation, Lemma \ref{lemma:Eisenstein:global_normalized_Fourier:global_normalization:global_functional_equation}).
\end{remark}

\begin{remark}\label{remark:arithmetic_Siegel-Weil:main_results:Faltings_height}
We explain how Theorem \ref{theorem:arithmetic_Siegel-Weil:main_results:main} may be reformulated in terms of Faltings heights. Assume $2$ is split in $\mc{O}_F$. Let $\widehat{\omega}$ be the metrized Hodge bundle on $\mc{M}$ as defined in \crefext{III:ssec:ab_var:integral_models} (also \crefext{I:ssec:part_I:arith_intersections:metrized_taut_bundle}). Take $T \in \mrm{Herm}_n(\Q)$ with $\rank(T) = n - 1$. By \crefext{III:equation:non-Arch_uniformization:horizontal:local_global_degrees:Faltings}, we have
    \begin{equation}
    \widehat{\deg}(\widehat{\omega} |_{\mc{Z}(T)_{\ms{H}}}) = (\deg_{\Z} \mc{Z}(T)_{\ms{H}}) \cdot n \cdot h_{\mrm{Fal}}^{\mrm{CM}} - 2 \sum_{p} \mrm{Int}_{\ms{H}, p, \mrm{global}}(T)
    \end{equation}
where $h_{\mrm{Fal}}^{\mrm{CM}}$ is the Faltings height of any elliptic curve with CM by $\mc{O}_F$ (as in \crefext{III:equation:part_II:arith_cycle_classes:Hodge_bundles:CM_Faltings_height}).
By definition of Faltings height, we have
    \begin{equation}\label{equation:arithmetic_Siegel-Weil:main_results:Faltings_height:bundle_degree}
    \widehat{\deg}(\widehat{\omega} |_{\mc{Z}(T)_{\ms{H}}}) = 2 \sum_{\a' \in \mc{Z}(T)(\C)} |\Aut(\a')|^{-1} h_{\mrm{Fal}}(A)
    \end{equation}
where $\a' = (A_0, \iota_0, \lambda_0, A, \iota, \lambda) \in \mc{Z}(T)(\C)$ (choose $F \ra \C$), and where $h_{\mrm{Fal}}(A)$ is the Faltings height of $A$ (as in \crefext{I:ssec:Faltings_and_taut:heights}) after descent to any number field, with metric normalized as in \crefext{I:equation:arith_cycle_classes:Hodge_bundles:Faltings_metric}. Alternatively, we could consider morphisms $\Spec \C \ra \mc{M}$ over $\Spec \Z$, which would remove the factor of $2$ in the previous formula.

Our main theorem (Theorem \ref{theorem:arithmetic_Siegel-Weil:main_results:main}) admits the equivalent formulation
    \begin{align}
    & \mrm{Int}_{\infty, \mrm{global}}(T, y) - \frac{1}{2} \widehat{\deg}(\widehat{\omega}|_{\mc{Z}(T)_{\ms{H}}}) + (\deg_{\Z} \mc{Z}(T)_{\ms{H}}) \cdot (h^{\mrm{CM}}_{\widehat{\mc{E}}^{\vee}} + \frac{n}{2} \cdot h^{\mrm{CM}}_{\mrm{Fal}}) + \sum_{p} \mrm{Int}_{\ms{V}, p, \mrm{global}}(T) \notag
    \\
    & = \frac{h_F}{w_F} \frac{d}{ds} \bigg|_{s = 0} E^*_T(y,s)^{\circ}_n
    \end{align}
via the decomposition in \eqref{equation:arithmetic_Siegel-Weil:main_results:finish_decomp}.
We remind the reader that $\deg_{\Z} \mc{Z}(T)_{\ms{H}}$ is essentially a special value of a $U(n - 1, n - 1)$ Eisenstein series \eqref{equation:arithmetic_Siegel-Weil:main_results:main:geometric_decomp:height_constant:global}.
For further discussion of the special case $n = 2$, see Section \ref{ssec:arithmetic_Siegel-Weil:Serre_tensor}.
\end{remark}

In the rest of Section \ref{ssec:arithmetic_Siegel-Weil:main_results}, we discuss some results which are applicable even if $L$ is not self-dual.
Allow possibly $2 \mid \Delta$, and let $L$ be any non-degenerate Hermitian $\mc{O}_F$-lattice of signature $(n - 1, 1)$ (with $n$ not necessarily even). Select any character $\chi \colon F^{\times} \backslash \A_F^{\times} \ra \C^{\times}$ such that $\chi|_{\A^{\times}} = \eta^n$, where $\eta$ is the quadratic character associated with $F / \Q$. Set $V = L \otimes_{\mc{O}_F} F$, with associated local Hermitian space $V_v$ for each place $v$ of $\Q$. Suppose $m^{\flat} \geq 0$ is an integer. For each prime $p$, let $\varphi_p^{\flat} = \pmb{1}_{L_p}^{m^{\flat}} \in \mc{S}(V_p^{m^{\flat}})$, form the local Siegel--Weil standard section $\Phi_{\varphi_v^{\flat}} \in I(\chi_v, s)$, and set
    \begin{equation}
    \Phi_L \coloneqq \Phi^{(n)}_{\infty} \bigotimes_p \Phi_{\varphi^{\flat}_p} \in I(\chi,s)
    \end{equation}
where the Archimedean component $\Phi^{(n)}_{\infty}$ is the standard (normalized) scalar weight section from Section \ref{ssec:Eisenstein:setup:adelic_v_classical}. Form the associated classical $U(m^{\flat}, m^{\flat})$ Eisenstein series $E(z^{\flat}, s, \Phi_L)_n$ for $z^{\flat} \in \mc{H}_{m^{\flat}}$, and consider the normalized Eisenstein series Fourier coefficients
    \begin{equation}
    E^*_{T^{\flat}}(y^{\flat}, s, \Phi_L)_n \coloneqq \left ( \prod_{p} \gamma_{\psi_p}(V_p)^{m^{\flat}} \mrm{vol}(L_p)^{-m^{\flat}} \right ) \Lambda_{m^{\flat}}(s)^{\circ}_n E_{T^{\flat}}(y^{\flat}, s, \Phi_L)_n
    \end{equation}
for $T^{\flat} \in \mrm{Herm}_{m^{\flat}}(\Q)$. We are not sure whether this is a ``good'' normalization if $L$ is not self-dual, so the preceding notation appears nowhere else in this work. As in Section \ref{ssec:Eisenstein:local_Whittaker:non-Arch}, $\gamma_{\psi_p}(V_p)$ is a Weil index and $\mrm{vol}(L_p)$ is the volume of $L_p$ with respect to a certain self-dual Haar measure on $V_p$ (these factors are $1$ for all but finitely many $p$).

Form the moduli stack $\mc{M} \ra \Spec \mc{O}_F[1/d_L]$ associated with $L$ as in \crefext{III:ssec:ab_var:integral_models}.

\begin{remark}\label{remark:arithmetic_Siegel-Weil:main_results:other_levels}
Since the proof of Theorem \ref{theorem:arithmetic_Siegel-Weil:main_results:main} is local in nature, it is possible to use our local main theorems to prove variants for non self-dual $L$, up to discarding finitely many primes.
This also allows the possibility of relaxing the $n \equiv 2 \pmod{4}$ assumption in
Theorem \ref{theorem:arithmetic_Siegel-Weil:main_results:main}
(at the cost of discrepancy up to $\Q \cdot \log p$ for the primes $p$
where $L$ is not self-dual), which was imposed
to ensure the existence of an everywhere self-dual lattice. 

Set $m^{\flat} = n - 1$.
Consider $T^{\flat} \in \mrm{Herm}_{n - 1}(\Q)$ with $\det T^{\flat} \neq 0$. Let $C \in \Q_{>0}$ be the volume constant from Lemma \ref{lemma:local_Siegel-Weil:uniformization_degree:main}(3), for the Hermitian space $V$ and with $v_0 = \infty$ etc. in the notation of loc. cit..
Consider $y^{\flat} \in \mrm{Herm}_{n - 1}(\R)_{>0}$. Form $[\widehat{\mc{Z}}(T^{\flat})]$ with Green current with respect to $y^{\flat}$. Arguing as in the proof of our main theorem (Theorem \ref{theorem:arithmetic_Siegel-Weil:main_results:main}) gives
    \begin{equation}\label{equation:arithmetic_Siegel-Weil:main_results:other_levels}
    \widehat{\deg}([\widehat{\mc{Z}}(T^{\flat})] \cdot \widehat{c}_1(\widehat{\mc{E}}^{\vee})) = C \cdot \frac{d}{d s} \bigg|_{s = 0} \left ( \frac{\Lambda_n(s)_n^{\circ}}{\Lambda_{n - 1}(s + 1/2)^{\circ}_n} E^*_{T^{\flat}}(y^{\flat}, s + 1/2, \Phi_L)_n \right ) \mod \sum_{p \mid 2 d_L} \Q \cdot \log p.
    \end{equation}
For proving \eqref{equation:arithmetic_Siegel-Weil:main_results:other_levels}, the diagonalization argument (Step 1) in the proof of Theorem \ref{theorem:arithmetic_Siegel-Weil:main_results:main} can be skipped. If $2$ is split in $\mc{O}_F$, the expression ``$2 d_L$'' in \eqref{equation:arithmetic_Siegel-Weil:main_results:other_levels} may be replaced by ``$d_L$''.

In the case $n = 1$, recall that $\mc{M}$ extends smoothly (and nontrivially) over all of $\Spec \mc{O}_F$ \crefext{III:remark:ab_var:integral_models:n=1}. In this case, we need not discard any primes in \eqref{equation:arithmetic_Siegel-Weil:main_results:other_levels}. As $m^{\flat} = 0$, the normalized $U(m^{\flat}, m^{\flat})$ Eisenstein series $E^*$ is the constant function $1$ in this case.
\end{remark}

Recall that our main Archimedean local result was valid in arbitrary ``codimension'' for empty local special cycles with possibly nontrivial Green current (``purely Archimedean intersection number''). This has the following global consequence.

\begin{theorem}\label{theorem:arithmetic_Siegel-Weil:main_results:Archimedean}
Let $m^{\flat}$ be any integer with $1 \leq m^{\flat} \leq n$. Consider $T^{\flat} \in \mrm{Herm}_{m^{\flat}}(\Q)$ which is nonsingular and not positive definite. Let $C \in \Q_{>0}$ be the volume constant from Lemma \ref{lemma:local_Siegel-Weil:uniformization_degree:main}(1), for the Hermitian space $V$, the lattice $L$, and $v_0 = \infty$ in the notation of loc. cit.. 

For any $y^{\flat} \in \mrm{Herm}_{m^{\flat}}(\R)_{>0}$, we have an equality of real numbers
    \begin{equation}\label{equation:arithmetic_Siegel-Weil:main_results:Archimedean}
    \widehat{\deg}([\widehat{\mc{Z}}(T^{\flat})] \cdot \widehat{c}_1(\widehat{\mc{E}}^{\vee})^{n - m^{\flat}}) \coloneqq \int_{\mc{M}_{\C}} g_{T^{\flat},y^{\flat}} \wedge c_1(\widehat{\mc{E}}^{\vee}_{\C})^{n - m^{\flat}} = (-1)^{n - m^{\flat}} C \cdot \frac{h_F}{w_F} \frac{d}{ds} \bigg|_{s = s_0^{\flat}} E^*_{T^{\flat}}(y^{\flat}, s, \Phi_L)_n
    \end{equation}
where $s_0^{\flat} \coloneqq (n - m^{\flat}) / 2$.
\end{theorem}
\begin{proof}
In the theorem statement, we set $\mc{M}_{\C} \coloneqq \mc{M} \times_{\Spec \mc{O}_F} \Spec \C$ for either choice of embedding $F \ra \C$. Recall that the special cycle $\mc{Z}(T^{\flat})$ is empty by the non positive definite-ness \crefext{III:ssec:ab_var:integral_models}. The current $g_{T^{\flat},y^{\flat}}$ associated with $[\widehat{\mc{Z}}(T^{\flat})]$ is formed with respect to $y^{\flat}$, as usual.

Using our main Archimedean result \crefext{II:theorem:local_identities:Archimedean_identity:statement:main_Archimedean} and local Siegel--Weil (Lemma \ref{lemma:local_Siegel-Weil:uniformization_degree:main}) for uniformization degrees, the theorem follows as in the proof of Theorem \ref{theorem:arithmetic_Siegel-Weil:main_results:main}, Step (3a). Since $\det T^{\flat} \neq 0$, the proof is simpler here as the diagonalization argument of loc. cit. plays no role. Recall $W^*_{T^{\flat}, \infty}(y^{\flat},s_0^{\flat})^{\circ}_n = 0$ \eqref{equation:Eisenstein:local_Whittaker:Archimedean:special_value}, so the derivatives of non-Archimedean Whittaker functions play no role. If $T^{\flat}$ has signature $(m^{\flat} - r, r)$ for $r \geq 2$, then both sides of \eqref{equation:arithmetic_Siegel-Weil:main_results:Archimedean} are zero. The sign $(-1)^{n - m^{\flat}}$ comes from the Archimedean local functional equation (Lemma \ref{lemma:Eisenstein:local_functional_equations:Archimedean:scalar_weight}), since \crefext{II:theorem:local_identities:Archimedean_identity:statement:main_Archimedean} was stated at $s = - s^{\flat}_0$.
\end{proof}

When $m^{\flat} = n$, the preceding result is due to Liu (see \cite[{Theorem 4.17, Proof of Theorem 4.20}]{Liu11} and also \cite[{Theorem 15.3.1}]{LZ22unitary}). We do not have a new proof of this case (we deduced our local result for arbitrary $m^{\flat}$ from Liu's result using our local limiting method).
    
            \subsection{Faltings heights of Hecke translates of CM elliptic curves}
            \label{ssec:arithmetic_Siegel-Weil:Serre_tensor}
                Using the Serre tensor construction, we restate part of the simplest case ($n = 2$) of our main theorem (Theorem \ref{theorem:arithmetic_Siegel-Weil:main_results:main}) in more elementary terms, via Faltings heights of Hecke translates of CM elliptic curves (Corollary \ref{corollary:arithmetic_Siegel-Weil:Serre_tensor}).

We assume $2 \nmid \Delta$, but allow $2$ inert or split in $\mc{O}_F$ for the moment. When $n = 2$ and $L$ is a self-dual Hermitian $\mc{O}_F$-lattice of signature $(1,1)$, recall
    \begin{equation}
    \mc{M} = \ms{M}_0 \times_{\Spec \mc{O}_F} \ms{M}(1,1)^{\circ}
    \end{equation}
in the notation of \crefext{I:ssec:part_I:arith_intersections:integral_models}. Recall that $\ms{M}_0$ is the moduli stack parameterizing $(A_0, \iota_0, \lambda_0)$ where $A_0$ is an elliptic curve with signature $(1,0)$ action $\iota_0$ by $\mc{O}_F$, and $\lambda_0$ the unique principal polarization. Recall that $\ms{M}(1,1)^{\circ}$ is the closure of the generic fiber in the moduli stack of signature $(1,1)$ Hermitian abelian schemes $(A, \iota, \lambda)$ where $|\Delta| \cdot \lambda$ is a polarization with $\ker (|\Delta| \cdot \lambda) = A[\sqrt{\Delta}]$.

For integers $j > 0$, we first recall how to relate the special cycles $\mc{Z}(j) \ra \mc{M}$ to Hecke translates of CM elliptic curves, as explained in \cite[{\S 14}]{KR14}. Our $|\Delta| \cdot \lambda$ is their $\lambda$.

Write $\ms{M}_{\text{ell}}$ for the moduli stack of elliptic curves base-changed to $\Spec \mc{O}_F$. If $\mc{O}_F^* \coloneqq \Hom_{\Z}(\mc{O}_F, \Z)$, we write $\lambda_{\mrm{tr}} \colon \mc{O}_F \ra \mc{O}_F^*$ for the $\s$-linear map corresponding to the symmetric $\Z$-bilinear pairing $\mrm{tr}_{F / \Q}(a^{\s} b)$ on $\mc{O}_F$. As in \cite[{\S 14}]{KR14}, there is a \emph{Serre tensor} morphism
    \begin{equation}\label{equation:arithmetic_Siegel-Weil:Serre_tensor}
    \begin{tikzcd}[row sep = tiny]
    \ms{M}_{\text{ell}} \arrow{r}{i_{\text{Serre}}} & \ms{M}(1,1)^{\circ} \\
    E \arrow[mapsto]{r} & E \otimes_{\Z} \mc{O}_F
    \end{tikzcd}
    \end{equation}
where $E \otimes_{\Z} \mc{O}_F$ is given the polarization $|\Delta|^{-1} (\lambda_E \otimes \lambda_{\mrm{tr}}) \colon E \otimes_{\Z} \mc{O}_F \ra E^{\vee} \otimes_{\Z} \mc{O}_F^*$. As we have seen previously, $E \otimes_{\Z} \mc{O}_F$ is (by definition) the functor given by $(E \otimes_{\Z} \mc{O}_F)(S') = E(S') \otimes_{\Z} \mc{O}_F$ for schemes $S'$ (over the understood base for $E$).

For the rest of Section \ref{ssec:arithmetic_Siegel-Weil:Serre_tensor}, we now assume $\mc{O}_F^{\times} = \{ \pm 1 \}$. In this case, the Serre tensor morphism is an open and closed immersion.\footnote{The hypothesis $\mc{O}_F^{\times} = \{ \pm 1\}$ should be added in \cite[{Proposition 14.4}]{KR14}, as otherwise $\Aut(E) \neq \Aut(E \otimes_{\Z} \mc{O}_F)$ (right-hand side means $\mc{O}_F$-linear automorphisms preserving the polarization) so $i_{\text{Serre}} \colon \ms{M}_{\mrm{ell}} \ra \ms{M}(1,1)^{\circ}$ is not a monomorphism and hence cannot be a closed immersion in the sense of \cite[\href{https://stacks.math.columbia.edu/tag/04YK}{Section 04YK}]{stacks-project}. The remaining arguments are the same at least if $2 \nmid \Delta$.} Indeed, $i_{\text{Serre}}$ is proper (valuative criterion) and a monomorphism of algebraic stacks, hence a closed immersion of algebraic stacks. Since the source and target are Deligne--Mumford, smooth, finite type, and separated over $\Spec \mc{O}_F$ of the same relative dimension, this implies that $i_{\text{Serre}}$ is also an open immersion. 

The class group $\mrm{Cl}(\mc{O}_F)$ acts on $\ms{M}(1,1)^{\circ}$ as follows. Given any fractional ideal $\mf{a} \subseteq F$, set $\mf{a}^{\vee} \coloneqq \Hom_{\mc{O}_F}(\mf{a}, \mc{O}_F)$, and consider the $\s$-linear map $\lambda_{\mf{a}} \colon \mf{a} \xra{\sim} \mf{a}^{\vee}$ given by the perfect positive-definite Hermitian pairing $a, b \mapsto N(\mf{a})^{-1} a^{\s} b$ on $\mf{a}$. There is an induced automorphism of $\ms{M}(1,1)^{\circ}$ sending 
    \begin{equation}
    (A, \iota, \lambda) \ra (A \otimes_{\mc{O}_F} \mf{a}, \iota, \lambda \otimes \lambda_{\mf{a}}).
    \end{equation}
The action of $\mrm{Cl}(\mc{O}_F)$ on $\ms{M}(1,1)^{\circ}$ is simply transitive on the set of connected components (see the proof of \cite[{Proposition 14.4}]{KR14}). 
There is a similar action of $\mrm{Cl}(\mc{O}_F)$ on $\ms{M}_0$ which sends $(A_0, \iota_0, \lambda_0) \mapsto (A_0 \otimes_{\mc{O}_F} \mf{a}, \iota_0, \lambda_0 \otimes \lambda_{\mf{a}})$. Given a fractional ideal $\mf{a} \subseteq F$, we write $f_{\mf{a}} \colon \mc{M} \ra \mc{M}$ for the induced automorphism just described.

Given any integer $j > 0$, the action of $\mrm{Cl}(\mc{O}_F)$ preserves $\mc{Z}(j)$, in the sense that there is a $2$-Cartesian diagram
    \begin{equation}
    \begin{tikzcd}
    \mc{Z}(j) \arrow{r}{\tilde{f}_{\mf{a}}} \arrow{d} & \mc{Z}(j) \arrow{d} \\
    \mc{M} \arrow{r}{f_{\mf{a}}} & \mc{M}
    \end{tikzcd}
    \end{equation}
for any fractional ideal $\mf{a}$, where $\tilde{f}_{\mf{a}}$ sends
    \begin{equation}
    (A_0, \iota_0, \lambda_0, A, \iota, \lambda, x) \mapsto (A_0 \otimes_{\mc{O}_F} \mf{a}, \iota_0, \lambda_0 \otimes \lambda_{\mf{a}}, A \otimes_{\mc{O}_F} \mf{a}, \iota, \lambda \otimes \lambda_{\mf{a}}, x \otimes 1)
    \end{equation}
for $x \in \Hom_{\mc{O}_F}(A_0, A)$ satisfying $x^{\dagger} x = j$.

Consider the \emph{$j$-th Hecke correspondence} $\mc{T}_j \ra \ms{M}_0 \times_{\Spec \mc{O}_F} \ms{M}_{\text{ell}}$, where $\mc{T}_j$ is the stack parameterizing tuples $(E_0, \iota_0, \lambda_0, E, w)$ for $(E_0, \iota_0, \lambda_0) \in \ms{M}_0$, for $E \in \ms{M}_{\text{ell}}$, and $w \colon E \ra E_0$ an isogeny of degree $j$. 

Consider the map $\ms{M}_0 \times \ms{M}_{\text{ell}} \ra \mc{M}$ induced by $i_{\text{Serre}}$ (and the identity on $\ms{M}_0$).
The Kudla--Rapoport cycle $\mc{Z}(j)$ pulls back to the Hecke correspondence $\mc{T}_j$, i.e. there is a $2$-Cartesian diagram
    \begin{equation}\label{equation:arithmetic_Siegel-Weil:Serre_tensor:pullback_to_Hecke}
    \begin{tikzcd}
    \mc{T}_j \arrow{d} \arrow{r} & \mc{Z}(j) \arrow{d} \\
    \ms{M}_0 \times_{\Spec \mc{O}_F} \ms{M}_{\text{ell}} \arrow{r} & \mc{M}
    \end{tikzcd}
    \end{equation}
where $\mc{T}_j \ra \mc{Z}(j)$ sends
    \begin{equation}
    (E_0, \iota_0, \lambda_0, E, w) \mapsto (E_0, \iota_0, \lambda_0, E \otimes_{\Z} \mc{O}_F, \iota, \lambda_E \otimes \lambda_{\mrm{tr}}, x_w)
    \end{equation}
(with $\lambda_E$ denoting the unique principal polarization of $E$) and where $x_w \colon E_0 \ra E \otimes_{\Z} \mc{O}_F$ is the $\mc{O}_F$-linear map such that $\sqrt{\Delta} x_w^{\dagger} \in \Hom_{\mc{O}_F}(E \otimes_{\Z} \mc{O}_F, E_0)$ corresponds to $w$ via the adjunction
    \begin{equation}
    \Hom_{\mc{O}_F}(E \otimes_{\Z} \mc{O}_F, E_0) = \Hom(E, E_0).
    \end{equation}
Here, we are implicitly claiming $\deg(w) = x^{\dagger}_w x_w$. The fact that \eqref{equation:arithmetic_Siegel-Weil:Serre_tensor:pullback_to_Hecke} is well-defined and $2$-Cartesian is proved in \cite[{Proposition 14.5}]{KR14}.

We next discuss the Eisenstein series of Theorem \ref{theorem:arithmetic_Siegel-Weil:main_results:main}(2) in more elementary terms when $n = 2$. In this case, the $U(1,1)$ Eisenstein series $E^*(z, s)^{\circ}_2$ (with $m = 1$ in our usual notation, and normalized as in Section \ref{ssec:Eisenstein:global_normalized_Fourier:global_normalization}) admits the classical expression
    \begin{equation}
    E^*(z, s)^{\circ}_2 = - \frac{\pi^{-s + 1/2}}{8 \pi^2} \Gamma(s + 3/2) \zeta(2s + 1) \sum_{\substack{c,d \in \Z \\ (c,d) = 1}} \frac{y^{s - 1/2}}{(c z + d)^2 |c z + d|^{2(s - 1/2)}}
    \end{equation}
for $z = x + i y \in \mc{H}$, where $\mc{H} \subseteq \C$ is the usual upper half-space (here $z$ corresponds to $z^{\flat}$ in Theorem \ref{theorem:arithmetic_Siegel-Weil:main_results:main}(2)).

For nonzero $j \in \Z$, the (normalized) $j$-th Fourier coefficient of $E^*(z,s)^{\circ}_2$ factorizes into (normalized) local Whittaker functions
    \begin{equation}
    E^*_j(y,s)^{\circ}_2 = W^*_{j, \infty}(y,s)^{\circ}_2 \prod_{p} W^*_{j,p}(s)^{\circ}_2
    \end{equation}
as in Section \ref{ssec:Eisenstein:global_normalized_Fourier:global_normalization}.
We have the formulas
    \begin{equation}\label{equation:arithmetic_Siegel-Weil:Serre_tensor:explicit_Whittaker}
    W^*_{j,p}(s)^{\circ}_2 = p^{v_p(j) (s + 1/2)} \s_{-2s}(p^{v_p(j)})  \quad \quad  \prod_{p} W^*_{j,p}(s)^{\circ}_2 = |j|^{s + 1/2} \s_{-2s}(|j|)
    \end{equation}
where $v_p(-)$ means $p$-adic valuation and
    \begin{equation}
    \s_s(|j|) \coloneqq \sum_{d \mid |j|} d^s
    \end{equation}
is the classical divisor function. These formulas for local Whittaker functions are likely classical, but they also follow from \crefext{I:equation:non-Arch_idnentity:horizontal_identity:rank_one_Den} on local densities (translation to local Whittaker functions via \eqref{equation:local_Whittaker:local_densities_spherical:Den_star}). An integral expression for $W^*_{j,\infty}(y,s)^{\circ}_2$ may be found in \crefext{II:ssec:Archimedean_identity:case_n_is_2}. For $j > 0$, recall $W^*_{j,\infty}(y,1/2)^{\circ}_2 = 1$ \eqref{equation:Eisenstein:local_Whittaker:Archimedean:special_value}.

We require $j > 0$ for the rest of Section \ref{ssec:arithmetic_Siegel-Weil:Serre_tensor}. Fix an embedding $F \ra \C$. Given a CM elliptic curve $(E_0, \iota_0, \lambda_0) \in \ms{M}_0(\C)$, we consider the set of \emph{$j$-th Hecke translates} of $E_0$ given by
    \begin{equation}
    \mc{T}_j(E_0) \coloneqq \{ (E_0, \iota_0, \lambda_0, E, w) \in \mc{T}_j(\C) \}.
    \end{equation}
Phrased alternatively, the fiber of $\mc{T}_j \ra \ms{M}_0$ over the point $\Spec \C \ra \ms{M}_0$ corresponding to $E_0$ is a finite scheme over $\Spec \C$,
and $\mc{T}_j(E_0)$ is its set of $\C$-points. We set
    \begin{equation}
    \deg \mc{T}_j(E_0) \coloneqq | \mc{T}_j(E_0) | \quad \quad h_{\mrm{Fal}}(\mc{T}_j(E_0)) \coloneqq \sum_{E \in \mc{T}_j(E_0)} h_{\mrm{Fal}}(E)
    \end{equation}
where $|-|$ denotes set cardinality, the sum runs over $(E_0, \iota_0, \lambda_0, E, w) \in \mc{T}_j(E_0)$, and $h_{\mrm{Fal}}(E)$ denotes the Faltings height of $E$ (with metric normalized as in \crefext{I:equation:arith_cycle_classes:Hodge_bundles:Faltings_metric}, see also \crefext{I:ssec:Faltings_and_taut:heights}) after descending from $\C$ to any number field.

The following lemma states that the (total) Faltings height of $j$-th Hecke translates of a chosen elliptic curve with CM by $\mc{O}_F$ does not depend on the choice of CM elliptic curve. It should admit a general formulation in terms of Hecke correspondences over $\ms{M}_0$. We give a more elementary treatment in the spirit of this section.

\begin{lemma}\label{lemma:arithmetic_Siegel-Weil:Serre_tensor:framing_object_independence}
Fix $j \in \Z_{>0}$. For any $(E_0, \iota_0, \lambda_0) \in \ms{M}_0(\C)$ and $(E'_0, \iota'_0, \lambda'_0) \in \ms{M}_0(\C)$, we have
    \begin{equation}
    \deg \mc{T}_j(E_0) = \deg \mc{T}_j(E'_0) \quad \quad h_{\mrm{Fal}} (\mc{T}_j(E_0)) = h_{\mrm{Fal}}(\mc{T}_j(E'_0)).
    \end{equation}
\end{lemma}
\begin{proof}
Given any $d \in \Z$, we claim that there exists an isogeny $\phi \colon E'_0 \ra E_0$ of degree prime to $d$. Consider
    \begin{equation}
    E_0(\C) = \C / \Lambda_0 \quad \quad E'_0(\C) = \C / \Lambda'_0
    \end{equation}
for lattices $\Lambda_0$ and $\Lambda'_0$. Without loss of generality, we may assume $\Lambda_0 = \mc{O}_F \subseteq \C$ and that $\Lambda'_0 = \mf{a}'_0$ for some fractional ideal $\mf{a}'_0 \subseteq \C$. By the Chinese remainder theorem, we can assume $\mf{a}'_0 \subseteq \mc{O}_F$ and that $\mf{a}'_0$ has norm prime to $d$ (without changing the ideal class of $\mf{a}'_0$). The inclusion $\mf{a}'_0 \subseteq \mc{O}_F$ gives an isogeny $E'_0 \ra E_0$ of degree prime to $d$.

Let $p$ be any prime. Let $\phi \colon E'_0 \ra E_0$ be an isogeny of degree prime to $p j$. As above, we view $\phi \colon E_0(\C) \ra E'_0(\C)$ as an inclusion of lattices $\Lambda'_0 \ra \Lambda_0$ of index prime to $p j$. There is an induced bijection
    \begin{equation}
    \begin{tikzcd}[row sep = tiny]
    \mc{T}_j(E_0) \arrow{r} & \mc{T}_j(E'_0) \\
    \Lambda \arrow[mapsto]{r} & \Lambda \cap \Lambda'_0.
    \end{tikzcd}
    \end{equation}
We are viewing $\Lambda$ as the element $\C / \Lambda \ra \C / \Lambda_0$ of $\mc{T}_j(E_0)$, and similarly for $\Lambda \cap \Lambda'_0$.

The isogeny $\C / (\Lambda \cap \Lambda'_0) \ra \C / \Lambda$ has degree $\deg \phi$, which is prime to $p$. As these elliptic curves are defined over $\overline{\Q}$, this isogeny also descends to $\overline{\Q}$. By the formula for change of Faltings height along an isogeny \crefext{I:equation:Faltings_and_taut:isogeny_change_global:Faltings:log_p}, we conclude $h_{\mrm{Fal}}(\mc{T}_j(E_0)) - h_{\mrm{Fal}}(\mc{T}_j(E'_0)) \in \sum_{\ell \mid \deg \phi} \Q \cdot \log \ell$. Varying $p$ shows $h_{\mrm{Fal}}(\mc{T}_j(E_0)) = h_{\mrm{Fal}}(\mc{T}_j(E'_0))$, as the real numbers $\log p$ are $\Q$-linearly independent for varying $p$.
\end{proof}

Consider any $(E_0, \iota_0, \lambda_0) \in \ms{M}_0(\C)$. Using \eqref{equation:arithmetic_Siegel-Weil:Serre_tensor:pullback_to_Hecke} (Kudla--Rapoport cycle pulls back to Hecke correspondence), the geometric Siegel--Weil statement in Remark \ref{remark:geometric_Siegel-Weil:degrees:self-dual} implies
    \begin{equation}\label{equation:arithmetic_Siegel-Weil:Serre_tensor:geometric_Siegel-Weil}
    \frac{h^2_F}{w_F} \deg \mc{T}_j(E_0) = 2 \frac{h^2_F}{w^2_F} E^*_{j}(y, 1/2)^{\circ}_2
    \end{equation}
for any $y \in \R_{>0}$. On the left, one factor of $h_F$ appears because the Serre tensor morphism $i_{\text{Serre}} \colon \ms{M}_{\text{ell}} \ra \ms{M}(1,1)^{\circ}$ is the inclusion of one connected component (and $\ms{M}(1,1)^{\circ}$ has $h_F$ connected components, by the action of $\mrm{Cl}(\mc{O}_F)$ discussed above; we discussed that this action is compatible with Kudla--Rapoport cycles). On the left, the additional factor $h_F / w_F$ appears via Lemma \ref{lemma:arithmetic_Siegel-Weil:Serre_tensor:framing_object_independence} (instead of summing over $\ms{M}_0(\C)$, it is enough to consider a fixed $E_0$ and multiply by $h_F / w_F = \deg_{\C} (\ms{M}_0 \times_{\Spec \mc{O}_F} \Spec \C)$).

By the formulas in \eqref{equation:arithmetic_Siegel-Weil:Serre_tensor:explicit_Whittaker} and surrounding discussion, this recovers the well-known identity $\deg \mc{T}_j(E_0) = \s_1(j)$ for degrees of Hecke correspondences (recall our running assumption $|\mc{O}_F^{\times}| = \{ \pm 1 \}$ for most of Section \ref{ssec:arithmetic_Siegel-Weil:Serre_tensor}, i.e. $w_F = 2$).

In the next lemma, $h^{\mrm{CM}}_{\mrm{Fal}} = h_{\mrm{Fal}}(E_0)$ is the Faltings height of any elliptic curve with CM by $\mc{O}_F$, normalized as in \crefext{III:equation:part_II:arith_cycle_classes:Hodge_bundles:CM_Faltings_height}.
It is well known that this does not depend on the choice of CM elliptic curve (also follows from Lemma \ref{lemma:arithmetic_Siegel-Weil:Serre_tensor:framing_object_independence}).

\begin{corollary}\label{corollary:arithmetic_Siegel-Weil:Serre_tensor}
Suppose $2$ is split in $\mc{O}_F$. For any integer $j > 0$ and any CM elliptic curve $(E_0, \iota_0, \lambda_0) \in \ms{M}_0(\C)$, we have
    \begin{equation}
    h_{\mrm{Fal}}(\mc{T}_j(E_0)) - \s_1(j) \cdot h^{\mrm{CM}}_{\mrm{Fal}} = \frac{1}{2} \frac{d}{ds} \bigg|_{s = 1/2} \left ( j^{s + 1/2} \s_{-2s}(j) \right ).
    \end{equation}
\end{corollary}
\begin{proof}
Set $n = 2$ and consider the $2 \times 2$ matrix $T = \mrm{diag}(0, j)$. 
Again using \eqref{equation:arithmetic_Siegel-Weil:Serre_tensor:pullback_to_Hecke} to pull back Kudla--Rapoport cycles to Hecke correspondences, we have
    \begin{equation}\label{equation:arithmetic_Siegel-Weil:Serre_tensor:height_difference}
    2 \frac{h^2_F}{w_F} \left ( 2 h_{\mrm{Fal}}(\mc{T}_j(E_0)) - 2 (\deg \mc{T}_j(E_0)) \cdot h_{\mrm{Fal}}^{\mrm{CM}} \right ) = 
    - 2 \sum_{p} \mrm{Int}_{\ms{H}, p, \mrm{global}}(T)
    \end{equation}
in our previous notation (Remark \ref{remark:arithmetic_Siegel-Weil:main_results:Faltings_height}). On the left, the outer factor of $2$ has the same explanation as in \eqref{equation:arithmetic_Siegel-Weil:main_results:Faltings_height:bundle_degree} (see following discussion). The factor $h^2_F/w_F$ has the same explanation as in \eqref{equation:arithmetic_Siegel-Weil:Serre_tensor:geometric_Siegel-Weil}, via Lemma \ref{lemma:arithmetic_Siegel-Weil:Serre_tensor:framing_object_independence} on Faltings height. The factor of $2$ in $2 h_{\mrm{Fal}}(\mc{T}_j(E_0))$ appears because $h_{\mrm{Fal}}(E \otimes_{\Z} \mc{O}_F) = h_{\mrm{Fal}}(E \times E) = 2 h_{\mrm{Fal}}(E)$. The factor of $2$ in $2 (\deg \mc{T}_j(E_0)) \cdot h_{\mrm{Fal}}^{\mrm{CM}}$ is the $n$ in Remark \ref{remark:arithmetic_Siegel-Weil:main_results:Faltings_height}. 

In our previous notation, we have $\mrm{Int}_{\ms{V},p,\mrm{global}}(T) = 0$ for all primes $p$ as the vertical special cycle class ${}^{\mbb{L}} \mc{Z}(T)_{\ms{V},p}$ is $0$ when $n = 2$ \crefext{III:lemma:non-Arch_uniformization:split_proper_quasi-finite_global}. Hence $\mrm{Int}_{p, \mrm{global}}(T) = \mrm{Int}_{\ms{H}, p, \mrm{global}}(T) + \mrm{Int}_{\ms{V},\ell,\mrm{global}}(T) = \mrm{Int}_{\ms{H}, p, \mrm{global}}(T)$.

Then \eqref{equation:arithmetic_Siegel-Weil:main_results:main:geometric_decomp:at_p:global} (``horizontal local part'' of our main result) implies
    \begin{equation}\label{equation:arithmetic_Siegel-Weil:Serre_tensor:Int_global}
    \mrm{Int}_{p, \mrm{global}}(T) = - \frac{2 h_F^2}{w^2_F} \left ( \frac{d}{ds} \bigg|_{s = 1/2} W^*_{j,p}(s)^{\circ}_2 \right ) \prod_{\ell \neq p} W_{j,\ell}^*(1/2)^{\circ}_2
    \end{equation}
for all $p$ (in the notation of loc. cit., take $T^{\flat} = j$, $a^{\flat}_{v} = 1$ for all $v < \infty$, and recall our notation $\tilde{W}^*_{T^{\flat},v}(1,s)^{\circ}_n = W^*_{T^{\flat},v}(1,s)^{\circ}_n \eqqcolon W^*_{T^{\flat},v}(s)^{\circ}_n$). Since $j > 0$, we have used $W^*_{j,\infty}(1/2)^{\circ}_2 = 1$ \eqref{equation:Eisenstein:local_Whittaker:Archimedean:special_value} as recalled above.

Combining \eqref{equation:arithmetic_Siegel-Weil:Serre_tensor:Int_global} and \eqref{equation:arithmetic_Siegel-Weil:Serre_tensor:height_difference} along with the formula $\deg \mc{T}_j(E_0) = \s_1(j)$, we obtain
    \begin{equation}
    h_{\mrm{Fal}}(\mc{T}_j(E_0)) - \s_1(j) \cdot h^{\mrm{CM}}_{\mrm{Fal}} = \frac{1}{2} \frac{d}{ds} \bigg|_{s = 1/2} \left ( \prod_{p} W^*_{j,p}(1/2)^{\circ}_2 \right )
    \end{equation}
where the product runs over all primes (not including the Archimedean place). The corollary now follows from the formulas in \eqref{equation:arithmetic_Siegel-Weil:Serre_tensor:explicit_Whittaker}.
\end{proof}
    
    \clearpage


    \phantomsection
    \addcontentsline{toc}{part}{References}
    \renewcommand{\addcontentsline}[3]{}
    \printbibliography

\end{document}